\documentclass[12pt]{amsart}
\usepackage{amssymb}
\usepackage{amsthm}
\usepackage{amscd}

\usepackage{amsmath}
\usepackage{epic,eepic}
\usepackage{graphicx}
\DeclareMathAlphabet{\mathpzc}{OT1}{pzc}{m}{it}

\usepackage{latexsym}

\usepackage{pifont}

\theoremstyle{plain}
\newtheorem{lemma}{Lemma}[subsection]
\newtheorem{prop}[lemma]{Proposition}
\newtheorem{thm}[lemma]{Theorem}
\newtheorem{cor}[lemma]{Corollary}
\newtheorem{aplemma}{Lemma~A.\hspace{-1.5mm}}
\newtheorem{approp}{Proposition~A.\hspace{-1.5mm}}
\newtheorem{apthm}{Theorem~A.\hspace{-1.5mm}}
\newtheorem{apcor}{Corollary~A.\hspace{-1.5mm}}
\newtheorem{intthm}{Theorem}

\usepackage[all]{xy}

\newtheorem{conj}[lemma]{Conjecture}

\usepackage{mathrsfs}

\theoremstyle{definition}

\newtheorem{rema}[lemma]{Remark}

\newtheorem{remb}{Remark}

\newtheorem{defi}[lemma]{Definition}
\newtheorem{exa}[lemma]{Example}
\newtheorem{aprem}{Remark~A.\hspace{-1.5mm}}
\newtheorem{apdefi}{Definition~A.\hspace{-1.5mm}}
\newcommand{\bde}{\begin{defi}}
\newcommand{\ede}{\end{defi}\vspace{1mm}}
\newcommand{\ble}{\begin{lemma}}
\newcommand{\ele}{\end{lemma}}
\newcommand{\bpr}{\begin{prop}}
\newcommand{\epr}{\end{prop}}
\newcommand{\bt}{\begin{thm}}
\newcommand{\et}{\end{thm}}
\newcommand{\bco}{\begin{cor}}
\newcommand{\eco}{\end{cor}}
\newcommand{\bre}{\begin{rema}}
\newcommand{\ere}{\end{rema}}
\newcommand{\brea}{\begin{rema}}
\newcommand{\erea}{\end{rema}\vspace{1mm}}
\newcommand{\breb}{\begin{remb}}
\newcommand{\ereb}{\end{remb}\vspace{1mm}}
\newcommand{\bex}{\begin{exa}}
\newcommand{\eex}{\end{exa}}
\newcommand{\bpf}{\begin{proof}}
\newcommand{\epf}{\end{proof}\vspace{1mm}}

\newcommand{\bade}{\begin{apdefi}}
\newcommand{\eade}{\end{apdefi}}
\newcommand{\bale}{\begin{aplemma}}
\newcommand{\eale}{\end{aplemma}}
\newcommand{\bapr}{\begin{approp}}
\newcommand{\eapr}{\end{approp}}
\newcommand{\bat}{\begin{apthm}}
\newcommand{\eat}{\end{apthm}}
\newcommand{\baco}{\begin{apcor}}
\newcommand{\eaco}{\end{apcor}}
\newcommand{\bare}{\begin{aprem}}
\newcommand{\eare}{\end{aprem}}

\newcommand{\be}{\begin{enumerate}}
\newcommand{\ee}{\end{enumerate}}
\newcommand{\bcd}{\[\begin{CD}}
\newcommand{\ecd}{\end{CD}\]}
\newcommand{\bit}{\begin{itemize}}
\newcommand{\eit}{\end{itemize}}
\newcommand{\bq}{\begin{quote}}
\newcommand{\eq}{\end{quote}}
\newcommand{\ba}{\begin{array}}
\newcommand{\ea}{\end{array}}
\newcommand{\mcD}{\mathcal{D}}
\newcommand{\mcE}{\mathcal{E}}
\newcommand{\mcF}{\mathcal{F}}
\newcommand{\mcG}{\mathcal{G}}
\newcommand{\mcH}{\mathcal{H}}

\newcommand{\mcK}{\mathcal{K}}
\newcommand{\mcL}{\mathcal{L}}

\newcommand{\mcO}{\mathcal{O}}
\newcommand{\mcP}{\mathcal{P}}
\newcommand{\mcQ}{\mathcal{Q}}

\newcommand{\mcS}{\mathcal{S}}
\newcommand{\mcT}{\mathcal{T}}
\newcommand{\mcU}{\mathcal{U}}
\newcommand{\mcV}{\mathcal{V}}
\newcommand{\mcW}{\mathcal{W}}

\newcommand{\mbA}{\mathbb{A}}

\newcommand{\mbC}{\mathbb{C}}
\newcommand{\mbD}{\mathbb{D}}

\newcommand{\mbF}{\mathbb{F}}
\newcommand{\mbG}{\mathbb{G}}
\newcommand{\mbH}{\mathbb{H}}

\newcommand{\mbP}{\mathbb{P}}

\newcommand{\mbR}{\mathbb{R}}

\newcommand{\mbZ}{\mathbb{Z}}
\newcommand{\mfA}{\mathfrak{A}}

\newcommand{\mfM}{\mathfrak{M}}

\newcommand{\mfS}{\mathfrak{S}}

\newcommand{\mfa}{\mathfrak{a}}
\newcommand{\mfb}{\mathfrak{b}}

\newcommand{\mfe}{\mathfrak{e}}

\newcommand{\mfg}{\mathfrak{g}}
\newcommand{\mfh}{\mathfrak{h}}

\newcommand{\mfm}{\mathfrak{m}}

\newcommand{\mfr}{\mathfrak{r}}

\newcommand{\mft}{\mathfrak{t}}

\newcommand{\migi}{\rightarrow}
\newcommand{\longmigi}{\longrightarrow}

\newcommand{\isom}{\stackrel{\sim}{\migi}}

\newcommand{\migiincl}{\hookrightarrow}

\newcommand{\migisurj}{\twoheadrightarrow}

\newcommand{\SSP}{\vspace{3mm}}
\newcommand{\LSP}{\vspace{5mm}}

\newcommand{\Y}{X}

\newcommand{\ST}{\blacklozenge}

\newcommand{\mr}{\mathrm}
\newcommand{\hidden}[1]{\,}

\newcommand{\N}{N}
\newcommand{\M}{m}

\newcommand{\PP}{P}

\newcommand{\hh}{\heartsuit}
\newcommand{\cc}{\clubsuit}
\newcommand{\sss}{\spadesuit}
\newcommand{\dd}{\diamondsuit}

\usepackage[all]{xy}
\usepackage {cancel}

\begin{document}

\title[Frobenius-Ehresmann structures and   Cartan geometries in char.\,$p>0$]{Frobenius-Ehresmann structures and \\  Cartan geometries  in positive characteristic}
\author{Yasuhiro Wakabayashi}
\date{}
\markboth{Yasuhiro Wakabayashi}{}
\maketitle
\footnotetext{Y. Wakabayashi: Department of Mathematics, Tokyo Institute of Technology, 2-12-1 Ookayama, Meguro-ku, Tokyo 152-8551, JAPAN;}
\footnotetext{e-mail: {\tt wkbysh@math.titech.ac.jp};}
\begin{abstract}
The aim of the present paper is to lay  the foundation for a theory of Ehresmann structures in positive characteristic, generalizing the Frobenius-projective and Frobenius-affine structures defined in the previous work. This theory deals with atlases  of  \'{e}tale coordinate charts  on varieties  modeled on homogenous spaces of algebraic groups, which we call   Frobenius-Ehresmann structures. These  structures are compared with Cartan geometries in positive characteristic, as well as with higher-dimensional  generalizations of dormant  indigenous bundles. In particular, we investigate the conditions under which   these geometric structures  are equivalent to each other. Also, we consider the classification problem of Frobenius-Ehresmann structures on algebraic curves. The latter half of the present paper discusses the deformation theory of  indigenous bundles in the algebraic setting. The tangent and obstruction spaces of various deformation functors are computed in terms of the hypercohomology groups of certain complexes. As a consequence, we formulate and prove the Ehresmann-Weil-Thurston principle for Frobenius-Ehresmann structures. This fact asserts that deformations of a  variety equipped with a Frobenius-Ehresmann structure are completely determined by their monodromy crystals.

\end{abstract}
\tableofcontents 
%\runninghead{}{}
%\markboth{Yasuhiro Wakabayashi}{An explicit formula for the number of dormant indigenous bundles}

%%%%%%%%%%%%%%%%%%%%%%%%%%%%%%%---[begin section]---%%%%%%%%%%%%%%
\section{Introduction}
%\vspace{5mm}

\LSP
%%%%%%%%%%%%%%%%%%%%%%%%%%%%%%%---[begin section]---%%%%%%%%%%%%%%
\subsection{Review of Ehresmann structures and Cartan geometries} \label{S01}

In the spirit of F. Klein's Erlangen program  of 1872,   the essence of classical geometry  is    the idea   
  that 
{\it geometry} is  the study of  properties on homogenous spaces that are invariant under  the action of transformations.
This idea includes, for instance, Euclidean geometry 
 by considering the isometry group on Euclidean space.
 On the basis of this idea  and those of H. Poincar\'{e}, E. Cartan,  and other mathematicians, C. Ehresmann 
  (cf. ~\cite{Ehr}) initiated a general study of geometric structures, which was later 
promoted by W. P. Thurston (cf. ~\cite{Thu}). 
 For  a homogenous space $P$ obtained as  the quotient of a finite-dimensional Lie group  $G$,
   we consider
 a maximal  atlas of coordinate charts on a manifold  into $P$  such that 
on any two overlapping coordinate patches, the change of  coordinates is given by  elements of 
 $G$;
this geometric structure is  called   an {\it Ehresmann $(G, P)$-structure} (following ~\cite{Gol1} and ~\cite{CG}), or simply a {\it $(G, P)$-structure}.
This notion allows us  to
 investigate  various  geometries in the unifying frame of the model geometry $(G, P)$.

A basic  example in the holomorphic setting is given by   
the complex projective line $\mbP^1_\mbC$ considered as 
a homogenous space for  
the M\"{o}bius group $\mr{PGL}_2 (\mbC)$.
A $(\mr{PGL}_2 (\mbC), \mbP^1_\mbC)$-structure on a Riemann surface is  known as a {\it (complex) projective structure},  and it  
appears in the study of second-order linear differential equations.
Also, flat $\mbP^1_\mbC$-bundles corresponding, via the Riemann-Hilbert correspondence,   to projective structures and their moduli spaces  have been studied under the name of {\it indigenous bundles} or {\it $\mr{PGL}_2$-opers} (cf., e.g.,  ~\cite{BD1}, ~\cite{Fr}, ~\cite{G2}).
They play a major role in understanding  the framework of the uniformization theorem.
As an important consequence of the uniformization theorem for Riemann surfaces, 
each  Riemann surface has a canonical projective structure that is obtained by collecting various local inverses of a universal covering map.
In the hyperbolic case, the monodromy of this projective structure  defines a Fuchsian subgroup of $\mr{PSL}(\mbR)$, which 
  reflects the geometry  of the underlying Riemann surface from the group-theoretic viewpoint.

The monodromy of projective structures can be generalized to  a general $(G, P)$, as above.
 Denote by $\mr{Hom}(\pi_1 (X), G)/G$  the representation variety 
    classifying   homomorphisms $\pi_1 (X) \migi G$ (up to an inner automorphism)  of the fundamental group $\pi_1 (X)$ of $X$ based at a fixed point. 
    It  may also be identified with  the moduli space of $G$-valued local systems on $X$.
Although there are numerous special cases  (cf. ~\cite{BG}, ~\cite{Hej}, ~\cite{Wei1}, ~\cite{Wei2}, ~\cite{Wei3}), 
  Thurston
  provided a guiding principle  for classifying  $(G, P)$-structures  by using $\mr{Hom}(\pi_1 (X), G)/G$   (cf. ~\cite[\S\,5.3.1]{Thu}).
 This  so-called  {\it Ehresmann-Weil-Thurston principle}  ensures that the map associating to
 each $(G, P)$-structure its monodromy representation $\pi_1 (X) \migi G$
  gives   a local homeomorphism between the deformation space of $(G, P)$-structures 
  and $\mr{Hom}(\pi_1 (X), G)/G$ (cf. ~\cite{CEG}, ~\cite{Gol1}, ~\cite{Lok}).

Also, Cartan generalized 
Klein's homogeneous model spaces to the corresponding infinitesimal notion (in exactly  the same way as Riemannian geometry is the infinitesimal version of Euclidean geometry), namely {\it Cartan geometries}.
Cartan geometries can be described by using  higher-dimensional generalizations of indigenous bundles (= $\mr{PGL}_2$-opers), which are  principal bundles equipped with a connection satisfying a certain  transversality condition.
A fundamental result in that theory is  that a Cartan geometry, or the corresponding  indigenous bundle, 
 admits a curvature tensor that is  zero exactly when it comes from an Ehresmann structure.
 See ~\cite{Sha} for a modern treatment of this theory.

\LSP
%%%%%%%%%%%%%%%%%%%%%%%%%%%%%%%---[begin section]---%%%%%%%%%%%%%%
\subsection{Previous work related to Frobenius-projective structures} \label{S0gotu1}

Then, is it possible to develop the geometry of Ehresmann structures and their  generalizations by Cartan  {\it  in positive characteristic}?
On the one hand, 
Cartan geometries and the corresponding  indigenous bundles can be defined   in any characteristic 
 because of the algebraic nature of their  formulations.
On the other hand, 
 due to its analytic formulation,
the definition of an Ehresmann structure cannot be adopted, at least  in positive characteristic, as it is.
Indeed, if we had defined it in a naive fashion, it would  be nothing other than a trivial example of a variety having such a structure.
By taking this unfortunate fact into account, we need to consider   appropriate replacements.
Regarding this matter, Y. Hoshi (cf. ~\cite{Hos2}) introduced {\it Frobenius-projective structures}  (of finite level) to deal with analogues  of complex projective  structures on a Riemann surface.
(He also introduced in ~\cite{Hos4} the notion of a {\it Frobenius-affine structure} as an analogue of an affine structure on a Riemann surface.)
Frobenius-projective structures  have provided a rich  and deep story under the identifications  with various  
equivalent realizations (in the case of  $\N=1$), e.g., {\it dormant indigenous bundles}   (=  {\it dormant $\mr{PGL}_2$-opers}) and projective connections with a full set of solutions.
See ~\cite{Mzk2}, ~\cite{Wak7},  and ~\cite{Wak8} for the study of 
 these objects and their moduli space
in the context of $p$-adic Teichm\"{u}ller theory developed by S. Mochizuki.

Also, in ~\cite{Wak6}, the author generalized Frobenius-projective  structures to higher-dimensional varieties (i.e., generalized the local model to the $n$-dimensional projective space $\mbP^n$ equipped with the natural  $\mr{PGL}_{n+1}$-action  for an arbitrary $n$) including the case of infinite level.
The  main subject of {\it loc.\,cit.} is  the classification problem, first studied by  S. Kobayashi and T. Ochiai (cf. ~\cite{KO1}, ~\cite{KO2}) for complex manifolds,  of varieties admitting   such structures.
 For example, 
 we proved 
 the positive characteristic version of Gunning's formula (cf. ~\cite[Theorem B]{Wak6}), which  gives a fundamental criterion  
   for the existence of  Frobenius-projective  structures on a variety. 
Various  other results proved there  are compared with  the classical  facts in the holomorphic setting; these comparisons show that
Frobenius-projective structures 
should be appropriate replacements for complex projective  structures.
 The present work is, in some sense, a continuation and generalization of those previous studies.
 In particular, it includes a very general construction of the correspondence 
   between Frobenius-projective structures  and dormant indigenous bundles.

\LSP
%%%%%%%%%%%%%%%%%%%%%%%%%%%%%%%---[begin section]---%%%%%%%%%%%%%%
\subsection{Frobenius-Ehresmann structures and Cartan geometries in char.$\,p>0$} \label{Sghhgotu1}

The aim of the present paper is 
 to lay
  the foundation of 
a theory of Ehresmann structures and Cartan geometries  in positive characteristic  for  arbitrary   local models not just $\mr{PGL}_{n+1}$. 
Let  $k$ be an algebraically close field  of characteristic $p >2$ and 
 $P$  
 a smooth homogeneous space  over  $k$ obtained as the quotient of  a connected smooth algebraic group $G$ by a closed subgroup $H$.
Also,  let $f: X \migi S$ be  a smooth morphism of schemes  over $k$.
For a positive integer $\N$, denote by $X^{(\N)}$ the 
$\N$-th relative Frobenius twist of $X$ over $S$ (cf. \S\,\ref{S1}), i.e., the base-change of $X$ via the $\N$-th iterate of the Frobenius endomorphism $F_S$ of $S$.
The sheaf $G_{X/S}^{(\N)}$  on $X^{(\N)}$ represented by $G$ can be  regarded  as a sheaf on $X$ via the underlying homeomorphism of the $\N$-th relative Frobenius morphism $F^{(\N)}_{X/S} : X \migi X^{(\N)}$.
Then, we define a {\it Frobenius-Ehresmann $(G, P)$-structure of level $\N$}, or an  {\it $F^\N$-$(G, P)$-structure} for short, on $X/S$  
 as a maximal collection $\mcS^\hh$ of \'{e}tale coordinate charts on $X$ forming a left $G_{X/S}^{(\N)}$-torsor (cf. Definition \ref{D023} for the precise definition).
The local sections of $G^{(\N)}_{X/S}$  acting on $\mcS^\hh$
play the same  role as 
 the coordinate changes appearing in the classical definition of an Ehresmann $(G, P)$-structure.
 Moreover,
   a compatible collection of $F^\N$-$(G, P)$-structures  for various $\N$ is called an {\it $F^\infty$-$(G, P)$-structure} (cf. Definition \ref{D027}). 
 For each $\N \in \mbZ_{> 0} \sqcup \{ \infty \}$, we shall write
 \begin{align}
 F^\N\text{-}\mr{Ehr}_{(G, P), X/S}
 \end{align}
(cf. (\ref{wer4456}), (\ref{wer4457}))  for  the set of  $F^\N$-$(G, P)$-structures   on $X/S$.
To simplify the notation slightly,  we will write $ F\text{-}\mr{Ehr}_{(G, P), X/S} :=  F^1\text{-}\mr{Ehr}_{(G, P), X/S}$,  and refer to any $F^1$-$(G, P)$-structure as an {\it $F$-$(G, P)$-structure}.
 Truncations to lower levels via Frobenius pull-back (cf. (\ref{gougoa})) give rise to a sequence  of sets
\begin{align}
 F^\infty\text{-}\mr{Ehr}_{(G, P), X/S} \migi \cdots \migi  F^\N\text{-}\mr{Ehr}_{(G, P), X/S} \migi \cdots  F^2\text{-}\mr{Ehr}_{(G, P), X/S} \migi  F\text{-}\mr{Ehr}_{(G, P), X/S}.
\end{align}
 Frobenius-Ehresmann structures, i.e., elements of sets in this sequence,  are the central characters of our discussion, and  they are specialized to  Frobenius-projective structures if 
 $G = \mr{PGL}_{n+1}$ ($n >0$) and $H$ is the subgroup $\mr{PGL}_{n+1}^\circledcirc$ (cf. (\ref{taae89})) of $\mr{PGL}_{n+1}$ consisting of elements stabilizing  the point $[1 : 0: 0: \cdots : 0]$  under an identification with the natural  actions on  the projective space $\mbP^n$.
 Some of the fundamental subjects in the theory  are  as follows:
 \begin{itemize}
 \item[(i)]
 Compare Frobenius-Ehresmann structures with Cartan geometries in positive characteristic in analogy with  the classical correspondence between  Ehresmann structures and flat Cartan geometries;
 \item[(ii)]
 Describe which model geometries $(G, P)$ a smooth scheme     $X/S$ may  have locally, and up to what level $\N$ it may  have an $F^\N$-$(G, P)$-structure;
 \item[(iii)]
 For a fixed model $(G, P)$ and $\N \in \mbZ_{>0}\sqcup \{ \infty \}$, describe all possible $F^\N$-$(G, P)$-structures on $X/S$ and all possible deformations  of a fixed $F^\N$-$(G, P)$-structure.
 \end{itemize}

 First, let us consider the relationship with Cartan geometries, i.e., (i) above.
We shall write $\mfg := \mr{Lie}(G)$, i.e., the Lie algebra of $G$.
Then,  just as in the classical case, 
 we can define a 
 Cartan geometry  on a scheme in positive characteristic.
 A {\it Cartan geometry}  with model $(G, H)$ on $X/S$  (cf. Definition \ref{D024} for the precise definition) is a pair
 \begin{align}
\mcE^\cc := (\mcE_H, \omega)
 \end{align}
 consisting of an $H$-bundle $\mcE_H$ on $X$ and a $\mfg$-valued $1$-form $\omega$ on $\mcE_H/S$ and 
 satisfying  the following conditions:
 \begin{itemize}
 \item
 The $\mcO_{\mcE_H}$-linear morphism $\omega^\triangleright  : \mcT_{\mcE_H/S} \migi \mcO_{\mcE_H} \otimes_k \mfg$  corresponding to $\omega$
 is an isomorphism and $H$-equivariant with $H$ acting on $\mfg$ via conjugation; 
 \item
 The restriction of $\omega$ to each fiber of the projection $\mcE_H \migi X$ coincides with the Maurer-Cartan form $\omega_H$ on $H$.
 \end{itemize}
 We say that the  Cartan geometry $\mcE^\cc$ is {\it flat}  if the curvature $\psi_\omega := d \omega + \frac{1}{2} \cdot [\omega, \omega]$ of $\omega$ (cf. (\ref{Erosmf}))  vanishes identically.
 The standard example  of a flat Cartan geometry is  
 defined on the homogeneous space $P$; it consists of 
  the natural projection $G \migisurj P$, considered as an $H$-bundle on $P$, together with the Maurer-Cartan form $\omega_G$ on $G$.
 More generally, 
 given an $F$-$(G, P)$-structure $\mcS^\hh$ on $X/S$, we can glue together  
   copies of the standard one $(G \migisurj P, \omega_G)$ defined on various 
\'{e}tale  charts defining  $\mcS^\hh$ and obtain a flat Cartan geometry $\mcS^{\hh \Rightarrow \cc}$.
 The resulting assignment $\mcS^\hh \mapsto \mcS^{\hh \Rightarrow \cc}$ 
 gives a map from $F\text{-}\mr{Ehr}_{(G, P), X/S}$ to the set of (isomorphism classes of) flat Cartan geometries  with model $(G, H)$ on $X/S$.
Unlike the classical case, this map is not surjective, meaning that
a flat Cartan geometry does not necessarily come from an $F$-$(G, P)$-structure.
Then, 
by what conditions on Cartan geometries can the image of this map be characterized?
The essential  ingredient in understanding  this matter is the $p$-curvature ${^p}\psi_\omega$ of  $\mfg$-valued $1$-forms  $\omega$ introduced in  Definition \ref{D00102}.
If $\mcE^\cc$ is a flat Cartan geometry as above, then 
 ${^p}\psi_\omega$ is defined as   the global section of $F_{X}^*(\Omega_{X/S} \otimes_k  \mfg)$ determined by 
\begin{align}
\langle F_X^{-1}(\partial), {^p}\psi_\omega \rangle = \omega^\triangleright (\partial)^{[p]} - \omega^\triangleright (\partial^{[p]}) + \sum_{i=1}^{p-1} \frac{s_i (\partial, \omega^\triangleright (\partial))}{i}
\end{align}
for each local section $\partial \in \mcT_{X/S}$, where $(-)^{[\partial]}$ denotes $p$-power operations and the $s_i$'s are defined at the beginning of \S\,\ref{SS0131}.
We shall say that  $\mcE^\cc$ is {\it $p$-flat} if ${^p}\psi_\omega$ vanishes identically.

Also, we will introduce a {\it dormant  indigenous $(G, H)$-bundle} (cf. Definitions \ref{D0013} and \ref{Ets395}), which generalizes  the classical notion of a dormant  indigenous bundle (= dormant $\mr{PGL}_2$-oper)  on an algebraic curve; this idea is, roughly speaking,  an $H$-bundle $\mcE_H$ together with  an $S$-connection $\nabla$ on the induced  $G$-bundle  $\mcE_G := \mcE_H \times^H G$ such that $\mcE_H$ is transverse to the horizontal distribution with respect to $\nabla$.
If $(G, H) = (\mr{PGL}_{n+1}, \mr{PGL}_{n+1}^\circledcirc)$,  it is as formulated in ~\cite[\S\,1.7, Definition 1.7.1]{Wak}, where it is called  an $F^1$-indigenous structure.
We  prove  (cf. Proposition \ref{P0234}) that 
  giving a $p$-flat Cartan geometry with model $(G, H)$ is equivalent to giving a dormant indigenous $(G, H)$-bundle.
Also, one of the fundamental results of the present paper  is 
  that
 the map 
 $\mcS^\hh \mapsto \mcS^{\hh \Rightarrow \cc}$
 mentioned above is  
 restricted to a bijective correspondence
 between $F$-$(G, P)$-structures and $p$-flat Cartan geometries with model $(G, H)$, or equivalently, dormant indigenous $(G, H)$-bundles  (cf. Theorem \ref{P0011}).

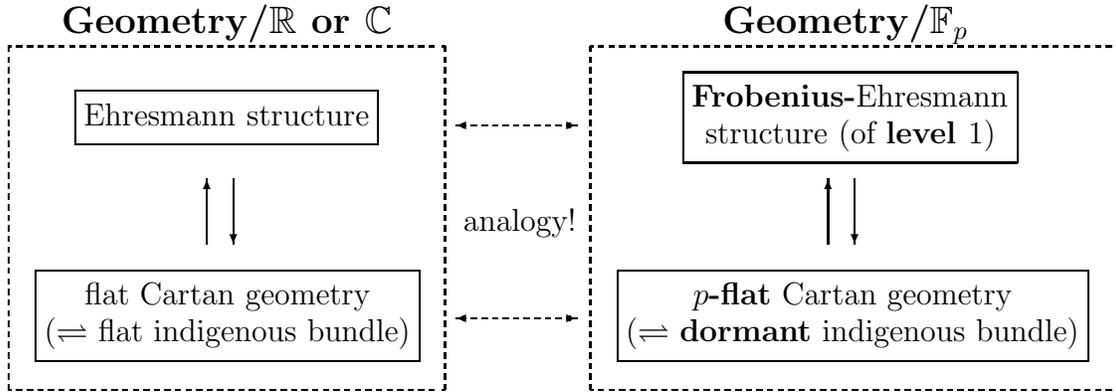
\vspace{2mm}
\begin{center}
\begin{picture}(400,150)

\put(0, 0){\dashbox{2.0}(165, 130){}}

\put(220, 0){\dashbox{2.0}(200, 130){}}

\put(172, 60){analogy!}

\put(25,100){\fbox{$\begin{matrix}
\text{Ehresmann structure}
% \\
%\text{structure of {\bf  level $1$}} 
%\text{\&} \\
%\text{decomp.\,arising from $\sigma$}
\end{matrix}$}}

\put(255,100){\fbox{$\begin{matrix}
\text{{\bf Frobenius-}Ehresmann} \\
\text{structure (of {\bf  level $1$})} 
\end{matrix}$}}

\put(10,25){\fbox{$\begin{matrix}
\text{flat Cartan geometry}
 \\
\text{($\rightleftharpoons$ flat indigenous bundle)} 
\end{matrix}$}}

\put(230,25){\fbox{$\begin{matrix}
\text{{\bf $p$-flat} Cartan geometry}
 \\
\text{($\rightleftharpoons$ {\bf dormant} indigenous bundle)} 
\end{matrix}$}}

\put(20,135){{\bf {\large Geometry/$\mbR$ or $\mbC$}}}
\put(270,135){{\bf {\large Geometry/$\mbF_p$}}}

\put(75,55){\vector(0,1){25}}
\put(85,80){\vector(0,-1){25}}

\put(310,55){\vector(0,1){25}}
\put(320,80){\vector(0,-1){25}}

\put(172,100){\vector(-1,0){3}}
\put(212,100){\vector(1,0){3}}
\put(172,27){\vector(-1,0){3}}
\put(212,27){\vector(1,0){3}}

\put(172, 100){\dashbox{2.0}(40, 0){}}

\put(172, 27){\dashbox{2.0}(40, 0){}}

\end{picture}
\vspace{2mm}
\end{center}
%-------------------------------------------------------------------

%To describe  these results in the following Theorem \ref{TheoremA1}, 
We shall denote by 
\begin{align}
\mr{Car}^{p\text{-}\mr{flat}}_{(G, H),  X/S}  \left(\text{resp.,} \  \mr{Ind}^{^\mr{Zzz...}}_{(G, H), X/S} \right)
\end{align}
(cf. (\ref{eoanv55677}), (\ref{e50q866})) the set of isomorphism classes of $p$-flat Cartan geometries with model $(G, H)$ (resp., dormant indigenous $(G, H)$-bundles) on $X/S$.
Then, 
what we  claimed above can be summarized below; 
this assertion is a positive characteristic version of Ehresmann's results  in ~\cite{Ehr} (cf. ~\cite[Theorem 2.2]{BD}).

\SSP
%---------------------------------------------------------------------[begin theorem]----------------------
\begin{intthm} [cf.  Theorem \ref{P0011}, Proposition \ref{P0234}] \label{TheoremA1}
Suppose that $G$ is affine.
Then, there exists  the following  commutative diagram in which  all the maps of sets  are bijective:
\begin{align}
\vcenter{\xymatrix@C=20pt@R=46pt{
&F\text{-}\mr{Ehr}_{(G, P), X/S} \ar[ld]^-{\sim}_-{\zeta^{\heartsuit \Rightarrow \clubsuit}}  \ar[rd]_-{\sim}^-{\zeta^{\heartsuit \Rightarrow \spadesuit}}&
\\
\mr{Car}^{p\text{-}\mr{flat}}_{(G, H), X/S} \ar[rr]^-{\sim}_-{\zeta^{\clubsuit \Rightarrow \spadesuit}}&& \mr{Ind}^{^\mr{Zzz...}}_{(G, H), X/S}.
}}
\end{align}
Moreover, the formations of these bijections commute with pull-back to \'{e}tale $X$-schemes, as well as with base-change to $S$-schemes.

 \end{intthm}
%\vspace{5mm}
%------------------------------------------------------------------------------[end theorem]-------------

\LSP
%%%%%%%%%%%%%%%%%%%%%%%%%%%%%%%---[begin section]---%%%%%%%%%%%%%%
\subsection{Generalizations to higher levels} \label{S099872}

 The next topic of the present paper will be  to extend 
%We will moreover discuss  $F^\N$-$(G, P)$-structures for large $\N$, and 
  the  bijections   in the above theorem   to higher level.
To do this, 
we equip a Cartan geometry $\mcE^\cc := (\mcE_H, \omega)$ with an additional  choice of
a $G$-bundle $\mcG$ on $X^{(\N)}$ together with an isomorphism $\left(\mcE_H \times^H G =:\right) \mcE_G \isom F^{(\N)*}_{X/S}(\mcG)$;
 this $G$-bundle $\mcG$ is taken to be compatible with $\omega$ in a certain sense.
We will call such an enriched Cartan geometry
an {\it  $F^\N$-Cartan geometry with model $(G, H)$} (cf. Definition \ref{Eotarmr90}).
Also, the notion of an $F^\infty$-Cartan geometry will be defined as a compatible system of $F^\N$-Cartan geometries for various $\N$.
 On the other hand, we introduce 
 (under the assumption that $G$ is affine)  {\it dormant indigenous $(G, H)$-bundles of level $\N$} (cf. Defnitions \ref{Eosi38392} and \ref{Ge990}) in terms of  both  Berthelot's PD-stratifications of higher level (cf. ~\cite{Ber1}, ~\cite{Ber2}) and the generalized notion of $p$-curvature (cf. ~\cite{LQ}).
 For each $\N \in \mbZ_{>0}\sqcup \{ \infty \}$,  denote by 
\begin{align}
F^\N\text{-}\mr{Car}_{(G, H),  X/S} \ \left(\text{resp.,} \  \mr{Ind}^{^\mr{Zzz...}}_{(G, H), X/S, \N}\right)
\end{align}
 the set of isomorphism classes of $F^\N$-Cartan geometries with model $(G, H)$  (resp., dormant indigenous $(G, H)$-bundles of level $\N$) on $X/S$.
Then, there exists  a  natural identification  $F^1\text{-}\mr{Car}_{(G, H),  X/S} = \mr{Car}_{(G, H),  X/S}^{p\text{-}\mr{flat}}$ (resp.,  $\mr{Ind}^{^\mr{Zzz...}}_{(G, H), X/S, 1} = \mr{Ind}^{^\mr{Zzz...}}_{(G, H), X/S}$).
In other words, the elements of $F^\N\text{-}\mr{Car}_{(G, H),  X/S}$ (resp., $\mr{Ind}^{^\mr{Zzz...}}_{(G, H), X/S, \N}$) may be regarded as higher-level  generalizations of  $p$-flat Cartan geometries  (resp., a dormant indigenous bundles).
By applying  a generalization of  Cartier's theorem (cf. Proposition \ref{Ertwq}),
we obtain the following   theorem as an extension of  Theorem \ref{TheoremA1}:

\SSP
%---------------------------------------------------------------------[begin theorem]----------------------
%\vspace{3mm}
\begin{intthm} [cf. Propositions \ref{P0201}, \ref{Efget}, \ref{Eqw209}, and Corollary \ref{Er281}] \label{TheoremA}
Suppose that $G$ is affine.
Then,
for each $N \in \mbZ_{>0}$, 
we obtain the following  commutative diagram in which  all the maps of sets  are bijective:
\begin{align}
\vcenter{\xymatrix@C=20pt@R=46pt{
&F^N\text{-}\mr{Ehr}_{(G, P), X/S} \ar[ld]^-{\sim}_-{\zeta_\N^{\heartsuit \Rightarrow \clubsuit}}  \ar[rd]_-{\sim}^-{\zeta_\N^{\heartsuit \Rightarrow \spadesuit}}&
\\
F^N\text{-}\mr{Car}_{(G, H). X/S} \ar[rr]^-{\sim}_-{\zeta_\N^{\clubsuit \Rightarrow \spadesuit}}&& \mr{Ind}^{^\mr{Zzz...}}_{(G, H), X/S}.
}}
\end{align}
Moreover, if  $S = \mr{Spec}(k)$ and $X$ is proper over $k$, then the above diagram is obtained even for  $\N = \infty$.
Finally, the formations of these bijections commute with pull-back to \'{e}tale  $X$-schemes, as well as with base-change to $S$-schemes.
 \end{intthm}
\SSP
%------------------------------------------------------------------------------[end theorem]-------------

Before proceeding,
we shall comment on the advantage of introducing various structures of infinite level, not only those of finite level.
As observed in ~\cite{Wak6}, 
the {\it stratified fundamental group} $\pi_1^\mr{str}(X)$ of $X$ can be effectively used in the study of 
 $F^\infty$-$(\mr{PGL}_{n+1}, \mbP^n)$-structures.
Here,  the  stratified fundamental group  of $X$ is
  the Tannakian fundamental group associated with  the category of
vector bundles equipped with an action of the ring of differential operators (in the sense of Grothendieck),
  or equivalently, the  category of $F$-divided sheaves.
 Each $F^\infty$-$(\mr{PGL}_{n+1}, \mbP^n)$-structure, or the corresponding dormant indigenous $(\mr{PGL}_{n+1}, \mr{PGL}_{n+1}^\circledcirc)$-bundle of level $\infty$,  specifies a $\mr{PGL}_{n+1}$-representation of
$\pi_1^\mr{str}(X)$ up to conjugation;
it may be thought of as the monodromy representation of this  structure, and it enables us to obtain a much better understanding of
related things from  the viewpoint of algebraic groups.
For example, by considering  $\pi_1^\mr{str}(X)$, we obtained some characterizations of $\mbP^n$ in terms of $F^\infty$-$(\mr{PGL}_{n+1}, \mbP^n)$-structures (cf. ~\cite[Theorem C, (ii)]{Wak6}).
Also,  it gives a complete description of the set of  $F^\infty$-$(\mr{PGL}_{n+1}, \mbP^n)$-structures on an Abelian variety  (cf. ~\cite[Theorem D, (ii)]{Wak6}).
For a general $(G, H)$, we will prove  some properties of a variety admitting an $F^\infty$-$(G, P)$-structure from the same viewpoint (cf. Propositions \ref{etoiuoa892} and \ref{etoiuoa893}).

After proving those properties,  we will focus on 
 Frobenius-Ehresmann structures on  smooth proper  curves in the context of subject (ii)   described in \S\,\ref{Sghhgotu1}.
According to a result by B. Laurent (cf. ~\cite[Theorem 1.1, (1)]{Lau}),  the list of homogenous curves, i.e., homogenous spaces of dimension $1$, may be described  as follows (where $e$ denotes the identity element of a group):
\begin{itemize}
\item[]
\begin{itemize}
\item[Case (a):] 
$G = \mr{PGL}_2$, $H= \mr{PGL}_2^\circledcirc$;
\item[Case (b):]
$G = \mr{Aff}_1$, $\mr{Aff}_1^\circledcirc$;
\item[Case (c):]
$G = \mbG_m$, $H  = \{ e \}$;
\item[Case (d):]
$G = \mbG_a$, $H = \{e \}$;
\item[Case (e):]
$G$ is an elliptic curve over $k$ and $H = \{e \}$.
\end{itemize}
\end{itemize}
Here, $\mr{Aff}_1$ denotes the group of affine transformations on the affine line $\mbA^1$ and $\mr{Aff}_1^\circledcirc$ denotes the subgroup of $\mr{Aff}_1$ consisting of affine transformations fixing the origin of $\mbA^1$.
Cases (a) and (b) were discussed in the previous studies (cf. ~\cite{Hos2}, ~\cite{Hos4}, ~\cite{Hos3}, ~\cite{Wak7}, and ~\cite{Wak6}).
By further examining the other cases, we will classify Frobenius-Ehresmann structures on curves  (cf. \S\,\ref{SSS0g1}).
At the same time, we will formulate a conjecture  (cf. Conjecture \ref{Efg09}) on computing  the explicit number of $F^\N$-$(\mr{PGL}_2, \mr{PGL}_2^\circledcirc)$-structures on a  general smooth proper curve; this is a generalization of 
Joshi's conjecture (cf. ~\cite[Conjecture 8.1]{Jo14}) for the rank $2$ case, which was proved in 
 ~\cite{Wak}.

\LSP
%%%%%%%%%%%%%%%%%%%%%%%%%%%%%%%---[begin section]---%%%%%%%%%%%%%%
\subsection{Deformation theory and the Ehresmann-Weil-Thurston principle} \label{S0eeghu1}

To consider subject (iii),  we  will develop the deformation  theory of (dormant) indigenous $(G, H)$-bundles.
The moduli stack of 
dormant indigenous $(\mr{PGL}_2, \mr{PGL}_2^\circledcirc)$-bundles on genus $g$ curves for  each integer $g>1$ (i.e., $\mfM_g^{^\mr{Zzz...}}\!$ mentioned in Example \ref{eso0w})
was substantially investigated by S. Mochizuki (cf.  ~\cite{Mzk2}) in the context of $p$-adic Teichm\"{u}ller theory.
In that work, it was proved that this moduli stack  is smooth over $k$ and generically \'{e}tale over the moduli stack $\mfM_{g}$ of genus $g$ curves (cf.  ~\cite[Chap.\,II, \S\,2.3, Theorem 2.8]{Mzk2}).
This fact is essential in the proof of Joshi's conjecture described in ~\cite{Wak}.
Indeed, we used generic \'{e}taleness  to lift  relevant  moduli stacks to characteristic $0$ and  then applied  a well-known formula for computing the Gromov-Witten invariant  of a certain type of Quot scheme.
Also,  a detailed understanding of
the moduli stack of dormant indigenous $(\mr{Aff}_1, \mr{Aff}_1^\circledcirc)$-bundles on curves was obtained  by investigating 
its local structure (cf. ~\cite{Wak7}).
As an application, we constructed  a family of  {\it pathological} algebraic surfaces (e.g., violating the Kodaira vanishing theorem)  in characteristic $p$ parametrized by  a higher dimensional base space (cf. ~\cite[Theorem C]{Wak7}).

In the present paper,  we study  the local structures of  various relevant  moduli spaces for a general $(G, H)$.
Part of the corresponding discussions in the holomorphic setting can be found in ~\cite{BD3}.
Here, we will briefly examine 
 one of the deformation functors 
  that we will  deal with.
  (See \S\,\ref{SSS03} for the other deformation functors and related results.)
Let $X$ be a smooth  scheme of finite type  over $k$  and 
 $\mcE^\sss := (\mcE_H, \nabla)$  a dormant indigenous $(G, H)$-bundle  on $X/k$.
Suppose that $G$ is affine, and consider the  functor 
  \begin{align}
  \mr{Def}^{^\mr{Zzz...}}_{(X, \mcE^\sss)}
  \end{align}
(cf. (\ref{esoriaa2}))  classifying deformations of the pair  $(X, \mcE^\sss)$ preserving  the {\it dormancy} condition.
A basic problem
regarding this  functor  is to determine its tangent and  obstruction spaces.
These spaces can  be computed 
by applying  an 
 explicit description of the hypercohomology  group of a certain complex associated with  $\mcE^\sss$.
In particular, we obtain a sufficient condition for the pro-representablity of $\mr{Def}^{^\mr{Zzz...}}_{(X, \mcE^\sss)}$.
The following Theorem \ref{TheoremB} states some  results on this functor (see Propositions \ref{P022}, \ref{P023}, and \ref{P051} for the results on the other deformation functors).

\SSP
%-----------------------[begin theorem]----------------------
%\vspace{3mm}
\begin{intthm}[cf. Propositions \ref{P050} and \ref{Etouaojm}]
\label{TheoremB}
Let us keep the above notation.
 \begin{itemize}
  \item[(i)]
Write  $\mfg_{\mcE_G}$  for the adjoint bundle   (cf. \S\,\ref{SS002})  of the $G$-bundle $\mcE_G := \mcE_H \times^H G$ and 
$\nabla^\mr{ad} : \mfg_{\mcE_G} \migi \Omega_{X/k} \otimes_{\mcO_X} \mfg_{\mcE_G}$ for  the $k$-connection on $\mfg_{\mcE_G}$   induced by $\nabla$.
 Then, the tangent space of 
 $\mr{Def}^{^\mr{Zzz...}}_{(X, \mcE^\sss)}$
 is isomorphic to $H^1 (X, \mr{Ker}(\nabla^\mr{ad}))$, and this functor has an obstruction theory by putting 
   $H^2 (X, \mr{Ker}(\nabla^\mr{ad}))$ as the obstruction space.
 \item[(ii)]
 Suppose that $X$ is proper over $k$.
 Then, 
 there exists a pro-representable hull for
 $\mr{Def}^{^\mr{Zzz...}}_{(X, \mcE^\sss)}$.
 \item[(iii)]
 Suppose  that $G$ is affine and 
 $\mcE^\sss$ is ordinary (cf. Definition \ref{ereao283}).
 Also, denote by $\mr{Def}_X$ (cf. (\ref{EWEW3})) the functor classifying deformations of $X$.
 Then, the morphism of functors 
 \begin{align}
 \mr{Def}^{^\mr{Zzz...}}_{(X, \mcE^\sss)} \migi \mr{Def}_X
 \end{align}
  given by forgetting the data of the deformations of $\mcE^\sss$ is a semi-isomorphism (cf. Definition \ref{pjapk0}).
 In particular,  if $X_R$ is a deformation of $X$ over
 a local Artinian $k$-algebra $R$ with residue field $k$ whose maximal ideal $\mfm_R$  is square-nilpotent,
then  there exists a unique (up to isomorphism) deformation of $\mcE^\sss$ over $X_R$ preserving the dormancy condition.
 \end{itemize}
 \end{intthm}
%\vspace{5mm}
%------------------------------[end theorem]-------------
\SSP

Finally, we formulate and prove  the {\it Ehresmann-Weil-Thurston principle for Frobenius-Ehresmann structures}.
Let $X$ be as above.
Notice that due to the absence of the topological fundamental group of $X$,
we have to consider
 other  approaches to define $G$-local systems (i.e., points of the representation variety $\mr{Hom}(\pi_1 (X), G)/G$)
   in an algebraic manner.
One candidate  is 
the category of 
 crystals  of $G$-bundles on 
% either  the infinitesimal or 
 the crystalline site of finite level.
 Indeed,  the notion of a crystal allows us to describe deformations of flat bundles 
  in an intrinsic way with respect to  $X$, just as the notion of a topological local system on a complex manifold does via the Riemann-Hilbert correspondence.

 Let  us fix
  a positive integer $\N$ and 
   an $F^\N$-$(G, P)$-structure
  $\mcS^\hh$   on $X/k$.
   Denote by
  $\mcE^\dd$  the {\it monodromy crystal}  of $\mcS^\hh$ (cf. Definition \ref{Efopwd}).
To state the following Theorem \ref{TheoremD},
we shall write
\begin{align}
\mr{Def}_{(X, \mcS^\heartsuit)} 
\hspace{3mm} \text{and} \hspace{3mm} 
%\ \left(\text{resp.,} \ 
\mr{Def}^{^\mr{Zzz...}}_{\mcE^\diamondsuit} 
%\right)
% : \mfA \mfr \mft_{/k}^\mr{PD} \migi \mfS \mfe \mft
\end{align}
(cf. (\ref{W23}), (\ref{W24})) for the deformation functors of $(X, \mcS^\heartsuit)$ and  $\mcE^\diamondsuit$ respectively.
(The latter functor  $\mr{Def}^{^\mr{Zzz...}}_{\mcE^\diamondsuit}\!$ classifies deformations of $\mcE^\diamondsuit$ satisfying the dormancy condition.)
Then, we will prove the  following assertion,   regarded as a positive characteristic version of the Ehresmann-Weil-Thurston principle; this result says that a local deformation of an  $F^\N$-$(G, P)$-structure is uniquely determined by its monodromy crystal.

\SSP
%-----------------------[begin theorem]----------------------
%\vspace{3mm}
\begin{intthm} [cf.  Theorem \ref{Wroauom}] \label{TheoremD}
Let us keep the above notation.
 Then, the morphism of deformation functors
 \begin{align}
\mu_{\N}^{\heartsuit \Rightarrow \diamondsuit} :  \mr{Def}_{(X, \mcS^\heartsuit)} \migi \mr{Def}_{\mcE^\diamondsuit}^{^\mr{Zzz...}}
  \end{align}
given by taking monodromy  crystals of deformations of $\mcS^\hh$
is a semi-isomorphism (cf. Definition \ref{pjapk0}).
 If, moreover,  $\mr{Def}_{\mcE^\dd}^{^\mr{Zzz...}}\!$ is rigid (cf. Definition \ref{Eotuamg}), then $\mu_{\N}^{\heartsuit \Rightarrow \diamondsuit}$ is  an isomorphism.
 \end{intthm}
%\vspace{5mm}
%------------------------------[end theorem]-------------

\LSP
%%%%%%%%%%%%%%%%%%%%%%%%%%%%%%%---[begin section]---%%%%%%%%%%%%%%
\subsection{Organization of the present paper} \label{S0eeghu1}

%We here outline  the organization of the present paper after this Introduction.
In \S\,2,  the notion of  a Frobenius-Ehresmann structure  is defined at the beginning (cf. Definition \ref{D023}).
We provide examples of this structure by  referring  to the previous work on   Frobenius-Ehresmann structures modeled on $(\mr{PGL}_{n+1}, \mr{PGL}_{n+1}^\circledcirc)$ and $(\mr{Aff}_n, \mr{Aff}_n^\circledcirc)$ (cf. Examples \ref{EX01} and \ref{Efgh8}).
In  \S\,3,  we start by  recalling  the definition of the Maurer-Cantan form and then  
discuss Cartan geometries  in positive characteristic, including the case of higher level (cf. Definitions \ref{D024} and  \ref{Eotarmr90}). 
A key ingredient of this study is the $p$-curvature of  a vector-valued $1$-form (cf.    Definition \ref{D00102}).
This invariant  leads us to define  $p$-flat  Cartan geometries (cf. Definition \ref{ea840q}).
The main parts of that section are the constructions of  the bijections between Frobenius-Ehresmann structures and certain  Cartan geometries asserted in Theorems \ref{TheoremA1} and \ref{TheoremA}  (cf.  
 Theorem \ref{P0011}, Propositions \ref{P0201} and \ref{Efget}).
In \S\,\ref{S0029}, we introduce   (dormant) indigenous $(G, H)$-bundles on varieties  of arbitrary dimension (cf. Definitions \ref{D0013}, \ref{Ets395}, \ref{Eosi38392}, and \ref{Ge990}).
The proofs of Theorems \ref{TheoremA1} and \ref{TheoremA} are completed   by constructing bijections between Cartan geometries and dormant indigenous bundles (cf. 
 Propositions \ref{P0234},  \ref{Eqw209}, and Corollary \ref{Er281}).
After that, we verify several properties  of varieties admitting an $F^\infty$-$(G, P)$-structure (cf. Propositions \ref{etoiuoa892} and \ref{etoiuoa893}).
The discussion at the end of the section  deals with $F^\N$-$(G, P)$-structures when $G$ is a product of $\mbG_m$'s and $\mbG_a$'s.
As a result,  it is shown  (cf. Proposition \ref{Erfgy}) that   the  $p$-rank and $a$-number  determine whether or not  such a Frobenius-Ehresmann structure exists on  a given  Abelian variety.
In \S\,\ref{SSS0g1}, we consider the classification problem of Frobenius-Ehresmann structures on  smooth proper curves while referring to the previous studies on $(G, H) = (\mr{PGL}_{2}, \mr{PGL}_{2}^\circledcirc), (\mr{Aff}_1, \mr{Aff}_1^\circledcirc)$. 
The content of \S\,\ref{SSS03} is  the deformation theory of (dormant) indigenous $(G, H)$-bundles.
First, we  introduce  several basic definitions in the general theory of deformation functors.
Then, we compute the tangent and obstruction spaces of various  deformation functors by  means of the hypercohomology  groups of certain complexes (cf. Proposition \ref{P022}, \ref{P023}, \ref{P050}, and \ref{P051}).
In particular, we   finish the proof of Theorem \ref{TheoremB}.
The Ehresmann-Weil-Thurston principle for Frobenius-Ehresmann structures, i.e., Theorem \ref{TheoremD}, is then  formulated  in terms of crystals and, at the same time,  proved by applying  the results obtained so far  (cf. Theorem \ref{Wroauom}).
In the middle of the section,
  we define the ordinariness of dormant indigenous $(G, H)$-bundles (cf. Definition \ref{ereao283}).
However, the discussion and results therein  will not be used after that.
In the Appendix, we briefly recall  some of the basic definitions concerning connections  on a principal bundle.
Moreover,  we develop the theory of higher-level PD stratifications on a principal bundle.
This includes  a generalization of  Cartier's theorem  (cf. Proposition \ref{Ertwq}).
Throughout the present paper, we will occasionally refer to definitions and facts discussed in the Appendix.

\LSP
%----------------------------------------------------------------------[begin subsection]-------------
%\subsection*{Notation and Convensions}
\subsection{Notation and Conventions} \label{S9025}
%\leavevmode\\ \vspace{-4mm}

Throughout the present paper, we fix
 an odd prime $p$
   and an algebraically closed field $k$ of characteristic $p$.

All schemes appearing in the present paper are assumed to be locally noetherian.
Given  a sheaf $\mcV$ on a scheme, 
 we use the notation ``\,$v \in \mcV$\,''
   for a local section $v$ of  $\mcV$.
If $S$ is a scheme and $X$ is a  scheme over $S$, then we denote by $\Omega_{X/S}$ (resp., $\mcT_{X/S}$) the sheaf of $1$-forms (resp., the sheaf of vector fields) on $X$ relative to $S$.
In \S\S\,\ref{S001}-\ref{SSS03},
 we will always be working  over $k$ unless stated otherwise; for example,  products  of $k$-schemes will be taken over $k$, i.e., $X_1 \times X_2 := X_1 \times_k X_2$.
Given a $k$-scheme  $X$,    we write  $\Omega_X := \Omega_{X/k}$ and  $\mcT_X := \mcT_{X/k}$ for simplicity.
Also, for  a $k$-vector space $M$ and an $\mcO_X$-module $\mcF$, we write $\mcF  \otimes M := \mcF  \otimes_k M$.

Let $S$ be a $k$-scheme,  $X$ a smooth $S$-scheme, $M$ a $k$-vector space, and  $q$ a nonnegative integer.
Then, 
an {\it $M$-valued $q$-form on  $X/S$}
 is  defined as   a global  section of the sheaf  $(\bigwedge^q \Omega_{X/S}) \otimes M$, where $\bigwedge^q \Omega_{X/S}$ denotes the $q$-th exterior power of $\Omega_{X/S}$ over $\mcO_X$.
For an $M$-valued $q$-form $\omega$ on $X/S$, we shall write
$\omega^\triangleright$ for the  $\mcO_X$-linear morphism $\bigwedge^q\mcT_{X/S} \migi \mcO_X \otimes M$ corresponding to $\omega$.

Given an algebraic group, we  will often  use the notation $e$ to denote   its identity element.
In \S\S\,\ref{S001}-\ref{SSS03}, we fix 
 a connected smooth  algebraic group $G$ over $k$ and a closed subgroup $H$ of $G$ such that 
   the homogenous space $P$  defined as the quotient $G/H$  is smooth over $k$.
 We will write 
 $\mfg \left(:= \mr{Lie}(G) \right)$ for   the Lie algebra of $G$, i.e., the tangent space $\mcT_{G}|_e \left(=e^*(\mcT_{G}) \right)$  at the identity element $e$ and also write $\mfh := \mr{Lie}(H)$, which defines a Lie subalgebra of $\mfg$.
Also, we write 
 $\pi_G$ for  the natural projection $G \migisurj P$.

Denote by $\mbG_m$  the multiplicative group over $k$ and by  $\mbG_a$ the additive group over $k$.
Also, denote by $\mu_p$ the subgroup of $\mbG_m$ consisting of $p$-th roots of unity and
by $\alpha_p$ the subgroup  of $\mbG_a$ consisting of elements $a$ with $a^p =0$.

Let $X$ be a scheme and $G'$ an algebraic group.
By a {\it $G'$-bundle} on $X$, we mean 
a principal $G'$-bundle (equipped with a right $G'$-action) on $X$ in the \'{e}tale topology.
 Let $g : G' \migi G''$ be a morphism of algebraic groups.
Then,  given a principal $G'$-bundle $\mcE'$ on $X$, we shall  write $\mcE' \times^{G'} G''$  for the  $G''$-bundle on $X$ obtained from $\mcE'$ via  change of structure group by $g$.
Also,  for  a $G''$-bundle $\mcE''$ on $X$, a {\it $G'$-reduction} of $\mcE''$
is defined as a $G'$-bundle $\mcE'$ on $X$
 together with an isomorphism $\mcE' \times^{G'} G'' \isom \mcE''$ of  $G''$-bundles.

Let $X$ and $G'$  be as above  and let
$h : X \migi G'$ be   a morphism of schemes.
Then, we write 
$\mr{R}_h$ (resp., $\mr{L}_h$)
  for the right-translation  (resp., the left-translation) by $h$, i.e., the automorphism of $X \times G$ given by $(x, a) \mapsto (x, a \cdot h)$ (resp., $(x, a) \mapsto (x, h \cdot a)$) for any $x \in X$ and  $a \in G$.
If  $\mcE$ is a $G$-bundle on $X$ and  $h$ is a point of $G$, then 
we  use the notation $\mr{R}_h$ to denote the automorphism of $\mcE$ induced from  the action by $h$.

\LSP
%%%%%%%%%%%%%%%%%%%%%%%%%%%%%%%---[begin section]---%%%%%%%%%%%%%%
\subsection*{Acknowledgements} 
 
 We would like to thank Professor Sorin Dumitrescu for his helpful comments concerning holomorphic Cartan geometries (cf. Remark \ref{ERqqw3}).
Also, we are grateful for the many constructive conversations we had with algebraic varieties having a Frobenius-Ehresmann structure, who live in the world of mathematics!\footnote{My conversations with mathematical entities usually go as follows: I'm walking around and see a black cat sleeping on a window sill. I say to ``Hey Math, look! A black cat!" Math replies ``No, no, there is at least one window sill with at least one sleeping cat on the sill such that the cat is black on at least one side of the cat during at least one period of time."}
Our work was partially supported by Grant-in-Aid for Scientific Research (KAKENHI No. 18K13385, 21K13770).

%%%%%%%---[begin section]---%%%%%%%%%%%%%%
\vspace{10mm} 
\section{Frobenius-Ehresmann structures} \label{S001} \vspace{3mm}

In this section, we introduce Frobenius-Ehresmann structures  of various levels,  which generalize the Frobenius-projective and Frobenius-affine structures discussed in ~\cite{Hos2}, ~\cite{Hos4}, and ~\cite{Wak6}.
Let us fix a $k$-scheme  $S$   and a smooth scheme $X$   over $S$.

\LSP
%----------------------------------------------------------------------[begin subsection]-------------
\subsection{Frobenius twists and relative Frobenius morphisms} \label{S1}
We shall  write $f : X \migi S$ for the structure morphism of $X/S$ and 
 $F_{X}$ (resp., $F_S$) for
the absolute Frobenius endomorphism of $X$ (resp., $S$).
For each positive integer $\N$, the  {\it $\N$-th Frobenius twist of $\Y$ over $S$} is, by definition,  the base-change $X^{(\N)} \left(:= \Y \times_{S, F_S^\N} S\right)$ of $\Y$ via
 the $\N$-th iterate $F_S^\N$ of $F_S$.
Denote by $f^{(\N)} : \Y^{(\N)} \migi S$ the structure morphism of 
$\Y^{(\N)}/S$, i.e., $f^{(\N)} := f \times \mr{id}_S$.
The {\it $\N$-th relative Frobenius morphism of $\Y$ over $S$} is  the unique morphism $F_{\Y/S}^{(\N)} : \Y \migi \Y^{(\N)}$ over $S$ that makes the following diagram commute:
\begin{align}
\vcenter{\xymatrix@C=36pt@R=36pt{
\Y \ar@/^10pt/[rrrrd]^{F_\Y^\N}\ar@/_10pt/[ddrr]_{f} \ar[rrd]_{  \ \ \ \  \  \   F_{\Y/S}^{(\N)}} & & &   &   \\
& & \Y^{(\N)}  \ar[rr]_{\mr{id}_\Y \times F_S^\N} \ar[d]^-{f^{(\N)}}  \ar@{}[rrd]|{\Box}  &  &  \Y \ar[d]^-{f} \\
&  & S \ar[rr]_{F_S^\N} & &  S.
}}
\end{align}
For simplicity, we shall write $F_{X/S} := F^{(1)}_{X/S}$.
Also, let $(X^{(0)}, F_X^{(0)}) := (X, \mr{id}_X)$ for convenience.

 Denote by 
$G_{\Y}$ the group-valued \'{e}tale sheaf on $\Y$ represented by $G$, i.e., $G_{\Y} := \mr{Hom}_k (-, G)$.
Moreover, given  a positive integer $\N$, we shall  write
\begin{align} \label{Eg60}
G_{X/S}^{(\N)} := (F^{(\N)}_{X/S})^{-1} (G_{\Y^{(\N)}}) \left(\subseteq G_{\Y} \right).
\end{align}
The sheaves $G_{X/S}^{(\N)}$ for various $\N$ constitute, via pull-back by  relative  Frobenius morphisms, 
 a sequence of inclusions
\begin{align} \label{E0078}
G_\Y \supseteq G_{X/S}^{(1)} \supseteq G_{X/S}^{(2)} \supseteq \cdots \supseteq G_{X/S}^{(\N)} \supseteq \cdots.
\end{align}

\LSP
%---------------------------[begin subsection]-------------
\subsection{$F^N$-$(G, P)$-structures  for $N < \infty$} \label{SS01}
We shall 
write 
\begin{align}\label{L014}
P_S := S \times P,
\end{align}
 and write $\mr{L}_{P_S}$ for the  left $G$-action 
$G \times P_S \migi P_S$ on the $S$-scheme $P_S$ induced by the natural left $G$-action  on $P$.
We shall set
\begin{align}
\mcP^{\text{\'{e}t}}_{P, X/S}
\end{align}
to be  the set-valued \'{e}tale  sheaf  on $X$  that assigns, to
each \'{e}tale scheme $U$ over  $X$, the set
of \'{e}tale $S$-morphisms from $U$ to $P_S$.
When there is no fear of confusion,  we will denote this sheaf  by $\mcP^{\text{\'{e}t}}$ for simplicity.
Now, let us fix a positive integer $\N$.

\SSP
%-----------------------------------------------------------------------------
 \ble
Let $U$ be an \'{e}tale  scheme over  $X$,
 $\phi$ 
 an $S$-morphism $U \migi \PP_S$, 
  and $h$ an element of $G_{X/S}^{(\N)}(U)  \left(= G (U^{(\N)}) \right)$.
Write $h^F := h \circ F_{U/S}^{(\N)} : U \migi G$ and define $h (\phi)$ to be 
  the composite
 \begin{align} \label{L01}
 h (\phi) : U \xrightarrow{(h^F, \phi)} G \times P_S  \xrightarrow{\mr{L}_{P_S}} P_S.
 \end{align}
Then,  $\phi \in \mcP^{\text{\'{e}t}} (U)$ implies 
$h (\phi) \in \mcP^{\text{\'{e}t}}(U)$.
\ele
%-----------------------------------------------------------------------------
\begin{proof}
Suppose that $\phi \in \mcP^{\text{\'{e}t}} (U)$, and 
let us choose  a point $q \in U (k)$.
The differential  of $(h^F, \phi): U \migi G \times P_S$ at $q$  defines a morphism  of $k$-vector spaces
\begin{align}
(dh^F, d\phi) : \mcT_{U/S} |_q  \migi \mcT_{G/k} |_{h^F (q)} \oplus \mcT_{P_S/S} |_{\phi (q)} \left(= \mcT_{(G \times P_S)/S} |_{(h^F (q), \phi (q))} \right).
\end{align}
Since $h^F$ factors through $F_{U/S}^{(\N)}:U \migi U^{(\N)}$,
the first component  $d h^F : \mcT_{U/S} |_q \migi \mcT_{G/k} |_{h^F (q)}$ of this morphism  becomes the zero map.
On the other hand, the \'{e}taleness assumption on $\phi$ implies the bijectivity of 
the second component  $d \phi : \mcT_{U/S} |_q \isom \mcT_{P_S/S} |_{\phi (q)}$.
Also, the differential $d \mr{L}_{P_S} : \mcT_{G/k} |_{h^F (q)} \oplus \mcT_{P_S/S} |_{\phi (q)}\migi \mcT_{P_S/S} |_{h(\phi)(q)}$ of $\mr{L}_{P_S}$
is restricted  to the identity morphism $\mcT_{P_S/S}|_{\phi (q)} \left(=0 \oplus \mcT_{P_S/S}|_{\phi (q)}\right) \isom \mcT_{P_S/S} |_{h(\phi)(q)}$ of $\mcT_{P_S/S}|_{\phi (q)}$.
It follows that the differential $d (h (\phi)) \left(=  d\mr{L}_{P_S} \circ (d h^F, d \phi)\right)$ of $h (\phi)$ at $q$ is  bijective.
Hence, since both $U$ and $P_S$ are smooth over $S$,  the morphism  $h (\phi)$ turns out to be  \'{e}tale.
This completes the proof of this lemma.
\end{proof}
%-----------------------------------------------------------------------------
\SSP

By the above lemma, the assignment $(h, \phi) \mapsto h (\phi)$  defines a left  $G_{X/S}^{(\N)}$-action
\begin{align}
\mr{L}_{\mcP}^{(\N)}: G_{X/S}^{(\N)} \times \mcP^{\text{\'{e}t}} \migi \mcP^{\text{\'{e}t}}
\end{align}
 on $\mcP^{\text{\'{e}t}}$.
The formation of this action for each $\N$  commutes  with the inclusions in (\ref{E0078}).

%-----------------------------------------------------------------------[begin definition]------------------
\SSP
\bde \label{D023}
 Let $\mcS^\hh$ be a subsheaf of $\mcP^{\text{\'{e}t}}$. 
We shall  say that
 $\mcS^\hh$ is 
a {\bf Frobenius-Ehresmann $(G, P)$-structure of level $\N$}, or simply  
 an {\bf $F^\N$-$(G, P)$-structure},  
  on $X/S$ if $\mcS^\hh$ is closed under the $G_{X/S}^{(\N)}$-action $\mr{L}_{\mcP}^{(\N)}$
    and forms a left $G_{X/S}^{(\N)}$-torsor on $X$ with respect to the resulting $G_{X/S}^{(\N)}$-action on $\mcS^\hh$.
   In the case of $\N=1$, it is called an $F$-$(G, P)$-structure.
 \ede
\SSP
%-----------------------------------------------------------------------[end definition]-------------------

We shall denote by
\begin{align} \label{wer4456}
F^N\text{-}\mr{Ehr}_{(G, P), X/S} 
\end{align}
the set of $F^N$-$(G, P)$-structures 
 on $X/S$.
 To simplify the notation slightly, we will write 
 $F\text{-}\mr{Ehr}_{(G, P), X/S}$ $ := F^1\text{-}\mr{Ehr}_{(G, P), X/S}$. 
Note that  $X/S$ admits  an $F^\N$-$(G, P)$-structure only when the relative dimension of $X/S$ coincides with the dimension  of $G$.

\LSP
%---------------------------[begin subsection]-------------
\subsection{$F^N$-$(G, P)$-structures for $N = \infty$} \label{SS0333}
Let $N$ and $N'$ be positive integers with $N \leq N'$ and $\mcS^\hh$ an $F^{N'}$-$(G, P)$-structure  on $X/S$.
Denote by 
\begin{align} \label{gougoa}
\mcS^\hh |^{\langle N \rangle}
\end{align}
the smallest subsheaf of $\mcP^{\text{\'{e}t}}$ which contains $\mcS^\hh$ and is closed under the $G_{X/S}^{(N)}$-action.
This sheaf can be verified to  form an $F^N$-$(G, P)$-structure  on $X/S$.
We refer to $\mcS^\hh |^{\langle N \rangle}$ as the {\bf $N$-th truncation} of $\mcS^\hh$.
The resulting assignments $\mcS^\hh \mapsto \mcS^\hh |^{\langle N \rangle}$ for various pairs $(\N, \N')$
  give a projective system of sets
\begin{align} \label{E009}
\cdots \migi F^N\text{-}\mr{Ehr}_{(G, P), X/S}
\migi \cdots 
 \migi 
 F^2\text{-}\mr{Ehr}_{(G, P), X/S}
\migi F^1\text{-}\mr{Ehr}_{(G, P), X/S}.
\end{align}

%-----------------------------------------------------------------------[begin definition]------------------
\SSP
\bde \label{D027}
 A {\bf Frobenius-Ehresmann $(G, P)$-structure of level $\infty$}, or simply
 an {\bf  $F^\infty$-$(G, P)$-structure},  on $X/S$ is a collection
 \begin{align}
 \mcS^\hh_\infty := \{ \mcS_N^\hh \}_{N \in \mbZ_{>0}},
 \end{align}
 where each $\mcS_N^\hh$ ($N \in \mbZ_{>0}$) is an  $F^N$-$(G, P)$-structure  on $X/S$ with $\mcS^\hh_{N+1} |^{\langle N \rangle} = \mcS_N^\hh$.
  \ede
\SSP
%-----------------------------------------------------------------------[end definition]-------------------

Denote by
\begin{align} \label{wer4457}
F^\infty\text{-}\mr{Ehr}_{(G, P), X/S} 
\end{align}
the set of $F^\infty$-$(G, P)$-structures  on $X/S$; this set 
 may be naturally identified with the limit of the projective system (\ref{E009}), i.e., 
\begin{align} \label{rosao3048}
F^\infty\text{-}\mr{Ehr}_{(G, P), X/S} = \varprojlim_{N \in \mbZ_{>0}}F^N\text{-}\mr{Ehr}_{(G, P), X/S}.
\end{align}

%-------------------------%-------------------------
\begin{exa}[Standard example]
For each positive integer $\N$, 
we shall construct a standard example of an $F^\N$-$(G, P)$-structure defined on the homogenous space $P$.
Let us set
\begin{align} \label{EE012}
\mcS^{\mr{triv}, \hh}_{N}
\end{align}
to be the subsheaf of $\mcP_{P, P_S/S}^{\text{\'{e}t}}$ which, to any \'{e}tale  scheme $U$ over  $P_S$, assigns the set
\begin{align}
\mcS^{\mr{triv}, \hh}_{N} (U) := \left\{h (\phi_U) \in \mcP_{P, P_S/S}^{\text{\'{e}t}} (U) \, \Big| \, h \in G_{X/S}^{(N)} (U)  \right\},
\end{align}
where $\phi_U$ denotes the structure morphism  $U \migi P_S$ of the $P_S$-scheme $U$.
That is to say, $\mcS^{\mr{triv}, \hh}_{N}$ is the smallest subsheaf of $\mcP_{P, P_S/S}^{\text{\'{e}t}}$ that is closed under the $G_{P_S/S}^{(N)}$-action and contains the global section determined by the identity morphism of $P_S$.
Then, $\mcS^{\mr{triv}, \hh}_{N}$ specifies an $F^N$-$(G, P)$-structure on $P_S/S$.
The formation of $\mcS^{\mr{triv}, \hh}_{N}$ commutes with truncation to lower levels, so the collection
\begin{align} \label{eptpss9}
\mcS^{\mr{triv}, \hh}_{\infty} := \{ \mcS^{\mr{triv}, \hh}_{N} \}_{N \in \mbZ_{>0}}
\end{align}
forms an  $F^\infty$-$(G, P)$-structure  on $P_S/S$.
\end{exa}
%-------------------------%-------------------------
\SSP

%------------------------------------------------------------------
\begin{exa} \label{EX01}
Recall the Frobenius-Ehresmann structures for the special  cases  discussed  in   ~\cite{Wak6}.
 Let $n$ be a positive integer and let  $\N \in \mbZ_{>0}\sqcup \{ \infty \}$.
 \begin{itemize}
\item[(i)]
Denote by  $\mbP^n$ the $n$-dimensional projective space over $k$.
In particular, the $k$-rational points of $\mbP^n$ correspond bijectively to the ratios $[a_0: a_1: \cdots: a_n]$ with $a_0, \cdots, a_n \in k$ and $(a_0, \cdots, a_n) \neq (0, \cdots, 0)$.
The rank $n+1$ projective linear group
$\mr{PGL}_{n+1}$ over $k$ 
may be identified with the automorphism group of $\mbP^n$ in such a way that
$\overline{A} [a_0:  \cdots : a_n] := [a_0: \cdots : a_n]\overline{{^t}A}$ for each $A \in \mr{GL}_{n+1}$, where $\overline{(-)}$ denotes the natural projection $\mr{GL}_{n+1} \migisurj \mr{PGL}_{n+1}$.
Denote  by 
\begin{align} \label{taae89}
\mr{PGL}_{n+1}^\circledcirc
\end{align}
  the subgroup of $\mr{PGL}_{n+1}$   consisting of automorphisms of $\mbP^n$
fixing $[1: 0: 0: \cdots : 0]$.
The homogenous space $\mr{PGL}_{n+1}/\mr{PGL}_{n+1}^\circledcirc$ is naturally isomorphic to $\mbP^n$.
Thus, we obtain  the notion of an $F^\N$-$(\mr{PGL}_{n+1}, \mbP^n)$-structure, which is equivalent to what we call an $F^\N$-projective structure (cf. ~\cite[\S\,1, Definition 1.2.1]{Wak6}) under the assumption that $X$ is a smooth quasi-projective $k$-variety   of dimension $n$ with  $p \nmid (n+1)$.
These notions  differ in   that, on the one hand 
an $F^\N$-$(\mr{PGL}_{n+1}, \mbP^n)$-structure   is  an {\it \'{e}tale}  sheaf of    locally defined  \'{e}tale morphisms to $\mbP^n$,
but on the other hand  an $F^\N$-projective structure  is   defined in {\it loc.\,cit.} as  a {\it Zariski} sheaf.
However,  the result of  ~\cite[Theorem A (and \S\,3, Corollary 3.4.4)]{Wak6} together with   Theorem \ref{TheoremA} in the present paper shows that these definitions turn out to be equivalent.
That is to say, under the assumption described above, the underlying  $\mr{PGL}_{n+1}$-torsor of any $F^\N$-$(\mr{PGL}_{n+1}, \mbP^n)$-structure   may be trivialized in the Zariski topology.

\item[(ii)]
Next, 
 denote by  $\mbA^n$ the $n$-dimensional  affine space over $k$.
 Also, denote by $\mr{Aff}_{n}$ the group of affine transformations on $\mbA^n$;
 it can be described concretely as the semidirect product  of  $\mbG_a^{\times n}$ (:= the group of translations on $\mbA^n$) by $\mr{GL}_n$ (:= the rank $n$ general linear group over $k$):
 \begin{align} \label{E03351}
 \mr{Aff}_n = \mbG_a^{\times n} \rtimes \mr{GL}_n,
 \end{align}
 where the $\mr{GL}_n$-action  on $\mbG_a^{\times n} \left(= \mbA^{n} \right)$ is taken to mean matrix multiplication of a vector.
  If 
  \begin{align} \label{rish34}
  \mr{Aff}_n^\circledcirc
  \end{align}
   denotes the subgroup of $\mr{Aff}_n$ consisting of affine transformations fixing the origin, then the homogenous space $\mr{Aff}_n/\mr{Aff}_n^\circledcirc$ is naturally isomorphic to $\mbA^n$.
 Thus, we obtain 
 the notion of an $F^\N$-$(\mr{Aff}_n, \mbA^n)$-structure.
 Just as in the above case,  this  is equivalent to what we call an {\it $F^\N$-affine structure} (cf. ~\cite[\S\,1, Definition 1.2.1,  and \S\,3, Corollary 3.4.4]{Wak}).

\end{itemize}
In ~\cite{Wak6},  we showed various fundamental  properties of a variety admitting either of these  structures   by establishing bijective correspondences with   other equivalent realizations  (cf. ~\cite[Theorem A]{Wak6}).
For example, we  proved formulas providing   necessary conditions on (crystalline) Chern classes for the existence of an $F^\N$-projective or $F^\N$-affine structure (cf. ~\cite[Theorem B]{Wak6}).
Moreover,  we investigated  the classification problems of curves,  surfaces, and Abelian varieties admitting 
such a structure
 (cf. ~\cite[Theorems E, F]{Wak6}).
\end{exa}
%------------------------------------------------------------------

\LSP
%---------------------------[begin subsection]-------------
\subsection{Pull-back and base-change of $F^\N$-$(G, P)$-structures} \label{SS066}
Consider the pull-back of a given Frobenius-Ehresmann $(G, P)$-structure by an \'{e}tale morphism between  the  underlying schemes.
Let $Y$ be a smooth scheme over  $S$ and $y : Y \migi X$ an \'{e}tale $S$-morphism.
Given each \'{e}tale $S$-morphism $\phi : U \migi P_S$ from an \'{e}tale $X$-scheme $U$  to $P_S$,
we obtain the \'{e}tale composite $\phi \circ y : y^{-1}(U) \left(:= Y \times_{X} U \right) \migi P_S$;
the assignment $\phi \mapsto \phi \circ y$ yields a morphism
\begin{align} \label{EE011}
y^{-1} (\mcP_{P, X/S}^{\text{\'{e}t}}) \migi \mcP_{P, Y/S}^{\text{\'{e}t}}
\end{align}
of \'{e}tale sheaves on $Y$.
For an $F^N$-$(G, P)$-structure  $\mcS^\hh$ on $X/S$ with  $N \in \mbZ_{>0}$,
we set
\begin{align}
y^*(\mcS^\hh)
\end{align}
to be the subsheaf of $\mcP_{P, Y/S}^{\text{\'{e}t}}$ defined to be the smallest subsheaf that is closed under the $G_{Y/S}^{(N)}$-action and contains the image of $y^{-1}(\mcS^\hh)$ via (\ref{EE011}).
One may verify  that $y^*(\mcS^\hh)$ specifies  an $F^\N$-$(G, P)$-structure on $Y/S$, 
 and the formation of $y^*(\mcS^\hh)$ commutes with truncation to lower levels.
This implies that, for each  $F^\infty$-$(G, P)$-structure  $\mcS^\hh := \{ \mcS^\hh_\N \}_{\N \in \mbZ_{>0}}$, 
the collection $y^*(\mcS^\hh) := \{ y^*(\mcS_\N^\hh) \}_{\N}$  forms an   $F^\infty$-$(G, P)$-structure on $Y/S$.
In each case of $N \in \mbZ_{>0}\sqcup \{ \infty \}$, $y^*(\mcS^\hh)$ is called the {\bf pull-back} of $\mcS^\hh$ via $y$, and the assignment  $\mcS^\hh \mapsto y^*(\mcS^\hh)$ defines a map of sets
\begin{align} \label{L0981}
F^\N\text{-}\mr{Ehr}_{(G, P), X/S} \migi F^\N\text{-}\mr{Ehr}_{(G, P), Y/S}.
\end{align}

Next,  consider the base-change of a Frobenius-Ehresmann $(G, P)$-structure.
Let $s: S' \migi S$ be a morphism of  $k$-schemes.
Given each \'{e}tale $S$-morphism $\phi : U \migi P_S$ from an \'{e}tale $X$-scheme $U$ to $P_S$, we obtain the base-change $s^*(\phi):  S' \times_S U \migi P_{S'} \left(=S' \times_S P_S \right)$, which is \'{e}tale;
the assignment $\phi \mapsto s^*(\phi)$ yields a morphism
\begin{align} \label{Ee211}
s_X^{-1}(\mcP_{P, X/S}^{\text{\'{e}t}}) \migi \mcP_{P, X'/S'}^{\text{\'{e}t}}
\end{align}
of \'{e}tale sheaves on $X' := S' \times_S X$, where $s_X$ denotes the natural projection $X' \migi X$.
For  an $F^\N$-$(G, P)$-structure $\mcS^\hh$ on $X/S$ with $\N \in \mbZ_{>0}$,
 we shall set
\begin{align}
s^*(\mcS^\hh)
\end{align}
to be 
 the smallest subsheaf of $\mcP^{\text{\'{e}t}}_{P, X'/S'}$  that is closed under the $G_{X'/S'}^{(\N)}$-action and contains the image of $s_X^{-1}(\mcS^\hh)$ via  (\ref{Ee211}).
One may verify that $s^*(\mcS^\hh)$ specifies  an  $F^\N$-$(G, P)$-structure on $X'/S'$, and 
the formation of $s^*(\mcS^\hh)$ commutes with truncation to lower levels.
This implies that, for  each  $F^\infty$-$(G, P)$-structure $\mcS^\hh := \{ \mcS_\N^\hh \}_{\N \in \mbZ_{>0}}$, the collection $s^*(\mcS^\hh) := \{ s^*(\mcS^\hh_\N) \}_{\N}$ forms an $F^\infty$-$(G, P)$-structure on $X'/S'$.
In each case of $\N \in \mbZ_{>0}\sqcup \{ \infty \}$, 
 $s^*(\mcS^\hh)$ is called the {\bf base-change} of $\mcS^\hh$, and the assignment  $\mcS^\hh \mapsto s^*(\mcS^\hh)$ defines a map of sets
\begin{align} \label{L010}
F^\N\text{-}\mr{Ehr}_{(G, P), X/S} \migi F^\N\text{-}\mr{Ehr}_{(G, P), X'/S'}.
\end{align}

\LSP
%---------------------------[begin subsection]-------------
\subsection{Change of structure group} \label{SS066g}
 Let $G'$ be another connected  smooth algebraic group over $k$ and $H'$ a   closed subgroup of $G'$ such that the quotient $P' := G'/H'$ is smooth over $k$.
Moreover, suppose that we are given a morphism of $k$-groups $g : G \migi G'$ such that $g (H) \subseteq H'$ and the $k$-morphism $\overline{g} : P \migi P'$ induced by $g$ is \'{e}tale.

 Now, let $\mcS^\hh$ be an $F^\N$-$(G, P)$-structure on $X/S$ with $\N \in \mbZ_{>0}$.
 Denote by 
 \begin{align} \label{L002}
 \mcS^\hh \times^P P'
 \end{align}
  the smallest subsheaf of $\mcP_{P', X/S}^{\text{\'{e}t}}$ that  is closed under the $G'^{(\N)}_{X/S}$-action on $\mcP_{P', X/S}^{\text{\'{e}t}}$  and contains the sections $\overline{g} \circ \phi$ for  various  local sections $\phi \in \mcS^\hh$.
 One  may verify that $\mcS^\hh \times^P P'$ defines  an $F^\N$-$(G', H')$-structure on $X/S$, and 
 the formation of  $\mcS^\hh \times^P P'$ commutes  with truncation to lower levels.
 This implies that,  for each  $F^\infty$-$(G, P)$-structure $\mcS^\hh := \{ \mcS^\hh_\N \}_{\N}$ on $X/S$,
 the collection $\mcS^\hh \times^P P' := \{ \mcS^\hh_N \times^P P' \}_{\N}$  forms
  an $F^\infty$-$(G', P')$-structure  on $X/S$.
 In each case of $\N \in \mbZ_{>0} \sqcup \{ \infty \}$, 
 the assignment $\mcS^\hh \mapsto \mcS^\hh \times^{P} P'$ defines a map of sets
  \begin{align} \label{L009}
  F^\N \text{-}\mr{Ehr}_{(G, P), X/S} \migi  F^\N \text{-}\mr{Ehr}_{(G', P'), X/S}.
  \end{align}
%  The formation of this map commutes with truncation to lower levels.
  
  \SSP
  %--------------------------%--------------------------
\begin{exa} \label{Efgh8}
We shall identify $\mbA^n$ with the open subscheme of $\mbP^n$ consisting of  the elements $\{ [1: a_1 : \cdots : a_n] \, | \, a_1, \cdots, a_n \in k \}$.
Consider  the morphism of $k$-groups $\mr{Aff}_n \migi \mr{PGL}_{n+1}$ given by $({\bf a}, A) \mapsto \overline{\begin{pmatrix} 1 & {\bf 0} \\ {^t}{\bf a} & A\end{pmatrix}}$ for each ${\bf a} \in \mbG_a^{\times n}$ and $A \in \mr{GL}_n$.
By this morphism, 
 $\mr{Aff}_n$ can be regarded as the subgroup of $\mr{PGL}_{n+1}$ consisting of automorphisms $h : \mbP^n \isom  \mbP^n$ with $h (\mbA^n) = \mbA^n$. 
 The inclusion $\mr{Aff}_n \migiincl \mr{PGL}_{n+1}$ induces a morphism  $\mr{Aff}_n / \mr{Aff}_n^\circledcirc \migi \mr{PGL}_{n+1}/\mr{PGL}_{n+1}^\circledcirc$, which   becomes the open immersion $\mbA^n \migiincl \mbP^n$ under the natural identifications $\mr{Aff}_n / \mr{Aff}_n^\circledcirc = \mbA^n$,  $\mr{PGL}_{n+1}/\mr{PGL}_{n+1}^\circledcirc = \mbP^n$.
 In particular,  for each $\N \in \mbZ_{>0}\sqcup \{ \infty \}$, we obtain a  map of sets
 \begin{align} \label{Fggth}
 F^\N\text{-}\mr{Ehr}_{(\mr{Aff}_n, \mbA^n), X/S} \migi  F^\N\text{-}\mr{Ehr}_{(\mr{PGL}_{n+1}, \mbP^n), X/S}.
 \end{align}
 This map is, in general,  neither surjective nor injective.
 In fact,   we  know that  it is not surjective (resp., not  injective)  when  $\N = \infty$ and $X/k$ is the projective space $\mbP^n$ (resp.,  an ordinary Abelian variety) (cf. ~\cite[Theorems C and  D]{Wak6}).
 \end{exa}
%--------------------------%--------------------------

 \SSP
  %--------------------------%--------------------------
\begin{exa} \label{es04xm}
Next,  let us consider the morphism of algebraic $k$-groups
\begin{align} \label{eos390}
\rho : \mbG_m^{\times n} \migi \mr{Aff}_n; (u_1, u_2, \cdots, u_n) \mapsto \overline{\begin{pmatrix} 1 & 0& 0 &  \cdots & 0 \\ 1-u_1 &u_1& 0& \cdots & 0 \\ 1-u_2 & 0 & u_2& \cdots &  0 \\ \vdots & \vdots & \vdots & \ddots & \vdots \\ 1-u_n & 0 & 0 & \cdots & u_n   \end{pmatrix}}.
\end{align}
Under the natural  identification $\mr{Aff}_n/\mr{Aff}_n^\circledcirc = \mbA^n$,
the $k$-morphism $\mbG_m^{\times n} \left(= \mbG_m^{\times n}/\{ e\} \right) \migi \mr{Aff}_n /\mr{Aff}_n^\circledcirc$ induced by $\rho$ coincides  with
the open immersion $\mbG_m^{\times n} \migiincl \left( (\mbA^1)^{\times n} =  \right) \mbA^n$ given by $(u_1, u_2, \cdots, u_n)$ $\mapsto (1-u_1, 1-u_2, \cdots, 1-u_n)$.
 In particular,  for each $\N \in \mbZ_{>0}\sqcup \{ \infty \}$,  we obtain the map of sets
 \begin{align} \label{Fgjjh}
 F^\N\text{-}\mr{Ehr}_{(\mbG_m^{\times n}, \mbG_m^{\times n}), X/S} \migi F^\N\text{-}\mr{Ehr}_{(\mr{Aff}_n, \mbA^n), X/S}.
 \end{align}
 
 Now, denote by $\mfS_n$ the symmetric group of $n$ letters.
 Each  element $\sigma \in \mfS_n$ determines an automorphism $\sigma^*$ of  $\mbG_m^{\times n}$ given by permuting the components according to $\sigma$; it induces a bijective endomap of $F^\N\text{-}\mr{Ehr}_{(\mbG_m^{\times n}, \mbG_m^{\times n}), X/S}$.
 In this way, we obtain an action of $\mfS_n$ on this set, and hence obtain the quotient set $F^\N\text{-}\mr{Ehr}_{(\mbG_m^{\times n}, \mbG_m^{\times n}), X/S}/\mfS_n$.
 It is verified  that (\ref{Fgjjh}) factors through the natural surjection 
 $F^\N\text{-}\mr{Ehr}_{(\mbG_m^{\times n}, \mbG_m^{\times n}), X/S} \migisurj F^\N\text{-}\mr{Ehr}_{(\mbG_m^{\times n}, \mbG_m^{\times n}), X/S}/\mfS_n$.
 Thus, (\ref{Fgjjh}) yields a map
  \begin{align} \label{Fgjjh8}
 F^\N\text{-}\mr{Ehr}_{(\mbG_m^{\times n}, \mbG_m^{\times n}), X/S}/\mfS_n \migi F^\N\text{-}\mr{Ehr}_{(\mr{Aff}_n, \mbA^n), X/S}.
 \end{align}
According to ~\cite[Theorem D, (ii)]{Wak6}, this map is bijective when  $\N = \infty$ and  $X$ is an ordinary Abelian variety over $\mr{Spec}(k) = S$ (cf. Proposition \ref{C0558}, (i)).
\end{exa}
%--------------------------%--------------------------

%%%%%%%---[begin section]---%%%%%%%%%%%%%%
\vspace{10mm}
\section{Cartan geometries in positive characteristic} \label{EosGG} \vspace{3mm}

In this section, we study   Cartan geometries in positive characteristic, including the case of higher level.
The main results of this section provide  bijective correspondences between Frobenius-Ehresmann structures and certain Cartan geometries (cf. Theorem \ref{P0011} and Proposition \ref{P0201}).
Let $S$ be a $k$-scheme and  $f: X \migi S$  a smooth scheme  over $S$.

\LSP
\subsection{The Maurer-Cartan form}
 In what follows, we  occasionally  identify the Lie algebra $\mfg$  with the space $\Gamma (G, \mcT_{G})^G$ of  vector fields on $G$    invariant under the left-translations $\mr{L}_{(-)}$.
That is to say, given any element $v \in \mfg = \mcT_{G} |_e$, we can use left-translation to produce a  global vector field $v^\dagger$ on $G$ with $v^\dagger |_h = d \mr{L}_h |_e (v)$ for every $h \in G$, where $d \mr{L}_h |_e$ denotes the $k$-linear morphism $\mcT_{G}|_e \left(=e^*(\mcT_{G}) \right)\migi \mcT_{G}|_h \left(= h^*(\mcT_{G}) \right)$  defined as the differential  of $\mr{L}_h$ at $e$.
Then, there exists 
 an $\mcO_G$-linear isomorphism
\begin{align} \label{E023}
\omega_{G}^\triangleright : \mcT_{G} \isom \mcO_G \otimes \mfg
\end{align}
determined uniquely 
by $\omega_{G}^\triangleright (v^\dagger) = 1 \otimes v$ for any $v \in \mfg$.
This morphism corresponds to a $\mfg$-valued $1$-form 
\begin{align}
\omega_G \in \Gamma (G, \Omega_{G} \otimes \mfg) \left(\cong \mr{Hom}_{\mcO_G}(\mcT_{G}, \mcO_G \otimes \mfg) \right),
\end{align}
which is called the {\bf  (left-invariant) Maurer-Cartan form} on $G$ (cf. ~\cite[Chap.\,3, \S\,1, Definition 1.3]{Sha}).
The  assertions in the following proposition  are well-known in the differential (and holomorphic) geometric setting.

\SSP
%------------------------------------------------------------
\bpr \label{L017}
\begin{itemize}
\item[(i)]
For any  $h \in G(k)$,
the following square diagrams are commutative:
 \begin{align} \label{E0057}
\vcenter{\xymatrix@C=46pt@R=36pt{
\mcT_{G}\ar[r]_-{\sim}^-{\omega^\triangleright_G} \ar[d]^{\wr}_{d \mr{L}_h}&  \mcO_{G} \otimes \mfg\ar[d]_{\wr}^{\mr{id}} \\
\mr{L}_h^*(\mcT_{G})\ar[r]^{\sim}_-{\mr{L}_h^*(\omega_G^\triangleright)}& \mcO_{G} \otimes \mfg,
 }}
 \hspace{10mm}
 \vcenter{\xymatrix@C=46pt@R=36pt{
\mcT_{G}\ar[r]_-{\sim}^-{\omega_G^\triangleright} \ar[d]^{\wr}_{d \mr{R}_h}&  \mcO_{G} \otimes \mfg\ar[d]_{\wr}^{\mr{id} \otimes \mr{Ad}_G(h^{-1})} \\
\mr{R}_h^*(\mcT_{G})\ar[r]^-{\sim}_-{\mr{R}_h^*(\omega^\triangleright_G)}& \mcO_{G} \otimes \mfg
 }}
\end{align}
(cf. (\ref{adjrep}) for the definition of $\mr{Ad}_G$).
\item[(ii)]
Let   $G'$ be  another  connected  smooth algebraic group over $k$ with Lie algebra $\mfg'$ and $g : G \migi G'$ a morphism of $k$-groups.
Then,  the Maurer-Cartan form $\omega_{G'}$ on  $G'$ makes the following square diagram commute:
 \begin{align} \label{E00234}
\vcenter{\xymatrix@C=46pt@R=36pt{
\mcT_{G}\ar[r]_-{\sim}^-{\omega_G^\triangleright} \ar[d]_{d g} & \mcO_G \otimes \mfg \ar[d]^{\mr{id} \otimes d g |_e} \\
g^*(\mcT_{G'}) \ar[r]^-{\sim}_-{g^*(\omega_{G'}^\triangleright)}& \mcO_G \otimes \mfg'.
 }}
\end{align}
\end{itemize}
\epr
%------------------------------------------------------------
\begin{proof}
Let $v$ and $a$ be arbitrary  elements of $\mfg$ and  $G (k)$ respectively.
Observe that
\begin{align} \label{eao3a930z0z}
(\mr{L}_h^*(\omega^\triangleright_G)\circ d \mr{L}_h) (v^\dagger |_a) = & \ 
(\mr{L}_h^*(\omega^\triangleright_G)\circ d \mr{L}_h) (d \mr{L}_a |_e (v)) \\
= & \  \omega^\triangleright_G (d \mr{L}_{h \cdot a} |_e (v)) \notag \\
= & \ \omega_G^\triangleright (v^\dagger |_{h \cdot a}) \notag \\
= & \  1 \otimes v \notag \\
= & \ \omega^\triangleright_G (v^\dagger |_a). \notag
\end{align}
Since $\mcT_G$ is generated by  the sections in $\Gamma (G, \mcT_G)$, 
(\ref{eao3a930z0z}) implies 
 the  commutativity of the left-hand square of (\ref{E0057}).
Moreover, for the same reason as above, the commutativity of the right-hand square of (\ref{E0057}) can be verified by the following sequence of equalities:
\begin{align}
(\mr{R}_h^*(\omega_G^\triangleright) \circ  d\mr{R}_h) (v^\dagger |_a) = & \ 
(\mr{R}_h^*(\omega_G^\triangleright) \circ  d\mr{R}_h) (d \mr{L}_{a} |_e (v)) \\
= & \ \omega_G^\triangleright ((d \mr{R}_h |_a \circ d \mr{L}_{h \cdot a} |_{h^{-1}} \circ d \mr{L}_{h^{-1}}|_e )(v))\notag \\
= & \ \omega_G^\triangleright ((d \mr{L}_{h \cdot a} |_e \circ d \mr{R}_{h} |_{h^{-1}} \circ d \mr{L}_{h^{-1}}|_e )(v)) \notag \\
= & \ \omega_G^\triangleright ((d \mr{R}_h |_{h^{-1}} \circ d \mr{L}_{h^{-1}}|_e (v))^\dagger |_{h \cdot a}) \notag  \\
= & \ \omega^\triangleright_G (\mr{Ad}_G (h^{-1})(v)^\dagger |_{h \cdot a}) \notag \\
= & \  1 \otimes \mr{Ad}_G (h^{-1})(v) \notag  \\
= & \ (\mr{id} \otimes \mr{Ad}_G (h^{-1}))(1 \otimes v) \notag \\
= & \ ((\mr{id} \otimes \mr{Ad}_G (h^{-1}))\circ \omega_G^\triangleright) (v^\dagger  |_a). \notag
\end{align}
Finally, we have 
\begin{align}
(g^*(\omega_{G'}^\triangleright) \circ dg) (v^\dagger |_a) = & \  (g^*(\omega_{G'}^\triangleright) \circ dg) (d \mr{L}_a |_e (v)) \\
= & \ \omega_{G'}^\triangleright ((d g |_a \circ d \mr{L}_a |_e )(v)) \notag \\
= & \ \omega_{G'}^\triangleright((d\mr{L}_{g (a)} |_e  \circ dg |_e) (v)) \notag \\
= & \ \omega_{G'}^\triangleright (dg|_e (v)^\dagger |_{g (a)}) \notag  \\
= &\ (1 \otimes dg|_e (v))   \notag \\
= & \ (\mr{id} \otimes dg |_e) (1 \otimes v)  \notag \\
= & \ ((\mr{id} \otimes dg|_e) \circ  \omega^\triangleright_G)(v^\dagger |_a). \notag
\end{align}
This proves the commutativity of (\ref{E00234}), and hence, completes 
the proof of this proposition.
\end{proof}
\SSP

 Next, let 
  $\mcE$  be a $G$-bundle on $X$.
The identification  $\mfg = \Gamma (G, \mcT_{G})^G$ gives  an $\mcO_X$-linear  isomorphism
\begin{align} \label{L0116}
\omega_{\mcE/X}^\triangleright : \mcT_{\mcE/X} \isom \mcO_{\mcE} \otimes \mfg.
\end{align}
 Indeed, if $\mcE = X \times G$ (i.e., the trivial $G$-bundle), then 
 we can  define $\omega_{X \times G/T}^\triangleright$  as the pull-back of $\omega^\triangleright_G$ via the second projection $\mr{pr}_2 : X \times G \migi G$, i.e., 
  the  composite isomorphism
 \begin{align} \label{L200}
\omega_{X \times G/T}^\triangleright : \mcT_{X \times G/X} \isom  \mr{pr}_2^*(\mcT_{G}) \xrightarrow{\mr{pr}_2^*(\omega_G^\triangleright)} \mr{pr}_2^*(\mcO_G \otimes \mfg) \isom \mcO_{X \times G} \otimes  \mfg,
 \end{align}
  where  the first arrow is obtained by differentiating  $\mr{pr}_2$.
 Moreover, suppose that we are given an \'{e}tale scheme $U$ over  $X$ and an automorphism 
   of  the trivial $G$-bundle  $U \times G$; 
  it    may be expressed as 
 the left-translation  by some $U$-rational point  of $G$.
The result of   Proposition \ref{L017}, (i), implies that 
 $\omega_{U \times G/G}^\triangleright$ is invariant under  pull-back via this automorphism.
 Thus,  the isomorphisms $\omega_{U \times G/G}^\triangleright$ defined  for various  $U$'s equipped  with a  local trivialization  $\mcE \cong U \times G$ of $\mcE$ may be glued together to form  a well-defined  isomorphism $\omega_{\mcE/X}^\triangleright : \mcT_{\mcE/X} \isom \mcO_{\mcE} \otimes \mfg$, as desired.
In particular,  we obtain the $\mfg$-valued $1$-form
\begin{align}
\omega_{\mcE/X} \in \Gamma (\mcE, \Omega_{\mcE/X} \otimes\mfg) 
\left(\cong \mr{Hom}_{\mcO_\mcE} (\mcT_{\mcE/X}, \mcO_\mcE \otimes \mfg) \right)
\end{align}
on $\mcE/X$ corresponding to $\omega_{\mcE/X}^\triangleright$.
By construction, we have $\omega_{G/k} = \omega_{G}$.

\LSP
%-----------------------------------------------------------------------------------------
\subsection{The curvature of a $1$-form} \label{SS090}

For $\mfg$-valued $1$-forms $\omega_1$, $\omega_2$ on $X/S$, 
we  define $[\omega_1, \omega_2]$ to be  the $\mfg$-valued $2$-form  on $X/S$ given by 
\begin{align} \label{QQ21}
[\omega_1, \omega_2]^\triangleright (\partial_1, \partial_2) := [\omega_1^\triangleright (\partial_1), \omega_2^\triangleright (\partial_2)] - [\omega_1^\triangleright (\partial_2), \omega_2^\triangleright (\partial_1)]
\end{align}
for any local sections $\partial_1, \partial_2 \in \mcT_{X/S}$.
In particular,  for a  $\mfg$-valued $1$-form $\omega$, we have 
$[\omega, \omega]^\triangleright (\partial_1, \partial_2) = 2 \cdot [\omega^\triangleright (\partial_1), \omega^\triangleright (\partial_2)]$.
Also, for a $\mfg$-valued $1$-form $\omega$ on $X/S$, 
 the {\it  exterior differential $d \omega$} of $\omega$ is defined as the $\mfg$-valued $2$-form on $X/S$ given by 
\begin{align}
(d \omega)^\triangleright (\partial_1, \partial_2) := \partial_1 (\omega^\triangleright (\partial_2)) -\partial_2 (\omega^\triangleright (\partial_1)) - \omega^\triangleright ([\partial_1, \partial_2])
\end{align}
for any local sections $\partial_1$, $\partial_2 \in \mcT_{X/S}$.
Then,
 the {\bf curvature} of a  $\mfg$-valued $1$-form $\omega$ on $X/S$  is defined as the  $\mfg$-valued  $2$-form  
 \begin{align} \label{Erosmf}
 \psi_\omega := d \omega + \frac{1}{2} \cdot [\omega, \omega] \in \Gamma (X, (\bigwedge^2 \Omega_{X/S}) \otimes \mfg).
 \end{align}
 That is to say,  it is  given by the formula 
 \begin{align}
 \psi_ \omega^\triangleright (\partial_1, \partial_2) =
 \partial_1 (\omega^\triangleright (\partial_2)) -\partial_2 (\omega^\triangleright (\partial_1)) + [\omega^\triangleright (\partial_1), \omega^\triangleright (\partial_2)] - \omega^\triangleright ([\partial_1, \partial_2])
 \end{align}
  for any local sections $\partial_1, \partial_2 \in \mcT_{X/S}$.
  We shall say that $\omega$ is {\bf flat} if $\psi_\omega =0$.
We know  (cf. ~\cite[\S\,1.1, Proposition 1.3]{Wak8}) that  the Maurer-Cartan form $\omega_G$ is flat, 
i.e., 
\begin{align} \label{E022}
\left(\psi_{\omega_G} :=  \right) d \omega_G  + \frac{1}{2} \cdot [\omega_G, \omega_G] =0.
\end{align}

\LSP
%-----------------------------------------------------------------------------------------
\subsection{Connections and $1$-forms on the trivial bundle} \label{QQ131}
Let us consider the relationship between the notion of curvature defined above and the classical notion of curvature of a connection (cf. (\ref{L054})).
Denote by $\mr{pr}_1$ (resp., $\mr{pr}_2$) the first projection $X \times G \migi X$ (resp., the second projection $X \times G \migi G$).
In particular, $\mr{pr}_1$ is nothing but the structure morphism of the trivial $G$-bundle $X \times  G$ on $X$.
The differentials of $\mr{pr}_1$ and $\mr{pr}_2$ give a direct sum decomposition  
\begin{align} \label{QQ13}
\mcT_{X \times G/S} \isom \mr{pr}_1^*(\mcT_{X/S})\oplus \mr{pr}_2^*(\mcT_{G}).
\end{align}
By applying the functor $\mr{pr}_{1*}(-)^G$ to this isomorphism, we obtain 
a decomposition 
\begin{align} \label{QQ136}
\tau : \widetilde{\mcT}_{X\times G/S} \isom \mcT_{X/S} \oplus \mcO_X \otimes \mfg
\index{$\tau$}
\end{align}
of $\widetilde{\mcT}_{X\times G/S}$, where the second factor  $\mcO_X \otimes \mfg$  in the codomain  is naturally isomorphic to $\mfg_{X \times G}$.
The morphism $d_{X\times G} : \widetilde{\mcT}_{X \times G/S} \migi \mcT_{X/S}$ (cf. (\ref{Ex0}))  coincides, via $\tau$, with the first projection $\mcT_{X/S} \oplus \mcO_X \otimes \mfg \migi \mcT_{X/S}$.
Under the identification given by  $\tau$, the Lie bracket operation on $\widetilde{\mcT}_{X \times G/S}$ is given by 
\begin{align} \label{QR922}
 [(\partial_1, a_1 \otimes v_1), (\partial_2, a_2 \otimes v_2)]  
=  ([\partial_1, \partial_2], a_1  a_2\otimes [v_1, v_2] + \partial_1 (a_2)\otimes v_2-\partial_2(a_1)\otimes v_1)
\end{align}
for any local sections $\partial_i \in \mcT_{X/S}$, $a_i \in \mcO_X$, and any  $v_i \in \mfg$ ($i=1,2$).

Now, let $\omega$ be a $\mfg$-valued $1$-form on $X/S$.
The assignment $\partial \mapsto (\partial, \omega^\triangleright (\partial))$ for each local section $\partial \in \mcT_{X/S}$ defines, via $\tau$,
an $S$-connection 
\begin{align} \label{L048}
\nabla_\omega : \mcT_{X/S} \migi \widetilde{\mcT}_{X \times G/S}
\end{align}
 on the trivial $G$-bundle $X \times G$.
 The connection 
 $\nabla_0$ (i.e., $\nabla_\omega$ with $\omega =0$) is called    the {\bf trivial connection}.
By taking account of (\ref{QR922}),
we can verify the equality 
\begin{align} \label{L049}
\psi_\omega = \psi_{\nabla_\omega}
\end{align}   
of elements in $\Gamma (X, (\bigwedge^2 \Omega_{X/S}) \otimes \mfg) \left(= \Gamma (X, (\bigwedge^2 \Omega_{X/S}) \otimes_{\mcO_X} \mfg_{X \times G}) \right)$.
In particular, $\omega$ is flat (in the above sense)  if and only if $\nabla_\omega$ is flat  (in the sense of \S\,\ref{SS002}).

\LSP
%-----------------------------------------------------------------------------------------
\subsection{Cartan geometries} \label{SS001}

We shall  recall the definition of a Cartan geometry.
A Cartan geometry has been formulated  mainly in the differential (and holomorphic) geometric setting, but it can also be defined for schemes.

\SSP
%----------------------------------[begin definition]------------------
%\vspace{3mm}
\bde  \label{D024}
\begin{itemize}
\item[(i)]
A {\bf Cartan geometry with model $(G, H)$} on $X/S$  is a pair
\begin{align}
\mcE^\cc := (\mcE_H, \omega)
\end{align}
consisting of an $H$-bundle $\mcE_H$ on $X$ and
a $\mfg$-valued $1$-form $\omega$ on $\mcE_H/S$ 
such that the corresponding $\mcO_{\mcE_H}$-linear morphism $\omega^\triangleright : \mcT_{\mcE_H/S} \migi \mcO_{\mcE_H} \otimes \mfg$
  satisfies  the following conditions:
\begin{itemize}
\item
The morphism  $\omega^\triangleright$ is an isomorphism and the following square diagram is commutative:
 \begin{align} \label{E0057}
\vcenter{\xymatrix@C=46pt@R=36pt{
\mcT_{\mcE_H/X}\ar[r]_-{\sim}^-{\omega^\triangleright_{\mcE_H/X}} \ar[d]_{\mr{inclusion}}&  \mcO_{\mcE_H} \otimes \mfh\ar[d]^{\mr{inclusion}} \\
\mcT_{\mcE_H/S}\ar[r]^-{\sim}_-{\omega^\triangleright}& \mcO_{\mcE_H} \otimes \mfg.
 }}
\end{align}
 \item
For any  $h \in H(k)$,
the following square diagram is commutative:
 \begin{align} \label{E0050}
\vcenter{\xymatrix@C=46pt@R=36pt{
\mcT_{\mcE_H/S}\ar[r]^-{\omega^\triangleright} \ar[d]^{\wr}_{d \mr{R}_h}&  \mcO_{\mcE_H} \otimes \mfg\ar[d]_{\wr}^{\mr{id} \otimes \mr{Ad}(h^{-1})} \\
\mr{R}_h^*(\mcT_{\mcE_H/S})\ar[r]_-{\mr{R}_h^*(\omega^\triangleright)}& \mcO_{\mcE_H} \otimes \mfg.
 }}
\end{align}
\end{itemize}
\item[(ii)]
Let $\mcE^\cc := (\mcE_H, \omega)$ and $\mcE'^\cc := (\mcE'_H, \omega')$ be  Cartan geometries with model $(G, H)$  on $X/S$.
An {\bf isomorphism of Cartan geometries with model $(G, H)$} from $\mcE^\cc$ to $\mcE'^\cc$ is defined as  an isomorphism $\eta : \mcE_H \isom \mcE'_H$ of $H$-bundles  with 
$\eta^* (\omega') = \omega$.
\end{itemize}
 \ede
%-------------------------[end definition]-------------------

\SSP

%----------------------------------[begin definition]------------------
\bde  \label{D027ff}
Let  $\mcE^\cc := (\mcE_H, \omega)$ be 
a Cartan geometry with model $(G, H)$ on $X/S$.
Then, 
we shall say that $(\mcE_H, \omega)$ is {\bf flat} if $\omega$ is flat.
\ede
%--------------------------------------------------------------------------------------------
\SSP

Denote by
\begin{align}
\mr{Car}_{(G, H), X/S} \ \left(\text{resp.,} \ \mr{Car}^{\mr{flat}}_{(G, H), X/S} \right)
\end{align}
the set of  isomorphism classes of  Cartan geometries  (resp., flat Cartan geometries) with model $(G, H)$ on $X/S$.

\SSP
%---------------------------------------------------------------------
\begin{exa}[Standard example] \label{Exam1}
Let us recall the standard (flat) Cartan geometry on the $S$-scheme  $P_S$ (cf. (\ref{L014})).
We shall write  $\mcE^\mr{triv}_H :=  S  \times G$; it  may be regarded as an $H$-bundle on $P_S$ by putting  the  morphism  $\mr{id}_S \times \pi_G : \mcE_H^\mr{triv}  \migisurj P_S$  as the structure morphism of $\mcE_H^\mr{triv}$ over $P_S$.
Note that  $\mcE_{H}^\mr{triv}$ specifies  an $H$-reduction on the trivial $G$-bundle $P_S \times G$ on $P_S$ via the morphism
\begin{align} \label{ghfjke}
(\mr{id}_S \times \pi_G, \mr{pr}_2) : \mcE_{H}^\mr{triv} \left(= S \times G \right) \migi P_S \times G,
\end{align}
where $\mr{pr}_2$ denotes the second projection $S \times G \migi G$.
Then, it follows from Proposition  \ref{L017}, (i), and (\ref{E022}) that
the pair 
\begin{align} \label{L045}
\mcE^{\mr{triv}, \cc} := (\mcE_H^{\mr{triv}}, \omega_{S \times G/S})
\end{align}
  forms a flat Cartan geometry with model $(G, H)$ on $P_S/S$.
\end{exa}
%---------------------------------------------------------------------

\LSP
%----------------------------------------------------------------------[begin subsection]-------------
\subsection{Pull-back and base-change of Cartan geometries} \label{SS02}

Let $Y$ and $y$ be as in the former half of \S\,\ref{SS066}.
Also, let $\mcE^\cc := (\mcE_H, \omega)$ be a  Cartan geometry with model $(G, H)$ on $X/S$.
The base-change   $y^*(\mcE_H)$ of $\mcE_H$ by $y$ defines an $H$-bundle on $Y$.
If $y_{\mcE_H}$ denotes the projection $y^*(\mcE_H) \migi \mcE_H$,
then its differential $d y_{\mcE_H} : \mcT_{y^*(\mcE_H)/S} \migi y_{\mcE_H}^*(\mcT_{\mcE_H /S})$ 
 is an isomorphism because of the \'{e}taleness of $y$.
 Hence, the pull-back $y^*(\omega) := y_{\mcE_H}^*(\omega)$
 specifies  a $\mfg$-valued $1$-form on $y^*(\mcE_H)/S$.
One may verify that
the pair 
 \begin{align} \label{L046}
 y^*(\mcE^\cc) := (y^*(\mcE_H), y^*(\omega))
 \end{align}
  forms  a Cartan geometry with model $(G, H)$ on $Y/S$.
 Moreover,  the equality $y^*_{\mcE_H} (\psi_\omega) = \psi_{y^*(\omega)}$ holds under the identification  $\mcT_{y^*(\mcE_H)/S} = y^*_{\mcE_H} (\mcT_{\mcE_H/S})$  given   by $d y_{\mcE_H}$, so
 $y^*(\mcE^\cc)$ is flat  if $\mcE^\cc$ is flat.
The Cartan geometry $y^*(\mcE^\cc)$ is called the {\bf pull-back} of $\mcE^\cc$. The resulting assignment
$\mcE^\cc \mapsto y^*(\mcE^\cc)$
defines maps
\begin{align} \label{L023}
\mr{Car}_{(G, H), X/S} \migi \mr{Car}_{(G, H), Y/S}, 
\hspace{5mm}
\mr{Car}_{(G, H), X/S}^{\mr{flat}} \migi \mr{Car}_{(G, H), Y/S}^{\mr{flat}}.
\end{align}
Note that  the formation of a (flat) Cartan geometry with model $(G, H)$ on $X/S$ has descent with respect to the \'{e}tale topology on $X$.

Next, 
let $s : S' \migi S$ and $X'$ be as in the latter half of  \S\,\ref{SS066}, and let $\mcE^\cc$ be as above.
The base-change $s^*(\mcE_H)$ defines an $H$-bundle on $X'$.
Since  the differential $d s_{\mcE_H} : \mcT_{s^*(\mcE_H)/S'} \isom s_{\mcE_H}^* (\mcT_{\mcE_H/S})$ of the natural  projection  $s_{\mcE_H} : s^*(\mcE_H) \migi \mcE_H$ is an isomorphism,
the base-change $s^*(\omega)$ of $\omega$ defines a $\mfg$-valued $1$-form on $s^*(\mcE_H)$.
Moreover, the pair
\begin{align}
s^*(\mcE^\cc):= (s^*(\mcE_H), s^*(\omega))
\end{align}
forms a Cartan geometry with model $(G, H)$ on $X'/S'$.
Under the identification $\mcT_{s^*(\mcE_H)/S'} = s_{\mcE_H}^* (\mcT_{\mcE_H/S})$ given by $d s_{\mcE_H}$, the equality $s^*_{\mcE_H}(\psi_\omega) = \psi_{s^*(\omega)}$ holds.
This implies that $s^*(\mcE^\cc)$ is flat if $\mcE^\cc$ is flat.
The Cartan geometry $s^*(\mcE^\cc)$ is called the {\bf base-change} of $\mcE^\cc$, and 
the assignment
$\mcE^\cc \mapsto s^*(\mcE^\cc)$
defines maps 
\begin{align} \label{L024}
\mr{Car}_{(G, H), X/S} \migi \mr{Car}_{(G, H), X'/S'},
\hspace{5mm}
 \mr{Car}_{(G, H), X/S}^{\mr{flat}} \migi \mr{Car}_{(G, H), X'/S'}^{\mr{flat}}.
\end{align}

\LSP
%---------------------------[begin subsection]-------------
\subsection{From Frobenius-Ehresmann structures  to Cartan geometries} \label{SS01}
In this subsection, we construct a Cartan geometry with model $(G, H)$ using an $F$-$(G, H)$-structure.

Let $\mcS^\hh$ be an $F$-$(G,H)$-structure   on $X/S$.
Since $G$ is quasi-projective (cf. ~\cite[\S\,8, Theorem 8.43]{Mil}),
we  obtain the  {\it right} $G$-bundle $\mcE_{\mcS^\hh, G}$ on $X$ representing the inverse of  the left $G_X$-torsor $\mcS^\hh \times^{G_{X/S}^{(1)}} G_X$ (cf. ~\cite[Chap.\,III, \S\,4, Theorem 4.3]{Mil4}). 
The  graphs $\Gamma_\phi : U \migi U \times P \left(= U \times_S P_S \right)$ for various local sections   $\phi : U \migi P_S$ in $\mcS^\hh$ 
may be glued together to construct
an $H$-reduction  $\mcE_{\mcS^\hh, H}$ of $\mcE_{\mcS^\hh, G}$.
In particular,  for such a section $\phi$, there exists a well-defined (up to a right translation by $H$)  isomorphism of $G$-bundles $\tau_{\phi} : \mcE_{\mcS^\hh, G} |_U \isom U \times G$ via which  the $H$-reduction $\mcE_{\mcS^\hh, H}|_U$ of $\mcE_{\mcS^\hh, G} |_U$ corresponds to  that of  $U \times G$ determined by  the graph $\Gamma_{\phi}$.  
Next, let us choose  a collection  $ \{ (U_\alpha, \phi_\alpha) \}_{\alpha \in I}$, where $I$ denotes  a finite set, $\{ U_\alpha \}_{\alpha \in I}$ is  an  \'{e}tale covering  of $X$,  and each $\phi_\alpha$ ($\alpha \in I$) denotes a morphism $U_\alpha \migi P_S$ belonging to   $\mcS^\hh (U_\alpha)$.
For each $\alpha \in I$, we shall denote by $\phi_{\mcE, \alpha}$ the $H$-equivariant morphism $\mcE_{\mcS^\hh, H} |_{U_\alpha} \migi  \mcE^\mr{triv}_H$ defined as  the following composite:
\begin{align}
\phi_{\mcE, \alpha} : \mcE_{\mcS^\hh, H} |_{U_\alpha}
 \xrightarrow{\mr{inclusion}}
 \mcE_{\mcS^\hh, G} |_{U_\alpha} \xrightarrow{\tau_{\phi_\alpha}}
U_\alpha \times G \left(= U_\alpha \times_S  \mcE^\mr{triv}_H\right) \xrightarrow{\mr{pr}_2}  \mcE^\mr{triv}_H.
\end{align}
By  the definition of $\tau_{\phi_\alpha}$, 
the following square diagram is commutative and hence cartesian:
 \begin{align} \label{E0021}
\vcenter{\xymatrix@C=46pt@R=36pt{
\mcE_{\mcS^\hh, H} |_{U_\alpha}\ar[r]^{\phi_{\mcE, \alpha}} \ar[d]_{\mr{projection}}&  \mcE^\mr{triv}_H \ar[d]^{\mr{id}_S \times \pi_G} \\
U_\alpha\ar[r]_{\phi_\alpha}& P_S.
 }}
\end{align}
Hence,  the \'{e}taleness of $\phi_\alpha$ implies that of $\phi_{\mcE, \alpha}$.
Under the identification $\mcE_{\mcS^\hh, H}|_{U_\alpha} = U_\alpha \times_{P_S} \mcE^{\mr{triv}}_H$ given by 
the above  diagram, 
the pull-back 
\begin{align}
(\mcE_{\mcS^\hh, H}|_{U_\alpha}, \phi_{\mcE, \alpha}^*(\omega_{S \times G/S}))
\end{align}
 of $(\mcE_H^{\mr{triv}}, \omega_{S \times G/S})$ forms a  flat  Cartan geometry with model $(G, H)$ on $U_\alpha /S$
(cf. (\ref{L045}), (\ref{L046})).
Note that $\phi_{\mcE, \alpha}$ does not depend on the choice of $\tau_{\phi_\alpha}$, so the $\mfg$-valued $1$-form $\phi_{\mcE, \alpha}^*(\omega_{S \times_k G/S})$ depends only on $(U_\alpha, \phi_\alpha)$.

 Next, denote by $I_2$ the subset of $I \times I$ consisting of pairs $(\alpha, \beta)$ with $U_{\alpha \beta} := U_\alpha \cap U_\beta \neq \emptyset$.
Write $\mcE_{H, \alpha \beta} := \mcE_{\mcS^\hh, H}|_{U_{\alpha \beta}}$ and $\mcE_{G, \alpha \beta} := \mcE_{\mcS^\hh, G}|_{U_{\alpha \beta}}$.
For each pair $(\alpha, \beta) \in I_2$,
there exists $h \in G_{X/S}^{(1)}(U_{\alpha \beta}) \left(\subseteq G (U_{\alpha \beta}) \right)$ satisfying the equality  $h (\phi_\alpha  |_{U_{\alpha \beta}}) = \phi_\beta |_{U_{\alpha \beta}}$.
We can find an element $g \in H (U_{\alpha \beta})$ with
$\mr{L}_h \circ  (\tau_{\phi_\alpha}  |_{\mcE_{G, \alpha \beta}}) = \mr{L}_g \circ (\tau_{\phi_\beta} |_{\mcE_{G, \alpha \beta}})$.
Then, we have
\begin{align}
\phi^*_{\mcE, \alpha} (\omega_{S \times G/S})  |_{\mcE_{H, \alpha \beta}}
& =  \tau_{\phi_\alpha}^*(\mr{pr}_2^*(\omega_{S \times G/S}))) |_{\mcE_{H, \alpha \beta}} \\
& =   
(\mr{L}_{h} \circ (\tau_{\phi_\alpha} |_{\mcE_{G, \alpha \beta}}))^*(\mr{L}_{h^{-1}}^*(\mr{pr}_2^*(\omega_{S \times G/S}))) |_{\mcE_{H, \alpha \beta}}  \notag \\
& =   
(\mr{L}_{h} \circ (\tau_{\phi_\alpha} |_{\mcE_{G, \alpha \beta}}))^*(\mr{pr}_2^*(\omega_{S \times G/S})) |_{\mcE_{H, \alpha \beta}}  \notag \\
& =  (\mr{L}_g \circ (\tau_{\phi_\beta} |_{\mcE_{G, \alpha \beta}}))^* (\mr{pr}_2^*(\omega_{S \times G/S}))|_{\mcE_{H, \alpha \beta}} \notag \\
& = 
 \tau_{\phi_\beta}^*(\mr{pr}_2^*(\omega_{S \times G/S}))) |_{\mcE_{H, \alpha \beta}}  \notag \\
 & = \phi^*_{\mcE, \beta} (\omega_{S \times G/S})  |_{\mcE_{H, \alpha \beta}}, \notag
\end{align}
where  the third  and fifth  equalities follow from Proposition  \ref{L017}, (i).
It follows that $\phi^*_{\mcE, \alpha}(\omega_{S \times G/S})$'s  may be glued together to obtain a $\mfg$-valued $1$-form $\omega_{\mcS^\hh}$ on $\mcE_{\mcS^\hh, H}/S$. 
 The resulting pair
\begin{align}
\mcS^{\hh \Rightarrow \cc} := (\mcE_{\mcS^\hh, H}, \omega_{\mcS^\hh})
\end{align}
 forms a flat Cartan geometry with model $(G, H)$ on $X/S$.
Since    $\mcS^{\hh \Rightarrow \cc}$ does not depend on the choice of $\{ (U_\alpha, \phi_\alpha) \}_\alpha$,  the assignment $\mcS^\hh \mapsto \mcS^{\hh \Rightarrow \cc}$ determines a well-defined  map of sets
 \begin{align} \label{L096}
 F\text{-} \mr{Ehr}_{(G, P), X/S} \migi \mr{Car}^\mr{flat}_{(G, H), X/S}.
 \end{align}
 In the case of $X = P_S$, 
  this map  sends the $F$-$(G, P)$-structure $\mcS_1^{\mr{triv}, \hh}$ (cf. (\ref{EE012}))  to the flat Cartan geometry $\mcE^{\mr{triv}, \cc}$ (cf. (\ref{L045})).

The formation of $\mcS^{\hh \Rightarrow \cc}$ commutes with pull-back 
to \'{e}tale $X$-schemes, as well as with base-change to $S$-schemes.
That is to say,  with the notation in \S\,\ref{SS066}, the following square diagrams are commutative:
\begin{align} \label{E00100}
\vcenter{\xymatrix@C=50pt@R=36pt{
F\text{-} \mr{Ehr}_{(G, H), X/S}  \ar[r]^{\text{(\ref{L096}) for $X/S$}}  \ar[d]_-{(\ref{L0981})} & \mr{Car}^\mr{flat}_{(G, H), X/S} \ar[d]^-{(\ref{L023})}
\\
F\text{-} \mr{Ehr}_{(G, H), Y/S} \ar[r]_{\text{(\ref{L096}) for $Y/S$}}   & \mr{Car}^\mr{flat}_{(G, H), Y/S},  
 }}\hspace{7mm}
 \vcenter{\xymatrix@C=50pt@R=36pt{
F\text{-} \mr{Ehr}_{(G, H), X/S}  \ar[r]^{\text{(\ref{L096}) for $X/S$}}  \ar[d]_-{(\ref{L009})} & \mr{Car}^\mr{flat}_{(G, H), X/S} \ar[d]^-{(\ref{L024})}
\\
F\text{-} \mr{Ehr}_{(G, H), X'/S'} \ar[r]_{\text{(\ref{L096}) for $X'/S'$}}   & \mr{Car}^\mr{flat}_{(G, H), X'/S'}. 
 }}
\end{align}

\LSP
%---------------------------[begin subsection]-------------
\subsection{The $p$-curvature of a $1$-form} \label{SS0131}

%In this section, we  introduce  the notion of $p$-curvature of a $\mfg$-valued $1$-form.

Recall that, for   a Lie algebra $\mfa$ over a $k$-algebra $R$,
a {\it $p$-power operation} on $\mfa$ is
 a map $(-)^{[p]} :\mfa \migi \mfa$ satisfying the following conditions:
\begin{itemize}
\item
$(\mr{ad}(v))^p = \mr{ad}(v^{[p]})$ for all $v \in \mfa$;
\item
$(a\cdot v)^{[p]} = a^p \cdot v^{[p]}$ for all $a \in R$, $v \in \mfa$;
\item
$(v + u)^{[p]} = v^{[p]}+ u^{[p]} + \sum_{i=1}^{p-1} \frac{s_i (v, u)}{i}$ for all $v, u \in \mfa$, where  $s_i (v, u)$ denotes the coefficient of $t^{i-1}$ in the  expression $\mr{ad}(vt + u)^{p-1}(v) \in \mfa [t] \left(:=\mfa \otimes_R R[t]\right)$.
\end{itemize}
A  Lie algebra over $R$ together with  a $p$-power operation is called a
 {\it  restricted Lie algebra} (or a {\it  $p$-Lie algebra}) over $R$.
For example,  
the space of $R$-linear endomorphisms $\mr{End}_R(V)$ of $R$-module $V$
 gives a typical example of  a restricted Lie algebra over $R$ by putting  
   $[v, u] := vu -uv$ as a bracket operation and putting $v^{[p]} := v^p$ as a  $p$-power operation.

In the natural fashion, we can define the notion of 
a {\it restricted $\mcO_S$-Lie algebra} on $X$ (i.e., a sheaf of restricted $\mcO_S$-Lie algebra in the sense of ~\cite{Lan}, \S\,4.1).
Just as in the case of $\mr{End}_R(V)$ above,
the sheaf   $\mcE nd_{f^{-1}(\mcO_S)} (\mcV)$  defined  for each $\mcO_X$-module  $\mcV$  has a natural structure of restricted $\mcO_S$-Lie algebra on $X$ with 
  $p$-power operation $v^{[p]} := v^p$ (= the $p$-th iteration of $v$) for any local section $v \in \mcE nd_{f^{-1}(\mcO_S)} (\mcV)$.
Recall that the tangent  bundle $\mcT_{X/S}$ may be identified with the subsheaf of 
$\mcE nd_{f^{-1}(\mcO_S)} (\mcO_X)$ consisting of 
locally defined derivations on $\mcO_X$ relative to $S$.
The $p$-power operation on $\mcE nd_{f^{-1}(\mcO_S)} (\mcO_X)$ is  restricted to  
 $\mcT_{X/S}$.
 In particular, 
 we obtain 
  the $p$-power operation
$(-)^{[p]} : \mfg \migi \mfg$
 on the Lie algebra $\mfg$ arising from that on $\mcT_{G}$ via  the  identification of  $\mfg  = \Gamma (G, \mcT_{G})^G$;
 it induces naturally a $p$-power operation on $\mcO_X \otimes \mfg$, by which $\mcO_X \otimes_k \mfg$ forms a restricted $\mcO_X$-Lie algebra.

Now, let 
 $\omega$ be  a flat $\mfg$-valued $1$-form on $X/S$.
For each local section $\partial \in \mcT_{X/S}$,
we shall write 
\begin{align} \label{eoa03-3x}
{^p}\psi_{\omega, \partial} := \omega^\triangleright (\partial)^{[p]} - \omega^\triangleright (\partial^{[p]}) +  \sum_{i=1}^{p-1} \frac{s_i (\partial, \omega^\triangleright (\partial))}{i} \in \mcO_X \otimes \mfg.
\end{align}

\SSP
%----------------------------------------------------------------
\ble \label{L0103}
Let 
$\nabla_\omega$ 
be  the flat $S$-connection  on the trivial $G$-bundle $X \times  G$ 
corresponding to $\omega$  (cf.  (\ref{L048})).
Then,  for each local section $\partial \in \mcT_{X/S}$, the following equality  holds:
\begin{align}
{^p}\psi_{\omega, \partial} =  \langle F^{-1}_X (\partial), {^p}\psi_{\nabla_\omega}\rangle
\end{align}
(cf. (\ref{L056}) for the definition of ${^p}\psi_{\nabla_\omega}$), where $\langle -, - \rangle$ denotes the natural pairing $F_X^*(\mcT_{X/S}) \times F^*_{X}(\Omega_{X/S}) \otimes_{\mcO_X} \mfg_{X \times G} \migi \mfg_{X \times G} \left(= \mcO_X \otimes \mfg \right)$.
\ele
%----------------------------------------------------------------
\begin{proof}
Under  the identification $\widetilde{\mcT}_{X \times G/S} = \mcT_{X/S} \oplus \mcO_X \otimes \mfg$ given by $\tau$ (cf. (\ref{QQ136})),
we have
\begin{align}
 & \ \ \  \ \langle F^{-1}_X(\partial), {^p}\psi_{\nabla_\omega}\rangle \\
 &=  (\partial, \omega^\triangleright (\partial))^{[p]} - (\partial^{[p]}, \omega^\triangleright (\partial^{[p]})) \notag \\
 & =  ((\partial, 0) + (0, \omega^\triangleright (\partial)))^{[p]}  -(\partial^{[p]}, \omega^\triangleright (\partial^{[p]})) \notag \\
 & = (\partial^{[p]}, 0) + (0, \omega^\triangleright (\partial)^{[p]}) + \left(0,  \sum_{i=1}^{p-1} \frac{s_i ((\partial, 0), (0, \omega^\triangleright (\partial)))}{i} \right)-(\partial^{[p]}, \omega^\triangleright (\partial^{[p]})) \notag \\
 & = \left(0, \omega^\triangleright (\partial)^{[p]} -\omega^\triangleright (\partial^{[p]}) + \sum_{i=1}^{p-1} \frac{s_i (\partial, \omega^\triangleright (\partial))}{i} \right) \notag \\
 & = \left(0, {^p}\psi_{\omega, \partial} \right). \notag
\end{align}
This completes the proof of the assertion.
\end{proof}
%-------------------------------------------------------------
\SSP

By  the above  lemma and a property of  ${^p}\psi_{\nabla_\omega}$,
we obtain the section
\begin{align}
{^p}\psi_\omega \in \Gamma (X, F^*_X (\Omega_{X/S}) \otimes \mfg)
\end{align}
determined by $\langle F^{-1}_X (\partial), {^p}\psi_\omega \rangle = {^p}\psi_{\omega, \partial}$ for any local section $\partial \in \mcT_{X/S}$.

\SSP
%------------------------------------------------------------------------
\bde \label{D00102}
We shall refer to  ${^p}\psi_\omega$ as the {\bf $p$-curvature} of $\omega$.
Also, we shall say that $\omega$ is {\bf $p$-flat} if ${^p}\psi_\omega =0$.
(Hence, the above lemma implies that $\omega$ is $p$-flat if and only if $\nabla_\omega$ is $p$-flat in the sense of \S\,\ref{SS002}).
\ede
%------------------------------------------------------------------
\SSP

%-------------------------------------------------------------
\begin{exa} \label{Eorkcr0}
Suppose that $G$ is abelian, which 
 implies that the Lie bracket operation on $\mfg$ is trivial.
Let us  write $\omega = \sum_{i=1}^l \omega_i \otimes e_i$ ($l \geq 1$)
for $\omega_i \in \Gamma (X, \Omega_{X/S})$ and $e_i \in \mfg$.
Then,  $\omega$ is flat if and only if, for every $i =1, \cdots, l$, the $1$-form  $\omega_i$ is flat, or equivalently, closed.
Moreover,  if $\mfg = k$ (i.e., $\mr{dim}(G) =1$) and  $\omega$ is flat, then
it follows from ~\cite[\S\,7, Proposition 7.1.2]{Katz2} that
the $p$-curvature ${^p}\psi_\omega$ of  $\omega$ satisfies 
\begin{align}
\langle F_X^{-1}(\partial), {^p}\psi_\omega \rangle \left(= \omega^\triangleright (\partial)^{[p]}- \omega^\triangleright (\partial^{[p]}) + \partial^{p-1} (\omega^\triangleright (\partial))\right) = \omega^\triangleright (\partial)^p - \langle F^*_S(\partial),   C (\omega) \rangle
\end{align}
 for any $\partial \in \mcT_{X/S}$, where $C$ denotes the usual  Cartier operator $F_{X/S*}(\mr{Ker}(\nabla_0^{(2)})) \migi \Omega_{X^{(1)}/S}$ (cf. (\ref{E0124}) for the definition of $(-)^{(2)}$) and $F^*_S (\partial)$ denotes the local section of $\mcT_{X^{(1)}/S}$ corresponding to $\partial$ via the isomorphism $(\mr{id}_X \times F_S)^*(\Omega_{X/S}) \isom \Omega_{X^{(1)}/S}$.
\end{exa}
%-------------------------------------------------------------

\LSP
%\vspace{5mm}
%---------------------------[begin subsection]-------------
\subsection{$p$-flatness of the Maurer-Cartan form} \label{SSG078}

In the differential geometric setting, 
the Maurer-Cartan equation (i.e., the flatness condition) associated with a $\mfg$-valued $1$-form  gives the local obstruction to its integrability
(cf.  ~\cite[Chap.\,3, \S\,6, Theorem 6.1]{Sha}).
By considering the $p$-flatness condition in addition, 
we can prove the  positive  characteristic version of this fact, as described in Proposition \ref{ETOFW2} below.
To begin with, let us observe the following proposition.

\SSP
%------------------------------------------------------------------------------------
\bpr \label{P0014}
The Maurer-Cartan form $\omega_G$ is $p$-flat.
\epr
%------------------------------------------------------------------------------------
\begin{proof}
Let us take an element $v$ of $\mfg$.
The isomorphism  $\omega_G^\triangleright : \mcT_G \isom \mcO_G \otimes \mfg$ (cf. (\ref{E023})) preserves the $p$-power operations, so we have
$\omega_G^\triangleright (v^{\dagger [p]}) = \omega_G^\triangleright (v^\dagger)^{[p]} = 1 \otimes v^{[p]}$.
Also, note that $\partial  (\omega_G^\triangleright (v^\dagger)) \left(=\partial (1 \otimes v) = \partial (1) \otimes v \right)=0$ for any $\partial \in \mcT_{G}$,
  in particular, the equality $s_i (v^\dagger, \omega^\triangleright_G (v^\dagger)) =0$ holds for every  $i=1, \cdots, p-1$.
Hence,
\begin{align} \label{e0r0lled}
{^p}\psi_{\omega_G, v^\dagger} = \omega^\triangleright_G (v^\dagger)^{[p]}-\omega^\triangleright_G (v^{\dagger [p]}) + \sum_{i=1}^{p-1} \frac{s_i (v^\dagger, \omega_G^\triangleright (v^\dagger))}{i} = 1 \otimes v^{[p]}- 1 \otimes v^{[p]} +0 = 0. 
\end{align}
Since $F^*_G (\mcT_G)$ is generated by the sections $F_G^{-1}(v^\dagger)$ for various $v \in \mfg$,  (\ref{e0r0lled}) implies  ${^p}\psi_{\omega_G} =0$.
This  completes  the proof of the  assertion.
\end{proof}
\SSP

%-------------------------------------------------------------
\bpr \label{ETOFW2}
Suppose that $G$ is affine.
\begin{itemize}
\item[(i)]
Let $\omega$ be a  $\mfg$-valued $1$-form on $X/S$.
Then, $\omega$  is $p$-flat if and only if each geometric point $x$ of $X$ admits  an \'{e}tale neighborhood $U \migi X$ together with   an $S$-morphism $h : U \migi G$ such that $h^*(\omega_G) = \omega$ (or equivalently, its differential $dh: \mcT_{U/S} \migi h^*(\mcT_G)$ composed with $h^*(\omega_G^\triangleright) : h^*(\mcT_G) \isom \mcO_U \otimes \mfg$ coincides with $\omega^\triangleright$).
\item[(ii)]
Let $h_1$ and $h_2$ be $S$-morphisms $X \migi G$ with $h_1^*(\omega_G) = h_2^*(\omega_G)$.
Then, there exists a unique global section   $g$ of $G_{X/S}^{(1)}$ with $h_1 = (g \circ F_{X/S}) \cdot h_2$, where the right-hand side of this equality denotes the morphism $X \migi G$ given by $a \mapsto (g \circ F_{X/S}) (a)\cdot h_2 (a)$.
\end{itemize}
\epr
%-------------------------------------------------------------
\begin{proof}
First, let us consider assertion (i).
The ``if" part of the required equivalence follows from Proposition \ref{P0014}.
To prove the ``only if" part, let us suppose that $\omega$ is $p$-flat.
Then, Lemma \ref{L0103} implies ${^p}\psi_{\nabla_\omega} =0$.
By  Proposition \ref{Ertwq}, (i), together with the affineness of $G$, there exist a $G$-bundle $\mcG$ on $X^{(1)}$  and  an isomorphism of flat $G$-bundles $\alpha : (X \times G, \nabla_\omega)  \isom (F_{X/S}^*(\mcG), \nabla^\mr{can}_\mcG)$ (cf. (\ref{QQwkko}) for the definition of $\nabla_\mcG^\mr{can}$).
For each geometric point $x$ of $X$, we can find an \'{e}tale  neighborhood $U \migi X$ of $x$ such that the $G$-bundle $\mcG |_{U^{(1)}}$ admits a trivialization $\beta : \mcG |_{U^{(1)}} \isom U^{(1)} \times G$.
The pull-back  
\begin{align}
\beta^F : F^*_{U/S}(\mcG |_{U^{(1)}}) \isom \left(F^*_{U/S} (U^{(1)} \times G) =  \right) U \times G
\end{align}
 of this trivialization   by $F_{U/S}$ 
 defines an isomorphism $(F^*_{U/S}(\mcG |_{U^{(1)}}), \nabla^\mr{can}_\mcG |_U) \isom (U \times G, \nabla_0)$.
 The composite automorphism  $\beta^F \circ \alpha$ of the trivial $G$-bundle $U \times G$ may be expressed as the left-translation $\mr{L}_h$ by some $h \in G (U)$.
In particular,  the gauge transformation  $\mr{L}_h^*(\nabla_0)$ of the trivial connection  $\nabla_0$ on $U \times G$ coincides with $\nabla_\omega$.
 On the other hand,  it follows from  ~\cite[Chap.\,1, Proposition 1.22]{Wak8} that  
 $\nabla_{\omega_G}$ coincides with $\mr{L}_{\mr{id}_G}^*(\nabla_0)$, i.e.,  the gauge transformation  of  $\nabla_0$  using  the left-translation $\mr{L}_{\mr{id}_G}$ by    $\mr{id}_G : G \migi  G$.
Since  the diagram
 \begin{align} \label{EriTGG}
\vcenter{\xymatrix@C=46pt@R=36pt{
U \times G \ar[r]^-{h \times \mr{id}_G} \ar[d]_-{\mr{L}_h : (a, b) \mapsto (a, h (a) \cdot b)} & G \times G \ar[d]^-{\mr{L}_{\mr{id}_G} : (a, b) \mapsto (a, a \cdot b)}
\\
U \times G\ar[r]_-{h \times \mr{id}_G} & G \times G
 }}
\end{align}
is commutative, the following sequence of equalities holds:
\begin{align}
h^*(\nabla_{\omega_G}) = h^*(\mr{L}_{\mr{id}_G}^*(\nabla_0)) = \mr{L}_h^*(h^*(\nabla_0)) = \mr{L}_h^*(\nabla_0) = \nabla_{\omega}
\end{align}
(cf. (\ref{E0068}) for the definition of pull-back  of a connection $h^*(-)$).
This implies $h^*(\omega_G) = \omega$, and hence, completes the proof of the ``only if" part.
 
 Next, we shall consider assertion (ii).
Let us write  $h := h_1 \cdot h_2^{-1}$, i.e.,  the morphism  $X \migi G$ given by $a \mapsto h_1 (a) \cdot h_2^{-1}(a)$ for each $a \in X$.
Then,  $h$ decomposes as 
\begin{align}
X \xrightarrow{(h_1, h_2)} G \times G 
\xrightarrow{\mr{id} \times \iota} G \times G \xrightarrow{\mu} G
\end{align}
(cf. Lemma \ref{E938dK} for the definitions of $\mu$ and $\iota$).
By  the following Lemma \ref{E938dK}, (i) and (ii), we have 
\begin{align} \label{Eoap029}
h^*(\omega_G) = h_2^*(\mr{Ad}_G)(h_1^*(\omega_G) - h_2^{*}(\omega_G)) = 0
\end{align}
(cf. ~\cite[Chap.\,3, \S\,4, Corollary 4.11]{Sha}).
Since $\omega_G^\triangleright : \mcT_{G} \migi \mcO_G \otimes \mfg$ is an isomorphism, (\ref{Eoap029}) shows that the differential $d h : \mcT_{X/S} \migi h^*(\mcT_G)$ vanishes identically.
Hence, it follows from Lemma \ref{EPWP934} below that
$h$ lies in $G_{X/S}^{(1)} \left(\subseteq G_X \right)$.
This completes the proof of assertion (ii).
\end{proof}
%-------------------------------------------------------------

\SSP

The following two lemmas were used in the proof of the above proposition.

\SSP

%-------------------------------------------------------------
\ble \label{E938dK}
We shall abuse notation by writing   $\mr{Ad}_G$ for the $\mcO_G$-linear  automorphism of $\mcO_G \otimes \mfg$ given by the adjoint representation $\mr{Ad}_G : G \migi \mr{GL}(\mfg)$.
\begin{itemize}
\item[(i)]
Let $\mu : G \times G \migi G$ denote the multiplication of $G$.
Then,
the following equality  holds:
\begin{align}
\mu^*(\omega_G) = \mr{pr}_2^*(\mr{Ad}_G^{-1})(\mr{pr}_1^*(\omega_G)) + \mr{pr}_2^*(\omega_G),
\end{align}
 where $\mr{pr}_i$ ($i=1,2$) denotes the $i$-th projection $G \times G \migi G$.
\item[(ii)]
Let $\iota : G \isom G$ denote the inverse of $G$.
Then, the  following equality  holds
\begin{align}
\iota^*(\omega_G) = -\mr{Ad}_G (\omega_G).
\end{align}
\end{itemize}
\ele
%-------------------------------------------------------------
\begin{proof}
The proofs are entirely similar to ~\cite[Chap.\,3, S\,4, Proposition 4.10]{Sha}.
\end{proof}
%-------------------------------------------------------------

\SSP
%-------------------------------------------------------------
\ble \label{EPWP934}
Let $Y$ be a scheme over $S$ and $g$ an $S$-morphism $X \migi Y$.
Then, the morphism $dg^\vee : g^* (\Omega_{Y/S}) \migi \Omega_{X/S}$
induced by $g$ (i.e., the dual of the differential of $g$)
  is the zero map if and only if $g$ factors through $F_{X/S} : X \migi X^{(1)}$.  
\ele
%-------------------------------------------------------------
\begin{proof}
The ``if" part of the required equivalence  is clear since $d F_{X/S} =0$.
Next, let us consider  the ``only if" part.
Observe that the following diagram is commutative:
\begin{align} \label{EriTY}
\vcenter{\xymatrix@C=46pt@R=36pt{
&& g^*(\mcO_Y) \ar[r]^-{g^*(d)} \ar[d]_-{g^\sharp}   & g^*(\Omega_{Y/S}) \ar[d]^-{dg^\vee }
\\
 0  \ar[r] & F_{X/S}^{-1}(\mcO_{X^{(1)}}) \ar[r]_-{\mr{inclusion}}& \mcO_X \ar[r]_-{d}& \Omega_{X/S},
 }}
\end{align}
where the lower horizontal sequence is exact because of 
the smoothness of $X/S$.
If $dg^\vee =0$, then  the morphism $g^\sharp  : g^*(\mcO_Y) \migi \mcO_X$ satisfies $d \circ g^\sharp \left(= dg^\vee \circ g^*(d) \right)=0$, so it factors through the inclusion $ F_{X/S}^{-1}(\mcO_{X^{(1)}}) \migiincl \mcO_X$.
The resulting morphism  $g^*(\mcO_Y) \migi F_{X/S}^{-1}(\mcO_{X^{(1)}})$ determines a factorization of $g$ as $X \xrightarrow{F_{X/S}}X^{(1)} \migi Y$.
This completes the proof of this lemma.
\end{proof}
%-------------------------------------------------------------

\LSP
%\vspace{5mm}
%---------------------------[begin subsection]-------------
\subsection{$p$-flat Cartan geometries} \label{SS078}

Let us  introduce $p$-flat Cartan geometries and give the bijective correspondence  with  $F$-$(G, P)$-structures (cf. Theorem \ref{P0011}); the construction of this correspondence  is a part of Theorem \ref{TheoremA1}.

\SSP
%------------------------------------------
\bde \label{ea840q}
Let $\mcE^\cc := (\mcE_H, \omega)$ be a $(G, H)$-Cartan geometry on $X/S$.  
We shall say that $\mcE^\cc$ is {\bf $p$-flat} if
 $\omega$ is flat and $p$-flat.
 \ede
%------------------------------------------
\SSP

%------------------------------------------
\begin{exa}\label{EQLS2312}
Recall  the flat Cartan geometry $\mcE^{\mr{triv}, \cc}$ on $P_S/S$
defined in Example \ref{Exam1}.
Then,  the result of  Proposition \ref{P0014} says that it is $p$-flat.
\end{exa}
%------------------------------------------

\SSP

Denote by
\begin{align} \label{eoanv55677}
\mr{Car}^{p\text{-}\mr{flat}}_{(G, H), X/S}
\end{align}
the subset of
$\mr{Car}^{\mr{flat}}_{(G, H), X/S}$ consisting of 
 $p$-flat  Cartan geometries.
 Here, we  shall observe the following assertion, which will be used in the proof of Theorem \ref{P0011} described next.

\SSP
%-------------------------------------------------------------
\bt \label{P0011}
Suppose that $G$ is affine.
Then, the  image of (\ref{L096}) is contained in $\mr{Car}_{(G, H), X/S}^{p\text{-}\mr{flat}}$.
Moreover, the resulting map 
\begin{align} \label{E0223}
\zeta^{\hh \Rightarrow \cc} :  F \text{-}\mr{Ehr}_{(G, P), X/S} \migi \mr{Car}^{p\text{-}\mr{flat}}_{(G, H), X/S}
\end{align}
is bijective.
In particular, a  Cartan geometry with model $(G, H)$ on $X/S$ comes from an  $F$-$(G, P)$-structure   if and only if it is $p$-flat.
\et
%-------------------------------------------------------------
\begin{proof}
The first  assertion follows from Corollary \ref{EQkji} and the construction of $\mcS^{\hh \Rightarrow \cc}$ for a $F$-$(G, P)$-structure $\mcS^\hh$.
Hence, the remaining portion is to prove  the bijectivity of $\zeta^{\hh \Rightarrow \cc}$.
To this end, we shall construct the inverse to that map.
Let $\mcE^\cc := (\mcE_H, \omega)$ be a $p$-flat Cartan geometry with model $(G, H)$ on $X/S$.
By  Lemma \ref{L0103} and  the $p$-flatness of $\omega$,  the $S$-connection $\nabla_\omega$ on the trivial $G$-bundle $\mcE_H \times G$ has vanishing $p$-curvature.
It follows that 
the $p$-flat $G$-bundle $(\mcE_H \times G, \nabla_\omega)$ corresponds to a $G$-bundle $\mcF$ on  $\mcE_H^{(1)} \left(:= \mcE_H \times_{F_S,  S} S \right)$ via 
 (\ref{EEjj32}).
That is to say, 
there exists  
 an isomorphism of $G$-bundles $\widetilde{\alpha} : \mcE_H \times G \isom F^*_{\mcE_H/S}(\mcF)$ compatible with the respective $S$-connections, i.e., $\nabla_\omega$ and $\nabla^\mr{can}_{\mcF}$ (cf. (\ref{QQwkko})).
For each $h \in H (k)$,
we denote by  $\widetilde{\mr{R}}_h^F$  the automorphism of $\mcE_H \times G$ defined as $\widetilde{\mr{R}}_h^F := \mr{R}_h \times \mr{L}_{h^{-1}}$.
It  is $G$-equivariant with respect to the right $G$-action arising  from the second factor $G$,  and the following square diagram is commutative: 
\begin{align} \label{E0104}
\vcenter{\xymatrix@C=46pt@R=36pt{
\mcE_H \times G \ar[r]^-{\widetilde{\mr{R}}_h^F}_-{\sim} \ar[d]_{\mr{pr}_1}& \mcE_H \times G \ar[d]^{\mr{pr}_1} \\
\mcE_H \ar[r]^-{\sim}_-{\mr{R}_h}& \mcE_H.
 }}
\end{align}
When  we regard, via $\widetilde{\mr{R}}_h^F$,   the base-change   $\mr{R}_h^*(\nabla_\omega)$ of $\nabla_\omega$ as
an $S$-connection on $\mcE_H \times_k G$ (i.e., the upper left corner of (\ref{E0104})), 
the equality   $\mr{R}_h^*(\nabla_\omega) = \nabla_\omega$  holds because of   the commutativity of  (\ref{E0050}) (cf. ~\cite[Chap.\,1, \S\,1.3 Lemma 1.21]{Wak8}).
Hence,  it follows from Corollary \ref{EQkji} again that
$\widetilde{\mr{R}}_h^F$ 
 induces a $G$-equivariant  automorphism $\widetilde{\mr{R}}_{h}$ of    $\mcF$ whose 
pull-back by  $F_{\mcE_H/S}$ coincides with $\widetilde{\mr{R}}_h^F$ under the identification $\mcE_H \times G = F^*_{\mcE_H/S}(\mcF)$ given by  $\widetilde{\alpha}$ and which makes the following square diagram  commute:
\begin{align} \label{E01087}
\vcenter{\xymatrix@C=46pt@R=36pt{
\mcF  \ar[r]^-{\widetilde{\mr{R}}_{h}}_-{\sim} \ar[d]_{\mr{projection}}& \mcF \ar[d]^{\mr{projection}} \\
\mcE_H^{(1)} \ar[r]^-{\sim}_-{\mr{R}_h^{(1)}}& \mcE_H^{(1)},
 }}
\end{align}
where $\mr{R}_h^{(1)}$ denotes the base-change of $\mr{R}_h$ by  $F_S$.
The quotient $\mcG := \mcF/H$ of $\mcF$ by the resulting $H$-action $\{ \widetilde{\mr{R}}_h \}_{h \in H}$ forms a $G$-bundle on $X^{(1)} \left(= \mcE^{(1)}_H/H \right)$, and $\widetilde{\alpha}$ induces  an isomorphism  $\alpha : \left((\mcE_H \times G)/H  =\right) \mcE_G \isom F^*_{X/S}(\mcG)$.

Now, let us take an \'{e}tale $X$-scheme $U$ and a trivialization $\beta : \mcG |_{U^{(1)}} \isom U^{(1)} \times G$ of the $G$-bundle $\mcG$ over $U^{(1)}$.
The pull-back  
\begin{align}
\beta^F : F^*_{U/S}(\mcG |_U) \isom \left(F^*_{U/S} (U^{(1)} \times G) =  \right) U \times G
\end{align}
 of this trivialization   by $F_{U/S}$  gives the $H$-equivariant  composite
\begin{align}
h: \mcE_H |_U \xrightarrow{\mr{inclusion}}\mcE_G |_U \xrightarrow{\alpha |_U} F^*_{U/S}(\mcG)  |_U \xrightarrow{\beta^F} U \times G \xrightarrow{\mr{pr}_2} G.
\end{align}
In what follows, we shall prove the claim that
the morphism $h_S : \mcE_H |_U \migi \mcE_H^{\mr{triv}} \left(:= S \times G\right)$ induced naturally by $h$ is \'{e}tale.
Let us consider the commutative square diagram 
\begin{align} \label{f7j}
\vcenter{\xymatrix@C=46pt@R=36pt{
\mcF \ar[r] \ar[d]_-{} & \mcE_H^{(1)}  \ar[d]^-{}
\\
\mcG \ar[r] & X^{(1)}
 }}
\end{align}
consisting of natural projections.
This diagram is cartesian because  the two horizontal (resp., vertical) arrows are the structure morphisms of $G$-bundles (resp., $H$-bundles).
In particular, we obtain an  isomorphism $\mcF |_{U^{(1)}} \isom \mcE_H^{(1)} |_{U^{(1)}} \times_{U^{(1)}} \mcG |_{U^{(1)}}$.
The base-change of $\beta$ by the projection $\mcE_H^{(1)}|_{U^{(1)}} \migi U^{(1)}$ defines a $G$-equivariant isomorphism $\widetilde{\beta} : \mcF |_{U^{(1)}} \isom \mcE^{(1)}_H |_{U^{(1)}} \times G$, making  the following diagram commute:
\begin{align} \label{f7}
\vcenter{\xymatrix@C=46pt@R=36pt{
\mcF |_{U^{(1)}}\ar[r]^-{\widetilde{\beta}}_-{\sim} \ar[d]_-{\mr{projection}} & \mcE_H^{(1)} |_{U^{(1)}} \times G\ar[d]^-{\mr{projection}}
\\
\mcG |_{U^{(1)}}\ar[r]_-{\beta}^-{\sim} &  U^{(1)} \times G.
 }}
\end{align}
By pulling-back this diagram by  relative Frobenius morphisms,
we obtain a commutative square diagram 
\begin{align} \label{f7g}
\vcenter{\xymatrix@C=46pt@R=36pt{
F^*_{\mcE_H /S}(\mcF) |_{U}\ar[r]^-{\widetilde{\beta}^F}_-{\sim} \ar[d]_-{\mr{projection}} & \mcE_H |_{U} \times G\ar[d]^-{\mr{projection}}
\\
F^*_{X/S}(\mcG) |_{U}\ar[r]_-{\beta^F} &  U \times G.
 }}
\end{align}
This square extends to   
 the following commutative diagram:
\begin{align} \label{E06747}
\vcenter{\xymatrix@C=46pt@R=36pt{
\mcE_H |_U \ar[r]^-{v \mapsto (v, e)} \ar[d]^-{\mr{id}} & \mcE_H |_U \times G \ar[r]_-{\sim}^-{\widetilde{\alpha}|_U} \ar[d] &  F^*_{\mcE_H/S} (\mcF)|_U \ar[r]_-{\sim}^{\widetilde{\beta}^F} \ar[d]^-{\mr{projection}} & \mcE_H |_U \times G  \ar[r]^-{\mr{pr}_2} \ar[d]^-{\mr{projection}} & G \ar[d]^-{\mr{id}}
\\
\mcE_H |_U \ar[r]_-{\mr{inclusion}} & \mcE_G |_U \ar[r]^-{\sim}_-{\alpha |_U} & F^*_{X/S}(\mcG)|_U \ar[r]^-{\sim}_-{\beta^F} &  U \times G \ar[r]_-{\mr{pr}_2} & G, 
 }}
\end{align}
where 
 the second vertical arrow from the left denotes   the natural projection   $\mcE_H |_U \times G \migisurj \left((\mcE_H |_U \times G)/H = \right) \mcE_G |_U$.
 The commutativity of this diagram implies that the composite of all the arrows in the upper horizontal sequence coincides with $h$.
In particular, $h$ is a morphism resulting from the ``only if" part of the equivalence asserted in Proposition \ref{ETOFW2}, (i)  (with $X$ replaced by $\mcE_H$),  and hence, satisfies $h^*(\omega_G) = \omega$.
That is to say,  the differential $d h : \mcT_{\mcE_H|_U/S} \migi h_S^*(\mcT_{\mcE_H^{\mr{triv}}/S})$
 of $h$  makes the diagram
\begin{align} \label{}
\vcenter{\xymatrix@C=46pt@R=36pt{
\mcT_{\mcE_H|_U/S}  \ar[rr]^-{dh} \ar[rd]^-{\sim}_-{\omega^\triangleright |_{U}}&& h_S^*(\mcT_{\mcE_H^{\mr{triv}}/S}) \left(=  h^*(\mcT_{G}) \right) \ar[ld]_-{\sim}^-{h^*(\omega^\triangleright_G)}
\\
&\mcO_{\mcE_H |_U} \otimes \mfg &
}}
\end{align}
commutes, so 
  $d h$ is an isomorphism.
This implies  that, since both $\mcE_H |_U$ and $\mcE^{\mr{triv}}_H$ are smooth over $S$,
the morphism $h_S :  \mcE_H |_U \migi \mcE_H^{\mr{triv}}$  turns out to be  \'{e}tale.
This completes the proof of the claim.

Moreover, 
 the \'{e}tale morphism $h_S$ induces, via   taking the quotients by the  $H$-actions, an \'{e}tale morphism $\phi : U \migi P_S$ over $S$, i.e., a section of $\mcP^{\text{\'{e}t}}_{P, X/S}$ over $U$.
The  collection of  $\phi$'s defined for various  pairs $(U, \beta)$ forms   a subsheaf  $\mcE^{\cc \Rightarrow \hh}$ of $\mcP_{P, X/S}^{\text{\'{e}t}}$,  which specifies   an   $F$-$(G, P)$-structure on $X/S$.
One may verify that the resulting assignment $\mcE^\cc \mapsto \mcE^{\cc\Rightarrow \hh}$ defines  the inverse to $\zeta^{\hh \Rightarrow \cc}$.
This completes the proof of this theorem. 
\end{proof}
%-------------------------------------------------------------

\LSP
%---------------------------[begin subsection]-------------
\subsection{$F^\N$-Cartan geometries for $\N < \infty$} \label{SS096}
We shall extend both  the notion of a Cartan geometry and the bijection $\zeta^{\hh \Rightarrow \cc}$ to  finite level.

In the rest of this section, we suppose that $G$ is affine.
Let $\N$ be a positive integer and 
$\mcG$  a $G$-bundle on $X^{(\N)}$.
In what follows,  we shall  construct a canonical $\mfg$-valued $1$-form on  the pull-back $F^{(\N)*}_{X/S}(\mcG)$.
Denote by $\pi$ the structure morphism of $F^{(\N)*}_{X/S}(\mcG)$ over $X$ and 
 by $F_{\mcG/X/S}^{(\N)}: F^{(\N)*}_{X/S}(\mcG) \migi \mcG$ the natural projection.
The differential 
\begin{align}
d F_{\mcG/X/S}^{(\N)} : \mcT_{F^{(\N)*}_{X/S}(\mcG)/S} \migi F_{\mcG/X/S}^{(\N)*}(\mcT_{\mcG/S})
\end{align}
  of $F_{\mcG/X/S}^{(\N)}$
fits into the following morphism of  short exact sequences:
 \begin{align} \label{E0157}
\vcenter{\xymatrix@C=26pt@R=36pt{
0 \ar[r] & \mcT_{F^{(\N)*}_{X/S}(\mcG)/X}\ar[r] \ar[d]^-{\wr} & \mcT_{F^{(\N)*}_{X/S}(\mcG)/S} \ar[r]^{d F_{\mcG/X/S}^{(\N)}} \ar[d]^-{d F_{\mcG/X/S}^{(\N)}} &  \pi^*(\mcT_{X/S})\ar[r] \ar[d] & 0
\\
0 \ar[r] & F_{\mcG/X/S}^{(\N)*}(\mcT_{\mcG/X^{(\N)}}) \ar[r]  & F_{\mcG/X/S}^{(\N)*}(\mcT_{\mcG/S})\ar[r]  & \pi^*(F^{(\N)*}_{X/S}(\mcT_{X^{(\N)}/S}))\ar[r] & 0,
 }}
\end{align}
where 
\begin{itemize}
\item
the upper horizontal sequence is obtained by differentiating $\pi$, and 
the lower horizontal sequence is obtained by differentiating 
the natural projection $\mcG \migi X^{(\N)}$ and successively pulling-back by  $F_{\mcG/X/S}^{(\N)}$;
\item
the right-hand  vertical arrow  denotes the pull-back by $\pi$ of the differential  of $F_{X/S}^{(\N)}$, and the left-hand vertical  arrow denotes the isomorphism  arising  from  the base-change by $F_{X/S}^{(\N)}$.
\end{itemize}
Since  the right-hand vertical arrow becomes the zero map  because of the definition of $F_{X/S}^{(\N)}$, 
  the morphism $d F_{\mcG/X/S}^{(\N)}$ factors through
the inclusion $F^{(\N)*}_{\mcG/X/S}(\mcT_{\mcG/X^{(\N)}}) \migiincl  F^{(\N)*}_{\mcG/X/S}(\mcT_{\mcG/S})$.
The resulting morphism $\mcT_{F^{(\N)*}_{X/S}(\mcG)/S} \migi F^{(\N)*}_{\mcG/X/S}(\mcT_{\mcG/X})$
induces the  $\mcO_{F^{(\N)*}_{X/S}(\mcG)}$-linear composite
 \begin{align} \label{39adids}
 \omega^{\mr{can}(\N), \triangleright}_\mcG :  \mcT_{F^{(\N)*}_{X/S}(\mcG)/S} \migi 
 F^{(\N)*}_{\mcG/X/S}(\mcT_{\mcG/X}) \isom  
 \mcT_{F^{(\N)*}_{X/S}(\mcG)/X} 
 \xrightarrow{\omega^\triangleright_{F^{(\N)*}_{X/S}(\mcG)/X}}
 \mcO_{F^{(\N)*}_{X/S}(\mcG)} \otimes \mfg,
 \end{align}
 where the second arrow denotes the inverse of the left-hand vertical arrow in (\ref{E0157}).
Hence, 
we obtain the $\mfg$-valued $1$-form
\begin{align} \label{Ewr342}
\omega^{\mr{can}(\N)}_\mcG \in \Gamma (F^{(\N)*}_{X/S}(\mcG), \Omega_{F^{(\N)*}_{X/S}(\mcG)/S} \otimes \mfg)
\end{align}
on $F^{(\N)*}_{X/S}(\mcG)/S$ corresponding to $\omega^{\mr{can}(\N), \triangleright}_\mcG$.
If $\N'$ is a positive integer with $\N' \leq  \N$, then  the equality  $\omega_\mcG^{\mr{can}(\N)} = \omega_{\mcG'}^{\mr{can}(\N')}$ holds, where $\mcG' := F^{(\N-\N')*}_{X^{(\N')}/S}(\mcG)$, under the natural identification $F^{(\N)*}_{X/S} (\mcG) = F^{(\N')*}_{X/S}(\mcG')$.

\SSP
%------------------------------------------
\ble \label{ea04394F}
The $\mfg$-valued $1$-form $\omega^{\mr{can}(\N)}_\mcG$ satisfies
$\psi_{\omega^{\mr{can}(\N)}_\mcG } =0$ and ${^p}\psi_{\omega^{\mr{can}(\N)}_\mcG } =0$.
\ele
%------------------------------------------
\begin{proof}
The problem is of local nature, 
so  it suffices to consider the case where $\mcG = X^{(\N)} \times G$.
In this case,  one may verify  that 
the $S$-connection $\nabla_{\omega^{\mr{can}(\N)}_\mcG }$ on $\mcO_{X \times G}$ (cf. (\ref{L048})) coincides with the trivial connection.
Hence, 
since both the curvature and $p$-curvature of $\nabla_{\omega^{\mr{can}(\N)}_\mcG }$ vanish identically, 
(\ref{L049}) implies  $\psi_{\omega^{\mr{can}(\N)}_\mcG }  = 0$  and Lemma \ref{L0103} implies ${^p}\psi_{\omega^\mr{can}_\mcG } =  0$.
\end{proof}
%------------------------------------------
\SSP

%------------------------------------------
\bde \label{Eotarmr90}
\begin{itemize}
\item[(i)]
 An {\bf  $F^\N$-Cartan geometry with model $(G, H)$}  on $X/S$ is a quadruple
 \begin{align}
 \mcE^\cc_\N := (\mcE_H, \omega, \mcG, \kappa)
 \end{align}
 consisting of
 a Cartan geometry $(\mcE_H, \omega)$ with model $(G, H)$ on $X/S$, a $G$-bundle $\mcG$ on $X^{(\N)}$, and an isomorphism  of $G$-bundles $\kappa : \mcE_G \left(:= \mcE_H \times^H G\right) \isom F_{X/S}^{(\N)*}(\mcG)$ with $\kappa^*(\omega^{\mr{can} (\N)}_\mcG) |_{\mcE_H} = \omega$. 
 \item[(ii)]
 Let $\mcE^\cc_\N := (\mcE_H, \omega, \mcG, \kappa)$ and $\mcE'^\cc_\N := (\mcE'_H, \omega', \mcG', \kappa')$ be $F^\N$-Cartan geometries with model $(G, H)$ on $X/S$.
 An {\bf isomorphism of $F^\N$-Cartan geometries with model $(G, H)$} from $\mcE^\cc_\N$ 
  to $\mcE'^\cc_\N$ is 
  a pair 
  \begin{align}
  (\eta_\mcE, \eta_\mcG)
  \end{align}
   consisting of an isomorphism $\eta_\mcE : (\mcE_H, \omega) \isom (\mcE'_H, \omega')$  of Cartan geometries with model $(G, H)$  
and   
   an isomorphism of $G$-bundles $\eta_\mcG : \mcG \isom \mcG'$ satisfying 
   $F_{X/S}^{(\N)*}(\eta_\mcG) \circ \kappa = \kappa' \circ (\eta_\mcE  \times^H G)$, where $\eta_\mcE  \times^H G$ denotes the isomorphism $\mcE_G \isom \mcE'_G$ induced by $\eta_\mcE$.
 \end{itemize}
\ede
%------------------------------------------
\SSP

Just as in the case of a Cartan geometry (cf. \S\,\ref{SS02}), we can define  the pull-back and the base-change  of an $F^\N$-Cartan geometry (but leave the details to the reader).

We shall denote by
\begin{align} \label{40su59}
F^\N\text{-}\mr{Car}_{(G, H), X/S}
\end{align}
the set of isomorphism classes of $F^\N$-Cartan geometries with model $(G, H)$ on $X/S$.

Let $\N$ and $\N'$ be positive integers with $\N \leq \N'$ and $\mcE^\cc_{\N'} := (\mcE_H, \omega, \mcG, \kappa)$ an $F^{\N'}$-Cartan geometry with model $(G, H)$ on $X/S$.
Then, under the canonical identification $F^{(\N')*}_{X/S} (\mcG) = F^{(\N)*}_{X/S}(F^{(\N' - \N)*}_{X^{(\N)}/S}(\mcG))$,
the quadruple
 \begin{align}
\mcE^\heartsuit_{\N'} |^{\langle \N\rangle} :=  (\mcE_H, \omega, F^{(\N' - \N)*}_{X^{(\N)}/S}(\mcG), \kappa)
\end{align}
specifies an $F^{(\N)}$-Cartan geometry with model $(G, H)$ on $X/S$.
We refer to $\mcE^\cc_{\N'} |^{\langle \N\rangle}$ as the {\bf $\N$-th truncation}   of $\mcE^\cc_{\N'}$.
The resulting assignments $\mcE^\cc_{\N'} \mapsto \mcE^\cc_{\N'} |^{\langle \N\rangle}$ 
for various pairs $(\N, \N')$ gives a projective system
\begin{align} \label{R4567}
\cdots \migi F^\N\text{-}\mr{Car}_{(G, H), X/S} \migi \cdots \migi 
F^2\text{-}\mr{Car}_{(G, H), X/S} \migi F^1\text{-}\mr{Car}_{(G, H), X/S}.
\end{align}

\SSP
%------------------------------------------
\begin{exa}[Standard example] \label{Efte123}
We shall consider  the standard example of an $F^\N$-$(G, H)$-Cartan geometry as an extension of  (\ref{L045}).
The identity morphism  $\mr{id}_{P \times G}$ of $P \times G$ may be considered   as an isomorphism $P \times G \isom F^{(\N)*}_{P_S/S} (P^{(\N)}_S \times G)$.
Then, 
the quadruple
\begin{align} \label{4598sj9}
\mcE_\N^{\mr{triv}, \cc} :=(\mcE_H^\mr{triv}, \omega_{S \times G/S}, P^{(N)}_S \times G, \mr{id}_{P \times G})
\end{align}
forms an $F^\N$-Cartan geometry with model $(G, H)$ on $P_S/S$.
To verify this, we shall  prove the equality 
\begin{align} \label{GH01}
\omega_{P^{(N)}_S \times G}^{\mr{can}(\N)} |_{\mcE_H^\mr{triv}} = \omega_{S \times G/S}.
\end{align}
%Recall from (\ref{QQ13}) that $\mcT_{P_S \times G/S}$ decomposes as $\mr{pr}_1^*(\mcT_{P_S/S}) \oplus \mr{pr}_2^*(\mcT_{G})$ ($= \mr{pr}_1^*(\mcT_{P_S/S}) \oplus  \mcO_{P_S \times G}\otimes \mfg$ via $\omega_G^\triangleright$).
Under the identification $\mcT_{P_S \times G/S} = \mr{pr}_1^*(\mcT_{P_S/S}) \oplus  \mcO_{P_S \times G}\otimes \mfg \left(=\mr{pr}_1^*(\mcT_{P_S/S}) \oplus \mr{pr}_2^*(\mcT_{G}) \right)$ given by (\ref{QQ13}),
the $\mcO_{P_S \times G}$-linear morphism $(\omega_{P^{(N)}_S \times G}^{\mr{can}(\N)})^\triangleright : \mcT_{P_S \times G/S} \migi \mcO_{P_S \times G} \otimes \mfg$ coincides with the second projection  $\mr{pr}_1^*(\mcT_{P_S/S}) \oplus  \mcO_{P_S \times G}\otimes \mfg \migisurj \mcO_{P_S \times G}\otimes \mfg$.
Hence, since $(\omega_{P^{(N)}_S \times G}^{\mr{can}(\N)} |_{\mcE_H^\mr{triv}})^\triangleright$ may be obtained  as
the restriction of $(\omega_{P^{(N)}_S \times G}^{\mr{can}(\N)})^\triangleright$ via  $(\mr{id}_S \times \pi_G, \mr{pr}_2): \mcE^\mr{triv}_H  \migiincl  P_S \times G$ (cf. (\ref{ghfjke})),
we see that it coincides with $\omega_{S \times G/S}^\triangleright$.
This implies (\ref{GH01}),
 as desired.
\end{exa}
%------------------------------------------
\SSP

%------------------------------------------
\bpr \label{P0200}
Let $\mcE_1^\cc := (\mcE_H, \omega, \mcG, \kappa)$ be an $F^1$-Cartan geometry with model $(G, H)$ on $X/S$.
Then, the underlying Cartan geometry  $(\mcE_H, \omega)$  is $p$-flat.
Moreover,  
the map of sets
\begin{align} \label{Gh045}
F^{1}\text{-}\mr{Car}_{(G, H), X/S} \migi  \mr{Car}^{p\text{-}\mr{flat}}_{(G, H), X/S}.
\end{align}
given by  $\mcE_1^\cc \mapsto (\mcE_H, \omega)$ is bijective.
\epr
%------------------------------------------
\begin{proof}
Let $\mcE_1^\cc$ be as in the statement.
Since the $p$-curvature of $\omega_{\mcG}^{\mr{can}(\N)}$ vanishes (cf. Lemma \ref{ea04394F}), 
the equality  $\kappa^*(\omega_{\mcG}^{\mr{can}(\N)}) |_{\mcE_H} = \omega$ implies  ${^p}\psi_\omega =0$.
Hence, $(\mcE_H, \omega)$ turns out to be $p$-flat.

Conversely,
let $\mcE^\cc := (\mcE_H, \omega)$ be a $p$-flat Cartan geometry with model $(G, H)$ on $X/S$.
Just as in the proof of Theorem \ref{P0011},   the data  $(\mcE_H, \omega)$ yields   a $G$-bundle $\mcG$ on $X^{(1)}$   and  an isomorphism of $G$-bundles $\alpha : \mcE_G \isom F^*_{X/S}(\mcG)$.
By  the construction of $\mcG$ together with  
 (\ref{GH01}), the equality  
$\alpha^*(\omega^{\mr{can}(1)}_\mcG) |_{\mcE_H}= \omega$ holds.
That is to say, the quadruple $(\mcE_H, \omega, \mcG, \alpha)$ forms an $F^1$-Cartan geometry.
One may verify that the resulting assignment $\mcE^\cc \mapsto (\mcE_H, \omega, \mcG, \alpha)$ determines  the inverse to (\ref{Gh045}).
This completes the proof of this proposition.
\end{proof}
%------------------------------------------

\SSP

%------------------------------------------
\bpr \label{P0201}
(Recall that $G$ is assumed to be affine.)
There exists  a canonical bijection
\begin{align} \label{Prt78}
\zeta_\N^{\hh \Rightarrow \cc} : F^{\N}\text{-}\mr{Ehr}_{(G, P), X/S} \isom F^{\N}\text{-}\mr{Car}_{(G, H), X/S},
\end{align}
and the formation of this bijection  commutes  with truncation to lower levels, pull-back to \'{e}tale $X$-schemes, and base-change to $S$-schemes.
If $\N=1$, then  this bijection coincides with (\ref{E0223}) under the identification $F^1\text{-}\mr{Car}_{(G, H), X/S} = \mr{Car}^{p\text{-}\mr{flat}}_{(G, H), X/S}$ given by  (\ref{Gh045}).
\epr
%------------------------------------------
\begin{proof}
First, we shall construct $\zeta_\N^{\hh \Rightarrow \cc}$.
Let $\mcS^\hh$ be an $F^{\N}$-$(G, P)$-structure  on $X/S$.
The $1$-st truncation $\mcS^\hh |^{\langle 1 \rangle}$ of $\mcS^\hh$ corresponds 
to
an $F^1$-Cartan geometry $(\mcE_H, \omega, \mcG, \kappa)$  with model $(G, H)$ via  $\zeta^{\hh \Rightarrow \cc}$ (cf.  (\ref{E0223})) and (\ref{Gh045}).
Let us set $\mcG_\N$ to be the $G$-bundle on $X^{(\N)}$ defined as the inverse of the left $G_{X/S}^{(\N)}$-torsor $\mcS^\hh$.
Since $\mcG \cong F^{(\N-1)*}_{X^{(1)}/S}(\mcG_\N)$ by construction,
  $\kappa$ may be regarded  as an isomorphism $\mcE_G \isom F^{(\N)*}_{X/S}(\mcG_\N) \left(=F^{*}_{X/S}(F^{(\N-1)*}_{X^{(1)}/S}(\mcG_\N)) \right)$.
Then, 
the quadruple $\mcS^{\hh \Rightarrow \cc}_\N := (\mcE_H, \omega, \mcG_\N, \kappa)$ specifies  an $F^\N$-Cartan geometry with model $(G, H)$ on $X/S$.
The resulting assignment $\mcS^\hh \mapsto \mcS^{\hh \Rightarrow \cc}_\N$ defines a map of sets 
$F^{\N}\text{-}\mr{Ehr}_{(G, P), X/S} \migi  F^{\N}\text{-}\mr{Car}_{(G, H), X/S}$, which we denote by $\zeta_\N^{\hh \Rightarrow \cc}$.

Next, we shall construct the inverse to this map.
Let  $\mcE^\cc_\N := (\mcE_H, \omega, \mcG, \kappa)$ be an $F^\N$-Cartan geometry
belonging to $ F^{\N}\text{-}\mr{Car}_{(G, H), X/S}$.
The $1$-st truncation $\mcE^\cc_\N |^{\langle 1 \rangle}$ of $\mcE_\N^\cc$ corresponds to an $F^1$-$(G, P)$-structure $\mcE^{\cc \Rightarrow \hh}_1$  via $\zeta^{\hh \Rightarrow \cc}$ and (\ref{Gh045}).
The $G$-bundle $\mcG$  
  induces  a left  $G_{X/S}^{(\N)}$-torsor $\breve{\mcE}^{\cc \Rightarrow \hh}_\N$ equipped with  an isomorphism   $\alpha_\mcS : \breve{\mcE}^{\cc \Rightarrow \hh}_\N \times^{G_{X/S}^{(\N)}} G_{X/S}^{(1)} \isom \mcE^{\cc \Rightarrow \hh}_1$.
Then,  the subsheaf $\mcE^{\cc \Rightarrow \hh}_\N$ of $\mcP_{P, X/S}^{\text{\'{e}t}}$ consisting of
the local sections of the form $\alpha_\mcS ((v, e))$ for various  $v \in \breve{\mcE}^{\cc \Rightarrow \hh}_\N$
   determines  an $F^\N$-$(G, P)$-structure on $X/S$.
The resulting assignment $\mcE^\cc_\N \mapsto  \mcE^{\cc \Rightarrow \hh}_\N$ forms the inverse to the map $\zeta_\N^{\hh \Rightarrow \cc}$, so 
this completes the proof of this proposition.
\end{proof}
%------------------------------------------

\LSP
%---------------------------[begin subsection]-------------
\subsection{$F^\N$-Cartan geometries for $\N = \infty$} \label{SSg096}

In this final subsection of \S\,\ref{EosGG}, we  generalize $F^\N$-Cartan geometries  to $\N = \infty$ in terms of $F$-divided $G$-bundles.

\SSP
%------------------------------------------
\bde \label{eia2FG}
An {\bf $F^\infty$-Cartan geometry with model $(G, H)$} on $X/S$ is a collection
\begin{align}
\mcE^\cc_\infty := (\mcE_H, \omega,  \{ (\mcG_l, \varpi_l ) \}_{l \in \mbZ_{>0}}, \kappa)
\end{align}
such that $(\mcE_H, \omega, \mcG_1, \kappa)$ forms an $F^1$-Cartan geometry with $(G, H)$  on $X/S$ and the collection $ \{ (\mcG_l, \varpi_l ) \}_{l \in \mbZ_{\geq 0}}$, where  $(\mcG_0, \varpi_0) := (\mcE_G, \kappa)$,  forms an $F$-divided $G$-bundle on $X/S$ (cf. Definition \ref{EIFO(79}).
The notion of an isomorphism between two $F^\infty$-Cartan geometries  can be defined in an evident manner.
\ede
%------------------------------------------
\SSP

Denote by
\begin{align}
F^\infty\text{-}\mr{Car}_{(G, H), X/S}
\end{align}
the set of isomorphism classes of $F^\infty$-Cartan geometries with model $(G, H)$ on $X/S$.
Let $\N$ be a positive integer and 
 $\mcE^\heartsuit_\infty := (\mcE_H, \omega,  \{ (\mcG_l, \varpi_l ) \}_l, \kappa)$  an $F^\infty$-Cartan geometry with model $(G, H)$ on $X/S$.
Then, $\mcE^\heartsuit_\N := (\mcE_H, \omega, \mcG_\N, \kappa_\N)$ forms an $F^\N$-Cartan geometry, where  $\kappa_\N := F^{(\N -1)*}_{X/S}(\varpi_{\N-1})\circ\cdots\circ F^{(1)*}_{X/S}(\varpi_1)\circ \kappa$.
The resulting assignment $\mcE^\heartsuit_\infty \mapsto \mcE^\heartsuit_\N$ defines a map of sets $F^\infty\text{-}\mr{Car}_{(G, H), X/S} \migi F^\N\text{-}\mr{Car}_{(G, H), X/S}$ compatible with  the morphisms in (\ref{R4567}).
Hence, we obtain a map of sets
\begin{align}\label{Ert2}
F^\infty\text{-}\mr{Car}_{(G, H), X/S} \migi \varprojlim_{\N \in \mbZ_{>0}}F^\N\text{-}\mr{Car}_{(G, H), X/S}.
\end{align}

\SSP
%-----------------------------------------------
\bpr \label{Efget} 
Suppose that $G$ is affine, $S = \mr{Spec}(k)$,  and  $X$ is proper over $k$.
Then, (\ref{Ert2}) is bijective.
In particular, we obtain a bijection of sets
\begin{align} \label{E(eoe9}
\zeta_\infty^{\hh \Rightarrow \cc} : F^\infty\text{-}\mr{Ehr}_{(G, P), X/S} \isom F^\infty\text{-}\mr{Car}_{(G, H), X/S}.
\end{align}
\epr
%-----------------------------------------------
\begin{proof}
The former assertion  is a direct consequence of  Lemma \ref{epGJJ} described below.
The latter assertion follows from the former assertion, (\ref{rosao3048}),  and  Proposition \ref{P0201}.
\end{proof}
\SSP
%-----------------------------------------------

The following lemma was already proved in  ~\cite[\S\,1, Proposition 1.7]{G}  when $G = \mr{GL}_n$.

\SSP

%-----------------------------------------------
\ble \label{epGJJ}
Let us keep the assumption in  Proposition \ref{Efget}.
For each nonnegative integer $l$,
denote by $G\text{-}\mr{Bun}_{X^{(l)}}$ the set of isomorphism classes of $G$-bundles on $X^{(l)}$.
Note that pulling-back $G$-bundles by relative Frobenius morphisms yields a projective system
\begin{align}
\cdots \migi G\text{-}\mr{Bun}_{X^{(l)}} \migi \cdots \migi G\text{-}\mr{Bun}_{X^{(2)}} \migi G\text{-}\mr{Bun}_{X^{(1)}} \migi G\text{-}\mr{Bun}_{X^{}}. 
\end{align}
Also, denote by $G\text{-}\mr{Bun}_{X^{(\infty)}}$ the set of isomorphism classes of $F$-divided $G$-bundles on $X/k$.
Then, the map of sets
\begin{align} \label{EOTIME)}
G\text{-}\mr{Bun}_{X^{(\infty)}} \migi \varprojlim_{l \in \mbZ_{\geq 0}}G\text{-}\mr{Bun}_{X^{(l)}}
\end{align}
given by  $\{ (\mcG_l, \varpi_l) \}_l \mapsto \{ \mcG_l \}_{l}$ is bijective.
\ele
%-----------------------------------------------
\begin{proof}
Since the surjectivity of (\ref{EOTIME)}) is clear, it suffices to consider the injectivity.
Let 
$\{ (\mcG_{i, l}, \varpi_{i, l})\}_l$ ($i=1,2$) be  $F$-divided $G$-bundles on $X/k$ with  $\mcG_{1, l} \cong \mcG_{2, l}$ for every $l$.
For each $l \in \mbZ_{\geq 0}$,  denote by $\mr{Aut} (\mcG_{2, l})$ 
the algebraic $k$-group classifying $G$-equivariant  automorphisms of   $\mcG_{2, l}$.
Notice that 
$\mr{Aut} (\mcG_{2, l})$ 
 is  of finite  type over $k$.
Indeed, let us  fix a closed immersion between algebraic $k$-groups  $G \migiincl \mr{GL}_m$ ($m \geq 1$).
(Such a morphism always exists because of the affineness assumption on $G$.)
  Write $\mcV_l$ for the vector bundle induced from  
$\mcG_{2, l}$ via change of structure group by this closed immersion.
By the properness of $X/k$,  the space  $\mr{End} (\mcV_l)$ of  endomorphisms of $\mcV_l$ may be represented by a finite type scheme over $k$.  Hence,   
$\mr{Aut} (\mcG_{2, l})$ 
turns out to be  of finite type  because it forms a subscheme of $\mr{End} (\mcV_l)$.
Now, let us consider 
the map of sets
$F_{X^{(l)}/k}^\sharp  : \mr{Aut} (\mcG_{2, l+1})(k) \migi \mr{Aut} (\mcG_{2, l})(k)$  arising from pull-back by $F_{X^{(l)}}$.
Since $X$ is reduced,  we may identify $\mcV_{l+1}$ as a subsheaf of $\mcV_{l}$, which implies
 that  $F_{X^{(l)}/k}^\sharp$ is  injective.
Moreover,  since $\mr{Aut} (\mcG_{2, l})$'s are  of finite type over $k$ as proved above,   there exists a positive  integer $L$ such that 
 $F_{X^{(l)}/k}^\sharp$ is bijective   for any $l \geq L$.
If we take  an isomorphism  of
$G$-bundles $\eta_l : \mcG_{1, l} \isom \mcG_{2, l}$ with $l \geq L$,
then
the bijectivity of $F_{X^{(l)}/k}^\sharp$  implies  that 
there exists an isomorphism $\eta_{l +1} : \mcG_{1, l+1} \isom \mcG_{2, l+1}$ with $\varpi_{2, l}\circ F_{X^{(l)}/k}^* (\eta_{l +1}) = \eta_l \circ \varpi_{1, l}$.
By starting from the choice of $\eta_L$ and  applying this argument inductively,  we obtain  a compatible collection of isomorphisms $\{ \eta_l \}_l$, specifying an isomorphism of $F$-divided $G$-bundles   $\{ (\mcG_{1, l}, \varpi_{1, l})\}_l \isom \{ (\mcG_{2, l}, \varpi_{2, l})\}_l$.
This concludes   the injectivity of (\ref{EOTIME)}), and hence, completes the proof of this lemma.
\end{proof}
%----------------------------------------------

%%%%%%%---[begin section]---%%%%%%%%%%%%%%
\vspace{10mm}
\section{Indigenous bundles on smooth schemes} \label{S0029}\vspace{3mm}

This section deals with   indigenous bundles (including the case of finite level) on a smooth scheme, generalizing the classical notion defined on a Riemann surface (cf.  ~\cite{G2}).
As one of the main results of this section,  we prove  that, for each level,  $p$-flat Cartan geometries correspond bijectively to dormant  indigenous bundles  (cf. Propositions \ref{P0234}, \ref{Eqw209}, and Corollary \ref{Er281}).
In particular,   the proofs of  Theorems \ref{TheoremA1} and \ref{TheoremA} will be completed in this section. 
Also,  the resulting correspondence enables us to prove  several properties on a smooth scheme having an $F^\infty$-$(G, P)$-structure (cf. Propositions  \ref{etoiuoa892}, \ref{etoiuoa893}).
The last two subsections will be devoted to studying  $F^\N$-$(G, P)$-structures in the case where $G$ is a product of $\mbG_m$'s and $\mbG_a$.

Let $S$ be a scheme over $k$ and  $X$  a  smooth scheme  over $S$.

\vspace{5mm}
%----------------------------------------------------------------------[begin subsection]-------------
\subsection{The Kodaira-Spencer map} \label{SS039}

 Let  $\mcE_H$ be an $H$-bundle on $X$, and 
 write $\mcE_G := \mcE_H \times^H G$.
The differential of  the natural inclusion $\mcE_H \migiincl \mcE_G$ induces
 an $\mcO_X$-linear injection
$\widetilde{\mcT}_{\mcE_H} \migiincl \widetilde{\mcT}_{\mcE_G}$ (cf. \S\,\ref{SS002} for the definition of $\widetilde{\mcT}_{(-)/S}$).
 We shall regard  $\widetilde{\mcT}_{\mcE_H/S}$ as an $\mcO_X$-submodule of $\widetilde{\mcT}_{\mcE_G/S}$ via this injection.

%-----------------------------------------------------------------------[begin definition]------------------
\vspace{3mm}
\bde \label{D0550}
Let $\nabla : \mcT_{X/S} \migiincl \widetilde{\mcT}_{\mcE_G/S}$ be an $S$-connection  on the $G$-bundle $\mcE_G$.
The $\mcO_X$-linear composite
 \begin{align} \label{E431}
\mr{KS}_{(\mcE_H, \nabla)} : \mcT_{X/S} \xrightarrow{\nabla} \widetilde{\mcT}_{\mcE_G/S} \migisurj \widetilde{\mcT}_{\mcE_G/S}/\widetilde{\mcT}_{\mcE_H/S}
\end{align}
is called the {\bf Kodaira-Spencer map} associated with $(\mcE_H, \nabla)$. 
    \ede
%-----------------------------------------------------------------------[end definition]-------------------
%\vspace{3mm}

%---------------------------------[begin remark]------------------
\begin{rema} \label{R13404}
 Let us keep the above notation.
 Differentiating  the inclusion $\mcE_H \migiincl  \mcE_G$ gives  the following morphism of short exact sequences:
 \begin{align} \label{E0ygg0}
\vcenter{\xymatrix@C=46pt@R=36pt{
0 \ar[r] & \mfh_{\mcE_H} \ar[r] \ar[d]^{\mr{inclusion}} & \widetilde{\mcT}_{\mcE_H /S} \ar[r]^{d_{\mcE_H}} \ar[d]^{\mr{inclusion}} & \mcT_{X/S} \ar[r] \ar[d]_{\wr}^{\mr{id}} & 0
\\
0  \ar[r] & \mfg_{\mcE_G} \ar[r] & \widetilde{\mcT}_{\mcE_G/S} \ar[r]_-{d_{\mcE_G}} & \mcT_{X/S} \ar[r]&0.
 }}
\end{align}
 This morphism   induces  an $\mcO_X$-linear  isomorphism  
 \begin{align} \label{Eoraplw2}
 \mfg_{\mcE_G}/\mfh_{\mcE_H} \isom \widetilde{\mcT}_{\mcE_G/S}/\widetilde{\mcT}_{\mcE_H/S}.
 \end{align}
By this isomorphism, 
%Hence, under the identification $\mfg_{\mcE_G}/\mfh_{\mcE_H} = \widetilde{\mcT}_{\mcE_G/S}/\widetilde{\mcT}_{\mcE_H/S}$ given by this isomorphism,
 the Kodaira-Spencer map $\mr{KS}_{(\mcE_H, \nabla)}$ may be regarded  as an $\mcO_X$-linear morphism $\mcT_{X/S} \migi \mfg_{\mcE_G}/\mfh_{\mcE_H}$.
 \end{rema}
%-----------------------------------[end remark]-------------------
%\vspace{3mm}

\SSP

\begin{exa}[The cases of $G = \mbG_m$ and $\mbG_a$] \label{EroDFG}
%Let us take a $1$-form $\omega$ on $X/S$, and denote by $\omega^\triangleright$
%the  $\mcO_X$-linear  morphism $\mcT_{X/S} \migi \mcO_X$ corresponding to $\omega$.
\begin{itemize}
\item[(i)]
First, we shall consider the case where  $G = \mbG_m$ and $H = \{ e\}$.
Let us take a $1$-form $\omega$ on $X/S$.
Recall that $\omega$ induces an $S$-connection $\nabla_\omega$  on the trivial $\mbG_m$-bundle $X \times \mbG_m$ (cf. (\ref{L048})).
Then,  the Kodaira-Spencer map $\mr{KS}_{(X \times \{e \}, \nabla_\omega)}$ associated with  the pair $(X \times \{ e \}, \nabla_\omega)$
 coincides with $\omega^\triangleright$ under the natural  identification 
 \begin{align} \label{QQEI2}
 \mcO_X \cong  \mr{Lie} (\mbG_m)_{X \times \mbG_m}/\mr{Lie}(\{e \})_{X \times \{ e \}}  \stackrel{(\ref{Eoraplw2})}{\cong}  \widetilde{\mcT}_{X \times \mbG_m/S}/\widetilde{\mcT}_{X \times \{ e \}/S}.
 \end{align}
Note that giving  a $\mbG_m$-bundle equipped with both  an $S$-connection and  an $\{ e \}$-reduction is  equivalent to giving an $S$-connection on the trivial $\mbG_m$-bundle, i.e., $\nabla_\omega$ for some $\omega$ as above.

\item[(ii)]
Next, suppose that $G = \mbG_a$ and $H = \{ e \}$.
We shall identify  $\mbG_a$ as the subgroup of $\mr{GL}_2$ consisting of matrices 
$\begin{pmatrix} 1 & a \\ 0 & 1\end{pmatrix}$ for some $a \in \mbG_a$.
Each $\mbG_a$-bundle may be thought of as a rank $2$ vector bundle $\mcV$ equipped with a structure of extension $0 \migi \mcO_X \migi \mcV \migi \mcO_X \migi 0$.
Moreover, a choice of an $\{e \}$-reduction  corresponds to a splitting $\mcV \cong \mcO_X^{\oplus 2}$ of this extension.
The vector bundle $\mcO_X^{\oplus 2}$  equipped with the  extension $0 \migi \mcO_X \xrightarrow{a \mapsto (a, 0)}\mcO_X^{\oplus 2} \xrightarrow{(a, b)\mapsto b} \mcO_X \migi 0$ corresponds to the trivial $\mbG_a$-bundle $X \times \mbG_a$.

For each $1$-form $\chi$ on $X/S$, 
 denote by 
\begin{align} \label{ZZAei2}
\chi^\dagger :\mcO_X^{\oplus 2} \migi \Omega_{X/S} \otimes_{\mcO_X} (\mcO_X^{\oplus 2}) \left(= \Omega_{X/S}^{\oplus 2} \right)
\end{align}
 the $\mcO_X$-linear morphism  defines as  $\begin{pmatrix} 0 & \chi \\ 0 & 0 \end{pmatrix}$.
Then,  the morphism 
\begin{align} \label{ZZAei}
\nabla_\chi^\dagger := d^{\oplus 2} +\chi^\dagger 
\end{align}
turns out to define  an $S$-connection on the trivial $\mbG_a$-bundle.
The  Kodaira-Spencer map $\mr{KS}_{(X \times \{ e \}, \nabla_\chi^\dagger)}$ associated with the pair $(X \times \{ e \}, \nabla^\dagger_\chi)$ coincides with
$\chi^\triangleright$ under the identification $\mcO_X \cong \widetilde{\mcT}_{X \times \mbG_a/S}/\widetilde{\mcT}_{X \times \{ e \} /S}$ similar to  (\ref{QQEI2}).  
Also, giving a $\mbG_a$-bundle  equipped with both an $S$-connection and an $\{e \}$-reduction is equivalent to giving an $S$-connection on $\mcO_X^{\oplus 2}$ of the form $\nabla_\chi^\dagger$ for some $\chi$ as above.
 \end{itemize}
\end{exa}

\SSP
%----------------------------------------------------------------------
\bde \label{D0013}
\begin{itemize}
\item[(i)]
An {\bf indigenous $(G, H)$-bundle} on $X/S$
is a pair
\begin{align}
\mcE^\sss := (\mcE_H, \nabla)
\end{align}
consisting of an $H$-bundle $\mcE_H$ on  $X$ and an $S$-connection $\nabla$ on $\mcE_G := \mcE_H \times^H G$ such that the Kodaira-Spencer map $\mr{KS}_{(\mcE_H, \nabla)}$ is an isomorphism.
\item[(ii)]
Let $\mcE^\sss :=(\mcE_H, \nabla)$ and $\mcE'^\sss := (\mcE'_H, \nabla')$ be indigenous $(G, H)$-bundles on $X/S$.
An {\bf isomorphism of indigenous $(G, H)$-bundles} from $\mcE^\sss$ to $\mcE'^\sss$ is defined as an isomorphism of $H$-bundles  $\mcE_H \isom \mcE'_H$  such  that the induced isomorphism of $G$-bundles $\mcE_G \isom \mcE'_G$  is compatible  with the respective connections, i.e.,  $\nabla$ and $\nabla'$.
\end{itemize}
\ede
%----------------------------------------------------------------------

%----------------------------------------------------------------------
\bde
Let $\mcE^\sss := (\mcE_H, \nabla)$ be an indigenous $(G, H)$-bundle on $X$.
We shall say that $\mcE^\sss$ is {\bf flat} if the curvature  $\psi_\nabla$ (cf. (\ref{L054})) of   $\nabla$ vanishes.
\ede
%----------------------------------------------------------------------
\SSP

Denote by 
\begin{align}
\mr{Ind}_{(G, H),  X/S} \ \left(\text{resp.,}  \ \mr{Ind}_{(G, H),  X/S}^\flat\right)
\end{align}
the set of isomorphism classes of indigenous $(G, H)$-bundles (resp., flat indigenous $(G, H)$-bundles) on $X/S$.

\SSP
%--------------------------------------------------------------------
\begin{exa}[Standard example] \label{J03}
Recall that $\mcE^{\mr{triv}}_H := S \times G$ (cf. Example \ref{Exam1}) specifies  an $H$-reduction  on $P_S \times G$ via $(\mr{id}_S \times \pi_G,  \mr{pr}_2) : \mcE^{\mr{triv}}_H \migiincl P_S \times G$.
In particular,  there exists a natural isomorphism  of  $G$-bundles $\mcE^{\mr{triv}}_G \left(:= \mcE^{\mr{triv}}_H \times^H G \right) \isom P_S \times G$ extending this inclusion.
Since the composite $G \xrightarrow{(\mr{id}_S \times \pi_G,  \mr{pr}_2)} P_S \times G \xrightarrow{\mr{pr}_2} G$ coincides with the identity morphism,
the composite 
\begin{align} \label{J02}
\widetilde{\mcT}_{\mcE_H^{\mr{triv}}/S} \migiincl \widetilde{\mcT}_{P_S \times G/S} \left(\stackrel{\tau}{\cong}  \mcT_{P_S/S} \oplus \mcO_{P_S} \otimes \mfg\right) \xrightarrow{\mr{pr}_2} \mcO_{P_S} \otimes \mfg 
\end{align}
is an isomorphism.
This implies that the Kodaira-Spencer map $\mr{KS}_{(\mcE_H^{\mr{triv}}, \nabla_0)}$ associated with 
$(\mcE_H^{\mr{triv}}, \nabla_0)$ (cf. (\ref{L048}) for the definition of $\nabla_{(-)}$) is an isomorphism.
That is to say,  the pair
\begin{align} \label{J01}
\mcE^{\mr{triv}, \sss} := (\mcE^{\mr{triv}}_H, \nabla_0)
\end{align}
forms a flat indigenous $(G, H)$-bundle on $P_S/S$.
It is clear that  $\mcE^{\mr{triv}, \sss}$ is also  dormant  (cf. Definition \ref{Ets395} described later).
\end{exa}
%--------------------------------------------------------------------

\LSP
%----------------------------------------------------------------------[begin subsection]-------------
\subsection{Pull-back and base-change of  Indigenous $(G, H)$-bundles} \label{SS0123}

Let $Y$ and $y$ be as in \S\,\ref{SS066}.
Also, let $\mcE^\sss := (\mcE_H, \nabla)$ be an indigenous $(G, H)$-bundle on $X/S$.
 The base-change $y^*(\mcE_H)$ (resp., $y^*(\mcE_G)$) of $\mcE_H$ (resp., $\mcE_G := \mcE_H \times^H G$) has  an \'{e}tale morphism $y^*(\mcE_H) \migi \mcE_H$ (resp., $y^*(\mcE_G) \migi \mcE_G$), whose differential  yields   an isomorphism $d\widetilde{y}_{H} : \widetilde{\mcT}_{y^*(\mcE_H)/S} \isom  y^* (\widetilde{\mcT}_{\mcE_H/S})$ (resp., $d \widetilde{y}_{G} : \widetilde{\mcT}_{y^*(\mcE_G)/S} \isom  y^* (\widetilde{\mcT}_{\mcE_G/S})$).
The pull-back $y^* (\nabla) : y^*(\mcT_{X/S}) \migi  y^*(\widetilde{\mcT}_{\mcE_G/S})$ of $\nabla$ specifies an $S$-connection  on $y^*(\mcE_G)$ under the identifications 
$\mcT_{Y/S} = y^*(\mcT_{X/S})$, $\widetilde{\mcT}_{y^*(\mcE_G)/S} =   y^* (\widetilde{\mcT}_{\mcE_G/S})$ given by $d y$, $d \widetilde{y}_{G}$ respectively.
Since $d\widetilde{y}_{G}$ is  restricted  to $d \widetilde{y}_{H}$,
we obtain the following commutative square diagram:
\begin{align} \label{E00100}
\vcenter{\xymatrix@C=76pt@R=36pt{
\mcT_{Y/S} \ar[r]^-{\mr{KS}_{(y^*(\mcE_H), y^*(\nabla))}} \ar[d]^-{\wr}_-{d y} & \widetilde{\mcT}_{y^*(\mcE_G)/S}/\widetilde{\mcT}_{y^*(\mcE_H)/S}\ar[d]_-{\wr}
\\
y^*(\mcT_{X/S})\ar[r]_-{y^*(\mr{KS}_{(\mcE_H, \nabla)})} & y^*(\widetilde{\mcT}_{\mcE_G/S})/y^*(\widetilde{\mcT}_{\mcE_H /S}).
 }}
\end{align}
This implies that the Kodaira-Spencer map $\mr{KS}_{(y^*(\mcE_H), y^*(\nabla))}$ is an isomorphism, and hence, the pair of pull-backs 
\begin{align} \label{L0112}
y^*(\mcE^\sss) := (y^*(\mcE_H), y^*(\nabla))
\end{align}
forms an indigenous  $(G, H)$-bundle on $Y/S$.
We shall refer  $y^*(\mcE^\sss)$ as  the {\bf pull-back} of $(\mcE_H, \nabla)$.
Because of  the equality $y^*(\psi_\nabla) = \psi_{y^*(\nabla)}$,
 we see   that $y^*(\mcE^\sss)$ is flat  if $\mcE^\sss$ is flat.
 Hence, the  assignment $\mcE^\sss \mapsto y^*(\mcE^\sss)$ defines  maps  
\begin{align}
\mr{Ind}_{(G, H), X/S} \migi \mr{Ind}_{(G, H), Y/S}, 
\hspace{5mm}
\mr{Ind}_{(G, H), X/S}^\flat \migi \mr{Ind}_{(G, H), Y/S}^\flat.
\end{align}
The formation of a(n) (flat) indigenous  $(G, H)$-bundle on $X/S$ has descent with respect to the \'{e}tale topology on $X$.

Next, let $s : S' \migi S$  and $X'$ be as in \S\,\ref{SS066} and let $\mcE^\sss$ be as above.
We shall use the notation $s^*(-)$ to denote the result of  base-change by $s$. 
The base-change $s^* (\mcE_H)$ defines an $H$-reduction of 
$s^*(\mcE_G)$,
and the differential of  the projection   $s^*(\mcE_G) \migi \mcE_G$ yields an isomorphism $d \widetilde{s}_{G} : \widetilde{\mcT}_{s^*(\mcE_G)/S'} \isom s^*(\widetilde{\mcT}_{\mcE_G/S})$.
Just as in the above discussion,
the base-change $s^*(\nabla)$ defines an $S'$-connection on $s^*(\mcE_G)$ under the identification $\widetilde{\mcT}_{s^*(\mcE_G)/S'} =  s^*(\widetilde{\mcT}_{\mcE_G/S})$   given by $d \widetilde{s}_{G}$.
Also,   
the Kodaira-Spencer map associated with  $(s^*(\mcE_H), s^*(\nabla))$ may be identified with the isomorphism $s^*(\mr{KS}_{(\mcE_H, \nabla)})$.
Thus, the pair of base-changes
\begin{align} \label{H009}
s^*(\mcE^\sss) \left(\text{or}  \ S' \times_S \mcE^\sss \right) := (s^*(\mcE_H), s^*(\nabla))
\end{align}
forms an indigenous $(G, H)$-bundle on $X'/S'$, called  the {\bf base-change} of $\mcE^\sss$.
Since 
the equality $s^*(\psi_{\nabla}) = \psi_{s^*(\nabla)}$ holds,
$s^*(\mcE^\sss)$ is flat if $\mcE^\sss$ is flat. 
Thus, the resulting assignment $\mcE^\sss \mapsto s^*(\mcE^\sss)$ determines  maps 
\begin{align}
\mr{Ind}_{(G,H), X/S} \migi \mr{Ind}_{(G, H), X'/S'}, 
\hspace{5mm}
\mr{Ind}_{(G,H), X/S}^{\flat} \migi \mr{Ind}_{(G, H), X'/S'}^{\flat}. 
\end{align}

\vspace{5mm}
%----------------------------------------------------------------------[begin subsection]-------------
\subsection{From $(G, H)$-Cartan geometries to Indigenous $(G, H)$-bundles} \label{SS089}

We shall construct an indigenous $(G, H)$-bundle by using a Cartan geometry with model $(G, H)$.

Let $\mcE^\cc := (\mcE_H, \omega)$ be a Cartan geometry with model $(G, H)$ on $X/S$.
Denote by $\pi : \mcE_H \migi X$
the structure morphism of $\mcE_H$ and by 
 $\jmath : \mcE_H \migiincl \mcE_G := \mcE_H \times^H G$  the natural inclusion.
The composite $\jmath^*(\omega^\triangleright_{\mcE_G/X})^{-1} \circ \omega^\triangleright$  (cf. (\ref{L0116}))  defines an  $\mcO_{\mcE_H}$-linear isomorphism $\mcT_{\mcE_H /S} \isom \jmath^*(\mcT_{\mcE_G/X})$.
The direct image  $\pi_*(\mcT_{\mcE_H /S}) \isom  \pi_*(\jmath^*(\mcT_{\mcE_G/X}))$ of this morphism    preserves  the $H$-actions,
so  it is  restricted  to an isomorphism $\widetilde{\mcT}_{\mcE_H/S} \isom  \pi_{*}(\jmath^*(\mcT_{\mcE_G/X}))^H$.
Since  $\mfg_{\mcE_G}$ is naturally isomorphic to $\pi_{*} (\jmath^*(\mcT_{\mcE_G/X}))^H$,
 it may be regarded as an isomorphism
\begin{align} \label{E0213}
 \widetilde{\omega}^\triangleright : \widetilde{\mcT}_{\mcE_H /S} \isom  \left(\mfg_{\mcE_H} = \right) \mfg_{\mcE_G}.
 \end{align}
  By the commutativity of (\ref{E0057}), 
 the restriction of $\widetilde{\omega}^\triangleright$ to   $\mfh_{\mcE_H}  \left(\subseteq \widetilde{\mcT}_{\mcE_H/S} \right)$ coincides with the inclusion $\mfh_{\mcE_H} \migiincl \mfg_{\mcE_G}$.
Observe  here that 
(\ref{E0ygg0})
 induces 
an  $\mcO_X$-linear  isomorphism  $\widetilde{\mcT}_{\mcE_H /S} \times^{\mfh_{\mcE_H}} \mfg_{\mcE_G} \isom \widetilde{\mcT}_{\mcE_G /S}$.
Under the identification $\widetilde{\mcT}_{\mcE_H /S} \times^{\mfh_{\mcE_H}} \mfg_{\mcE_G} = \widetilde{\mcT}_{\mcE_G /S}$ given by this isomorphism,  
the morphism $(\widetilde{\omega}^\triangleright, \mr{id}_{\mfg_{\mcE_G}}) : \widetilde{\mcT}_{\mcE_H} \times^{\mfh_{\mcE_H}} \mfg_{\mcE_G} \migi \mfg_{\mcE_G}$
determines  a split 
surjection $\widetilde{\mcT}_{\mcE_G} \migisurj \mfg_{\mcE_G}$
of (\ref{Ex0})   for $\mcE_G$.
The corresponding   split injection 
\begin{align}
\nabla_{\mcE^\cc} : \mcT_{X/S} \migiincl \widetilde{\mcT}_{\mcE_G /S}
\end{align}
specifies   
an $S$-connection on the $G$-bundle $\mcE_G$.

 The direct image $\pi_{*}(\bigwedge^2 \Omega_{\mcE_H/S} \otimes \mfg) \left(= \pi_{*}(\bigwedge^2 \Omega_{\mcE_H/S}) \otimes \mfg \right)$
has a right $H$-action  given by $(a \otimes v) \cdot h := d\mr{R}_h^\vee (a) \otimes \mr{Ad}_G (h^{-1})(v)$ for $a \in \pi_*(\bigwedge^2 \Omega_{\mcE_H/S})$, $v \in \mfg$, and $h \in H$.
Because of  the second condition in Definition \ref{D024}, (i), 
 the curvature  $\psi_\omega$ of the $1$-form $\nabla_{\omega}$, being an element of $\Gamma (\mcE_H, \bigwedge^2 \Omega_{\mcE_H/S} \otimes \mfg) \left(= \Gamma (X, \pi_* (\bigwedge^2 \Omega_{\mcE_H/S} \otimes \mfg)) \right)$,  is invariant under this $H$-action.
 Hence,  $\psi_\omega$ defines an element of
 $\Gamma (X, \pi_* (\bigwedge^2 \Omega_{\mcE_H/S} \otimes \mfg)^H)$.

 \SSP
 %--------------------------------------------------------------
 \bpr \label{P033}
 The pair  
 \begin{align} \label{Eorr093}
 \mcE^{\cc \Rightarrow \sss} := (\mcE_H, \nabla_{\mcE^\cc})
 \end{align}
  defines an indigenous $(G, H)$-bundle on $X/S$.
Moreover, $ \mcE^{\cc \Rightarrow \sss}$ is flat if and only if
$\mcE^\cc$ is flat.
 \epr
 %--------------------------------------------------------------
 \begin{proof}
For simplicity, we shall write $\nabla := \nabla_{\mcE^\cc}$.
First, let us consider 
 the former assertion.
By  the definition of $\mr{KS}_{(\mcE_H, \nabla)}$, 
  the following diagram forms a morphism of short exact sequences:
\begin{align} \label{E0020}
\vcenter{\xymatrix@C=46pt@R=36pt{
0 \ar[r] &  \mfh_{\mcE_H} 
 \ar[r]^{} \ar[d]^{\mr{id}}& \widetilde{\mcT}_{\mcE_H}   \ar[r]^{d_{\mcE_H}} \ar[d]^{\widetilde{\omega}^\triangleright}& \mcT_{X/k}\ar[d]^{\mr{KS}_{(\mcE_H, \nabla)}} \ar[r] & 0\\
0 \ar[r] & \mfh_{\mcE_H}
 \ar[r]_-{\mr{inclusion}}& \mfg_{\mcE_G}  \ar[r]_-{\mr{quotient}}& \mfg_{\mcE_G}/\mfh_{\mcE_H}  \ar[r] &0.
 }}
\end{align}
Hence, since $\widetilde{\omega}^\triangleright$ is an isomorphism, 
the Kodaira-Spencer map $\mr{KS}_{(\mcE_H, \nabla)}$  is an isomorphism.
That is to say, the pair $\mcE^{\cc \Rightarrow \sss} := (\mcE_H, \nabla)$ specifies   an indigenous $(G, H)$-bundle on $X/S$, which 
 completes the proof of the former assertion.

Next, let us consider the latter assertion.
Note that
the assignments $\mcV \mapsto \pi_*(\mcV)^H$ and  $\mcU \mapsto \pi^*(\mcU)$ define an equivalence of categories
\begin{align} \label{Er11g}
\left(\begin{matrix} \text{the category of $H$-equivariant } \\ \text{vector bundles on $\mcE_H$}\end{matrix} \right)
\isom \left(\begin{matrix} \text{the category of} \\ \text{vector bundles on $X$}\end{matrix} \right).
\end{align}
In particular, 
since  $\widetilde{\mcT}_{\mcE_H/S} = \pi_*(\mcT_{\mcE_H/S})^H$ and $\mfg_{\mcE_G} = \pi_*(\mcO_{\mcE_H} \otimes \mfg)^H$,
the natural morphisms
 \begin{align} \label{Efg01}
 \pi^*(\widetilde{\mcT}_{\mcE_H/S}) \migi \mcT_{\mcE_H/S} \hspace{3mm}  \text{and}
 \hspace{3mm}
 \pi^*(\mfg_{\mcE_G}) \migi \mcO_{\mcE_H} \otimes \mfg
 \end{align}
  are isomorphisms.
We shall define  $\psi_{\widetilde{\omega}^\triangleright}$ to be   the $\mcO_X$-linear morphism $\bigwedge^2 \widetilde{\mcT}_{\mcE_H/S} \migi \mfg_{\mcE_G}$ determined by 
\begin{align} \label{ER01}
\psi_{\widetilde{\omega}^\triangleright} (\partial_1, \partial_2) = \partial_1 (\widetilde{\omega}^\triangleright (\partial_2)) - \partial_2 (\widetilde{\omega}^\triangleright (\partial_1)) + [\widetilde{\omega}^\triangleright (\partial_1), \widetilde{\omega}^\triangleright (\partial_2)] - \widetilde{\omega}^\triangleright ([\partial_1, \partial_2])
\end{align}
 for any local sections $\partial_1, \partial_2 \in \widetilde{\mcT}_{\mcE_H /S}$.
 Under the identifications $\pi^*(\widetilde{\mcT}_{\mcE_H/S}) = \mcT_{\mcE_H/S}$,  $\pi^*(\mfg_{\mcE_G}) = \mcO_{\mcE_H} \otimes \mfg$ given by (\ref{Efg01}),
 the equality $\pi^*(\psi_{\widetilde{\omega}^\triangleright}) = \psi_\omega$ holds (cf. the last paragraph in the discussion preceding this proposition).
 Hence, $\psi_{\widetilde{\omega}^\triangleright} =0$ if and only if $\psi_\omega =0$.
 On the other hand, 
 $\widetilde{\omega}^\triangleright$ can be expressed as $\widetilde{\omega}^\triangleright  =  \nabla \circ  d_{\mcE_H} -\mr{id}_{\widetilde{\mcT}_{\mcE_H/S}}$.
If we take arbitrary  local sections $\partial_1, \partial_2 \in \widetilde{\mcT}_{\mcE_H/S}$ and write $\overline{\partial}_1 := d_{\mcE_H}(\partial_1)$, $\overline{\partial}_2 := d_{\mcE_H} (\partial_2)  \left(\in \mcT_{X/S}\right)$, then   (\ref{ER01}) induces  the following sequence of equalities:
 \begin{align}
 & \ \psi_{\widetilde{\omega}^\triangleright} (\partial_1, \partial_2)\\ 
 = & \ 
 [\partial_1, \nabla (\overline{\partial}_2) - \partial_2] - [\partial_2, \nabla (\overline{\partial}_1) - \partial_1] 
 \notag \\
  & \   + [\nabla (\overline{\partial}_1) - \partial_1, \nabla (\overline{\partial}_2) - \partial_2] - \left(\nabla (d_{\mcE_H} ([\partial_1, \partial_2])) - [\partial_1, \partial_2]\right) 
  \notag \\
  = & \ \left(\cancel{[\partial_1, \nabla (\overline{\partial}_2)]} - \cancel{[\partial_1, \partial_2]} \right) - \left(\cancel{[\partial_2, \nabla (\overline{\partial}_1)]} -\cancel{[\partial_2, \partial_1]}\right) 
  \notag \\
   & \ +  \left([\nabla (\overline{\partial}_1), \nabla (\overline{\partial}_2)] - \cancel{[\nabla (\overline{\partial}_1), \partial_2]} - \cancel{[\partial_1, \nabla (\overline{\partial}_2) ]} + \cancel{[\partial_1, \partial_2]}\right) - \left(\nabla ([\overline{\partial}_1, \overline{\partial}_2]) - \cancel{[\partial_1, \partial_2]}\right) \notag \\
    = &  \ \psi_{\nabla} (\overline{\partial}_1, \overline{\partial}_2). \notag
  \end{align}
 Thus, since $d_{\mcE_H} : \widetilde{\mcT}_{\mcE_H/S} \migi \mcT_{X/S}$ is surjective, we see that $\psi_\nabla = 0$  if and only if  $\psi_{\widetilde{\omega}^\triangleright} = 0$, or equivalently,  $\psi_\omega = 0$ by 
  the above discussion.
  This completes the proof of the latter assertion.
 \end{proof}
 %--------------------------------------------------------------
\SSP

 %--------------------------------------------------------------------
 \bco \label{C090}
 The assignment $\mcE^\cc \mapsto \mcE^{\cc \Rightarrow \sss}$ resulting from Proposition \ref{P033}  defines a bijection of sets
 \begin{align} \label{E096}
 \mr{Car}_{(G, H), X/S} \isom  \mr{Ind}_{(G, H), X/S} \ \left(\text{resp.,} \ \mr{Car}_{(G, H), X/S}^\flat \isom  \mr{Ind}_{(G, H), X/S}^\flat  \right).
 \end{align}
 Moreover, the formation of this bijection is compatible  with pull-back to \'{e}tale schemes over $X$, as well as with base-change to schemes over $S$.
 \eco
 %--------------------------------------------------------------------
\begin{proof}
Since the resp'd assertion follows from the non-resp'd assertion and the latter assertion of Proposition \ref{P033}, it suffices to consider the non-resp'd assertion.
To prove this,
we shall construct the inverse to the  assignment $\mcE^\cc \mapsto \mcE^{\cc \Rightarrow \sss}$.
Let $\mcE^\sss := (\mcE_H, \nabla)$ be an indigenous $(G, H)$-bundle on $X/S$.
Denote by $\widetilde{\omega}^\triangleright_{\mcE^\sss}$ the composite
of the natural inclusion $\widetilde{\mcT}_{\mcE_H} \migiincl \widetilde{\mcT}_{\mcE_G}$ and the split surjection $\widetilde{\mcT}_{\mcE_G/S} \migisurj \mfg_{\mcE_G/S}$ of (\ref{Ex0}) corresponding to $\nabla$.
This morphism fits into the following morphism of short exact sequences:
\begin{align} \label{E0029}
\vcenter{\xymatrix@C=46pt@R=36pt{
0 \ar[r] &  \mfh_{\mcE_H} 
 \ar[r]^{} \ar[d]^-{\mr{id}}& \widetilde{\mcT}_{\mcE_H/S}   \ar[r]^-{d_{\mcE_H}} \ar[d]^{\widetilde{\omega}_{\mcE^\sss}^\triangleright}& \mcT_{X/S}\ar[d]^-{\mr{KS}_{(\mcE_H, \nabla)}} \ar[r] & 0\\
0 \ar[r] & \mfh_{\mcE_H}
 \ar[r]_-{\mr{inclusion}}& \mfg_{\mcE_G}  \ar[r]_-{\mr{quotient}}& \mfg_{\mcE_G}/\mfh_{\mcE_H}  \ar[r] &0.
 }}
\end{align}
The Kodaira-Spencer map $\mr{KS}_{(\mcE_H, \nabla)}$ is an isomorphism,
so the five lemma shows that 
 $\widetilde{\omega}_{\mcE^\sss}^\triangleright$  is  an isomorphism.
Under the identifications given by (\ref{Efg01}),
the pull-back of $\widetilde{\omega}_{\mcE^\sss}^\triangleright$ by $\pi$
defines an $\mcO_{\mcE_H}$-linear isomorphism $\mcT_{\mcE_H/S} \isom \mcO_{\mcE_H} \otimes \mfg$.
If $\omega_{\mcE^\sss}$  denotes the corresponding  $\mfg$-valued $1$-form on $\mcE_H/S$, then   the pair $\mcE^{\sss \Rightarrow \cc} := (\mcE_H, \omega_{\mcE^\sss})$ forms a Cartan geometry with model $(G, H)$ on $X/S$.
(The commutativity of the left-hand square in (\ref{E0029}) implies, via pull-back by $\pi$,  the first condition in Definition \ref{D024}, (i).)
Moreover,  the resulting assignment $\mcE^\sss \mapsto \mcE^{\sss \Rightarrow \cc}$  is verified to define   the inverse to the assignment $\mcE^\cc \mapsto \mcE^{\cc \Rightarrow \sss}$.
This implies  the  bijectivity of   (\ref{E096}).
 The required compatibility of the formation 
  follows immediately from its  construction, and hence, 
we finish  the proof of this corollary.
\end{proof}
%---------------------------------------------------------------------

\LSP
%---------------------------[begin subsection]-------------
\subsection{Dormant indigenous $(G, H)$-bundles} \label{SS088}

Next, we make the definition of a dormant indigenous bundle.
If  the underlying space is a curve, then this notion was already defined and has been studied  in  previous works (cf., e.g., ~\cite{Mzk2}).

\SSP
%-----------------------------------------------------
\bde \label{Ets395}
We shall say that  an indigenous $(G, H)$-bundle  $\mcE^\sss := (\mcE_H, \nabla)$ is {\bf dormant} if  it is flat and the $p$-curvature ${^p}\psi_\nabla$ of $\nabla$ (cf. (\ref{L056})) vanishes identically on $X$.
\ede
%----------------------------------------------------
\SSP

Denote by 
\begin{align} \label{e50q866}
\mr{Ind}_{(G, H), X/S}^{^\mr{Zzz...}}
\end{align}
the subset of $\mr{Ind}_{(G, H), X}^\flat$ consisting of dormant indigenous $(G, H)$-bundles.
Then, we obtain 
 the remaining portion of Theorem \ref{TheoremA1} described below.
(The map $\zeta^{\hh \Rightarrow \sss}$, which is one of the required bijections  in that theorem,   is defined as the composite of  $\zeta^{\cc \Rightarrow \sss}$ constructed below and $\zeta^{\hh \Rightarrow \cc}$.)

\SSP
%-----------------------------------------------------
\bpr \label{P0234}
The bijection $\mr{Car}_{(G, H), X/S} \isom \mr{Ind}_{(G, H), X/S}$ resulting from Corollary \ref{C090} is restricted  to a bijection of sets
\begin{align}
\zeta^{\cc \Rightarrow \sss} : \mr{Car}^{p\text{-}\mr{flat}}_{(G, H), X/S} \isom \mr{Ind}_{(G, H), X/S}^{^\mr{Zzz...}}.
\end{align}
Moreover, the formation of this bijection commutes with pull-back to \'{e}tale schemes over $X$, as well as with base-change to schemes over $S$.
\epr
%----------------------------------------------------
\begin{proof}
We shall consider the former assertion.
Let $\mcE^\cc := (\mcE_H, \omega)$ be a Cartan  geometry with model $(G, H)$ on $X/S$.
Denote by  $\pi: \mcE_H \migi X$ the structure morphism of $\mcE_H$ and by  $\mcE^{\cc \Rightarrow \sss} := (\mcE_H, \nabla)$ 
  the  indigenous $(G, H)$-bundle corresponding to $\mcE^\cc$ via (\ref{E096}), i.e., $\nabla := \nabla_{\mcE^\cc}$.
  Under the identifications $\pi^*(\widetilde{\mcT}_{\mcE_H/S}) = \mcT_{\mcE_H/S}$, $\pi^*(\mfg_{\mcE_G}) = \mcO_{\mcE_H} \otimes \mfg$ given by (\ref{Efg01}),
  the $p$-curvature ${^p}\psi_\omega$  of $\omega$ may be regarded as an $\mcO_{\mcE_H} \left(= \pi^*(\mcO_X) \right)$-linear morphism $\left(F_{\mcE_H}^*(\pi^*(\widetilde{\mcT}_{\mcE_H/S}))= \right) \pi^*(F^*_{X}(\widetilde{\mcT}_{\mcE_H/S})) \migi \pi^*(\mfg_{\mcE_G})$.
  This morphism comes  from an $\mcO_X$-linear morphism ${^p}\psi_{\widetilde{\omega}^\triangleright} : F^*_{X}(\widetilde{\mcT}_{\mcE_H/S}) \migi \mfg_{\mcE_G}$   via pull-back by $\pi$.
 If  $\widetilde{\omega}^\triangleright : \widetilde{\mcT}_{\mcE_H/S} \isom \mfg_{\mcE_G}$ denotes the isomorphism constructed from $\mcE^\cc$   as in  (\ref{E0213}), 
 then  ${^p}\psi_{\widetilde{\omega}^\triangleright}$ is defined
  by 
  \begin{align}
  {^p}\psi_{\widetilde{\omega}^\triangleright} (F^{-1}_X (\partial)) = 
  \widetilde{\omega}^\triangleright (\partial)^{[p]} - \widetilde{\omega}^\triangleright (\partial^{[p]}) + \sum_{i=1}^{p-1} \frac{s_i (\partial, \widetilde{\omega}^\triangleright (\partial))}{i}
  \end{align}
  for any local section $\partial \in \mcT_{X/S}$.
 Let us take an arbitrary local section $\partial \in \widetilde{\mcT}_{\mcE_H/S}$  and write  $\overline{\partial} := d_{\mcE_H} (\partial) \in \mcT_{X/S}$.
 Since $\widetilde{\omega}^\triangleright = \nabla \circ d_{\mcE_H} - \mr{id}_{\widetilde{\mcT}_{\mcE_H/S}}$, we obtain the following sequence of  equalities:
  \begin{align}
& \   {^p}\psi_{\widetilde{\omega}^\triangleright} (F^{-1}_X (\partial)) \\
   = & \ 
  (\nabla \circ d_{\mcE_H}-\mr{id})(\partial)^{[p]} - (\nabla \circ d_{\mcE_H} -\mr{id})(\partial^{[p]})  + \sum_{i=1}^{p-1} \frac{s_i (\partial, (\nabla \circ d_{\mcE_H}-\mr{id})(\partial))}{i} \notag \\
  = & \ (\nabla (\overline{\partial})-\partial)^{[p]} - (\nabla (\overline{\partial}^{[p]})-\partial^{[p]}) + \sum_{i=1}^{p-1} \frac{s_i (\partial, \nabla (\overline{\partial})-\partial)}{i} \notag \\
  = & \   \cancel{(-1)^{p} \cdot \partial^{[p]}} + \nabla (\overline{\partial})^{[p]}  + 
  \cancel{\sum_{i=1}^{p-1} \frac{s_i (-\partial, \nabla (\overline{\partial}))}{i}} - \nabla (\overline{\partial}^{[p]}) + \cancel{\partial^{[p]}} +   \cancel{\sum_{i=1}^{p-1} \frac{s_i (\partial, \nabla (\overline{\partial})-\partial)}{i}} \notag \\
  = & \ {^p}\psi_{\nabla} (F^{-1}_X (\overline{\partial})),\notag
  \end{align}
  where the last equality follows from Lemma \ref{L1123} described below.
  Hence, the surjectivity of $d_{\mcE_H} : \widetilde{\mcT}_{\mcE_H/S} \migi \mcT_{X/S}$ implies  that  ${^p}\psi_{\nabla} =0$  if and only if  ${^p}\psi_{\widetilde{\omega}^\triangleright} =0$,  or equivalently, ${^p}\psi_\omega = 0$ by the above discussion.
 This completes the proof of the former assertion.
  The  latter assertion  follows immediately from its  construction.
 \end{proof}
 %---------------------------
\SSP

The following lemma was used  in the proof of the above proposition.

\SSP
%---------------------------
\ble \label{L1123}
Let $\mfa$ be a Lie algebra over a $k$-algebra $R$.
Then, for any elements $v, u \in \mfa$, 
the following equality holds:
\begin{align}
\sum_{i=1}^{p-1} \frac{s_i (v, u-v)}{i} = - \sum_{i=1}^{p-1} \frac{s_i (-v, u)}{i}.
\end{align}
\ele
%---------------------------
\begin{proof}
By the definition of $s_i$,  the following sequence of equalities holds:
\begin{align}
\sum_{i=1}^{p-1}s_i (v, u-v) t^{i-1} = & \  \mr{ad} (vt + u-v)^{p-1} (v)  \\
= &  \ \mr{ad} (v(t-1) + u)^{p-1} (v) \notag \\
 = & \  \sum_{j=1}^{p-1} s_j (v, u) (t-1)^{j-1} \notag \\
 = & \ \sum_{i=1}^{p-1} \left(\sum_{j=i}^{p-1} (-1)^{j-i} \cdot  \binom{j-1}{i-1} \cdot s_j (v, u)  \right)  t^{i-1}. \notag 
\end{align}
This implies $s_i (v, u- v) = \sum\limits_{j=i}^{p-1} (-1)^{j-i} \cdot  \binom{j-1}{i-1} \cdot s_j (v, u) $ for $i=1, \cdots , p-1$.
Thus, we have 
\begin{align}
\sum_{i=1}^{p-1} \frac{s_i (v, u-v)}{i} =
& \ \sum_{i=1}^{p-1}  \sum\limits_{j=i}^{p-1} \frac{1}{i} \cdot (-1)^{j-i} \cdot  \binom{j-1}{i-1} \cdot s_j (v, u)\\
  = & \  \sum_{j=1}^{p-1} \frac{s_j (v, u)}{j} \cdot \left(\sum_{i=1}^{j}  (-1)^{j-i} \cdot \binom{j}{i} \right) \notag \\
  = & \  \sum_{j=1}^{p-1} \frac{s_j (v, u)}{j} \cdot (-1)^{j+1}  \notag \\
  = & \  - \sum_{j=1}^{p-1} \frac{s_j (-v, u)}{j}. \notag
\end{align}
This completes the proof of this lemma.
\end{proof}
%---------------------------

%-------------------------------------------------------------
\LSP
\subsection{Indigenous $(G, H)$-bundles of higher level} \label{S0210}\vspace{3mm}

We generalize both  the notion of an indigenous $(G, H)$-bundle and the bijection $\zeta^{\cc \Rightarrow \sss}$ to higher level.
In the rest of this section, suppose that $G$ is affine.
Let $\N$ be an element of $\mbZ_{> 0} \sqcup \{ \infty \}$.
(Whenever we deal with $(\N-1)$-PD stratifications on a $G$-bundle for $\N \in \mbZ_{>0}$, we always suppose that the scheme $S$ is equipped with  an $(\N-1)$-PD ideal $(\mfa,  \mfb, \gamma)$, as in \S\,\ref{SS040}, such that $\gamma$ extends to $X$.)

\SSP
%-------------------------------------------------------------
\bde \label{Eosi38392}
\begin{itemize}
\item[(i)]
An {\bf indigenous $(G, H)$-bundle of level $\N$} on $X/S$ is a pair
\begin{align}
\mcE^\sss_\N := (\mcE_H,  \varepsilon^\ST)
\end{align}
consisting of an $H$-bundle $\mcE_H$ on $X$ and an $(\N-1)$-PD stratification $\varepsilon^\ST$ (cf. Definition \ref{D057}) on $\mcE_G := \mcE_H \times^H G$ such that if $\nabla$ denotes the flat $S$-connection on $\mcE_G$ induced by the $0$-PD stratification $\varepsilon^\ST |^{\langle 0 \rangle}$, i.e., $\nabla := \nabla_{\varepsilon^\ST |^{\langle 0 \rangle}}$ (cf. (\ref{Rt567}) and (\ref{eowip0})), then $(\mcE_H, \nabla)$ forms an indigenous $(G, H)$-bundle on $X/S$.
\item[(ii)]
Let $\mcE^\sss_\N := (\mcE_H, \varepsilon^\ST)$ and $\mcE'^\sss_\N := (\mcE'_H, \varepsilon'^\ST)$
be indigenous $(G, H)$-bundles of level $\N$ on  $X/S$.
An {\bf isomorphism of indigenous $(G, H)$-bundles of level $\N$} from $\mcE^\sss_\N$ to $\mcE'^\sss_\N$ is an isomorphism $\mcE_H \isom \mcE'_H$ of $H$-bundles such that the induced isomorphism of $G$-bundles $\mcE_G \isom \mcE'_G$  is compatible with the respective $(\N-1)$-PD stratifications, i.e., $\varepsilon^\ST$ and $\varepsilon'^\ST$.
\end{itemize}
\ede
%-------------------------------------------------------------
\SSP

Given an  indigenous $(G, H)$-bundle $\mcE_\N^\ST$ of level $\N$,  we can define,  in an obvious way,  the pull-back of $\mcE_\N^\ST$  to   an \'{e}tale $X$-scheme, as well as  with the base-change of $\mcE_\N^\ST$ to an $S$-scheme. (But we leave the details to the reader.)

\SSP
%-------------------------------------------------------------
\bde \label{Ge990}
Let $\mcE^\sss_\N := (\mcE_H,  \varepsilon^\ST)$ be an indigenous $(G, H)$-bundle of level $\N$.
If $\N \in \mbZ_{> 0}$ (i.e., $\N \neq \infty$), then we shall say that $\mcE^\sss_\N$ is {\bf dormant}  if the $p$-$(\N-1)$-curvature ${^p}\psi_{\varepsilon^\ST}$ of $\varepsilon^\ST$ (cf. Definition \ref{rtyu5}) vanishes identically.
Also, if $\N = \infty$, then we shall refer to any such $\mcE^\sss_\N$ as being {\bf dormant}, for convenience.
\ede
%-------------------------------------------------------------
\SSP

Denote by
\begin{align}
\mr{Ind}_{(G, H), X/S, \N}^{^\mr{Zzz...}}
\end{align}
the set of isomorphism classes of dormant indigenous $(G, H)$-bundles of level $\N$ on $X/S$.
It follows from Proposition \ref{C071} that
the notion of dormant indigenous $(G, H)$-bundle is equivalent to that of level $1$.
That is to say, we have a natural identification
\begin{align} \label{QQWv4}
\mr{Ind}_{(G, H), X/S,1}^{^\mr{Zzz...}} = \mr{Ind}_{(G, H), X/S}^{^\mr{Zzz...}}.
\end{align}
Let 
$\N'$ be another element of $\mbZ_{>0} \sqcup \{ \infty \}$ with $\N \leq \N'$ and 
and $\mcE^\sss_\N := (\mcE_H, \varepsilon^\ST)$ a(n) (dormant) indigenous $(G, H)$-bundle  of level $\N'$.
Then,  the pair 
\begin{align} \label{gRT4}
\mcE^\sss_\N |^{\langle \N \rangle} := (\mcE_H, \varepsilon^\ST |^{(\N-1)})
\end{align} 
(cf. (\ref{Rt567})) specifies a(n) (dormant) indigenous $(G, H)$-bundle  of level $\N$, which is called the {\bf $\N$-th truncation} of $\mcE^\sss$.
The resulting assignments $\mcE^\sss \mapsto \mcE^\sss |^{\langle \N \rangle}$ for various pairs $(\N, \N')$ with $\N \leq \N'$ give a projective system
\begin{align} \label{Ertq1}
\mr{Ind}_{(G, H), X/S, \infty}^{^\mr{Zzz...}} \migi \cdots \migi \mr{Ind}_{(G, H), X/S, \N}^{^\mr{Zzz...}} \migi \cdots
\migi \mr{Ind}_{(G, H), X/S, 2}^{^\mr{Zzz...}} \migi \mr{Ind}_{(G, H), X/S, 1}^{^\mr{Zzz...}}.
\end{align}

\SSP
%-------------------------------------------------------------
\begin{exa}[Standard example] \label{Erqw2311}
Let $\N \in \mbZ_{>0}$.
Then, the trivial $G$-bundle $P_S^{(\N)} \times G$ on $P_S^{(\N)}$ determines an $(\N-1)$-PD
 stratification $\varepsilon^{\mr{triv}, \ST}_\N$ on $P_S \times G$ via (\ref{E3erf002}).
 Then, the pair
 \begin{align}
 \mcE_\N^{\mr{triv}, \sss} := (\mcE_H^\mr{triv}, \varepsilon^{\mr{triv}, \ST}_\N)
 \end{align}
 (cf. (\ref{J01})) forms a dormant indigenous $(G, P)$-bundle of level $\N$ on $P_S/S$.
 Since the formation of $\varepsilon^{\mr{triv}, \ST}_\N$ commutes with truncation to lower levels, the collection $\{\varepsilon^{\mr{triv}, \ST}_\N \}_\N$ specifies an $\infty$-PD stratification $\varepsilon_\infty^{\mr{triv}, \ST}$.
 Thus, we obtain a dormant indigenous $(G, H)$-bundle 
 $\mcE_\infty^{\mr{triv}, \sss} := (\mcE_H^\mr{triv}, \varepsilon^{\mr{triv}, \ST}_\infty)$   of  level $\infty$.
 It is immediately verified that $ \mcE_\N^{\mr{triv}, \sss}$ corresponds to the standard example  of $F^\N$-Cartan geometry $\mcE^{\mr{triv}, \clubsuit}_\N$  (cf. (\ref{4598sj9})) via the bijection asserted in Proposition \ref{Eqw209} later.
 \end{exa}
%-------------------------------------------------------------
\SSP

The remaining portions of Theorem \ref{TheoremA} are completed by the following Proposition \ref{Eqw209} and Corollary \ref{Er281}.
(The map $\zeta_\N^{\hh \Rightarrow \sss}$, which is one of the required bijections in that theorem, is defined as the composite of $\zeta_\N^{\cc \Rightarrow \sss}$ below and  the bijection $\zeta^{\hh \Rightarrow \cc}_\N$ asserted in Propositions \ref{P0201} and \ref{Efget}.)

\SSP
%-------------------------------------------------------------
\bpr \label{Eqw209}
(Recall that   $G$ is assumed to be  affine.)
For any $\N \in \mbZ_{> 0} \sqcup \{ \infty \}$, 
there exists a canonical  bijection
of sets
\begin{align} \label{TRe46}
\zeta_\N^{\cc \Rightarrow \sss} : F^\N\text{-}\mr{Car}_{(G, H), X/S} \isom \mr{Ind}_{(G, H), X/S, \N}^{^\mr{Zzz...}},
\end{align}
and the formation of this bijection commutes with truncation to lower levels, pull-back to \'{e}tale $X$-schemes, and base-change to $S$-schemes.
\epr
%-------------------------------------------------------------
\begin{proof}
If $\N \in \mbZ_{>0}$ and  $\mcE_\N^\cc := (\mcE_H, \omega, \mcG, \kappa)$ is an $F^\N$-Cartan geometry with model  $(G, H)$,
then the  pair $\mcE_\N^{\cc \Rightarrow \sss} := (\mcE_H, \varepsilon^\ST)$, where $ \varepsilon^\ST$ denotes the $(\N-1)$-PD stratification on $\mcE_G := \mcE_H \times^H G$ corresponding to $\mcG$ via (\ref{E3erf002}), forms a dormant indigenous $(G, H)$-bundle of level $\N$;
the resulting  assignment $\mcE_\N^\cc \mapsto \mcE_\N^{\cc \Rightarrow \sss}$ defines the required bijection for $\N \in \mbZ_{>0}$ because of Proposition \ref{P0234} and  Proposition \ref{Ertwq}, (i).
Moreover, since this bijection commutes with truncation to lower levels (cf. Proposition \ref{Ertwq}, (ii)),  the case of $\N = \infty$ follows from Corollary \ref{C9865}.
 \end{proof}
%-------------------------------------------------------------
\SSP

%-------------------------------------------------------------
%\begin{rema}
\bco \label{Er281}
Suppose that $G$ is affine, $S = \mr{Spec}(k)$, and $X$ is  proper over $k$.
Then, 
the composite
\begin{align} \label{Rgjout}
\zeta_\infty^{\hh \Rightarrow \sss} := \zeta_{\infty}^{\cc \Rightarrow \sss} \circ \zeta_\infty^{\hh \Rightarrow \cc}  : F^\infty \text{-}\mr{Ehr}_{(G, P), X/S} \isom \mr{Ind}_{(G, H), X/S, \infty}^{^\mr{Zzz...}}
\end{align}
(cf. (\ref{E(eoe9}) for the definition of $\zeta_\infty^{\hh \Rightarrow \cc}$) is bijective.
In particular,
  the map of sets
\begin{align} \label{Rt980}
\mr{Ind}_{(G, H), X/S, \infty}^{^\mr{Zzz...}} \migi \varprojlim_{\N \geq \in \mbZ_{>0}} \mr{Ind}_{(G, H), X/S, \N}^{^\mr{Zzz...}}
\end{align}
induced by (\ref{Ertq1}) is bijective.
\eco
%\end{rema}
%-------------------------------------------------------------
\begin{proof}
The assertion follows   from Propositions \ref{Efget} and \ref{Eqw209}.
\end{proof}

\SSP
%-------------------------------------------------------------
\begin{exa} \label{R099}
We shall briefly recall the study of ~\cite{Wak6}, in which indigenous $(G, H)$-bundles for  $(G, H)= (\mr{PGL}_{n+1}, \mr{PGL}_{n+1}^\circledcirc)$ and $(\mr{PGL}_{n+1}, \mr{PGL}_{n+1}^\circledcirc)$ are discussed in terms of vector bundles.
Suppose that
  $X$ is a smooth projective scheme over $S = \mr{Spec}(k)$  of dimension $n$, and denote by 
$\omega_X$ the canonical line  bundle  on $X$.
 Then,  we obtain the composite isomorphism
\begin{align} \label{Ewq1}
\mr{det}(\mcD_{X/k, 1}^{(\N-1)}) \isom \mr{det} (\mcD_{X/k, 0}^{(\N-1)})
\otimes_{\mcO_X} \mr{det}(\mcD_{X/k, 1}^{(\N-1)}/\mcD_{X/k, 0}^{(\N-1)})
\isom \mr{det}(\mcT_{X/k}) \isom \omega_X^\vee
\end{align}
(cf. (\ref{MROD2}) for the definition of $\mcD_{X/k, (-)}^{(\N-1)}$).
\begin{itemize}
\item[(i)]
First, let us consider  the case where  $(G, H)= (\mr{PGL}_{n+1}, \mr{PGL}_{n+1}^\circledcirc)$ (cf. Example \ref{EX01}).
Suppose further that $p \nmid (n+1)$  and $X$ admits an  $F^\N$-theta characteristic (cf. ~\cite[\S\,3, Definition 3.5.1]{Wak6}), i.e., a pair $(\varTheta, \nabla^0)$ consisting of   a line bundle $\varTheta$ on $X$ and a $\mcD_{X/k}^{(\N-1)}$-module structure $\nabla^0$ on $\omega_X^\vee \otimes_{\mcO_X} \varTheta^{\otimes (n+1)}$  with vanishing $p$-$(\N-1)$-curvature (cf. (\ref{MSKO3})).
According to ~\cite[\S\,3, Theorem 3.6.4, (i)]{Wak6},  each element of $\mr{Ind}_{(\mr{PGL}_{n+1}, \mr{PGL}_{n+1}^\circledcirc), X/k, \N}^{^\mr{Zzz...}}$ may be expressed as a $\mcD_{X/k}^{(\N-1)}$-module structure 
\begin{align} \label{ew094szk}
\nabla : \mcD_{X/k}^{(\N-1)} \otimes_{\mcO_X} (\mcD_{X/k, 1}^{(\N-1)} \otimes_{\mcO_X} \varTheta) \migi \mcD_{X/k, 1}^{(\N-1)} \otimes_{\mcO_X} \varTheta
\end{align}
 on the rank $(n+1)$ vector bundle $\mcD_{X/k, 1}^{(\N-1)} \otimes_{\mcO_X} \varTheta$
   with vanishing $p$-$(\N-1)$-curvature satisfying the following conditions:
\begin{itemize}
\item[(a)]
The following  composite  defines an automorphism of  $\mcD_{X/k, 1}^{(\N-1)} \otimes_{\mcO_X} \varTheta$ (cf. ~\cite[\S\,2, Proposition 2.3.3]{Wak6}):
\begin{align}
\mcD_{X/k, 1}^{(\N-1)} \otimes_{\mcO_X} \varTheta \migiincl \mcD_{X/k}^{(\N-1)} \otimes_{\mcO_X} (\mcD_{X/k, 1}^{(\N-1)} \otimes_{\mcO_X} \varTheta) \xrightarrow{\nabla}\mcD_{X/k, 1}^{(\N-1)} \otimes_{\mcO_X} \varTheta,
\end{align}
 where the first arrow  arises from the natural inclusions $\mcD_{X/k, 1}^{(\N-1)} \migiincl \mcD_{X/k}^{(\N-1)}$ and $\varTheta \left(=  \mcD_{X/k, 0}^{(\N-1)}  \otimes_{\mcO_X} \varTheta\right) \migiincl \mcD_{X/k, 1}^{(\N-1)} \otimes_{\mcO_X} \varTheta$.
\item[(b)]
The determinant of $\nabla$ corresponds to $\nabla^0$ via the composite  isomorphism
\begin{align}
\mr{det}(\mcD_{X/k, 1}^{(\N-1)} \otimes_{\mcO_X} \varTheta) \isom \mr{det}(\mcD_{X/k, 1}^{(\N-1)}) \otimes_{\mcO_X} \varTheta^{\otimes (n+1)} \xrightarrow{ (\ref{Ewq1}) \otimes \mr{id}} \omega_X^\vee \otimes_{\mcO_X} \varTheta^{\otimes (n+1)}.
\end{align} 
\end{itemize}
If $X$ is a curve, then  we can take  $\varTheta$  as a theta characteristic in the usual sense and  $\nabla^0$ as the trivial connection $\nabla_0 \left(= d \right)$ under a fixed identification $\omega_X^\vee \otimes_{\mcO_X} \varTheta^{\otimes 2} = \mcO_X$.
In this situation, 
  a connection $\nabla$ as above  is nothing but  what we call  a  {\it  dormant $\mr{SL}_2$-oper}.
\item[(ii)]
Next,  let us recall  the case of $(G, H) = (\mr{Aff}_n, \mr{Aff}_n^\circledcirc)$.
According to ~\cite[\S\,3, Theorem 3.6.4, (i)]{Wak6},  each element of $\mr{Ind}_{(\mr{Aff}_n, \mr{Aff}_n^\circledcirc), X/k, \N}^{^\mr{Zzz...}}$ may be expressed as a   pair 
\begin{align} \label{e0s9z3}
(\nabla, \delta)
\end{align}
 consisting of a $\mcD_{X/k}^{(\N-1)}$-module structure $\nabla$   on the vector bundle $\mcD_{X/k, 1}^{(\N-1)}$ with vanishing $p$-$(\N-1)$-curvature  and a left inverse morphism   $\delta : \mcD_{X/k}^{(\N -1)} \migisurj \mcO_X$ of the natural inclusion $\mcO_X \migiincl \mcD_{X/k, 1}^{(\N-1)}$
satisfying the following conditions:
\begin{itemize}
\item[(c)]
The following  composite  defines an  automorphism of $\mcD_{X/k, 1}^{(\N-1)}$:
\begin{align}
\mcD_{X/k, 1}^{(\N-1)} \left(= \mcD_{X/k, 1}^{(\N-1)} \otimes_{\mcO_X} \mcO_X \right)
\migiincl  \mcD_{X/k}^{(\N-1)} \otimes_{\mcO_X} \mcD_{X/k, 1}^{(\N-1)} \xrightarrow{\nabla} \mcD_{X/k, 1}^{(\N-1)},
\end{align}
 where the first arrow arises from the natural inclusions
$\mcD_{X/k, 1}^{(\N-1)} \migiincl \mcD_{X/k}^{(\N-1)}$ and $\mcO_X \left(= \mcD_{X/k, 0}^{(\N-1)} \right)\migiincl \mcD_{X/k, 1}^{(\N-1)}$.
\item[(d)]
The subbundle $\mr{Ker}(\delta)$ of $\mcD_{X/k, 1}^{(\N-1)}$ is closed under the $\mcD_{X/k}^{(\N-1)}$-action $\nabla$ and 
the $\mcD_{X/k}^{(\N-1)}$-module structure on $\mcO_X$ induced from $\nabla$ via $\delta$ coincides with the trivial  one.
\end{itemize}
Now,  suppose that $p \nmid (n+1)$ and  there exists  an $F^\N$-theta characteristic $(\varTheta, \nabla^0)$ on $X$.
Let us take a pair $(\nabla, \delta)$ as above.
Since $\mr{det}((\mcD_{X/k, 1}^{(\N-1)})^\vee) \cong \omega_X$ (cf. (\ref{Ewq1})), we have the  $\mcD_{X/k}^{(\N-1)}$-module structure $\nabla_{\omega_X}$ on $\omega_X$ induced by (the determinant of the dual of) $\nabla$.
Then, it follows from 
~\cite[\S\,3, Lemma 3.6.3]{Wak6} that there exists uniquely
a $\mcD_{X/k}^{(\N-1)}$-module structure $\nabla_\varTheta$ on $\varTheta$
such that its $(n+1)$-st tensor product coincides, via $\varTheta^{\otimes (n+1)} \cong \omega_X \otimes_{\mcO_X} (\omega_X^\vee \otimes_{\mcO_X} \varTheta^{\otimes (n+1)})$,  with the tensor product of $\nabla_{\omega_X}$ and $\nabla^0$ (cf. \cite[\S\,2, (70)]{Wak6}).
Then, the tensor product of $\nabla$ and $\nabla_\varTheta$ defines a $\mcD_{X/k}^{(\N-1)}$-module structure $\nabla'$ on the vector bundle $\mcD_{X/k, 1}^{(\N-1)} \otimes_{\mcO_X} \varTheta$, satisfying the two conditions (a) and (b) in (i).
One may verify that the resulting assignment $\nabla \mapsto \nabla'$  coincides with 
(\ref{Fggth}) under the identification $F^\N\text{-}\mr{Ehr}_{(G, H), X/S} = \mr{Ind}_{(G, H), X/S, \N}^{^\mr{Zzz...}}$ given by  the bijection $\zeta_\N^{\hh \Rightarrow \sss} :=  \zeta_\N^{\cc \Rightarrow \sss}\circ \zeta_\N^{\hh \Rightarrow \cc}$.
\end{itemize}
\end{exa}
%-------------------------------------------------------------
\SSP

Let us keep the assumption in \S\,\ref{SS040}
  and take 
 an $F^\N$-$(G, P)$-structure $\mcS^\hh$  on $X/S$ with $\N \in \mbZ_{>0}$.
Then, we obtain   the indigenous $(G, H)$-bundle $\mcE^\hh_\N := (\mcE_H, \varepsilon^\ST)$   of level $\N$ corresponding to  $\mcS^\hh$ via $\zeta_\N^{\hh \Rightarrow \sss}$.
Then, the $(N-1)$-PD stratified $G$-bundle $(\mcE_G, \varepsilon^\ST)$, where $\mcE_G := \mcE_H \times^H G$,  detemines  an $(\N-1)$-crystal of $G$-bundles $\mcE^\diamondsuit$ on  $X/S$ (cf. Proposition \ref{P019}).

\SSP
%-------------------------------------------------------------
\bde \label{Efopwd}
We shall refer to $\mcE^\diamondsuit$ as the {\bf monodromy crystal} associated with $\mcS^\hh$.
\ede
%-------------------------------------------------------------

  \LSP
%-------------------------------------------------------------
%\vspace{5mm}
\subsection{Basic properties on $F^\infty$-$(G, P)$-structures} \label{S090}

By using the result of Theorem \ref{TheoremA}, we prove several properties concerning    smooth varieties admitting an $F^\infty$-$(G, P)$-structure.

 Recall that a {\it rational curve} in a $k$-scheme $X$ is a non-constant $k$-morphism $f : \mbP^1 \migi X$.
 We shall say that a rational curve $f: \mbP^1 \migi X$ is 
  {\it free} if $f^*(\mcT_{X})$ is spanned by its global sections.
 Also, a $k$-scheme  is called {\it convex} if every rational curve in it is free.
 
 \SSP
 %------------------------------------------------------------
 \bpr \label{etoiuoa892}
 Suppose that $G$ is affine.
 Then,
 any proper smooth $k$-scheme  admitting an $F^\infty$-$(G, P)$-structure  is convex.
 \epr
  %------------------------------------------------------------
  \begin{proof}
  Let  $X$ be a proper smooth $k$-scheme admitting an $F^\infty$-$(G, P)$-structure  $\mcS^\hh$.
  Denote by $\mcE^\sss := (\mcE_H, \varepsilon^\ST)$ the indigenous $(G, H)$-bundle of level $\infty$ corresponding to $\mcS^\hh$ via  $\zeta_{\infty}^{\hh \Rightarrow \sss}$ asserted in  Theorem \ref{TheoremA} and by 
   $\mcG_\infty := \{ (\mcG_l, \varpi_l ) \}_l$  the $F$-divided $G$-bundle  (cf. Definition \ref{EIFO(79}) corresponding to the $\infty$-PD stratified $G$-bundle $(\mcE_G, \varepsilon^\ST)$ via (\ref{E3erEr}).
 Now,  let  us take an arbitrary rational curve $f : \mbP^1 \migi X$  in $X$.
By pulling-back  $\mcG_\infty$ by $f$, we obtain an $F$-divided $G$-bundle $f^*(\mcG_\infty)$ on $\mbP^1$.
 Recall from  ~\cite[\S\,2, Theorem 2.2]{G} that there are no nontrivial $F$-divided  vector bundles on $\mbP^1$.
 Hence, for each  finite-dimensional representation $V$ of $G$ over $k$,
 the associated $F$-divided vector  bundle $f^*(\mcG_\infty) \times^G V$ may be trivialized.
 By   ~\cite[\S\,3, Lemma 3.5,  (ii)]{BS}, 
we see that  the $F$-divided $G$-bundle $f^*(\mcG_\infty)$ itself  may be trivialized.
 Let us fix an isomorphism $f^*(\mcE_G) \isom \mbP^1 \times G$ preserving the $F$-divided structures, where  $\mbP^1 \times G$ is equipped  with the trivial $F$-divided structure.
  It gives  
  an  $H$-equivariant composite 
  \begin{align} \label{Ef01}
 f^*(\mcE_H) \migiincl f^*(\mcE_G) \isom \mbP^1 \times G \xrightarrow{\mr{pr}_2} G.
  \end{align}
This morphism induces, via taking the quotients with respect to  the $H$-actions,  a morphism $h : \mbP^1 \left(:= f^*(\mcE_H)/H \right) \migi P \left(:= G/H \right)$.
Since the  commutative diagram 
 \begin{align} \label{E00217}
\vcenter{\xymatrix@C=46pt@R=36pt{
f^*(\mcE_H) \ar[r]^{(\ref{Ef01})} \ar[d]_-{\mr{projection}}& G \ar[d]^-{\pi_G} 
 \\
\mbP^1 \ar[r]_-{h} & P
 }}
\end{align}
 is  cartesian,   we obtain 
  an isomorphism  of $H$-bundles $f^*(\mcE_H) \isom h^*(G) \left(= \mbP^1 \times_{h, P, \pi_G} G \right)$.
Thus, we have
\begin{align} \label{Ed002}
f^*(\mcT_{X}) &\cong f^*(\mfg_{\mcE_G}/\mfh_{\mcE_H}) \\
& \cong f^*(\mcE_H \times^H (\mfg/\mfh)) \notag  \\
& \cong f^*(\mcE_H)  \times^H (\mfg/\mfh) \notag \\
&\cong h^*(G) \times^H (\mfg/\mfh) \notag \\
& \cong h^*(\mcE_H^\mr{triv} \times^H (\mfg/\mfh)) \notag \\
& \cong h^*(\mcT_{P}), \notag
\end{align}
where the first  and last ``$\cong$" follow from the Kodaira-Spencer maps $\mr{KS}_{(\mcE_H, \nabla)}$ and $\mr{KS}_{(\mcE_H^\mr{triv}, \nabla_0)}$ (cf. Example \ref{J03}) respectively.
 Observe  here that
 $\mcT_{P}$, as well as its pull-back $h^*(\mcT_{P})$,  
 is  spanned by the global sections because 
 the following composite is surjective:
 \begin{align}
 \mcO_{P} \otimes \mfg \left(\cong \pi_{G*}(\mcO_{G} \otimes \mfg)^H\right)  
 \isom
  \pi_{G*} (\mcT_{G})^H \xrightarrow{d_{\mcE_H^\mr{triv}}} \mcT_{P},
 \end{align}
 where the first arrow arises from $\omega_G^\triangleright$.
  Hence, by (\ref{Ed002}),  $f^*(\mcT_{X})$ turns out to be  spanned by the global sections, i.e., 
  $f$ is free.
  This completes the proof of the assertion.
   \end{proof}
  %------------------------------------------------------------
  \SSP

 Next,  a connected $k$-scheme $X$ is called {\it  simply connected} if  the \'{e}tale fundamental group 
   of $X$ is trivial.
   Then, we obtain  the following assertion, generalizing  ~\cite[Theorem C, (ii)]{Wak6}.

    \SSP
   %------------------------------------------------------------
 \bpr \label{etoiuoa893}
 Let $X$ be a simply connected, projective, and  smooth scheme  over $k$.
 Suppose that $G$ is affine.
  Then, the following assertions hold:
 \begin{itemize}
 \item[(i)]
 Suppose that $X/k$ admits  an $F^\infty$-$(G, P)$-structure.
 Then, there exists  a finite \'{e}tale covering  $X \migi P$ over $k$.
 In particular, if, moreover, $P$ is  simply connected, 
 then $X$ is isomorphic to $P$.
 \item[(ii)]
Suppose that $P$ is affine
  and not isomorphic to $\mr{Spec}(k)$.
  Then, there are  no   $F^\infty$-$(G, P)$-structures on $X/k$.
 \end{itemize}
  \epr
  %------------------------------------------------------------
  \begin{proof}
  We shall consider the former assertion of (i).
  Let $\mcS^\hh := \{ \mcS_\N^\hh \}_\N$ be an $F^\infty$-$(G, P)$-structure  on $X/k$. 
  Denote by $\mcG_\infty := \{ (\mcG_l, \varpi_l) \}_{l \in \mbZ_{\geq 0}}$ the $F$-divided $G$-bundle induced by $\mcS^\hh$, i.e., the underlying $F$-divided $G$-bundle of the $F^\infty$-Cartan geometry corresponding to $\mcS^{\hh}$.
Once one chooses  a morphism of algebraic $k$-groups $G \migi \mr{GL}_n$ ($n \geq 1$),
 $\mcG_\infty$ induces 
  an $F$-divided vector bundle  $\mcV_\infty$ via change of structure group by this morphism.
But, according to ~\cite[Theorem 1.1]{EM},   the simple connectedness of $X$  implies that   there are no nontrivial $F$-divided vector bundles on $X$, so $\mcV_\infty$ must be  trivial.
Hence, by ~\cite[\S\,3.4, Lemma 3.5, (ii)]{BS}, 
the $G$-bundle $\mcG_1$ is trivial, or equivalently,
there exists 
 a global section of $\mcS_1^\hh$. 
 (Note that the result in {\it loc.\,cit.} is asserted under the assumption that $X$ is separably rationally connected.
 But, this assumption was not used in its proof.)
 This global section specifies an \'{e}tale $k$-morphism $f : X \migi P$; it is finite since 
$X$ is projective.
This completes the proof of  the former assertion (i).
Moreover, both the latter assertion of (i)  and  assertion (ii)  follow immediately from this fact.
       \end{proof}
  %------------------------------------------------------------

%-------------------------------------------------------------
\LSP
\subsection{$F^\N$-$(\mbG_m^{\times r_m}\times \mbG_a^{\times r_a},  \{e \})$-structures}
 \label{S0910}

This subsection is devoted to studying
the case where $G$ is  the product of  finite copies of $\mbG_m$ and $\mbG_a$ (and $H = \{ e \}$). 

%-------------------------------------------------------------
\subsubsection{} \label{SSS09}
Let $\N$ be an element of $\mbZ_{>0}$.
First, we discuss  $(\N-1)$-PD stratified $\mbG_m$-bundles equipped with an $\{ e\}$-reduction.
Denote by $C$ the Cartier operator $C: F_{X/S*}(\mr{Ker}(\nabla_0^{(2)})) \migi \Omega_{X^{(1)}/S}$ (cf. (\ref{E0124}) for the definition of $(-)^{(2)}$).
Also, denote by 
\begin{align} \label{EQWS23}
\Gamma (X, \Omega_{X/S})^{C-1}
\end{align}
 the submodule of $\Gamma (X, \mr{Ker}(\nabla_0^{(2)})) \left(\subseteq \Gamma (X, \Omega_{X/S}) \right)$ consisting of closed $1$-forms $\omega$ with $C (\omega) = F^{*}_S(\omega) \left(:= (\mr{id}_X \times F_S)^*(\omega) \right)$.
Recall from ~\cite[\S\,7, Corollary 7.1.3]{Katz2} that, for a $1$-form $\omega$ on $X/S$,  the $S$-connection  $\nabla_\omega := d + \omega$ is $p$-flat if and only if $\omega \in \Gamma (X, \Omega_{X/S})^{C-1}$.

Now, let us consider a triple
\begin{align} \label{Er560}
\omega^+ := (\omega, \mcL, \eta),
\end{align} 
where $\omega$ denotes a $1$-form on $X/S$, $\mcL$ denotes a line bundle on $X^{(\N)}$, and $\eta$ denotes an isomorphism  of flat line bundles $(\mcO_X, \nabla_\omega) \isom (F^{(\N)}_{X/S}(\mcL), \nabla^{\mr{can} (\N)}_\mcL)$ (cf. (\ref{Ew9090})).
The $p$-flatness of $\nabla_\mcL^{\mr{can}(\N)}$ implies that of $\nabla_\omega$.
Hence, the above discussion shows that $\omega$ automatically belongs to $\Gamma (X, \Omega_{X/S})^{C-1}$.
Given two such triples $\omega^+_i := (\omega_i, \mcL_i, \eta_i)$ ($i=1,2$),
we define an isomorphism $\omega_1^+ \isom \omega_2^+$ to be an isomorphism $\iota : \mcL_1 \isom \mcL_2$ with $F^{(\N)}_{X/S}(\iota) \circ \eta_1 = \eta_2$.
Since the morphism
$\mcO_{X^{(\N)}}^\times  \migi F_{X/S*}^{(\N)}(\mcO_X^\times)$ induced by $F_{X/S}^{(\N)}$ is injective,  
an isomorphism between $\omega_1^+$ and $\omega_2^+$ is, if it exists,  uniquely determined.
We shall denote by 
\begin{align} \label{eghiti22}
B_{m, X/S}^{(\N)}
\end{align}
the set consisting of isomorphism classes  of   triples $(\omega, \mcL, \eta)$ as above.
By Proposition \ref{Ertwq}, (i), and the last comment in Example \ref{EroDFG}, (i), the set $B_{m, X/S}^{(\N)}$
 is naturally identified with  the set of isomorphism classes of $(\N-1)$-PD stratified $\mbG_m$-bundles on $X/S$ with vanishing $p$-$(\N-1)$-curvature whose underlying $\mbG_m$-bundle  is trivial.
If $\N =1$, then
we obtain the  canonical bijection
\begin{align} \label{Efgt3}
B_{m, X/S}^{(1)} \isom \mr{Coker} \left( \Gamma (X, \mcO_X^\times) \xrightarrow{u \mapsto d \mr{log}(u)}\Gamma (X, \Omega_{X/S})^{C-1}\right)
\end{align}
given by $(\omega, \mcL, \eta) \mapsto \omega$.

If $\N'$ is another positive integer with $\N \leq \N'$ and $\omega^+ := (\omega, \mcL, \eta) \in B_{m, X/S}^{(\N')}$, then the $\N$-th truncation of $\omega^+$ is defined as the element
\begin{align} \label{ehos930}
\omega^+ |^{\langle \N \rangle} := (\omega, F^{(\N'-\N)*}_{X^{(\N)}/S}(\mcL), \eta)
\end{align}
of $B_{m, X/S}^{(\N)}$.
Truncations to lower levels  gives rise to  a projective system 
\begin{align} \label{r04so}
\cdots  \migi B_{m, X/S}^{(\N)} \migi \cdots \migi B_{m, X/S}^{(2)} \migi B_{m, X/S}^{(1)}.
\end{align}

 %-------------------------------------------------------------
\subsubsection{}\label{SSS10}

Next, 
we shall  consider $(\N-1)$-PD stratified $\mbG_a$-bundles equipped with an $\{ e\}$-reduction.
Recall  from  the discussion in Example \ref{EroDFG} that $\mbG_a$ is identified with  the subgroup of $\mr{GL}_2$ consisting of matrices of the form 
$\begin{pmatrix} 1 & a \\ 0 & 1 \end{pmatrix}$ for  $a \in \mbG_a$.
Denote by 
\begin{align} \label{EQAV34}
\Gamma (X, \Omega_{X/S})^C
\end{align}
 the submodule of $\Gamma (X, \mr{Ker}(\nabla_0^{(2)})) \left( \subseteq \Gamma (X, \Omega_{X/S})\right)$ consisting of closed $1$-forms $\chi$ with $C (\chi) = 0$.
Here, let us observe the following assertion.

\SSP
%-------------------------------------------------------------
\ble \label{L03365}
Let   $\chi$  be a closed $1$-form on $X/S$, i.e., an element of $\Gamma (X, \mr{Ker}(\nabla_0^{(2)}))$.
Then, the $p$-curvature ${^p}\psi_{\nabla_\chi^\dagger}$ of $\nabla^\dagger_\chi$ (cf. (\ref{ZZAei})) is given by 
\begin{align}
{^p}\psi_{\nabla^\dagger_\chi} (F_X^{-1}(\partial)) = -\begin{pmatrix} 0 & F_{X/S}^{-1}(\langle F_S^*(\partial), C  (\chi)\rangle)   \\ 0 & 0 \end{pmatrix} 
\end{align} 
for any local section $\partial \in \mcT_{X/S}$.
In particular, $\nabla^\dagger_\chi$ has vanishing $p$-curvature if and only if $\chi \in \Gamma (X, \Omega_{X/S})^C$.
\ele
%-------------------------------------------------------------
\begin{proof}
Let us take a local section $\partial \in \mcT_{X/S}$.
Then, we have the following sequence of equalities:
\begin{align}
{^p}\psi_{\nabla^\dagger_\chi} (F_X^{-1}(\partial)) &= \nabla (\partial)^p - \nabla (\partial^p) \\
& = \left(\partial^{\oplus 2}+ \begin{pmatrix} 0 & \langle \partial, \chi\rangle \\ 0 & 0 \end{pmatrix} \right)^p - (\partial^p)^{\oplus 2} - \begin{pmatrix} 0 & \langle \partial^p,  \chi\rangle \\ 0 & 0 \end{pmatrix}  \notag \\
& =
\left((\partial^p)^{\oplus 2}  + \begin{pmatrix} 0 &  \partial^{p-1}(\langle  \partial, \chi\rangle) \\ 0 & 0 \end{pmatrix} \right)  - (\partial^p)^{\oplus 2} - \begin{pmatrix} 0 & \langle  \partial^p, \chi\rangle \\ 0 & 0 \end{pmatrix} \notag  \\
& =  \begin{pmatrix} 0 & \partial^{p-1}(\langle \partial, \chi \rangle) - \langle \partial^p, \chi \rangle \\ 0& 0 \end{pmatrix} \notag \\
& = - \begin{pmatrix} 0 & F_{X/S}^{-1}(\langle F_S^*(\partial), C (\chi) \rangle) \\ 0 & 0 \end{pmatrix}, \notag
\end{align}
where the third equality follows from Jacobson's formula and the last equality follows from ~\cite[\S\,7, (7.1.2.6)]{Katz2}.
This completes the proof of this lemma.
\end{proof}
%-------------------------------------------------------------
\SSP

Now, we shall consider a triple
\begin{align}
\chi^+ := (\chi, \mcV, \lambda)
\end{align}
consisting of 
 a $1$-form $\chi$ on $X/S$,  a rank $2$ vector bundle $\mcV$ on $X^{(\N)}$ equipped with a structure of extension $0 \migi \mcO_{X^{(\N)}}\migi \mcV \migi \mcO_{X^{(\N)}} \migi 0$, and an isomorphism  of flat bundles
$\lambda : (\mcO_X^{\oplus 2}, \nabla_\chi^\dagger) \isom (F_{X/S}^{(\N)*}(\mcV), \nabla^{\mr{can}(\N)}_\mcV)$ preserving the extension structures.
The $p$-flatness of $\nabla_\mcV^{\mr{can}(\N)}$ implies that of $\nabla_\omega^\dagger$, so the above lemma  shows  that $\chi$  belongs to $\Gamma  (X, \Omega_{X/S})^C$.
An isomorphism between such triples can be defined naturally and is, if it exists,  uniquely determined for the same reason as the case of $(G, H) = (\mbG_m, \{ e\})$.
Denote by 
\begin{align}
B_{a, X/S}^{(\N)}
\end{align}
 the set of isomorphism classes of triples $(\chi, \mcV, \lambda)$ as above.
 By Proposition \ref{Ertwq}, (i), and the last comment in Example \ref{EroDFG}, (ii),  this set is naturally identified with the set of isomorphism classes of $(\N-1)$-PD stratified $\mbG_a$-bundles on $X/S$ with vanishing $p$-$(\N-1)$-curvature whose underlying $\mbG_a$-bundle is trivial.
 If $\N =1$, then we obtain the  canonical bijection
\begin{align} \label{Efgt4}
B_{a, X/S}^{(1)} \isom \mr{Coker} \left(\Gamma (X, \mcO_X) \xrightarrow{u \mapsto du} \Gamma (X, \Omega_{X/S}) \right)
\end{align}
given by $(\chi, \mcV, \lambda) \mapsto \chi$.
 
 Just as in the case of $B_{m, X/S}^{(\N)}$,
 the $\N$-th truncation $\chi^+ |^{\langle \N \rangle}$ of $\chi^+ \in  B_{a, X/S}^{(\N')}$ with $\N \leq \N'$ can be  defined. 
In particular, truncations to various lower levels yield a projective system consisting of $B_{a, X/S}^{(\N)}$'s.

%-------------------------------------------------------------
\subsubsection{}\label{SSS11}

Finally, we combine the above discussions to study $(\N-1)$-PD stratified bundles whose structure group is the product of finite copies of  $\mbG_m$ and $\mbG_a$.
Let $s$ and $t$ be nonnegative integers with $s + t = n \left(:= \mr{dim}(X/S) \right)$.
Write $G:= \mbG_m^{\times s} \times \mbG_a^{\times t}$.
The product of sets $(B_{m, X/S}^{(\N)})^{\times s} \times (B_{a, X/S}^{(\N)})^{\times t}$ is in bijection  with 
the set of isomorphism classes of $(\N-1)$-PD stratified $G$-bundles with vanishing $p$-$(\N-1)$-curvature   equipped  with   an $\{e \}$-reduction (i.e., whose underlying $G$-bundles are trivial).
Denote by
\begin{align} \label{eo4pfdc}
B_{X/S}^{(\N), s, t}
\end{align}
the elements  of 
$(B_{m, X/S}^{(\N)})^{\times s} \times (B_{a, X/S}^{(\N)})^{\times t}$
represented by  collections 
$(\omega^+_1, \cdots, \omega_s^+, \chi_1^+, \cdots, \chi_t^+)$, where $\omega_i^+ := (\omega_i, \mcL_i, \eta_i)$ and $\chi_j^+ := (\chi_j, \mcV_j, \lambda_j)$ ($i=1, \cdots, s$ and $j = 1, \cdots, t$),
 such that the global sections  $\omega_1, \cdots, \omega_{s}, \chi_1, \cdots, \chi_{t}$ forms a basis of  the $\Gamma (X, \mcO_X)$-module $\Gamma (X, \Omega_{X/S})$.
Truncations to lower levels defined before give rise to  a projective system
consisting of $B_{X/S}^{(\N), s, t}$'s for various $\N$.
Hence, we obtain the set
\begin{align} \label{r4k39485}
B_{X/S}^{(\infty), s, t} := \varprojlim_{\N \in \mbZ_{>0}} B_{X/S}^{(\N), s, t}.
\end{align}

\SSP
%-------------------------------------------------------------
\bpr \label{L0045}
Let us keep the above notation.
\begin{itemize}
\item[(i)]
For any $\N \in \mbZ_{>0}\sqcup \{ \infty \}$, there exists a canonical bijection   of sets
\begin{align} \label{Gheit}
F^\N\text{-}\mr{Ehr}_{(G, G), X/S} \isom B_{X/S}^{(\N), s, t}.
\end{align}
%whose  formation  commutes with truncation to lower levels.
In particular, 
$X/S$ admits  no $F^\N$-$(G, G)$-structures
 unless $\Gamma (X, \Omega_{X/S}) \cong \Gamma (X, \mcO_X)^{\oplus n}$.
\item[(ii)]
The following two conditions are equivalent to each other:
\begin{itemize}
\item[(a)]
There exists an $F^1$-$(G, G)$-structure on $X/S$;
\item[(b)]
There exists a basis  $\{ \omega_1, \cdots, \omega_s, \chi_1, \cdots, \chi_t \}$ of $\Gamma (X, \Omega_{X/S})$ such that
each $\omega_i$ ($i=1, \cdots, s$) belongs to $\Gamma (X, \Omega_{X/S})^{C-1}$
and each $\chi_j$ ($j =1, \cdots, t$) belongs to $\Gamma (X, \Omega_{X/S})^C$.
\end{itemize}
Moreover, if $S = \mr{Spec}(k)$ and  $X$ is proper over $k$, then these conditions are equivalent to the following condition:
\begin{itemize}
\item[(c)]
The following equalities hold:
 \begin{align}
 n = \mr{dim}_k (\Gamma (X, \Omega_{X/k})), \hspace{3mm}  
 s = \mr{dim}_{\mbF_p}(H^1_{\mr{fl}} (X, \mu_p)),  \hspace{3mm}
 t = \mr{dim}_k (H^1_{\mr{fl}}(X, \alpha_p)),
 \end{align}
where $H^1_{\mr{fl}}$ denotes the  $1$-st flat cohomology.
\end{itemize}

\end{itemize}
\epr
%-------------------------------------------------------------
\begin{proof}
Let us consider assertion (i).
The problem is reduced to the case of $\N \in \mbZ_{>0}$ because the assertion  for $\N = \infty$ follows immediately  from this case together with (\ref{rosao3048}) and    (\ref{r4k39485}).
First, we shall construct a map of sets
$F^\N\text{-}\mr{Ehr}_{(G, G), X/S} \migi B_{X/S}^{(\N), s, t}$.
Take an $F^\N$-$(G, G)$-structure  $\mcS^\hh$ on $X/S$, and denote by $\mcE^\sss_\N$
the dormant indigenous $(G, \{ e \})$-bundle  of level $\N$ 
corresponding to $\mcS^\hh$ (cf. Theorem \ref{TheoremA}).
As mentioned above,  the data $\mcE^\sss_\N$ gives 
  a collection of data
\begin{align} \label{L0909}
(\omega_1^+, \cdots, \omega_s^+, \chi_1^+, \cdots, \chi^+_t),
\end{align}
 where each $\omega^+_i := (\omega_i, \mcL_i, \eta_i)$
 ($i=1, \cdots, s$) denotes an element of $B_{m, X/S}^{(\N)}$ and each $\chi^+_j := (\chi_j, \mcV_j, \lambda_j)$ ($j=1, \cdots, t$) denotes an element of $B_{a, X/S}^{(\N)}$. 
Under the natural identifications $\mr{Lie}(G) \cong  \mr{Lie}(\mbG_m)^{\oplus s} \oplus \mr{Lie}(\mbG_a)^{\oplus t} \cong k^{\oplus (s +t)}$, 
the vector bundle $\mr{Lie}(G)_{X \times G}/\mr{Lie}(\{ e \})_{X \times \{e \}}$ may be identified with $\mcO_X^{\oplus (s +t)}$.
In particular, if $\nabla$ denotes the $S$-connection on the trivial $G$-bundle 
defining $\mcE^\sss_\N |^{\langle 1 \rangle}$ (cf. (\ref{QQWv4})), then
 the  Kodaira-Spencer map $\mr{KS}_{(X \times \{ e \}, \nabla)}$
 may be regarded as an $\mcO_X$-linear morphism $\mcT_{X/S} \migi \mcO_X^{\oplus (s +t)}$, or equivalently, an element of $\Gamma (X, \Omega_{X/S})^{\oplus (s+t)}$ (cf. Remark \ref{R13404}).
 By taking account of 
 the discussions in Example \ref{EroDFG}, (i) and (ii),
 we see that
 $\mr{KS}_{(X \times \{ e \}, \nabla)}$ coincides with $(\omega_1, \cdots, \omega_s, \chi_1, \cdots, \chi_t)$.
Since 
 $\mr{KS}_{(X \times \{e \}, \nabla)}$ is an isomorphism,
 the global sections  $\omega_1, \cdots, \omega_s, \chi_1, \cdots, \chi_t$
 form a basis of $\Gamma (X, \Omega_{X/S})$.
Thus, the collection of data
(\ref{L0909})
   specifies an element of $B_{X/S}^{(\N), s, t}$.
The resulting assignment $\mcS^\hh \mapsto (\omega_1^+, \cdots, \omega_s^+, \chi_1^+, \cdots, \chi^+_t)$ defines a map of sets
$F^\N\text{-}\mr{Ehr}_{(G, G), X/S} \migi B_{X/S}^{(\N), s, t}$.
Moreover,   by  the above discussion,  
 this map is verified to be bijective, so  this completes the proof  of assertion (i).
 
 Next, let us consider assertion (ii).
 The equivalence (a) $\Leftrightarrow$ (b) 
  follows immediately from assertion (i) together with (\ref{Efgt3}) and (\ref{Efgt4}).
 The equivalence (b) $\Leftrightarrow$ (c) under the assumption mentioned in the statement  follows from the fact that $\Gamma (X, \Omega_{X/k})^{C-1}$,  $\Gamma (X, \Omega_{X/k})^C$ are naturally  isomorphic to the flat cohomology groups $H^1_\mr{fl}(X, \mu_p)$,  $H^1_\mr{fl}(X, \alpha_p)$ respectively  (cf. ~\cite[Chap.\,III, \S\,4, Proposition 4.14]{Mil4}).
 (Note that the natural $k$-linear morphism $k \otimes_{\mbF_p} \Gamma(X, \Omega_{X/k})^{C-1} \migi \Gamma (X, \Omega_{X/k})$ is verified to be injective.)
 This completes the proof of the assertion.
\end{proof}
%-------------------------------------------------------------
\SSP

In the following examples,  suppose that $\N$ is a positive integer.

\SSP
%-------------------------------------------------------------
\begin{exa}[]
Let us consider $F^\N$-$(\mbG_m, \mbG_m)$-structures on  the $k$-curve $\mbG_m := \mr{Spec}(k[T^{\pm 1}])$.
For each  $u \in \Gamma (\mbG_m, \mcO^\times_{\mbG_m}) \left(= \{ a \cdot T^{n} \, | \, n \in \mbZ, a \in k^\times\} \right)$,
we shall write $\mr{mult}_u$ for the automorphism of $\mcO_{\mbG_m}$ given by multiplication by $u$.
The pull-back  $\mr{mult}_u^*(\nabla_0)$ of the trivial connection $\nabla_0$ on $\mcO_{\mbG_m}$ via $\mr{mult}_u$ is the $k$-connection expressed as
 $\mr{mult}_u^*(\nabla_0) = \nabla_{d\mr{log}(u)}$, where $d \mr{log}(u) := u^{-1} \cdot du$.
Hence, 
since any line bundle on $\mbG_m^{(\N)}$ is trivial,
each element of $B_{m, \mbG_m/k}^{(\N)}$ may be represented by
$(d\mr{log}(u), \mcO_{\mbG_m^{(\N)}}, \mr{mult}_u)$ for some   $u \in \Gamma (\mbG_m, \mcO^\times_{\mbG_m})$.
Such a triple  belongs to  
$B_{\mbG_m/k}^{(\N), 1, 0} \left(\subseteq B_{m, \mbG_m/k}^{(\N)} \right)$
if and only if $d \mr{log}(u)$ is nowhere vanishing, i.e.,
$u \in  \{ a \cdot T^{n} \, | \, n \in \mbZ, (p, n)=1, a \in k^\times\}$.
Also, an isomorphism from $(d\mr{log}(u), \mcO_{\mbG_m^{(\N)}}, \mr{mult}_u)$ to $(d\mr{log}(v), \mcO_{\mbG_m^{(\N)}}, \mr{mult}_v)$ is given by the automorphism of $\mcO_{\mbG_m^{(\N)}}$ 
defined as multiplication by an element 
$w \in \Gamma (\mbG_m^{(\N)}, \mcO_{\mbG_m^{(\N)}}^\times)$ 
$\left(= \{ a \cdot T^{n \cdot p^\N} \, | \, n \in \mbZ, a \in k^\times\} \right)$ satisfying  $u = w \cdot v$.   
This implies that  each isomorphism class in $B_{\mbG_m/k}^{(\N), 1, 0}$
may be represented by a triple $(d \mr{log}(u), \mcO_{\mbG_m^{(\N)}}, \mr{mult}_u)$  for  $u  = T^m$ with  $1 \leq m \leq p^{\N}-1$ and $(p, m) =1$.
Consequently,
the assignment $m \mapsto (d\mr{log}(T^m), \mcO_{\mbG_m^{(\N)}}, \mr{mult}_{T^m})$ defines a canonical bijection 
\begin{align} \label{L0223}
\left\{ m \in \mbZ \, | \, 1 \leq m \leq p^\N-1, (p, m) =1  \right\} 
\isom  B_{\mbG_m/k}^{(\N), 1,0}.
\end{align}
In particular, $\mbG_m/k$ admits an $F^\N$-$(\mbG_m, \mbG_m)$-structure for any $\N$.
%, including the case of $\N = \infty$.
\end{exa}
%-------------------------------------------------------------
\SSP

%-------------------------------------------------------------
\begin{exa} \label{R0034}
Next,  we shall consider $F^\N$-$(\mbG_a, \mbG_a)$-structures on  the affine line  $\mbA^1 := \mr{Spec}(k[T])$.
For any $u \in \Gamma (\mbA^1, \mcO_{\mbA^1})$, we shall write $\mr{add}_u$
for the automorphism of $\mcO_{\mbA^1}^{\oplus 2}$ given by multiplication by 
$\begin{pmatrix} 1 & u \\ 0 & 1\end{pmatrix}$.
The pull-back  $\mr{add}^{*}_u (\nabla_0)$ of the trivial connection $\nabla_0$ on $\mcO_{\mbA^1}$  via $\mr{add}_u$ is the $k$-connection expressed as 
 $\mr{add}_u^*(\nabla_0) = \nabla_{d u}^\dagger$.
Hence, since any rank $2$ vector bundle on $\mbA^{1(\N)}$ is trivial, each element of $B_{a, \mbA^1/k}^{(\N)}$ may be represented by
$(du, \mcO^{\oplus 2}_{\mbA^{1(\N)}}, \mr{add}_u)$ for $u \in \Gamma (\mbA^1, \mcO_{\mbA^1})$.
Such a triple 
 belongs to  $B_{\mbA^1/k}^{(\N), 0,1} \left(\subseteq B_{a, \mbA^1/k}^{(\N)} \right)$
if and only if $d u$ is nowhere vanishing, i.e., $u \in \left\{ a \cdot T +  f(T) \, | \, a \in k^\times, f(T) \in k[T^p] \right\}$.
Also, an isomorphism from $(du, \mcO^{\oplus 2}_{\mbA^{1(\N)}}, \mr{add}_u)$ to $(dv, \mcO^{\oplus 2}_{\mbA^{1(\N)}}, \mr{add}_v)$
is given by the automorphism of $\mcO_{\mbA^{1 (\N)}}^{\oplus 2}$
defined as multiplication by $\begin{pmatrix} 1 & w \\ 0 & 1\end{pmatrix}$
for some $w \in \Gamma (\mbA^{1 (\N)}, \mcO_{\mbA^{1 (\N)}})  \left(= k [T^{p^\N}] \right)$ satisfying 
 $u = v + w$.
This implies that each isomorphism class in $B_{\mbA^1/k}^{(\N), 0,1}$
may be represented by a triple $(du, \mcO^{\oplus 2}_{\mbA^{1(\N)}}, \mr{add}_u)$ for 
$u =  a \cdot T + \sum_{i \geq 0 } a_i \cdot  T^{ i}  \in k[T]$ such that  $a \in k^\times$ and 
$a_{i} \neq  0$  only if $i \in p \mbZ \setminus p^\N \mbZ$.
%$u \in \{ a \cdot T + \sum_{i \geq 0 } a_i \cdot  T^{ i}  \in k[T] \, | \, a \in k^\times \ \text{and} \  a_{i} \neq  0 \ \text{only if} \  i \in p \mbZ \setminus p^\N \mbZ \}$.
Consequently, the assignment $u \mapsto (du, \mcO_{\mbA^{1(\N)}}, \mr{add}_u)$ defines a canonical bijection
\begin{align}
 \left\{ a \cdot T + \sum_{i \geq 0 } a_i \cdot  T^{i}  \in k[T] \, \Bigg| \,  a \in k^\times \ \text{and} \  a_{i} \neq  0 \ \text{only if} \  i \in p \mbZ \setminus p^\N \mbZ \right\}
\isom B_{\mbA^1/k}^{(\N), 0,1}.
\end{align}
In particular,  the affine line $\mbA^1/k$ admits an $F^\N$-$(\mbG_a, \mbG_a)$-structure for any $\N$.
%, including the case of $\N = \infty$.
\end{exa}
%-------------------------------------------------------------

\LSP
%-------------------------------------------------------------
%\vspace{5mm}
\subsection{$F^\N$-$(\mbG_m^{\times s}\times \mbG_a^{\times t},  \{e \})$-structures on Abelian varieties}
 \label{S091k0}

%In this final subsection, we examine  $F^\N$-$(G_m^{\times s} \times \mbG_a^{\times t}, \{ e \})$-structures on Abelian varieties.
Let $X$ be an Abelian variety over $k$ of dimension $n$.
For each positive integer $l$, denote by $l_X$ the endomorphism of $X$ given by multiplication by $l$.
Let $\widehat{X}$ denote the dual Abelian variety of $X$.
Recall that the {\it $p$-rank}  of $X$ is  the nonnegative integer  $\mu (X) := \mr{dim}_{\mbF_p}(\mr{Hom}(\mu_p, X))\left(= \mr{dim}_{\mbF_p}(\mr{Hom}(\mu_p, \mr{Ker}(p_X))) \right)$, which coincides with the value $\mr{dim}_{\mbF_p}(\mbZ/p\mbZ, X)$.
Also, 
the {\it $a$-number} of $X$  is defined as  $\alpha (X) := \mr{dim}_k (\mr{Hom}(\alpha_p, X)) \left(=  \mr{dim}_k (\mr{Hom}(\alpha_p, \mr{Ker}(p_{X}))) \right)$.
These values satisfy the inequalities $0 \leq \mu (X)  \leq g$ and $1 \leq \alpha (X) + \mu (X) \leq g$. 

Regarding  $F^\N$-$(\mbG_m^{\times s} \times \mbG_a^{\times t}, \{ e\})$-structures for $\N =1$,   Proposition \ref{L0045} implies the following assertion.

\SSP
%-----------------------------------------------------------------------
\bpr \label{Erfgy}
Let us keep the above notation.
Also, 
let $(s, t)$ be a pair of nonnegative integers with $n = s + t$.
Then, 
the Abelian variety  $X$ admits an $F^1$-$(\mbG_m^{\times s} \times \mbG_a^{\times t}, \{ e\})$-structure if and only if $s = \mu (\widehat{X})$ and $t = \alpha (\widehat{X})$.
\epr
%-----------------------------------------------------------------------
\begin{proof}
Observe that
\begin{align}
\mr{dim}_{\mbF_p}(H^1_{\mr{fl}} (X, \mu_p)) & = \mr{dim}_{\mbF_p}(\mr{Hom}(\varprojlim_{l} \mr{Ker}(l_X), \mu_p))
\\
&= \mr{dim}_{\mbF_p}(\mr{Hom}(\mr{Ker}(p_X), \mu_p))
\notag \\
& = \mr{dim}_{\mbF_p} (\mr{Hom} (\widehat{\mu_p}, \widehat{\mr{Ker}(p_X)}))
 \notag  \\
& = \mr{dim}_{\mbF_p}(\mbZ/p\mbZ, \mr{Ker}(p_{\widehat{X}}))\notag \\
&= \mu (\widehat{X}), \notag
\end{align}
where $\widehat{(-)}$ denotes Cartier dual and  the first equality follows from ~\cite[Chap.\,III, \S\,4, Corollary 4.19, (a)]{Mil4}.
By a similar discussion together with $\widehat{\alpha_p} = \alpha_p$, we have $\mr{dim}_{k}(H^1_{\mr{fl}} (X, \alpha_p)) = \alpha (\widehat{X})$. 
Thus,  
the assertion follows from the equivalence (a) $\Leftrightarrow$ (c) asserted in Proposition \ref{L0045}, (ii).
\end{proof}
%-----------------------------------------------------------------------
\SSP

  Recall that  the Abelian variety $X$ is called {\it ordinary} (resp., {\it superspecial}) if $\mu (X) =g$ (resp.,  $\alpha (X) = g$).
It is verified that $X$ is ordinary (resp., superspecial) if and only if the dual $\widehat{X}$ is ordinary (resp., superspecial).
The following assertion concerns   the cases where  $G = \mbG_m^{\times n}$ and $G = \mbG_a^{\times n}$. 
   
   \SSP
 %------------------------------------------------------------
 \bpr \label{C0558}
\begin{itemize}
\item[(i)]
Let $\N \in \mbZ_{>0}\sqcup \{ \infty \}$.
Then, 
$X$ admits an $F^\N$-$(\mbG_m^{\times n}, \mbG_m^{\times n })$-structure  
if and only if $X$ is ordinary.
Moreover, when $X$ is  ordinary,
the map of sets
\begin{align} \label{rpap459}
F^\N\text{-}\mr{Ehr}_{(\mbG_m^{\times n}, \mbG_m^{\times n}), X/k}/\mfS_n \migi F^\N \text{-}\mr{Ehr}_{(\mr{Aff}_n, \mbA^n), X/k}
\end{align}
defined in Example \ref{es04xm} is injective, and moreover  bijective if $\N = \infty$.
\item[(ii)]
There exists  an  $F^1$-$(\mbG_a^{\times n}, \mbG_a^{\times n})$-structure  on $X/k$
if and only if $X$ is superspecial.
Moreover, for any $\N \geq 2$, there are no $F^\N$-$(\mbG_a^{\times n}, \mbG_a^{\times n})$-structures on any Abelian variety.
\end{itemize}
 \epr
 %------------------------------------------------------------
 \begin{proof}
 Let us consider the former  assertion of (i).
 The ``only if" part of the required equivalence  is a direct consequence of Proposition \ref{Erfgy} for $(s, t)=(n, 0)$ together with the equivalence ``$X$ is ordinary" $\Leftrightarrow$ ``$\widehat{X}$ is ordinary".
  To prove the ``if" part, 
  let us 
 suppose that $X$, or equivalently, $\widehat{X}$ is ordinary.
 Also,  for a positive integer $\N$, suppose that 
 we are given an $F^\N$-$(\mbG_m^{\times n}, \mbG_m^{\times n})$-structure on $X/k$;   it specifies an element  $(\omega_1^+, \cdots, \omega_n^+)$  of
 $B_{X/k}^{(\N), n, 0}$, where $\omega_i^+ := (\omega_i, \mcL_i, \eta_i)$ ($i=1, \cdots, n$),  via (\ref{Gheit}).
 Since the map  $\widehat{X^{(\N +1)}} (k) \migi \widehat{X^{(\N)}} (k)$
arising from  pull-back by $F_{X^{(\N)}/k} : X^{(\N)} \migi X^{(\N +1)}$  is surjective,
there exists a collection of line bundles $\mcL_1^\nabla, \cdots, \mcL_n^\nabla$ on $X^{(\N +1)}$ with $F^{*}_{X^{(\N)}/k}(\mcL_i^\nabla) \cong \mcL_i$ ($i=1, \cdots, n$).
The collection 
$(\widetilde{\omega}_1^+, \cdots, \widetilde{\omega}_n^+)$ consisting of 
triples $\widetilde{\omega}_i^+ := (\omega_i, \mcL_i^\nabla, \eta_i)$ ($i=1, \cdots, n$)
 defines an element of $B_{X/k}^{(\N +1), n, 0}$ whose $\N$-th truncation coincides with $(\omega_1^+, \cdots, \omega_n^+)$.
 This proves the surjectivity of  
 the truncation map  $B_{X/k}^{(\N +1), n, 0} \migi B_{X/k}^{(\N), n, 0}$.  
 Observe here that 
$B_{X/k}^{(1), n, 0}$ is nonempty 
because of the result of   Proposition \ref{Erfgy} and the ordinariness assumption on $X$.
Hence, we can apply the above argument inductively to conclude that   $B_{X/k}^{(\N), n, 0} \neq \emptyset$  for any $\N \in \mbZ_{>0} \sqcup \{ \infty \}$. 
That is to say,  $X/k$ admits  an  $F^\N$-$(\mbG_m^{\times n}, \mbG_m^{\times n})$-structure on $X/k$  via  (\ref{Gheit}).
This completes the former assertion of (i).

Next, let us consider the latter assertion of (i) under the assumption that $X$ is ordinary.
Suppose here that $\N \in \mbZ_{>0}$, and 
let us take an element
$(\omega_1^+, \cdots, \omega_n^+)$ of $B_{X/k}^{(\N), n, 0}$, where $\omega_{i}^+ := (\omega_i, \mcL_{i}, \eta_{i})$.
Write $\nabla := \nabla_0 \oplus  \bigoplus_{i=1}^{n} \nabla^{\mr{can}(\N-1)}_{\mcL_i}$, which specifies a $\mcD_{X/k}^{(\N-1)}$-module structure  on $\mcO_X^{\oplus (n+1)}$ via $\mr{id}_{\mcO_X} \oplus \bigoplus_{i=1}^{n}\eta_i$.
Since $\omega_1, \cdots, \omega_n$ forms a basis of $\Gamma (X, \Omega_X)$,
the following composite turns out to be an isomorphism:
\begin{align}
\mcD_{X/k, 1}^{(\N-1)} \left(= \mcD_{X/k, 1}^{(\N-1)} \otimes_{\mcO_X} \mcO_X \right) \migiincl \mcD_{X/k}^{(\N-1)} \otimes_{\mcO_X} \mcO^{\oplus (n+1)}_X \xrightarrow{\nabla} \mcO^{\oplus (n+1)}_{X},
\end{align}
 where the first arrow arises from both the natural inclusion $\mcD_{X/k, 1}^{(\N -1)} \migiincl \mcD_{X/k}^{(\N-1)}$ and the diagonal embedding $\mcO_X \migiincl \mcO_X^{\oplus (n+1)}$.
Under the identification $\mcD_{X/k, 1}^{(\N-1)} = \mcO^{\oplus (n+1)}_{X}$ given by this isomorphism, we shall regard $\nabla$ as a $\mcD_{X/k}^{(\N-1)}$-module structure on $\mcD_{X/k, 1}^{(\N)}$.
Also, we regard the $(n+1)$-st projection $\delta:= \mr{pr}_{n+1} : \mcO_X^{\oplus (n+1)} \migi \mcO_X$  as a left inverse morphism $\mcD_{X/k}^{(\N-1)} \migisurj \mcO_X$ of the natural inclusion $\mcO_X \left(:= \mcD_{X/k, 0}^{(\N-1)} \right) \migiincl \mcD_{X/k, 1}^{(\N-1)}$.
By construction, the pair $(\nabla, \delta)$ satisfies the conditions (c) and (d) described in Example \ref{R099}.
That is to say, 
$(\nabla, \delta)$ specifies an elmement of $F^\N\text{-}\mr{Ehr}_{(\mr{Aff}_n, \mbA^n), X/k}$.
If $B_{X/k}^{(\N), n, 0}/\mfS_n$ denotes the quotient set of $B_{X/k}^{(\N), n, 0}$ by the $\mfS_n$-action determined by permutation of factors, then
the resulting assignment $(\omega_1^+, \cdots, \omega_n^+) \mapsto (\nabla, \delta)$ induces a map of sets
 $\varsigma_\N : B_{X/k}^{(\N), n, 0}/\mfS_n \migi F^\N\text{-}\mr{Ehr}_{(\mr{Aff}_n, \mbA^n), X/k}$.
 One may verify that $\varsigma_\N$ coincides with  (\ref{rpap459}) under the identification 
 $B_{X/k}^{(\N), n, 0}/\mfS_n = F^\N\text{-}\mr{Ehr}_{(\mbG_m^{\times n}, \mbG_m^{\times n}), X/k}/\mfS_n$ induced by (\ref{Gheit}).
Then, the required  assertion for $\N \in \mbZ_{>0}$ follows from the injectivity of $\varsigma_\N$ resulting from    ~\cite[\S\,6, Theorem 6.4.1]{Wak6}.
Moreover, the case of  $\N = \infty$  follows from the bijectivity of $\varprojlim_{\N \in \mbZ_{>0}}\varsigma_\N$ asserted in  the same theorem.

Finally, we shall consider assertion (ii).
 The former assertion follows directly from Proposition \ref{Erfgy} for $(s, t) = (0, n)$
 together with the equivalence ``$X$ is superspecial" $\Leftrightarrow$ ``$\widehat{X}$ is superspecial".
 To prove  the latter assertion, it suffices, by the former assertion,  to consider the case where $X$  is superspecial.
 Suppose, on the contrary, that there exists an $F^\N$-$(\mbG_a^{\times n}, \mbG_a^{\times n})$-structure $\mcS^\hh$ on $X/k$ for some $\N \geq 2$; it specifies,  via  truncation to level $2$,   
 an element 
 $(\chi^+_1, \cdots, \chi_n^+)$, where $\chi_j^+ := (\chi_j, \mcW_j, \lambda_j)$ ($j =1, \cdots, n$), 
 of $B_{X/k}^{(2),  0, n}$   via (\ref{Gheit}).
Let us fix $j \in \{ 1, \cdots, n \}$ and  write $\mcV_j := F_{X^{(1)}/k}^*(\mcW_j)$.
 Also, denote by $v$ (resp., $w$)   the element of $H^1 (X^{(1)}, \mcO_{X^{(1)}})$ (resp., $H^1 (X^{(2)}, \mcO_{X^{(2)}})$) determined by the extension structure of $\mcV_j$ induced by that of $\mcW_j$ (resp.,  the extension structure of $\mcW_j$) under the natural identification 
 $H^1 (X^{(1)}, \mcO_{X^{(1)}}) \cong \mr{Ext}^1_{\mcO_{X^{(1)}}} (\mcO_{X^{(1)}}, \mcO_{X^{(1)}})$ (resp., $H^1 (X^{(2)}, \mcO_{X^{(2)}})\cong \mr{Ext}^1_{\mcO_{X^{(2)}}} (\mcO_{X^{(2)}}, \mcO_{X^{(2)}})$).
 Since the $1$-form  $\chi_j$ is nozero,  there  are no split injections $\mcO_X \migiincl \mcO_X^{\oplus 2}$ of the trivial  extension $0 \migi \mcO_X \migi \mcO_X^{\oplus 2} \migi \mcO_X \migi 0$ whose image is closed under $\nabla_{\chi_j}^\dagger$.
 It follows  that the extension structure of  $\mcV_j$ is nontrivial, or equivalently, $v \neq 0$.
 On the other hand, 
 since $X$, as well as $X^{(1)}$,  is superspecial,
 the morphism $F^*_{X^{(1)}/k} : H^1 (X^{(2)}, \mcO_{X^{(2)}}) \migi H^1 (X^{(1)}, \mcO_{X^{(1)}})$ induced by $F_{X^{(1)}/k}$ is zero.
 This implies  $0 = F^*_{X^{(1)}/k} (w) = v$, so we obtain a contradiction.
 This completes the proof of the latter assertion.
 \end{proof}
%------------------------------------------------------------

\SSP

%%%%%%%%%%%%%--[ begin  section1]---%%%%%%
\vspace{10mm}
\section{Classification for  algebraic curves} \label{SSS0g1}\vspace{3mm}

The purpose of this section is to examine the classification problem  of 
  Frobenius-Ehresmann structures on smooth proper  curves.
Recall   a result by  B. Laurent (cf. ~\cite[Theorem 1.1, (1)]{Lau}), asserting that   the list of homogenous curves, i.e., homogenous space of dimension $1$, may be  described as follows:

\begin{itemize}
\item[]
\begin{itemize}
\item[Case (a):]
$G = \mr{PGL}_2$, $H= \mr{PGL}_2^\circledcirc$;
\item[Case (b):]
$G = \mr{Aff}_1$, $H = \mr{Aff}_1^\circledcirc$;
\item[Case (c):]
$G = \mbG_m$, $H  = \{e \}$;
\item[Case (d):]
$G = \mbG_a$, $H = \{e \}$;
\item[Case (e):]
$G$ is an elliptic curve over $k$ and $H = \{e \}$.
\end{itemize}
\end{itemize}

In what follows, 
we determine which curves have Frobenius-Ehresmann structures with a fixed model $(G, H)$ according to this classification.
Let $X$ be a smooth proper  curve over $k$ of genus $g \geq 0$.

%-------------------------------------------------------------
\vspace{5mm}
\subsection{Case (a): } \label{S0f910}

 The first case is the pair of the projective linear group  $G = \mr{PGL}_2$ and its subgroup  $H =\mr{PGL}^\circledcirc_2$.
In several  literatures,  $F^\N$-$(\mr{PGL}_2, \mbP^1)$-structures  on curves have been studied  
 under the name of  {\it Frobenius-projective structure}  (cf. ~\cite{Hos2}, ~\cite{Wak6}).
If $\N=1$, then 
 the corresponding dormant indigenous $(\mr{PGL}_2, \mr{PGL}_2^\circledcirc)$-bundles are  called 
 {\it dormant indigenous bundles}  or {\it dormant $\mr{PGL}_2$-opers} (cf. ~\cite{Mzk2}, ~\cite{Wak}, ~\cite{Wak8}).
Moreover, these objects are closely  related   to second order differential operators   having a full set of roots functions (cf. ~\cite[Chap.\,4, \S\,4.10, Proposition 4.65]{Wak8}).
Here, let us pick up a previous result of ~\cite[Theorem E]{Wak6},  concerning the existence of  $F^\N$-$(\mr{PGL}_2, \mbP^1)$-structures.

\SSP
%-------------------------------------------------------------
\bpr \label{EorGHE}
For any $\N \in \mbZ_{>0}\sqcup \{ \infty \}$, 
$X$ admits an $F^\N$-$(\mr{PGL}_2, \mbP^1)$-structure.
 Moreover, if $g >1$,  $\N \in \mbZ_{>0}$, and $\Theta$ is a theta characteristic of $X$ (i.e., a line bundle on $X$ equipped with an isomorphism $\Theta^{\otimes 2} \isom \Omega_{X}$), then there exists a canonical bijection of sets
 \begin{align} \label{Er45}
 F^\N\text{-}\mr{Ehr}_{(\mr{PGL}_2, \mbP^1), X/k} \isom \mcQ_{\N, \Theta}^{2, \mr{triv}}(k),
 \end{align}
 where $\mcQ_{\N, \Theta}^{2, \mr{triv}}$ denotes the Quot-scheme over $k$ classifying  $\mcO_{X^{(\N)}}$-submodules $\mcF$ of $F_{X/k}^{(\N)*}(\Theta^\vee)$ with $\mr{rank}(\mcF) = 2$ and $\mr{det} (\mcF) \cong \mcO_{X^{(\N)}}$.
\epr
%-------------------------------------------------------------
\SSP

There are some results related to  the counting problem of $F^\N$-$(\mr{PGL}_2, \mbP^1)$-structures.
S. Mochizuki showed that 
there are only finitely many $F^1$-$(\mr{PGL}_2, \mbP^1)$-structures  on  $X$ (cf. ~\cite[Chap.\,II, \S\,2.3, Theorem 2.8]{Mzk2}).
Moreover,  
the main  result of  ~\cite[Theorem A]{Wak} explicitly computes   the number of such structures   by means of the genus ``$g$'' of $X$ and the characteristic  ``$p$''  of the base field $k$.
To be precise, if $X$ is sufficiently general in the moduli stack of smooth proper genus $g$ curves, then
we have the equality
\begin{align} \label{eraq345}
F\text{-}\mr{Ehr}_{(\mr{PGL}_2, \mbP^1), X/k} = \frac{p^{g-1}}{2^{2g-1}} \cdot \sum_{\theta =1}^{p-1} \frac{1}{\sin^{2g-2} \big( \frac{\pi \cdot \theta}{p}\big)}.
\end{align}
This formula is proved  by using (\ref{Er45}) for $\N =1$ and  applying a formula computing the degree of 
certain Quot-schemes (i.e., a special case  of the Vafa-Intriligator formula).
However, for  large  $\N$'s, 
the author  does not know much about  the finiteness of $F^\N$-$(\mr{PGL}_2, \mbP^1)$-structures at the time of writing the present paper.
 So far, 
  the finiteness of  the set $F^\N\text{-}\mr{Ehr}_{(\mr{PGL}_2, \mbP^1), X/k}$ with $\N \geq 2$ has been shown only  for genus $2$ curves  (cf. ~\cite[\S\,7,  Theorem 7.4.2]{Wak6}).
By taking account of 
  (\ref{eraq345}) and its proof, we expect that the following conjectural formula will hold 
for any large $\N$; this conjecture may be thought of as a generalization of Joshi's conjecture  for the rank $2$ case (cf. ~\cite[Conjecture 8.1]{Jo14}).

\SSP
%-------------------------------------------------------------
\begin{conj} \label{Efg09}
Suppose that $X$ is a sufficiently general   smooth proper curve of genus $g$ over $k$.
Then, for each  $\N \in \mbZ_{>0}$, the cardinality  of the set  $F^\N\text{-}\mr{Ehr}_{(\mr{PGL}_2, \mbP^1), X/k}$ is given by the following formula:
\begin{align}
F^\N\text{-}\mr{Ehr}_{(\mr{PGL}_2, \mbP^1), X/k} = \frac{p^{n(g-1)}}{2^{2g-1}} \cdot \sum_{\theta =1}^{p^n-1} \frac{1}{\sin^{2g-2} \big( \frac{\pi \cdot \theta}{p^n}\big)}.
\end{align}
\end{conj}
%-------------------------------------------------------------

%-------------------------------------------------------------
\LSP
\subsection{Case (b):} \label{S0j910}

Next,
we  briefly review some results of the case where $G = \mr{Aff}_1$ and $H = \mr{Aff}_1^\circledcirc$.
Recall that $F^\N$-$(\mr{Aff}_1, \mr{Aff}_1^\circledcirc)$-structures on curves have been  studied in ~\cite{Hos4} and ~\cite{Wak6} under the name of  {\it Frobenius-affine structure}.
If $\N=1$, then these objects are equivalent to what we call {\it dormant (generic) Miura $\mr{PGL}_2$-opers} or {\it Tango structures} (cf. ~\cite{Muk}, ~\cite{Wak7}, ~\cite{Wak9}).

The existence of such an object  is very restrictive because 
if  there is a  Tango curve (= a curve admitting a Tango structure) of genus $g$, then 
$p$ must divide $g-1$.
M. Raynaud provided   an explicit  example of a
 Tango curve (cf. ~\cite{Ray} or ~\cite[\S\,1.1, Example 1.3]{Muk}).
Also,  Tango curves can be constructed   by means of solutions  to the Bethe ansatz equations for Gaudin model over finite fields (cf. ~\cite{Wak9}).
According to  ~\cite[Theorem B]{Wak7},  
the moduli stack 
 parametrizing Tango curves of genus $g$  forms, if it is nonempty,  
 an equidimensional smooth Deligne-Mumford  stack of dimension $2 \cdot g-2+ \frac{2 \cdot g-2}{p}$.
Regarding   higher-level cases,
 Y. Hoshi   proved (cf. ~\cite[Theorem 3]{Hos3}) that if an integer $g>1$   satisfies   $p^\N \mid (g-1)$ for a positive integer $\N$, then there exists a smooth proper curve of genus $g$ admitting  an $F^\N$-$(\mr{Aff}_1, \mbA^1)$-structure.
The following two propositions  describe some  results 
   proved previously   (cf.  ~\cite{Hos4}, ~\cite[Theorem D and \S\,7, Propositions 7.5.1,  7.5.2]{Wak6}).

\SSP
%-------------------------------------------------------------
\bpr \label{EQPPVL}
\begin{itemize}
\item[(i)]
For $\N =1$,
$X/k$ admits an $F^1$-$(\mr{Aff}_1, \mbA^1)$-structure  if and only if
$X$ is either a genus $1$ curve or a Tango curve.
\item[(ii)]
For an  integer $\N$ with $\N \geq 2$,
$X/k$ admits an $F^\N$-$(\mr{Aff}_1, \mbA^1)$-structure if and only if 
$X$ is either an ordinary genus $1$ curve  or a Tango curve with $p^\N \mid g-1$.
\item[(iii)]
For $\N =\infty$,
$X/k$ admits an $F^\infty$-$(\mr{Aff}_1, \mbA^1)$-structure if and only if
$X$ is an ordinary genus $1$ curve.
\end{itemize}
\epr
%-------------------------------------------------------------
\SSP

Here, for each ring $R$ and a free $R$-module $M$ of rank $1$, we shall denote by 
$\mbG_m (R; M)$ the set of generators in $M$.

\SSP
%-------------------------------------------------------------
\bpr \label{QAAV}
Suppose that $X$ is an ordinary genus $1$ curve, and 
denote by  $T_p X$   the $p$-adic Tate module of $X$.
Then, there exists a bijection of sets
\begin{align}
F^\N \text{-}\mr{Ehr}_{(\mr{Aff}_1, \mbA^1), X/k} \isom \mbG_m (R_\N ; T_p X \otimes_{\mbZ_p} R_\N),
\end{align}
where $R_\N := \mbZ/p^\N \mbZ$ if $\N \in \mbZ_{>0}$ and $R_\N := \mbZ_p$ if $\N = \infty$.
In particular,  if $\N \in \mbZ_{>0}$, then  the set $F^\N \text{-}\mr{Ehr}_{(\mr{Aff}_1, \mbA^1), X/k}$ is finite and its cardinality can be given by the following formula:
 \begin{align} \label{EWqddc}
 \sharp (F^\N \text{-}\mr{Ehr}_{(\mr{Aff}_1, \mbA^1), X/k}) = p^{\N -1}(p-1).
 \end{align}

\epr
%-------------------------------------------------------------

\LSP
%-------------------------------------------------------------
%\vspace{5mm}
\subsection{Cases (c) and (d): } \label{S092f1}

The classification problems  for the cases (c) and (d) were  already finished in the discussion of  \S\,\ref{S091k0}.
The following assertions are the results for  curves obtained there.

\SSP
%---------------------------------------------------------------
\bpr 
Let
 $\N\in \mbZ_{>0} \sqcup \{ \infty \}$.
 Then,    $X$ admits an $F^\N$-$(\mbG_m, \mbG_m)$-structure if and only if $X$ is  an  ordinary genus $1$ curve.
 If $X$ is  an ordinary  genus $1$ curve, then 
 the map of sets 
 \begin{align} \label{QLKLK}
 F^\N\text{-}\mr{Ehr}_{(\mbG_m, \mbG_m), X/k} \migi
 F^\N\text{-}\mr{Ehr}_{(\mr{Aff}_1, \mbA^1), X/k}
 \end{align}
 defined in Example \ref{es04xm}  is bijective.
 \epr
%---------------------------------------------------------------
\begin{proof}
The assertion  follow from
 Proposition \ref{C0558}, (i),  for $n=1$. 
Note that the surjectivity of  (\ref{QLKLK}) for $\N \in \mbZ_{>0}$ can be verified from  the discussion in the proof of Proposition \ref{C0558}, (i).
\end{proof}
%---------------------------------------------------------------
\SSP

 %------------------------------------------------------------
 \bpr \label{C055hk98}
 \begin{itemize}
 \item[(i)]
There exists  an  $F^1$-$(\mbG_a, \mbG_a)$-structure  on $X/k$
if and only if $X$ is  a  supersingular  genus $1$ curve.
Moreover, for any $\N \geq 2$, there are no $F^\N$-$(\mbG_a, \mbG_a)$-structures on $X/k$.
\item[(ii)]
Suppose that $X$ is a supersingular genus $1$ curve.
Then, there exists a canonical bijection
\begin{align}
F\text{-}\mr{Ehr}_{(\mbG_a, \mbG_a), X/k} \isom \mbG_m (k;  \Gamma (X, \Omega_{X/k})) \left(= \Gamma (X, \Omega_{X/k}) \setminus \{ 0 \} \right).
\end{align}
\end{itemize}
 \epr
 %------------------------------------------------------------
 \begin{proof}
 The assertions  follow from (\ref{Efgt4}),  the latter assertion of  Proposition \ref{L0045}, (i), and Proposition \ref{C0558}, (ii). 
 \end{proof}
%------------------------------------------------------------

%-------------------------------------------------------------
\vspace{5mm}
\subsection{Case (e):} \label{S0926}

Finally, let us consider the case where $G$ is an elliptic curve and $H =\{e \}$.

%-------------------------------------------------------------
\bpr
For each $\N \in \mbZ_{>0} \sqcup \{ \infty \}$, 
there exists   an $F^\N$-$(G, G)$-structure on $X/k$ if and only if $X$ is of genus $1$ and  isogenous to  $G$.
In particular, there are no $F^\N$-$(G, G)$-structures on any smooth proper curve of genus $>1$.
\epr
%-------------------------------------------------------------
\begin{proof}
The latter assertion is a direct consequence of the former assertion, so we only   discuss  the former assertion.
First, let us consider the ``if" part of the required equivalence.
Suppose that $X$ is of genus $1$ and there exists  an isogeny morphism $X \migi G$.
This morphism  decomposes as $X \xrightarrow{F_{X/k}^{(m)}} X^{(m)} \xrightarrow{w} G$ for some $m \geq 0$, where $w$ denotes a finite \'{e}tale $k$-covering.
Then, the pull-back $w^*(\mcS^{\mr{triv},\hh}_\N)$ of the standard $F^\N$-$(G, G)$-structure  $\mcS^{\mr{triv}, \hh}_\N$ on $G/k$ (cf. (\ref{EE012}))  forms 
an $F^\N$-$(G, G)$-structure  on $X^{(m)}$.
Moreover, it induces 
an $F^\N$-$(G, G)$-structure on $X/k$ via base-change by  $(F_k^{-1})^m : \mr{Spec}(k)\isom \mr{Spec}(k)$.
This completes the proof of the  ``if" part.

Next, 
we shall  consider  the ``only if" part.
To prove this, 
 it suffices to consider the case of $\N=1$.
Suppose that  $X$ is not isogenous to $G$ and  admits an $F^1$-$(G, G)$-structure $\mcS^\hh$.
Denote by $\mcE^\sss := (X \times \{ e\}, \nabla)$ the dormant indigenous $(G, \{e \})$-bundle corresponding to $\mcS^\hh$.
Under the natural identification   $\mr{Lie}(G)_{X \times G}/\mr{Lie}(\{ e\})_{X \times \{ e\}}\isom  \left(\mcO_X \otimes \mr{Lie}(G) =\right) \mcO_X$, 
the Kodaira-Spencer map  $\mr{KS}_{(X \times \{ e\}, \nabla)}$ may be regarded as an $\mcO_X$-linear  morphism  $\mcT_{X} \migi \mcO_X$.
Since $\mr{KS}_{(X \times \{ e\}, \nabla)}$  is assumed to be  an isomorphism, 
the equality  $\mr{deg}(\mcT_{X}) =0$ holds, in particular, $X$ must be  of genus $1$.
Next, denote by $\mcG$ the (right) $G$-bundle  on $X^{(1)}$ corresponding to 
the left $G_{X/S}^{(1)}$-torsor $\mcS^\heartsuit$.
Let us fix an isomorphism of $G$-bundles $\alpha : X \times G \isom F^*_{X/k} (\mcG)$.
Note that  there exist  a smooth proper curve $Y$  and a finite \'{e}tale covering $y : Y \migi X$ such that $\mcG$ can be  trivialized after pull-back   by the  finite \'{e}tale covering  $y^{(1)} : Y^{(1)} \migi X^{(1)}$  induced by $y$ (cf. ~\cite[Chap\,II, \S\,5, Lemma 5.5]{Mil2}).
That is to say, we can find  an isomorphism $\beta : Y^{(1)} \times_{X^{(1)}}\mcG  \isom Y^{(1)} \times G$.
Both $\alpha$ and $\beta$  yield the composite
\begin{align}
\gamma : Y \xrightarrow{a \mapsto (a, e)}    Y \times G \xrightarrow{y^*(\alpha)}  y^*(F^*_{X/k}(\mcG)) \left(= F^*_{Y/k}(y^{(1)*}(\mcG)) \right) \xrightarrow{F^*_{Y/k}(\beta)} Y \times G \xrightarrow{\mr{pr}_2} G.
\end{align}
The curve  $Y$ is not isogenous to $G$  because  $X$ is assumed  not to be  isogenous to $G$.
Hence, the image of $\gamma$ consists exactly of  one point.
It follows that    its  differential $d \gamma : \mcT_{Y} \migi \gamma^*(\mcT_{G}) \left(=  \mcO_Y \otimes \mr{Lie}(G) \right)$  is the zero map.
This contradicts  the fact that  
$d \gamma$
coincides with  the pull-back  of the  isomorphism  $\mr{KS}_{(X \times \{ e \}, \nabla)}$ by $y$.
This completes the proof of the ``only if" part.
\end{proof}
%-------------------------------------------------------------

%%%%%%%%%%%%%--[ begin  section1]---%%%%%%
\vspace{10mm}
\section{Deformation theory of (dormant) indigenous bundles} \label{SSS03}\vspace{3mm}

The topic of this section is  the deformation theory of  (dormant) indigenous $(G, H)$-bundles.
After proving basic facts on various deformation functors, we formulate and prove the Ehresmann-Weil-Thurston principle for Frobenius-Ehresmann structures (cf. Theorem \ref{Wroauom}), which is the main result of this section.

 Denote by
 $\mfS \mfe \mft$ the category of (small) sets and by 
  $\mfA \mfr \mft_{/k}$ the category of local Artinian $k$-algebras with residue field $k$.
 We shall write  $k_\epsilon := k[\epsilon]/(\epsilon^2)$, which is an object of $\mfA \mfr \mft_{/k}$.
For each $R \in \mr{Ob}(\mfA \mfr \mft_{/k})$, we write   $\mfm_R$ for  the maximal ideal of $R$.
Throughout this section, we fix   a smooth scheme $X$ of finite type over $k$.
Note that each morphism of sheaves $\nabla : \mcK^0 \migi \mcK^1$ on a scheme may be regarded as a complex concentrated at degrees $0$ and $1$; we shall denote this complex by $\mcK^\bullet [\nabla]$.
Also, for each sheaf $\mcF$, we shall write $\mcF [0]$ for the complex defined to be $\mcF$ concentrated at degree $0$.

\LSP
%---------------------------[begin subsection]-------------
\subsection{Deformation functors} \label{SSr086}

By a {\it  functor on $\mfA \mfr \mft_{/k}$} (or sometimes just a {\it functor} when there is no fear of confusion), we shall mean a covariant functor
 $D : \mfA \mfr \mft_{/k} \migi \mfS \mfe \mft$ such that  $D (k)$ is a singleton.
Also, following ~\cite[\S\, Definition 2.8]{FaMa}, we  refer to a functor $D$ on   $\mfA \mfr \mft_{/k}$ as 
a {\bf functor with good deformation theory} (or simply, a {\bf functor with gdt}) if it satisfies the following two of the Schlessinger conditions:
 \begin{itemize}
 \item[(H1)]
 For any morphism $\sigma' : R' \migi R$ and any surjective morphism  $\sigma'' : R'' \migisurj  R$ both of which belong  to  $\mfA \mfr \mft_{/k}$, 
 the map 
 \begin{align} \label{eaori}
 D (R' \times_R R'') \migi D (R') \times_{D (R)} D (R'')
 \end{align}
induced naturally by $\sigma'$ and $\sigma''$ is surjective.
(This condition remains equivalent even if $\sigma''$ is moreover assumed to be a principal small extension.)
 \item[(H2)]
 With the above notation, 
 (\ref{eaori}) becomes bijective when $R = k$.
 (This condition remains equivalent even if $R''$ is moreover assumed to be  $k_\epsilon$.)
 \end{itemize}

Let  $D : \mfA \mfr \mft_{/k} \migi \mfS \mfe \mft$ be  a deformation functor with gdt  and write $T D := D (k_\epsilon)$.
 Then, according to ~\cite[\S\,2, Lemma 2.12]{FaMa},
 the following properties  hold:
 \begin{itemize}
 \item
 $TD$ has a canonical structure of $k$-vector space;
 \item
For any small extension in $\mfA \mfr \mft_{/k}$ of the form $e := (0 \migi M \migi R \migi k \migi 0)$,
there exists a canonical  bijection $M \otimes TD \isom D (R)$;
 \item
 For any small extension $e := (0 \migi M \migi R \xrightarrow{\sigma} R_0 \migi 0)$ in  $\mfA \mfr \mft_{/k}$, there exists an $M \otimes TD$-action on $D (R)$ which is  restricted  to a transitive action on the fibers of $D (\sigma) : D(R) \migi D(R_0)$.
 Moreover, the formation of this action is compatible, in an evident sense,  with morphisms of small extensions.
 \end{itemize}
\SSP

Here, we recall several notions concerning  deformation functors used in our discussion.
 A morphism of  functors  $\gamma : D \migi D'$ is called {\it  formally smooth}  if for every small extension $e := (0 \migi M \migi R \xrightarrow{\sigma} R_0 \migi  0)$ in $\mfA \mfr \mft_{/k}$, the  map of sets $D (R) \migi D (R_0) \times_{D' (R_0)} D' (R)$ induced by $D (\sigma)$ and $\gamma (R)$ is surjective.
 If $\gamma : D \migi D'$ is a morphism between   functors with gdt,
 then  the map $d\gamma := \gamma (k_\epsilon) : TD \migi TD'$ induced by  it  becomes a $k$-linear morphism of vector spaces. We shall refer to $d \gamma$ as the {\it  differential} of $\gamma$.
  
  Next,  an {\it  obstruction theory} for a functor $D$ is  defined as a collection $(O, \{ o_e \}_e)$  consisting of a $k$-vector space $O$ and a map $o_e : D(R_0) \migi M \otimes O$ for every small extension $e := (0 \migi M \migi R \migi R_0 \migi 0)$ such that an element $a \in D (R_0)$ satisfies $o_e (a) = 0$ if and only if $a$ lifts to $D(R)$.
 The $k$-vector space $O$ is called an {\it  obstruction space} of $D$.
 We shall often abuse notation by just  writing  $O$ for $(O, \{ o_e \}_e)$.

 To every complete local $k$-algebra $R$  with residue field $k$,
 we can associate a functor with gdt $h_R : \mfA \mfr \mft_{/k} \migi \mfS \mfe \mft$ defined by $h_R (R') := \mr{Hom}_k (R, R')$.
 Recall that a functor on $\mfA \mfr \mft_{/k}$ is called {\it pro-representable} if it is isomorphic to $h_R$ for such an $R$.

\SSP
%---------------------------------------------------------------------------
\bde \label{e3099580}
\begin{itemize}
\item[(i)]
A {\bf deformation triple} is a collection
\begin{align}
\mbD := (D, T, O),
\end{align}
consisting of  a functor $D$  with gdt,  a $k$-vector space $T$  equipped with an identification  $T \isom T D$, and    an obstruction theory $O$ for $D$.
\item[(ii)]
Let $\mbD := (D, T, O)$ and $\mbD' := (D', T', O')$ be deformation triples.
A {\bf morphism of deformation triples} from $\mbD$ to $\mbD'$ is defined as a triple
\begin{align}
(\gamma_D, \gamma_T, \gamma_O),
\end{align}
where
\begin{itemize}
\item
$\gamma_D$ denotes  a morphism of  functors  $\gamma_D : D \migi D'$;
\item
$\gamma_T$ denotes  a $k$-linear morphism  $T \migi T'$
such that for any small extension $e := (0 \migi M \migi R \xrightarrow{\sigma} R_0 \migi 0)$ in $\mfA \mfr \mft_{/k}$, the $M \otimes T$-action on  $D(R)$ is compatible with
the $M \otimes T'$-action on $D' (R)$ via $\gamma_T$ and $\gamma_D$;
\item
$\gamma_O$ denotes a $k$-linear morphism $O \migi O'$ such that for any $e$ as above, the following square diagram is commutative:
\begin{align} \label{roi39}
\vcenter{\xymatrix@C=26pt@R=36pt{
D(R_0) \ar[r]^-{o_e} \ar[d]_-{\gamma_D (R_0) } & M \otimes O \ar[d]^-{\mr{id}_M \otimes \gamma_O} 
\\
D' (R_0)\ar[r]_-{o'_e} & M \otimes O'.
}}
\end{align}
 \end{itemize}
\end{itemize}
\ede
%---------------------------------------------------------------------------

\SSP
%---------------------------------------------------------------------------
\bde \label{Eotuamg}
Let $D$ be a functor with gdt.
We shall say that $D$ is {\bf rigid}
if  for any small extension $e:= (0 \migi M \migi R \xrightarrow{\sigma}R_0\migi 0)$ in $\mfA \mfr \mft_{/k}$,
the transitive $M \otimes TD$-action on every  fiber  of $D (\sigma) : D(R) \migi D(R_0)$ is free, i.e., the fibers form $M \otimes TD$-torsors.
\ede
%---------------------------------------------------------------------------

\SSP
%---------------------------------------------------------------------------
\begin{rema} \label{Rtiwisk}
Let $D : \mfA \mfr \mft_{/k} \migi \mfS \mfe \mft$ be a functor with gdt.
\begin{itemize}
\item[(i)]
It follows from ~\cite[\S\,2, Theorem 2.11]{Sch} 
 that 
$D$ is pro-representable  if and only if $D$ is rigid and satisfies $\mr{dim}_k (TD) < \infty$.
\item[(ii)]
Suppose that $D$ is associated with  a deformation category $\mcD$ in the sense of ~\cite[\S\,2, Definition 2.16]{TV}.
For each $R \in \mr{Ob}(\mfA \mfr \mft_{/k})$ and $\xi \in \mcD (R)$, we shall write $\mr{Aut}_R(\xi)$ for the automorphism of $\xi$ over  the identity morphism of $R$.
Then, 
 the rigidity of $D$ holds if, for any small extension $0 \migi M \migi R \xrightarrow{\sigma} R_0 \migi 0$ in $\mfA \mfr \mft_{/k}$, the group homomorphism $\mr{Aut}_{R}(\xi) \migi \mr{Aut}_{R_0} (\sigma^*(\xi))$ is surjective. 
\end{itemize}
\end{rema}
%---------------------------------------------------------------------------
\SSP

%-------------------------------------------------------------------------
\bde \label{pjapk0}
Let $\gamma : D \migi D'$ be a morphism of functors with gdt.
We shall say that $\gamma$ is a {\bf semi-isomorphism}
if it is formally smooth and its differential  $d \gamma : TD  \migi  TD'$ is bijective.
(Recall that a {\it pro-representable hull} for a functor $D$  is a semi-isomorphism $h_R \migi D$, where $R$ is  a completes local $k$-algebra $R$ with residue field $k$.)
\ede
%-------------------------------------------------------------------------

\SSP

%-------------------------------------------------------------------------
\bpr \label{pyp90s4}
Let $\gamma : D \migi D'$ be a morphism of functors with gdt.
\begin{itemize}
\item[(i)]
Suppose that  $D'$ is  rigid.
Then, $\gamma$ is a semi-isomorphism if and only if it is an isomorphism.
\item[(ii)]
Suppose that $\gamma$ is a semi-isomorphism
 and  there exists 
  a pro-representable hull $\nu : h_{R} \migi D$ for $D$.
Then, the composite $\gamma \circ \nu : h_R \migi D'$ is a pro-representable hull for $D'$. 
\item[(iii)]
Suppose that there exist  pro-representable hulls $\nu : h_{R} \migi D$, 
 $\nu' : h_{R'} \migi D'$
    for
$D$, $D'$ respectively.    
Then, the following conditions are equivalent to each other:
  \begin{itemize}
\item
The morphism $\gamma$ is a semi-isomorphism.
\item
There exists an isomorphism $\sigma : R \isom  R'$ of complete local $k$-algebras such that the induced isomorphism $h_\sigma : h_R \isom  h_{R'}$  makes the following diagram commute:
 \begin{align}\label{Diugh4}
\vcenter{\xymatrix@C=46pt@R=36pt{
h_R \ar[r]^{\nu} \ar[d]^-{\wr}_-{h_\sigma} & D \ar[d]^-{\gamma}
\\
h_{R'} \ar[r]_-{\nu'} & D'.
}}
\end{align}
\end{itemize}
\end{itemize}
\epr
%-------------------------------------------------------------------------
\begin{proof}
The assertions follow immediately  from various definitions involved together with  routine arguments of deformation functors.
\end{proof}
%-------------------------------------------------------------------------

\SSP

%-------------------------------------------------------------------------
\begin{rema} \label{E2Rgju}
Let $\gamma : D \migi D'$ be a morphism of  functors. 
Then, without the rigidity assumption on $D'$,
the morphism $\gamma$ can be a semi-isomorphism, yet $\gamma$ need not be an isomorphism.
It is because of this subtle difference that a semiuniversal family of a deformation functor  can fail to be universal (see ~\cite[\S\,3, Example 12]{Nit} for an example).
\end{rema}

\SSP

%---------------------------------------------------------------------------
\bpr \label{E2Rgju2}
Let $(\gamma_D, \gamma_T, \gamma_O) : (D, T, O) \migi (D', T', O')$ be a morphism of deformation triples.
Suppose that  $\gamma_T : T \migi T'$  is an isomorphism and $\gamma_O: O \migi O'$ is injective.
Then, $\gamma_D : D \migi D'$ is a semi-isomorphism.
\epr
%---------------------------------------------------------------------------
\begin{proof}
The formal smoothness of $\gamma_D$
   can be proved by a routine diagram-chasing argument  using the commutative  diagram
 \begin{align}
\vcenter{\xymatrix@C=46pt@R=36pt{
\hspace{-25mm}\left(M \otimes T D \curvearrowright \right) D(R)\ar[r]^-{D (\sigma)} \ar[d]^-{\gamma_D (R)} & D(R_0) \ar[d]^-{\gamma_D (R_0)} \ar[r]^-{o_e}& M\otimes O\ar[d]^-{\mr{id}_M \otimes \gamma_O}
\\
\hspace{-25mm} \left(M \otimes T D' \curvearrowright \right) D' (R) \ar[r]_-{D' (\sigma)} & D'(R_0) \ar[r]_-{o'_e} & M \otimes O'
}} \hspace{-15mm}
\end{align}
defined for every small extension $e:= (0 \migi M \migi R \xrightarrow{\sigma} R_0 \migi 0)$ in $\mfA \mfr \mft_{/k}$.
\end{proof}
%---------------------------------------------------------------------------

\LSP
%---------------------------[begin subsection]-------------
\subsection{Deformation of smooth schemes} \label{SS08df6}

For  $R \in \mr{Ob}(\mfA \mfr \mft_{/k})$, 
a {\bf  deformation of $X$} over $R$ is a smooth scheme $X_R$ over $R$ together with a $k$-isomorphism $X \isom k \times_{R}X_R$.
The notion of an isomorphism between two such deformations can be defined in an evident manner.
The {\bf  trivial deformation of $X$} over $R$ is the smooth $R$-scheme 
\begin{align} \label{ER456}
X_R^\circ := R \times X 
\end{align}
 together with the natural isomorphism $X \isom k \times_R (R \times X)$.
The assignment  from each $R \in \mr{Ob}(\mfA \mfr \mft_{/k})$ to the set $\mr{Def}_X (R)$ of isomorphism classes of deformations of $X$ over $R$ defines a functor  with gdt
\begin{align} \label{EWEW3}
\mr{Def}_{X} : \mfA \mfr \mft_{/k} \migi \mfS \mfe \mft.
\end{align}
For each $i \geq 0$, we shall write
\begin{align}
H^i := H^i (X, \mcT_{X}).
\end{align}
Then, it is well-known that 
$(\mr{Def}_X, H^1, H^2)$ forms a deformation triple, and 
%the tangent space of $\mr{Def}_X$  is isomorphic to $H^1$,  and $\mr{Def}_X$ has an obstruction theory by putting $H^2$ as the obstruction space.
%Moreover,
 the group of infinitesimal automorphisms of  $X$ is isomorphic to  $H^0$.

\LSP
%\vspace{5mm}
%---------------------------[begin subsection]-------------
\subsection{Deformation of indigenous $(G, H)$-bundles} \label{SS086}

%In this section, we shall discuss  the deformation functor of an indigenous $(G, H)$-bundle.
We shall fix an indigenous $(G, H)$-bundle $\mcE^\sss := (\mcE_H, \nabla)$ on $X/k$.

\SSP
%---------------------------------------------------
\bde
\begin{itemize}
\item[(i)]
Let $R \in \mr{Ob}(\mfA \mfr \mft_{/k})$.
A {\bf deformation of $(X, \mcE^\sss)$} over $R$ is 
a collection of data 
\begin{align}
(X_R, \mcE_{R}^\sss, \varrho),
\end{align}
where
$X_R$ denotes a deformation of $X$ over $R$, $\mcE_R^\sss$ denotes an indigenous $(G, H)$-bundle on $X_R/R$, and $\varrho$ denotes  an isomorphism  of indigenous $(G, H)$-bundles
$\mcE^\sss \isom k \times_R \mcE_R^\sss$.
A {\bf deformation of $\mcE^\sss$} over $R$ is defined as
a pair 
\begin{align}
(\mcE_{R}^\sss, \varrho)
\end{align}
 so that 
 $(X_R^\circ, \mcE_{R}^\sss, \varrho)$ forms a deformation  of $(X, \mcE^\sss)$ over $R$. 
 We  shall often abuse notation by just  writing $(X_R, \mcE_{R}^\sss)$ (resp., $\mcE_{R}^\sss$)  for $(X_R, \mcE_{R}^\sss, \varrho)$ (resp., $(\mcE_{R}^\sss, \varrho)$).  
\item[(ii)]
Let $\sigma : R' \migi R$ be a morphism in $\mfA \mfr \mft_{/k}$ and denote by
 $\sigma^\sharp$ the morphism $\mr{Spec}(R) \migi \mr{Spec}(R')$ induced by $\sigma$.
Also, let $(X_R, \mcE_{R}^\sss, \varrho)$ (resp., $(X_{R'}, \mcE_{R'}^\sss, \varrho')$) be a deformation of $(X, \mcE^\sss)$ over $R$ (resp., $R'$).
A {\bf morphism over $\sigma$} from $(X_{R}, \mcE_R^\sss, \varrho)$
to $(X_{R'}, \mcE^\sss_{R'}, \varrho')$ is a pair
\begin{align}
(\xi_X, \xi_\mcE)
\end{align}
consisting of   a morphism $\xi_X : X_{R} \migi X_{R'}$ over  $\sigma^\sharp$ 
 and  an isomorphism  of indigenous $(G, H)$-bundles $\xi_\mcE : \mcE_R^\sss \isom \xi^*_X (\mcE_{R'}^\sss)$ whose   reduction modulo $\mfm_R$ coincides with   the identity morphism  of $(X, \mcE^\sss)$ via $\varrho$ and $\varrho'$.
 If $\sigma = \mr{id}_R$, then we shall refer to such a   pair as a {\bf morphism over $R$}. 
Finally, we can define the notion of a morphism between deformations of
$\mcE^\sss$ in an evident manner. 
\end{itemize}
\ede
 %---------------------------------------------------

 For each $R \in \mr{Ob}(\mfA \mfr \mft_{/k})$,
 we shall write
 \begin{align}
 \mr{Def}_{\mcE^\sss} (R) := 
 \left(\begin{matrix} \text{the set of isomorphism classes} \\ \text{of deformations of $\mcE^\sss$ over $R$} 
 \end{matrix} \right) \hspace{12mm} \\
 \left(\text{resp.,} \  \mr{Def}_{(X, \mcE^\sss)} (R) := 
 \left(\begin{matrix} \text{the set of isomorphism classes} \\ \text{of deformations of $(X, \mcE^\sss)$ over $R$} 
 \end{matrix} \right)  \right). \notag 
 \end{align}
 If $\sigma  : R' \migi R$ is a morphism in $\mfA \mfr \mft_{/k}$,
 then changing the base by the induced morphism $\mr{Spec}(R) \migi \mr{Spec}(R')$ gives a map of sets $ \sigma^* : \mr{Def}_{\mcE^\sss} (R') \migi  \mr{Def}_{\mcE^\sss} (R)$ (resp., $ \mr{Def}_{(X, \mcE^\sss)} (R') \migi  \mr{Def}_{(X, \mcE^\sss)} (R)$).
 Thus, 
 the assignments $R \mapsto \mr{Def}_{\mcE^\sss}(R)$ (resp., $R \mapsto \mr{Def}_{(X, \mcE^\sss)}(R)$) and $\sigma \mapsto \sigma^*$   determine
a  functor 
  \begin{align} \label{Ew3}
 \mr{Def}_{\mcE^\sss}
 \  \left(\text{resp.,} \ \mr{Def}_{(X, \mcE^\sss)} \right) : \mfA \mfr \mft_{/k} \migi \mfS \mfe \mft.
 \end{align}
One may verify  from various  definitions  involved  that
 this functor specifies 
 a functor with gdt.
 The assignment
$\mcE^\sss_R \mapsto (X_{R}^{\circ}, \mcE_R^\sss)$ (resp., $(X_R, \mcE_R^\sss) \mapsto X_R$) 
 determines  a morphism of  functors
\begin{align} \label{L093f5}
\gamma_{\mcE^\sss} : \mr{Def}_{\mcE^\sss} \migi  \mr{Def}_{(X, \mcE^\sss)} \ \left(\text{resp.,} \  \delta_{\mcE^\sss} : \mr{Def}_{(X, \mcE^\sss)} \migi \mr{Def}_X \right).
\end{align}

\LSP
%---------------------------[begin subsection]-------------
\subsection{Cohomological descriptions of $\mr{Def}_{\mcE^\sss}$ and $\mr{Def}_{(X, \mcE^\sss)}$} \label{SS0100}
 In what follows, let us describe the tangent  and obstruction spaces of the deformation functors introduced above  in terms of cohomology groups.
 The corresponding discussion
 for holomorphic Cartan geometries 
   can be found in ~\cite{BD3} (cf. Remark \ref{ERqqw3}).
 
 Let $\mcE^\sss := (\mcE_H, \nabla)$ be as above
    and write $\mcE_G := \mcE_H \times^H G$.
 Denote by $\nabla^\mr{ad}_G$ the $k$-connection $\mfg_{\mcE_G} \migi \Omega_{X} \otimes \mfg_{\mcE_G}$ on the adjoint bundle $\mfg_{\mcE_G}$ induced by $\nabla$ via change of structure group by  $\mr{Ad}_G : G \migi \mr{GL}(\mfg)$.
 Denote by
 \begin{align}
 \nabla^\mr{ad}_H : \mfh_{\mcE_H} \migi \Omega_{X} \otimes \mfg_{\mcE_G}
 \end{align}
 the morphism obtained by restricting the domain of $\nabla_G^\mr{ad}$ to $\mfh_{\mcE_H} \left(\subseteq \mfg_{\mcE_G} \right)$.
   For any local sections
$v \in  \widetilde{\mcT}_{\mcE_G} \left(:= \widetilde{\mcT}_{\mcE_G/k} \right)$ and
$\partial \in \mcT_{X}$, 
we shall write 
\begin{align}
\widetilde{\nabla}^\mr{ad}_{\partial} (v):= [\nabla (\partial), v] - \nabla ([\partial, d_{\mcE_G}(v)]) \in \widetilde{\mcT}_{\mcE_G}.
\end{align}

\SSP
%-----------------------------------------------------------------------[begin definition]------------------
\ble \label{y0145}
The local section $\widetilde{\nabla}^\mr{ad}_{\partial} (v)$
lies  in 
$\mfg_{\mcE_G} \left(\subseteq \widetilde{\mcT}_{\mcE_G}\right)$.
\ele
%-----------------------------------------------------------------------[end definition]-------------------
\bpf
Note that $d_{\mcE_G}$ preserves the Lie bracket operations.
Hence, for any local sections $v \in  \widetilde{\mcT}_{\mcE_G}$ and  $\partial \in \mcT_{X}$,    the following sequence of equalities holds:
 \begin{align}   d_{\mcE_G}  ( \widetilde{\nabla}_{\partial}^{\mr{ad}} (v)) 
 & = d_{\mcE_G} ( [\nabla (\partial), v] - \nabla ([\partial,  d_{\mcE_G}(v)]) )  \\
 & = [d_{\mcE_G} \circ \nabla (\partial), d_{\mcE_G} (v)]-  d_{\mcE_G} \circ \nabla ([\partial, d_{\mcE_G}(v)])  \notag \\
&
 =  [\partial, d_{\mcE_G} (v)] -[\partial, d_{\mcE_G} (v)] \notag \\
 &  = 0. \notag
 \end{align} 
This means  that the image of $\widetilde{\nabla}^\mr{ad}_{\partial} (v)$ lies  in 
$\mr{Ker}(d_{\mcE_G})  \left(= \mfg_{\mcE_G}\right)$.
\epf
%-----------------------------------------------------------------------[end definition]-------------------
%\vspace{3mm}
\SSP

Let $\Box$ denote either $G$ or $H$.
By the above lemma, the assignment
$v \mapsto  \widetilde{\nabla}^{\mr{ad}}_{(-)} (v)$ defines 
a  $k$-linear morphism
\begin{align} 
 \widetilde{\nabla}_{\Box}^{\mr{ad}} :  \widetilde{\mcT}_{\mcE_\Box} \migi  \left(\mcH om_{\mcO_X} (\mcT_{X},   \mfg_{\mcE_G})= \right)\Omega_{X}\otimes  \mfg_{\mcE_G}.
 \label{circledast}
\end{align}
This morphism fits into the following morphism of short exact sequences:
\begin{align} \label{QW7010}
\vcenter{\xymatrix@C=40pt@R=36pt{
0 \ar[r] & \mr{Lie} (\Box)_{\mcE_\Box} \ar[r]^-{\mr{inclusion}}\ar[d]^{\nabla^\mr{ad}_\Box}& \widetilde{\mcT}_{\mcE_\Box} \ar[r]^-{d_{\mcE_\Box}}\ar[d]^{\widetilde{\nabla}^\mr{ad}_\Box}& \mcT_{X}  \ar[r] \ar[d] & 0 \\
0 \ar[r]& \Omega_{X} \otimes \mfg_{\mcE_G} \ar[r]_-{\mr{id}}& \Omega_{X} \otimes \mfg_{\mcE_G} \ar[r] & 0  \ar[r] & 0. 
}}
\end{align}
In other words,
we have  the following  short exact sequence of complexes:
\begin{align} \label{QW7034}
0 \longmigi \mcK^\bullet [\nabla^\mr{ad}_\Box] \longmigi \mcK^\bullet [\widetilde{\nabla}^\mr{ad}_\Box] \longmigi \mcT_{X}[0] \longmigi 0. 
\end{align}
For each $i \geq 0$,  we shall write
\begin{align}
\mbH^i_\Box := \mbH^i (X, \mcK^\bullet [\nabla_\Box^\mr{ad}]),
\hspace{5mm}
\widetilde{\mbH}_\Box^i := \mbH^i (X, \mcK^\bullet [\widetilde{\nabla}_\Box^\mr{ad}]).
\end{align}
By applying the hypercohomology functor $\mbH^i (X, -)$ ($i \geq 0$) to (\ref{QW7034}), we obtain  morphisms
\begin{align} \label{QW7011}
\mbH^i_\Box \xrightarrow{\gamma^i_\Box} \widetilde{\mbH}_\Box^i \xrightarrow{\delta^i_\Box}  H^i.
\end{align}
If $X$ is proper over $k$, then
 the $k$-vector spaces appearing in (\ref{QW7011}) are all finite-dimensional.
 
 \SSP
 %------------------------------------------------------------------------
 \ble \label{L0345s}
 Let us consider the $\mcO_X$-linear endomorphism $\eta$ of  $\widetilde{\mcT}_{\mcE_G}$ given by $v \mapsto v - (\nabla \circ  d_{\mcE_G})(v)$ for each local section $v \in \widetilde{\mcT}_{\mcE_G}$.
Then, this endomorphism is restricted to a morphism  (resp.,  an isomorphism)
\begin{align} \label{SSKKp0}
\eta_G : \widetilde{\mcT}_{\mcE_G} \migi \mfg_{\mcE_G} \ \left(\text{resp.,} \ 
\eta_H : \widetilde{\mcT}_{\mcE_H} \isom \mfg_{\mcE_G}\right).
\end{align}
 Moreover,  for  $\Box \in \{ G, H \}$,   the following  diagram  is commutative:
\begin{align} \label{Dig01}
\vcenter{\xymatrix@C=40pt@R=36pt{
\widetilde{\mcT}_{\mcE_\Box/k} \ar[rr]_-{}^-{\eta_\Box} \ar[rd]_-{\widetilde{\nabla}^\mr{ad}_\Box} && \mfg_{\mcE_G}\ar[ld]^-{\nabla^\mr{ad}_G}
\\
& \Omega_{X} \otimes \mfg_{\mcE_G}.&
}}
\end{align}
 In particular,  
 the isomorphism $\eta^\bullet_H : \mcK^\bullet [\widetilde{\nabla}^\mr{ad}_H] \isom \mcK^\bullet [\nabla^\mr{ad}_G]$ determined by 
   the pair of morphisms $(\eta_H, \mr{id}_{\Omega_{X} \otimes \mfg_{\mcE_G}})$ induces, for each $i \geq 0$,   an isomorphism of $k$-vector spaces $\widetilde{\mbH}_H^i \isom \mbH^i_G$.
 \ele
 %------------------------------------------------------------------------
 \begin{proof}
  First,  let us prove  the former assertion.
 To this end,  we  only  consider the case of $\eta_H$ because the other case is entirely similar.
 If we  take a local section $v \in \widetilde{\mcT}_{\mcE_G}$,
  then we have 
 \begin{align}
 d_{\mcE_G} (\eta (v)) &= d_{\mcE_G} (v - (\nabla \circ d_{\mcE_G})(v)) \\
 & = 
 d_{\mcE_G}(v) -  (d_{\mcE_G}\circ \nabla) (d_{\mcE_G}(v)) \notag  \\
 & = d_{\mcE_G}(v) - d_{\mcE_G}(v) \notag  \\
 & =0. \notag
 \end{align}
 This implies that $\eta (v) \in \mr{Ker}(d_{\mcE_G}) = \mfg_{\mcE_G}$, so  the image of $\eta$ is contained in $\mfg_{\mcE_G}$.
 By restricting the domain and codomain of $\eta$,
 we obtain an $\mcO_X$-linear morphism $\eta_H : \widetilde{\mcT}_{\mcE_H} \migi \mfg_{\mcE_G}$.
 One may verify that $\eta_H$ fits into the following morphism of short exact sequences:
 \begin{align} \label{QW70fg0}
\vcenter{\xymatrix@C=46pt@R=36pt{
0 \ar[r] & \mfh_{\mcE_H} \ar[r]^-{\mr{inclusion}}\ar[d]^{\mr{id}_{}}& \widetilde{\mcT}_{\mcE_H} \ar[r]^-{d_{\mcE_H}}\ar[d]^{\eta_H}& \mcT_{X}  \ar[r] \ar[d]^-{\mr{KS}_{(\mcE_H, \nabla)}} & 0 \\
0 \ar[r]& \mfh_{\mcE_H} \ar[r]_-{\mr{inclusion}}&  \mfg_{\mcE_G} \ar[r]_-{\mr{quotient}} & \mfg_{\mcE_G}/\mfh_{\mcE_H} \left(\cong \widetilde{\mcT}_{\mcE_G}/\widetilde{\mcT}_{\mcE_H} \right)  \ar[r] & 0. 
}}
\end{align}
Since the Kodaira-Spencer map $\mr{KS}_{(\mcE_H, \nabla)}$ is an isomorphism,
$\eta_H$ turns out to be an isomorphism by the five lemma.
  
Next, we shall prove the latter assertion, i.e., the commutativity of (\ref{Dig01}).
Clearly, it suffices to consider the case of $\eta_G$.
Denote by $\langle -, - \rangle$ the pairing $\mcT_{X} \times (\Omega_{X} \otimes \mfg_{\mcE_G})  \migi \mfg_{\mcE_G}$ arising from the natural pairing $\mcT_{X} \times \Omega_{X}  \migi \mcO_X$.
Then,
for each local section $v \in \widetilde{\mcT}_{\mcE_G}$, we have
\begin{align}
\langle \partial, (\nabla_G^\mr{ad} \circ \eta_G) (v) \rangle 
&= \langle \partial, \nabla_G^\mr{ad} (v - (\nabla \circ d_{\mcE_H})(v)) \rangle \\
&= 
[\nabla (\partial),  v - (\nabla \circ d_{\mcE_H})(v)]  \notag \\
& =  [\nabla (\partial), v]-[\nabla (\partial), (\nabla \circ d_{\mcE_G})(v)] \notag \\
& = [\nabla (\partial), v]- \nabla ([\partial, d_{\mcE_G}(v)]) \notag \\
& = \langle \partial, \widetilde{\nabla}^\mr{ad}_G (v) \rangle. \notag
\end{align}
This implies $\nabla_G^\mr{ad} \circ \eta_G = \widetilde{\nabla}_G^\mr{ad}$,   completing  the proof of this lemma. 
 \end{proof}
 %------------------------------------------------------------------------
 \SSP

 %------------------------------------------------------------
 \bpr \label{P022}
\begin{itemize}
\item[(i)]
The collection  $(\mr{Def}_{\mcE^\sss}, \mbH^1_H, \mbH^2_H)$ forms a deformation triple.
%The tangent space of $\mr{Def}_{\mcE^\sss}$ is isomorphic to
%$\mbH^1_H$,  and 
%this functor
% has an obstruction theory by putting 
%$\mbH^2_H$ as the obstruction space.
 Moreover, the  group of infinitesimal automorphisms  of $\mcE^\sss$ (i.e., the group of  automorphisms of $k_\epsilon \times_k \mcE^\sss$ inducing the identity morphism of $\mcE^\sss$ via reduction modulo $(\epsilon)$) is isomorphic to $\mbH^0_H$.
 \item[(ii)]
 The collection  $(\mr{Def}_{(X, \mcE^\sss)}, \widetilde{\mbH}^1_H, \widetilde{\mbH}^2_H)$ ($\cong (\mr{Def}_{(X, \mcE^\sss)}, \mbH^1_G, \mbH^2_G)$ by Lemma \ref{L0345s}) forms a deformation triple.
% The tangent space of $\mr{Def}_{(X, \mcE^\sss)}$ 
% is isomorphic to $\widetilde{\mbH}^1_H$ ($\cong \mbH^1_G$ by Lemma \ref{L0345s}),
% and 
% this functor 
%  has an obstruction theory by putting 
%  $\widetilde{\mbH}^2_H$ ($\cong \mbH^2_G$ by  the same lemma)  as the obstruction space.
   Moreover, the group of infinitesimal  automorphisms of $(X, \mcE^\sss)$ is isomorphic to $\widetilde{\mbH}_H^0$ ($\cong \mbH^0_G$).
  \item[(iii)]
  Suppose that $X$ is proper over $k$.
  Then,   there exist  pro-representable hulls for   $\mr{Def}_{\mcE^\sss}$ and  $\mr{Def}_{(X, \mcE^\sss)}$.
  \end{itemize}
 \epr
 %------------------------------------------------------------
  \begin{proof}
  We will omit the proof of assertion (i)
     because it is entirely similar to 
      the proof of assertion (ii).
 
  First, we shall  consider the tangent space of $\mr{Def}_{(X, \mcE^\sss)}$.
 Let $e := (0 \migi M \migi R \xrightarrow{\sigma} R_0 \migi 0)$ be a small extension in $\mfA \mfr \mft_{/k}$ and    $(X_R, \mcE^\sss_R)$  (where $\mcE_R^\sss := (\mcE_{R, H}, \nabla_R)$)
  a deformation of $(X, \mcE^\sss)$  over $R$.
 Also, let us take
  an element $v$ of $\mbH^1 (X, \mcK^\bullet [ \mr{id}_M \otimes \widetilde{\nabla}_H^\mr{ad}]) \left( \cong M \otimes \widetilde{\mbH}^1_H \right)$, where $\mr{id}_M \otimes \widetilde{\nabla}_H^\mr{ad}$ denotes the tensor product  $M \otimes \widetilde{\mcT}_{\mcE_H} \migi M \otimes \mfg_{\mcE_G}$ of $\mr{id}_M$ and $\widetilde{\nabla}_H^\mr{ad}$.
Note here that  $\mbH^1 (X, \mcK^\bullet [ \mr{id}_M \otimes \widetilde{\nabla}_H^\mr{ad}])$ can be realized as the total cohomology of   the \v{C}ech double complex   of 
  $\mcK^\bullet [\mr{id}_M \otimes \widetilde{\nabla}_H^\mr{ad}]$ associated with  an affine covering.
That is to say, we can find an affine open covering $\mcU := \{ U_\alpha \}_{\alpha \in I}$ (where 
$I$ is a finite index set) of $X$ and  $v$ may be represented
by
 a collection of data
\begin{align}\label{W105}
 v = (\{ \widetilde{\partial}_{\alpha \beta} \}_{(\alpha, \beta) \in I_2}, \{ \delta_\alpha \}_{\alpha \in I}),
\end{align}
where
\begin{itemize}
\item
 $I_2$ denotes  the subset of $I \times I$ consisting of pairs $(\alpha, \beta)$ with $U_{\alpha \beta}:= U_\alpha \cap U_\beta \neq \emptyset$;
\item
$\{ \widetilde{\partial}_{\alpha\beta} \}_{(\alpha, \beta) \in I_2}$ is 
 a \v{C}ech $1$-cocycle  in 
$\check{C}^1 (\mcU, M \otimes \widetilde{\mcT}_{\mcE_H})$, where $\widetilde{\partial}_{\alpha\beta} \in \Gamma (U_{\alpha \beta},  M \otimes \widetilde{\mcT}_{\mcE_H})$;
 \item
$\{ \delta_\alpha \}_{\alpha \in I}$ is  a \v{C}ech $0$-cochain in   $\check{C}^0 (\mcU, M \otimes (\Omega_{X}  \otimes \mfg_{\mcE_G}))$, where $\delta_\alpha \in \Gamma (U_\alpha, M \otimes (\Omega_{X} \otimes \mfg_{\mcE_G})) = \mr{Hom}_{\mcO_{U_\alpha}}(\mcT_{X} |_{U_\alpha},  M \otimes \mfg_{\mcE_G} |_{U_\alpha})$, 
\end{itemize}
such that
 $\{ \widetilde{\partial}_{\alpha\beta} \}_{(\alpha, \beta) \in I_2}$ and $\{ \delta_\alpha \}_{\alpha \in I}$
agree under $\mr{id}_M \otimes \widetilde{\nabla}_H^\mr{ad}$ and the \v{C}ech coboundary map.
 For each $(\alpha, \beta) \in I_2$,
 we set $\partial_{\alpha \beta} := d_{\mcE_H} (\widetilde{\partial}_{\alpha \beta})$.
We shall write $U_{R, \alpha} := X_{R} |_{U_\alpha}$ for each $\alpha \in I$ and write
$U_{R, \alpha \beta} := X_{R} |_{U_{\alpha\beta}}$ for each $(\alpha, \beta) \in I_2$.
 The element $\partial_{\alpha \beta}$ (resp.,  $\widetilde{\partial}_{\alpha \beta}$)  may be regarded as  a $k$-linear morphism 
 $\partial_{\alpha \beta} : \mcO_{U_{\alpha \beta}} \migi M \otimes \mcO_{U_{\alpha \beta}}$ 
    (resp., $\widetilde{\partial}_{\alpha \beta} : \mcO_{\mcE_H |_{U_{\alpha \beta}}} \migi M \otimes \mcO_{\mcE_H |_{U_{\alpha \beta}}}$)  satisfying
    the Leibniz rule.
 Then, 
 $\mr{id}_{\mcO_{U_{R, \alpha \beta}}} + \partial_{\alpha \beta}$  (resp., $\mr{id}_{\mcO_{\mcE_{R, H}|_{U_{R, \alpha \beta}}}} + \widetilde{\partial}_{\alpha \beta}$) defines a well-defined  automorphism $\partial^\sharp_{\alpha \beta}$ (resp., $\widetilde{\partial}^\sharp_{\alpha \beta}$) of  $U_{R, \alpha \beta}$ (resp., $\mcE_{R, H}|_{U_{R, \alpha \beta}}$).
The automorphism  $\widetilde{\partial}^\sharp_{\alpha \beta}$ is $H$-equivariant and 
 the following square diagram is commutative and cartesian:
  \begin{align} \label{DD010}
\vcenter{\xymatrix@C=40pt@R=36pt{
\mcE_{R, H}|_{U_{R, \alpha \beta}} \ar[r]_-{\sim}^{\widetilde{\partial}^\sharp_{\alpha \beta}} \ar[d]_-{\mr{projection}}
& \mcE_{R, H}|_{U_{R, \alpha \beta}} \ar[d]^-{\mr{projection}}
\\
U_{R, \alpha \beta} \ar[r]^-{\sim}_{\partial^\sharp_{\alpha \beta}}& U_{R, \alpha \beta}.
}}
\end{align}
By means of $\partial_{\alpha \beta}^\sharp$'s ($(\alpha, \beta) \in I_2$), we can glue together  $U_{R, \alpha}$'s ($\alpha \in I$) to
  obtain a smooth scheme  $X^v_{R}$ over $R$.
  Moreover, 
 the pairs $(\mcE_{R, H}|_{U_{R, \alpha}}, \nabla |_{U_{R, \alpha}} + \delta_\alpha)$ 
  may be glued together by means of 
 $\widetilde{\partial}^\sharp_{\alpha \beta}$'s.
 Thus, we obtain an $H$-bundle $\mcE_{R, H}^v$ on $X_{R}^v$ and an $R$-connection $\nabla^v_{R}$ on $\mcE_{R, G}^v := \mcE_{R, H}^v \times^H G$.
The Kodaira-Spencer map $\mr{KS}_{(\mcE_{R, H}^v, \nabla^v_{R})}$
is an isomorphism because its reduction modulo $\mfm_R$ coincides with the isomorphism $\mr{KS}_{(\mcE_H, \nabla)}$.
It follows that the  resulting pair $\mcE_{R}^{\sss, v} := (\mcE_{R, H}^v, \nabla^v_{R})$ forms an indigenous $(G, H)$-bundle on $X^v_{R}$ and the pair  $(X_{R}^v, \mcE_{R}^{\sss, v})$ specifies an element of $\mr{Def}_{(X, \mcE^\sss)}(R)$.
The assignment  $((X_{R}, \mcE_{R}^\sss), v) \mapsto (X_{R}^v, \mcE_{R}^{\sss, v})$ defines a transitive  action of $M \otimes  \widetilde{\mbH}_H^1$ on 
the fibers of $\mr{Def}_{(X, \mcE^\sss)}(\sigma) : \mr{Def}_{(X, \mcE^\sss)} (R) \migi \mr{Def}_{(X, \mcE^\sss)} (R_0)$.
(To verify the transitivity of this action, we may need the fact that $M \otimes \widetilde{\mbH}_H^1$ is canonically isomorphic to the $1$-st hypercohomology group of the complex of  {\it \'{e}tale} sheaves  associated with $\mcK^\bullet [\mr{id}_M \otimes \widetilde{\nabla}^\mr{ad}_H]$, each of whose elements may be represented by a \v{C}ech $1$-cocycle with respect to an \'{e}tale covering; this is because $H$-bundles are locally trivial in the \'{e}tale topology, but not  necessarily in the Zariski topology.)
This action is, by  construction,  verified to be free if $R = k$.
Consequently, $\widetilde{\mbH}_H^1$  is isomorphic to the tangent space of $\mr{Def}_{(X, \mcE^\sss)}$.

Next,  we shall consider the assertion about an obstruction space.
Let $e$ be as above and 
$a$ an element of $\mr{Def}_{(X, \mcE^\sss)}(R_0)$.
Choose a deformation  $(X_{R_0}, \mcE_{R_0}^\sss)$ (where 
$\mcE^\sss_{R_0} := (\mcE_{R_0, H}, \nabla_{R_0})$) of $(X, \mcE^\sss)$ representing $a$.
Observe that
there exists an affine open covering 
$\mcV := \{ V_\alpha \}_{\alpha \in J}$ of $X$  such that, for each $\alpha \in J$,
the $R_0$-scheme $V_{R_0, \alpha} := X_{R_0} |_{V_\alpha}$ has  a deformation   $V_{R, \alpha}$ over $R$ and 
 $(\mcE_{R_0, H} |_{V_{\alpha}}, \nabla_{R_0}  |_{V_{\alpha}})$ may be lifted to a pair  $\mcE_{R, \alpha}^\sss :=(\mcE_{R,  H, \alpha}, \nabla_{R, \alpha})$ defined over $V_{R, \alpha}$.
After possibly replacing $\mcV$ by its refinement,
we can find a collection $\{ \eta_{\alpha \beta} \}_{(\alpha, \beta) \in J_2}$, where each $\eta_{\alpha \beta}$ denotes an $H$-equivariant  isomorphism $\mcE_{R, H,  \alpha} |_{V_{\alpha \beta}} \isom \mcE_{R, H,  \beta} |_{V_{\alpha \beta}}$ inducing the identity morphism of $\mcE_{R_0, H} |_{V_{\alpha \beta}}$ via reduction modulo $M$.
There exists an element $\delta_{\alpha \beta} \in \Gamma (U_{\alpha \beta}, M \otimes (\Omega_{X} \otimes \mfg_{\mcE_G}))$ with $\eta^*_{\alpha \beta} (\nabla_{R, \alpha}) = \nabla_{R, \beta} + \delta_{\alpha \beta}$.
Moreover, for each triple $(\alpha, \beta, \gamma) \in I \times I \times I$ with $U_{\alpha \beta \gamma} := U_\alpha \cap U_\beta \cap U_\gamma \neq \emptyset$, there exists an element $\widetilde{\partial}_{\alpha \beta \gamma}  \in \Gamma (U_{\alpha \beta \gamma}, M \otimes \widetilde{\mcT}_{\mcE_H})$ with 
 $\eta_{\gamma \alpha} \circ \eta_{\beta \gamma} \circ \eta_{\alpha \beta} = \mr{id} + \widetilde{\partial}_{\alpha \beta \gamma}$.
The collection $(\{ \widetilde{\partial}_{\alpha \beta \gamma} \}_{(\alpha, \beta, \gamma)}, \{ \delta_{\alpha \beta} \}_{(\alpha, \beta)})$ represents an element  $o_e (a)$ of $M \otimes \widetilde{\mbH}_H^2$.
In this way, to each small extension $e$ as above, we   associate a map of sets $o_e : \mr{Def}_{(X, \mcE^\sss)}(R_0) \migi M \otimes \widetilde{\mbH}^2_H$.
By construction, the equality  $o_e (a) =0$ holds  if and only if the pairs $(V_{R, \alpha}, \mcE_{R, \alpha}^\sss)$   may be glued together, by means of $\{ \eta_{\alpha \beta} \}_{(\alpha, \beta) \in J_2}$, to obtain a deformation of $(X, \mcE^\sss)$ over $R$ lifting  $(X_R, \mcE_{R}^\sss)$.
That is to say, $(\widetilde{\mbH}^2_H, \{ o_e\}_e)$ forms an obstruction theory for $\mr{Def}_{(X, \mcE^\sss)}$.

The last assertion of (ii) follows immediately from the \v{C}ech cohomological description of deformations discussed above.
Finally, assertion (iii) follows from ~\cite[\S\,2, Theorem 2.11]{Sch} together with   assertions (i), (ii), and  the fact that $\mr{dim}_k (\mbH_H^1) < \infty$, $\mr{dim}_k (\widetilde{\mbH}_H^1) < \infty$ under the properness assumption on $X$.
 \end{proof}
%------------------------------------------------------------
\SSP

 %------------------------------------------------------------
 \bco \label{Cf02k2}
 Suppose that $\mbH^0_H = 0$ (resp., $\widetilde{\mbH}_H^0 =0$) and $X$ is proper over $k$.
 Then,  $\mr{Def}_{\mcE^\sss}$ (resp., $\mr{Def}_{(X, \mcE^\sss)}$) is  pro-representable. 
 If, moreover, $\mbH^2_H =0$ (resp., $\widetilde{\mbH}_H^2 =0$), then this functor 
  may be pro-representable by $k [\! [x_1, \cdots, x_d ]\! ]$, where $d := \mr{dim}_k (\mbH^1_H)$ (resp., $d := \mr{dim}_k (\widetilde{\mbH}^1_H)$).
 \eco
 %------------------------------------------------------------
  \begin{proof}
 The assertion follows from Proposition \ref{P022} and ~\cite[Chap.\,6, \S\,6.2, Corollary 6.2.5, and  \S\,6.3, Corollary 6.3.5]{FGA}.
 \end{proof}
%------------------------------------------------------------

\SSP

%------------------------------------------------------------
\begin{rema} \label{ERqqw3}
The space of infinitesimal automorphisms, as well as 
the infinitesimal deformations,  of a {\it holomorphic} Cartan geometry
 are computed  in ~\cite[\S\,3, Propositions 3.1, 3.3, and Theorem 3.4]{BD3}.
As S. Dumitrescu commented in our personal communication,  the statement of Theorem 3.4 in {\it loc.\,cit.} will need to be corrected.
(The updated version may be available at a later date in some way.)
 These results after correction are consistent with the assertion in  
 Proposition \ref{P022}, (i), under the correspondence between Cartan geometries and indigenous bundles by (\ref{E096}). 
\end{rema}
%------------------------------------------------------------

\SSP

%------------------------------------------------------------
\begin{rema} \label{Rwwqa}
We here recall well-known facts concerning the deformation theory of $G$-bundles and connections.  (These facts will be used in the proofs of some propositions described  later and hold without the assumption that $\mcE_G$ has an $H$-reduction defining the structure of indigenous $(G, H)$-bundle, i.e., $\mcE_H$.)
Denote by  
\begin{align} \label{RwwqaW}
\mr{Def}_{(X, \mcE_G)} \  (\text{resp.,}   \ \mr{Def}_{(\mcE_G, \nabla)}; \text{resp.,}  \ \mr{Def}_{(X, \mcE_G, \nabla)})
\end{align}
the deformation functor classifying deformations of  the pair $(X, \mcE_G)$ (resp.,  the pair   $(\mcE_G, \nabla)$; resp., the triple $(X, \mcE_G, \nabla)$), viewed as
a smooth $k$-scheme and a $G$-bundle (resp., 
  a $G$-bundle  equipped with a  $k$-connection; resp., a smooth $k$-scheme  and a $G$-bundle on it equipped with a $k$-connection).
Then,  by a similar argument in the proof of Proposition \ref{P022},
we see that
$(\mr{Def}_{(X, \mcE_G)}, H^1 (X, \widetilde{\mcT}_{\mcE_G}), H^2 (X, \widetilde{\mcT}_{\mcE_G}))$
(resp., $(\mr{Def}_{(\mcE_G, \nabla)}, \mbH^1_G, \mbH^2_G)$; resp., $(\mr{Def}_{(X, \mcE_G, \nabla)}, \widetilde{\mbH}_G^1, \widetilde{\mbH}_G^2)$)
forms a deformation triple.
%$\mr{Def}_{(X, \mcE_G)}$ (resp.,  $\mr{Def}_{(\mcE_G, \nabla)}$; resp., $\mr{Def}_{(X, \mcE_G, \nabla)}$) defines   a functor with gdt whose tangent space
%is isomorphic to 
%$H^1 (X, \widetilde{\mcT}_{\mcE_G})$ (resp., 
%$\mbH^1_G$; resp., $\widetilde{\mbH}_G^1$),
%and  has an obstruction theory by putting  $H^2 (X, \widetilde{\mcT}_{\mcE_G})$ (resp., $\mbH^2_G$; resp., $\widetilde{\mbH}_G^2$)
%as the obstruction space.
Moreover, the group of infinitesimal  automorphisms of $(X, \mcE_G)$ (resp., $(\mcE_G, \nabla)$; resp., $(X, \mcE_G, \nabla)$) is isomorphic to    $H^0 (X, \widetilde{\mcT}_{\mcE_G})$ (resp.,  $\mbH^0_G$; resp., $\widetilde{\mbH}^0_G$).
\end{rema}
\SSP
%------------------------------------------------------------

By the above proposition and  various constructions involved, 
the triple 
\begin{align} \label{Erty}
(\gamma_{\mcE^\sss}, \gamma^1_H, \gamma_H^2) : (\mr{Def}_{\mcE^\sss}, \mbH^1_H, \mbH_H^2) \migi (\mr{Def}_{(X, \mcE^\sss)}, \widetilde{\mbH}_H^1,  \widetilde{\mbH}_H^2) \\
\left(\text{resp.,} \  (\delta_{\mcE^\sss}, \delta_H^1, \delta_H^2) : 
(\mr{Def}_{(X, \mcE^\sss)}, \widetilde{\mbH}_H^1, \widetilde{\mbH}_H^2) \migi (\mr{Def}_X, H^1, H^2)
 \right) \notag
\end{align}
(cf. (\ref{L093f5}),  (\ref{QW7011})) forms  a morphism of deformation triples.

\LSP
%---------------------------[begin subsection]-------------
\subsection{Deformation of flat indigenous $(G, H)$-bundles} \label{SS077}
Next, suppose that the indigenous $(G, H)$-bundle  $\mcE^\sss := (\mcE_H, \nabla)$ is flat.

\SSP
%------------------------------------------------------------
\bde \label{D0100}
A deformation  $\mcE_R^\sss$ (resp., $(X_R, \mcE_R^\sss)$) of $\mcE^\sss$ (resp., $(X, \mcE^\sss)$) over $R \in \mr{Ob}(\mfA \mfr \mft_{/k})$
is called  {\bf flat} if the indigenous $(G, H)$-bundle $\mcE_R^\sss$ is flat.
\ede
%------------------------------------------------------------
\SSP
 
 For each $R \in \mr{Ob}(\mfA \mfr \mft_{/k})$,
 we shall 
 denote by
 \begin{align}
 \mr{Def}_{\mcE^\sss}^\flat (R) \ \left(\text{resp.,} \   \mr{Def}_{(X, \mcE^\sss)}^\flat (R) \right)
 \end{align}
 the subset of $\mr{Def}_{\mcE^\sss}(R)$ (resp., $\mr{Def}_{(X, \mcE^\sss)}(R)$) consisting of flat deformations.
 The assignment $R \mapsto \mr{Def}_{\mcE^\sss}^\flat (R)$ (resp., $R \mapsto \mr{Def}_{(X, \mcE^\sss)}^\flat (R)$)
 defines 
  a  functor with gdt
 \begin{align}
 \mr{Def}_{\mcE^\sss}^\flat \ \left(\text{resp.,} \  \mr{Def}_{(X, \mcE^\sss)}^\flat \right) : \mfA \mfr \mft_{/k} \migi \mfS \mfe \mft.
 \end{align}
Also, we have  a natural inclusion of  functors $\mr{Def}_{\mcE^\sss}^\flat \migi \mr{Def}_{\mcE^\sss}$ (resp., $\mr{Def}_{(X, \mcE^\sss)}^\flat \migi \mr{Def}_{(X, \mcE^\sss)}$).

In what follows, we compute  the tangent and obstruction spaces of  $\mr{Def}_{\mcE^\sss}^\flat$,  $\mr{Def}_{(X, \mcE^\sss)}^\flat$.
Let  $\Box \in \{ G, H \}$.
The flatness of $\nabla$ implies that of $\nabla^\mr{ad}_G$, so the equality 
 $\nabla^{\mr{ad}(2)}_G \circ \nabla^\mr{ad}_\Box =0$ holds (cf. (\ref{E0124})).
It follows that we can  restrict the codomain (resp., the domain and codomain) of $\nabla^{\mr{ad}}_G$ to obtain a $k$-linear morphism
\begin{align}
\nabla^{\mr{ad}}_{\Box, Z} : \mr{Lie}(\Box)_{\mcE_\Box} \migi  \mr{Ker} (\nabla_G^{\mr{ad}(2)})
\end{align}
The commutativity of (\ref{Dig01})
implies $\mr{Im}(\widetilde{\nabla}^\mr{ad}_\Box) \subseteq \mr{Im}(\nabla_G^\mr{ad}) \subseteq \mr{Ker} (\nabla_G^{\mr{ad}(2)})$.
Hence,
the morphism $\widetilde{\nabla}_\Box^{\mr{ad}} : \widetilde{\mcT}_{\mcE_\Box} \migi \Omega_{X} \otimes \mfg_{\mcE_G}$ induces 
 a $k$-linear morphism
 \begin{align} \label{FR45}
 \widetilde{\nabla}_{\Box, Z}^{\mr{ad}} : \widetilde{\mcT}_{\mcE_\Box} \migi \mr{Ker}(\nabla_G^{\mr{ad}(2)}).
 \end{align}
 For each $i \geq 0$, we shall write 
 \begin{align}
   \mbH^i_{\Box, Z} := \mbH^i (X, \mcK^\bullet [\nabla_{\Box, Z}^\mr{ad}]),
  \hspace{5mm}
  \widetilde{\mbH}^i_{\Box, Z} := \mbH^i (X, \mcK^\bullet [\widetilde{\nabla}_{\Box, Z}^\mr{ad}]).
 \end{align}
 Then, 
 the pair of morphisms $(\eta_\Box, \mr{id}_{\mr{Ker}(\nabla_G^{\mr{ad}(2)})})$ defines
 a morphism $\eta^\bullet_{\Box, Z} : \mcK^\bullet [\widetilde{\nabla}^{\mr{ad}}_{\Box, Z}] \migi  \mcK^\bullet [\nabla_{G, Z}^{\mr{ad}}]$,  and hence,  induces  a morphism of $k$-vector spaces
\begin{align} \label{Fer4}
\mbH^i (\eta^\bullet_{\Box, Z}) : \widetilde{\mbH}^i_{\Box, Z} \migi \mbH^i_{G, Z}.
\end{align}
 In the case of $\Box = H$, this morphism in fact is an isomorphism.
 \SSP
 %------------------------------------------------------------
 \bpr \label{P023}
 Recall that  $\mcE^\sss$
 is assumed to be flat.
\begin{itemize}
\item[(ii)]
The collections $(\mr{Def}^\flat_{\mcE^\sss}, \mbH^1_{H, Z}, \mbH^2_{H, Z})$
and $(\mr{Def}^\flat_{(X, \mcE^\sss)}, \widetilde{\mbH}_{H, Z}^1, \widetilde{\mbH}_{H, Z}^2)$ ($\cong (\mr{Def}^\flat_{(X, \mcE^\sss)}, \mbH^1_{G, Z}, \mbH^2_{G, Z})$ by (\ref{Fer4})) form deformation triples.
%The tangent space of $\mr{Def}^\flat_{\mcE^\sss}$ is isomorphic to
%$\mbH^1_{H, Z}$, and this functor has an obstruction theory by putting 
%$\mbH^2_{H, Z}$ as the obstruction space.
%\item[(ii)]
% The tangent space  of  $\mr{Def}^\flat_{(X, \mcE^\sss)}$
% is isomorphic to $\widetilde{\mbH}_{H, Z}^1$ ($\cong \mbH^1_{G, Z}$ by (\ref{Fer4})),
% and this functor has an obstruction theory by putting 
%  $\widetilde{\mbH}_{H, Z}^2$ ($\cong \mbH^2_{G, Z}$) as the  obstruction space.
  \item[(ii)]
  Suppose that $X$ is proper over $k$.
  Then, 
  there exist pro-representable hulls for 
   $\mr{Def}^\flat_{\mcE^\sss}$ and $\mr{Def}^\flat_{(X, \mcE^\sss)}$.
  \item[(iii)]
  Regarding the pro-representability of $\mr{Def}^\flat_{\mcE^\sss}$ (resp., $\mr{Def}^\flat_{(X, \mcE^\sss)}$),
   the statement of  Corollary \ref{Cf02k2} with $\mbH^i_H$'s (resp., $\widetilde{\mbH}^i_{H}$'s) replaced by $\mbH_{H, Z}^i$'s (resp., $\widetilde{\mbH}_{H, Z}^i$'s) holds  for $\mr{Def}^\flat_{\mcE^\sss}$ (resp., $\mr{Def}^\flat_{(X, \mcE^\sss)}$).
\end{itemize}
 \epr
 %------------------------------------------------------------
 \begin{proof}
 The assertions follow from  Proposition \ref{P022} (and its proof) together with the following lemma.
 \end{proof}
%------------------------------------------------------------

\SSP

%---------------------------------------------------------
\ble \label{L010}
Let $0 \migi M \migi R \migi R_0 \migi 0$ be a small extension in $\mfA \mfr \mft_{/k}$ and 
$(X_{R}, \mcE_{R}^\sss)$ (where $\mcE_{R}^\sss := (\mcE_{R, H}, \nabla_{R})$)   a flat deformation of $(X, \mcE^\sss)$ over $R$. 
We shall write   $\mcE_{R, G} := \mcE_{R, H} \times^H G$.
Also, let
$a$  be an element of
$\Gamma (X, M \otimes \Omega_{X} \otimes \mfg_{\mcE_G}) \left(\subseteq \mr{Hom}_{\mcO_X} (\mcT_{X}, M \otimes \widetilde{\mcT}_{\mcE_G}) \right)$, 
 and we regard it  as   an element of $\mr{Hom}_{\mcO_{X_{R}}}(\mcT_{X_{R}/R}, \widetilde{\mcT}_{\mcE_{R, G}/R})$ in a natural manner.
Then, the curvature $\psi_{\nabla_{R}+ a}$ of the sum $\nabla_{R} + a$ is given by
\begin{align} \label{Ehj2}
\psi_{\nabla_{R}+a} = (\mr{id}_M \otimes \nabla_G^{\mr{ad}(2)}) (a)
\end{align}
via the natural   inclusion 
\begin{align} 
\Gamma (X, M \otimes (\bigwedge^2 \Omega_{X}) \otimes \mfg_{\mcE_G}) \subseteq \Gamma (X_{R}, (\bigwedge^2 \Omega_{X_{R}/R}) \otimes_{R} \mfg_{\mcE_{R, G}}).
\end{align}
In particular, $\psi_{\nabla_{R} + a} =0$ if and only if $a \in \mr{id}_M \otimes \mr{Ker}(\nabla_G^{\mr{ad}(2)}) \left(= \mr{Ker}(\mr{id}_M \otimes \nabla_G^{\mr{ad}(2)}) \right)$.
\ele
%---------------------------------------------------------
\begin{proof}
For any local sections $\partial_1, \partial_2 \in \mcT_{X_{R}/R}$,
we have
\begin{align} \label{Eddw}
\psi_{\nabla_{R}+a}^\triangleright (\partial_1, \partial_2) = & \  [(\nabla_{R}+a)(\partial_1), (\nabla_{R}+ a)(\partial_2)] - (\nabla_{R}+a)([\partial_1, \partial_2]) \\
= & \ \cancel{[\nabla_{R}(\partial_1), \nabla_{R}(\partial_2)]} + [\nabla_{R}(\partial_1), a (\partial_2)] + [a (\partial_1), \nabla_{R}(\partial_2)]\notag \\
& \ + [a (\partial_1), a (\partial_2)] - \cancel{\nabla_{R}([\partial_1, \partial_2])} - a ([\partial_1, \partial_2]) \notag \\
= & \ [\nabla_{R}(\partial_1), a (\partial_2)] + [a (\partial_1), \nabla_{R}(\partial_2)]+ \cancel{[a (\partial_1), a (\partial_2)]}- a ([\partial_1, \partial_2]) \notag \\
= & \ [\nabla_{R}(\partial_1), a (\partial_2)] + [a (\partial_1), \nabla_{R}(\partial_2)]- a ([\partial_1, \partial_2]), \notag
\end{align}
 where the third equality follows from the flatness  assumption on $\nabla_{R}$
 and the last equality follows from the fact that $\mr{Im}(a) \subseteq M   \otimes \widetilde{\mcT}_{\mcE_{R, G}/R}$ and $M^2 =0$.
 In particular, the assignment $a \mapsto \psi_{\nabla_{R}+a}$ is additive.
 Hence, in order to prove (\ref{Ehj2}), it suffices to consider the case where $a = \omega \otimes v$ for  $\omega \in \Omega_{X}$, $v \in  M \otimes \mfg_{\mcE}$.
For $\partial_1$ and  $\partial_2$ as above (and we occasionally  regard them as local sections of $\mcT_{X}$ via reduction via $\mfm_{R}$), 
the following sequence of equalities holds:
\begin{align}
& \ (\mr{id}_M \otimes \nabla_G^{\mr{ad}(2)}) (\omega \otimes v)^\triangleright (\partial_1, \partial_2) \\
= & \ (d \omega \otimes v)^\triangleright (\partial_1, \partial_2) - (\omega \wedge (\mr{id}_M \otimes \nabla_G^{\mr{ad}})(v))^\triangleright(\partial_1, \partial_2) \notag  \\
= & \  \left( \partial_1 (\omega^\triangleright (\partial_2))-\partial_2 (\omega^\triangleright (\partial_1))- \omega^\triangleright ([\partial_1, \partial_2])\right) \otimes v
\notag \\
& \ - \omega^\triangleright (\partial_1) \cdot (\mr{id}_M \otimes \nabla_G^\mr{ad})(v)^\triangleright (\partial_2) + \omega^\triangleright (\partial_2) \cdot (\mr{id}_M \otimes\nabla_G^\mr{ad})(v)^\triangleright (\partial_1) \notag \\
= & \ \partial_1 (\omega^\triangleright (\partial_2)) \otimes v - \partial_2 (\omega^\triangleright (\partial_1)) \otimes v  - \omega^\triangleright ([\partial_1, \partial_2]) \otimes v - \omega^\triangleright (\partial_1) \cdot [\nabla (\partial_2), v]   \notag \\
 & \ + \omega^\triangleright (\partial_2) \cdot [\nabla (\partial_2), v]\notag \\
=& \ \left(\omega^\triangleright (\partial_2)\cdot [\nabla_{R}(\partial_1), v] + d_{\mcE_{R', G}}(\nabla_{R}(\partial_1))(\omega^\triangleright (\partial_2) \otimes v)\right) \notag \\
& \ - \left(\omega^\triangleright (\partial_1)\cdot [\nabla_{R}(\partial_2), v] + d_{\mcE_{R', G}}(\nabla_{R}(\partial_2))(\omega^\triangleright (\partial_1) \otimes v)\right) - \omega^\triangleright ([\partial_1, \partial_2]) \otimes v \notag \\
= & \ [\nabla_{R}(\partial_1), \omega^\triangleright (\partial_2) \otimes v] + [\omega^\triangleright (\partial_1) \otimes v, \nabla_{R}(\partial_2)] - \omega^\triangleright ([\partial_1, \partial_2]) \otimes v\notag \\
\stackrel{(\ref{Eddw})}{=} & \ \psi_{\nabla_{R}+\omega \otimes v}^\triangleright (\partial_1, \partial_2), \notag
\end{align}
where the fifth equality follows from the fact that $(\widetilde{\mcT}_{\mcE_{R, G}/R}, d_{\mcE_{R, G}})$ forms a Lie algebroid.
This completes the proof of this lemma.
\end{proof}
 %---------------------------------------------------------
 \SSP
 
  %---------------------------------------------------------
 \begin{rema} \label{FFmK0}
 We shall denote by
 \begin{align} \label{EJJKL}
 \mr{Def}_{(\mcE_G, \nabla)}^\flat \ \left(\text{resp.,} \ \mr{Def}_{(X, \mcE_G, \nabla)}^\flat \right)
 \end{align}
 the subfunctor of $\mr{Def}_{(\mcE_G, \nabla)}$ (resp., $\mr{Def}_{(X, \mcE_G, \nabla)}^\flat$) introduced in (\ref{RwwqaW}) classifying flat deformations.
 Note that the result of Lemma \ref{L010} holds without the assumption that $\mcE_G$ has an $H$-reduction defining the structure of indigenous $(G, H)$-bundle, i.e., $\mcE_H$.
 Hence, this lemma implies that
 $(\mr{Def}_{(\mcE_G, \nabla)}^\flat, \mbH_{G, Z}^1, \mbH_{G, Z}^2)$
 (resp., $(\mr{Def}_{(X, \mcE_G, \nabla)}^\flat, \widetilde{\mbH}_{G, Z}^1, \widetilde{\mbH}_{G, Z}^2)$) forms a deformation triple.
 % $ \mr{Def}_{(\mcE_G, \nabla)}^\flat$
 %(resp., $\mr{Def}_{(X, \mcE_G, \nabla)}^\flat$)
 %defines a functor with gdt whose tangent space is isomorphic to $\mbH_{G, Z}^1$ (resp., $\widetilde{\mbH}_{G, Z}^1$) and has an obstruction theory by putting $\mbH^2_{G, Z}$ (resp., $\widetilde{\mbH}_{G, Z}^2$) as the obstruction space.
 \end{rema}
  %---------------------------------------------------------

\LSP
%---------------------------[begin subsection]-------------
\subsection{Deformation of dormant indigenous $(G, H)$-bundles} \label{SS0100}

Next, suppose further  that the indigenous $(G, H)$-bundle $\mcE^\sss$ is dormant.

 \SSP
 %------------------------------------------------------------
\bde \label{D011}
A deformation $\mcE_R^\sss$ (resp., $(X_R, \mcE^\sss_R)$) of $\mcE^\sss$ (resp., $(X, \mcE^\sss)$) over $R \in \mr{Ob}(\mfA \mfr \mft_{/k})$ is called {\bf dormant} if $\mcE_R^\sss$ is dormant.
\ede
%------------------------------------------------------------
\SSP
 
 For each $R \in \mr{Ob}(\mfA \mfr \mft_{/k})$, we shall denote by
 \begin{align} \label{ei4s038}
  \mr{Def}^{^\mr{Zzz...}}_{\mcE^\sss}\! (R) \
  \left(\text{resp.,} \ \mr{Def}^{^\mr{Zzz...}}_{(X, \mcE^\sss)}(R) \right)
 \end{align}
 the subset of $\mr{Def}_{\mcE^\sss}^\flat (R)$ (resp., $\mr{Def}_{(X, \mcE^\sss)}^\flat (R)$) consisting of dormant deformations.
 One may verify that the assignment  $R \mapsto  \mr{Def}^{^\mr{Zzz...}}_{\mcE^\sss}\! (R)$ (resp., $R \mapsto \mr{Def}^{^\mr{Zzz...}}_{(X, \mcE^\sss)}(R)$) defines a  functor with gdt
 \begin{align} \label{esoriaa2}
  \mr{Def}^{^\mr{Zzz...}}_{\mcE^\sss} \  \left(\text{resp.,} \  \mr{Def}^{^\mr{Zzz...}}_{(X, \mcE^\sss)}  \right) : \mfA \mfr \mft_{/k} \migi \mfS \mfe \mft.
 \end{align}
  Also, there exists a natural inclusion  of  functors $ \mr{Def}^{^\mr{Zzz...}}_{\mcE^\sss} \! \migi  \mr{Def}^\flat_{\mcE^\sss}$
  (resp., $ \mr{Def}^{^\mr{Zzz...}}_{(X, \mcE^\sss)} \migi  \mr{Def}^\flat_{(X, \mcE^\sss)}$).

  In what follows, let us describe the tangent and obstruction spaces of $\mr{Def}_{\mcE^\sss}^{^\mr{Zzz...}}\!$ and $\mr{Def}_{(X, \mcE^\sss)}^{^\mr{Zzz...}}$.
  Note that the morphism $\nabla_{G, Z}^\mr{ad}$  (resp.,  $\nabla_{H, Z}^{\mr{ad}}$; resp., $\widetilde{\nabla}_{H, Z}^{\mr{ad}}$) induces, via restricting the codomain, a morphism
  \begin{align}
   \nabla_{G, B}^{\mr{ad}} : \mfg_{\mcE_G}  \migi  \mr{Im}(\nabla_G^\mr{ad})
   \ 
  \left(\text{resp.,} \ \nabla_{H, B}^{\mr{ad}} : \mfh_{\mcE_H}  \migi  \mr{Im}(\nabla_G^\mr{ad}); \text{resp.,} \    \widetilde{\nabla}_{H, B}^{\mr{ad}} : \widetilde{\mcT}_{\mcE_H}  \migi  \mr{Im}(\nabla^\mr{ad}_G)  \right). \hspace{-5mm}
  \end{align}
For each $i \geq 0$, we shall write
\begin{align}
\mbH^i_{G, B} := \mbH^i (X, \mcK^\bullet [\nabla_{G, B}^\mr{ad}]),
\hspace{5mm}
\mbH^i_{H, B} := \mbH^i (X, \mcK^\bullet [\nabla_{H, B}^\mr{ad}]),
\hspace{5mm}
\widetilde{\mbH}^i_{H, B} := \mbH^i (X, \mcK^\bullet [\widetilde{\nabla}_{H, B}^\mr{ad}]).
\end{align}
  The isomorphism $\eta^\bullet_{H, Z} : \mcK^\bullet [\widetilde{\nabla}_{H, Z}^{\mr{ad}}] \isom \mcK^\bullet [\nabla_{G, Z}^{\mr{ad}}]$ is restricted to an isomorphism of complexes
  $\eta^\bullet_{H, B} : \mcK^\bullet [\widetilde{\nabla}_{H, B}^{\mr{ad}}] \isom \mcK^\bullet [\nabla_{G, B}^{\mr{ad}}]$.
  By applying the functor $\mbH^i (X, -)$ to this isomorphism for each $i \geq 0$, 
  we obtain  an isomorphism of $k$-vector spaces
  \begin{align} \label{Ty67kk}
  \mbH^1 (\eta^\bullet_{H, B}) : \widetilde{\mbH}^i_{H, B} \isom \mbH^i_{G, B}.
  \end{align}
   Moreover,
    the natural inclusion $\mr{Ker}(\nabla_G^\mr{ad})[0] \migiincl  \mcK^\bullet [\nabla^{\mr{ad}}_{G, B}]$ is a quasi-isomorphism, so   it yields an  isomorphism 
 \begin{align} \label{Ty67}
 H^i (X, \mr{Ker}(\nabla_G^\mr{ad})) \isom \mbH^i_{G, B}.
 \end{align}
   
\SSP 
 %------------------------------------------------------------
 \bpr \label{P050}
 Recall that  $\mcE^\sss$ is  assumed to be dormant.
 \begin{itemize}
 \item[(i)]
 The collections 
 $(\mr{Def}_{\mcE^\sss}^{^\mr{Zzz...}}\!, \mbH^1_{H, B}, \mbH^2_{H, B})$ and
 $(\mr{Def}_{(X, \mcE^\sss)}^{^\mr{Zzz...}}, \widetilde{\mbH}^1_{H, B}, \widetilde{\mbH}^2_{H, B})$ (which is isomorphic to $(\mr{Def}_{(X, \mcE^\sss)}^{^\mr{Zzz...}}, H^1 (X, \mr{Ker}(\nabla_G^\mr{ad})), H^2 (X, \mr{Ker}(\nabla_G^\mr{ad})))$ by (\ref{Ty67kk}) and (\ref{Ty67}))
 form deformation triples.
  %The tangent space of $\mr{Def}_{\mcE^\sss}^{^\mr{Zzz...}}\!$ 
 %is isomorphic to $\mbH^1_{H, B}$, and this functor has an obstruction theory
 %by putting 
 %$\mbH^2_{H, B}$  as the obstruction space.
 %\item[(ii)]
% The tangent space of  $\mr{Def}_{(X, \mcE^\sss)}^{^\mr{Zzz...}}$
 %is isomorphic to 
  %$\widetilde{\mbH}^1_{H, B}$ ($\cong H^1 (X, \mr{Ker}(\nabla_G^\mr{ad}))$ by  (\ref{Ty67kk}) and (\ref{Ty67})), and this functor has an obstruction theory by putting 
 %$\widetilde{\mbH}^2_{H, B}$ ($\cong H^2 (X, \mr{Ker}(\nabla_G^\mr{ad}))$) as the obstruction space.
  \item[(ii)]
 Suppose that $X$ is proper over $k$.
 Then, 
 there exist pro-representable hulls for 
 $\mr{Def}_{\mcE^\sss}^{^\mr{Zzz...}}\!$ and $\mr{Def}_{(X, \mcE^\sss)}^{^\mr{Zzz...}}$. 
 \item[(iii)]
  Regarding the pro-representability of $\mr{Def}^{^\mr{Zzz...}}_{\mcE^\sss}\!$ (resp., $\mr{Def}^{^{\mr{Zzz...}}}_{(X, \mcE^\sss)}$),
   the statement of  Corollary \ref{Cf02k2} with $\mbH^i_H$'s (resp., $\widetilde{\mbH}^i_{H}$'s) replaced by $\mbH_{H, B}^i$'s (resp., $\widetilde{\mbH}_{H, B}^i$'s) holds  for $\mr{Def}^{^\mr{Zzz...}}_{\mcE^\sss}\!$ (resp., $\mr{Def}^{^\mr{Zzz...}}_{(X, \mcE^\sss)}$).
 \end{itemize}
 \epr
 %------------------------------------------------------------
 \begin{proof}
 The assertions follow from Proposition \ref{P022} (and its proof) together with the following lemma.
 \end{proof}
%------------------------------------------------------------

 \SSP

%-----------------------------------------------------------------------[begin lemma]------------------
%\vspace{3mm}
\ble \label{y0157}
 Let us keep the notation in Lemma \ref{L010} and suppose further that the deformation $(X_{R}, \mcE_{R}^\sss)$ is dormant and the $R$-connection $\nabla_{R} + a$ is flat (or equivalently, $a \in \mr{id}_M \otimes \mr{Ker}(\nabla_G^{\mr{ad}(2)})$ by Lemma \ref{L010}).
 Then, the  $p$-curvature ${^p}\psi_{\nabla_{R}+a}$ of $\nabla_{R}+a$ is given by 
\begin{align} \label{910}
 {^p \psi}_{\nabla_{R}+a}
  = - (\mr{id}_M \otimes  C_{(\mfg_{\mcE_G}, \nabla_G^\mr{ad})} )(a)
\end{align}
(cf. (\ref{Efty}) for the definition of $C_{(-, -)}$) under the natural inclusion relation
\begin{align}
\Gamma (X, M \otimes F^*_X(\Omega_{X})\otimes \mfg_{\mcE_G}) \subseteq \Gamma (X_{R}, F^*_{X_{R}}(\Omega_{X_{R}/R}) \otimes \mfg_{\mcE_{R, G}}).
\end{align}
In particular, 
${^p}\psi_{\nabla_{R}+a} = 0$ if and only if $a$ may be expressed, Zariski locally on $X$, as
$a = \nabla_G^\mr{ad}(v)$ for some local section $v \in \mfg_{\mcE_G}$. 
  \ele
%-----------------------------------------------------------------------[begin proof]-------------------
\begin{proof}
Let us take a local section $\partial \in \mcT_{X_{R}/R}$, which will be 
 occasionally regarded  as a local section of $\mcT_{X}$ via reduction modulo $\mfm_{R}$.
 Observe 
 the following  sequence of equalities   in the  enveloping algebra of the Lie algebroid $(\widetilde{\mcT}_{\mcE_{G}}, d_{\mcE_{G}})$:
\begin{align}  \label{xi1}
 & \  \ \  \  \langle F_{X}^{-1}(\partial), - (\mr{id}_M \otimes  C_{(\mfg_{\mcE_G}, \nabla_G^\mr{ad})} )(a) \rangle \\
  & =  -  \langle  F_k^*(\partial), (\mr{id}_M \otimes C_{(\mfg_{\mcE_G}, \nabla_G^\mr{ad})})(a)\rangle \notag  \\
 & =   -  \langle \partial^{[p]}, a\rangle + (\mr{id}_M \otimes  (\nabla_{G, \partial}^\mr{ad})^{p-1})(\langle \partial, a  \rangle)  \notag \\
&=   -   \langle \partial^{[p]}, a \rangle +  (\mr{id}_M \otimes  \mr{ad}(\nabla (\partial))^{p-1})(\langle \partial, a \rangle) \notag \\
& = -  \langle \partial^{[p]}, a\rangle + \sum_{i=0}^{p-1} (-1)^i \cdot \begin{pmatrix} p-1 \\ i \end{pmatrix} \cdot \nabla (\partial)^{p-1 -i} \cdot \langle \partial, a \rangle  \cdot \nabla (\partial)^{i} \notag \\ 
& =    - \langle \partial^{[p]}, a\rangle +  \sum_{i=0}^{p-1} \nabla (\partial)^{p-1 -i} \cdot \langle  \partial, a\rangle  \cdot \nabla (\partial)^{i}, \notag
\end{align}
where the second  equality follows from  (\ref{Wer4}).
On the other hand,  we have 
\begin{align}  \label{xi2}
 & \  \langle F_{X}^{-1}(\partial), {^p}\psi_{\nabla_{R} + a} \rangle  \\
= & \  (\nabla_{R} + a)(\partial)^{[p]}-(\nabla_{R} + a)(\partial^{[p]}) \notag \\
= &  \  (\nabla_{R} (\partial) +  \langle \partial, a \rangle)^{[p]} -  \nabla_{R} (\partial^{[p]}) - \langle \partial^{[p]}, a \rangle\notag \\
=  & \  \sum_{i=0}^{p-1} \nabla (\partial)^{ p-1-i} \cdot  \langle \partial, a \rangle \cdot  \nabla (\partial)^{i}  
 +\cancel{\nabla (\partial)^{[p]} - \nabla  (\partial^{[p]})}  - \langle \partial^{[p]}, a \rangle \notag \\
= & \ -    \langle \partial^{[p]}, a  \rangle +    \sum_{i=0}^{p-1} \nabla (\partial)^{ p-1-i} \cdot  \langle \partial, a \rangle \cdot  \nabla (\partial)^{i}, \notag
\end{align} 
where  the last equality follows from   the assumption ${^p}\psi_\nabla  =0$.
Thus,    (\ref{xi1}) and  (\ref{xi2}) together show that
\begin{align}
\langle F_{X}^{-1}(\partial), {^p}\psi_{\nabla_{R} + a} \rangle =  \langle F_{X}^{-1}(\partial), - (\mr{id}_M \otimes C_{(\mfg_{\mcE_G}, \nabla_G^\mr{ad})})(a) \rangle.
\end{align}
By considering this equality for every  $\partial \in \mcT_{X}$, 
 we obtain  (\ref{910}).
 This   completes  the proof of the former assertion.
The latter assertion follows from
 the former assertion  and the commutativity of (\ref{Cartier3F}).
\end{proof}
%-----------------------------------------------------------------------[end lemma]-------------------

\SSP

\begin{rema} \label{QAAEoe8}
If the underlying space $X$ is a curve, then the above lemma was (essentially) proved in 
  ~\cite[Chap.\,6, \S\,6.4, Proposition 6.10]{Wak8}.
  \end{rema}

\SSP
%-----------------------------------------------------------------
\begin{exa}[Dormant indigenous $(\mr{PGL}_2, \mr{PGL}_2^\circledcirc)$-bundles on a curve] \label{E099}
Let $X$ be a smooth proper  curve over $k$ of genus $g>1$ and $\mcE^\sss := (\mcE_H, \nabla)$ be a dormant indigenous $(\mr{PGL}_2, \mr{PGL}_2^\circledcirc)$-bundle on $X/k$.
We know  that   $H^i (X, \mr{Ker}(\nabla_G^\mr{ad})) =0$ for $i=0,2$ and $\mr{dim}_k(H^1 (X, \mr{Ker}(\nabla_G^\mr{ad}))) =3g-3$ (cf. ~\cite[Chap.\,6, \S\,6.2, Proposition 6.4, and \S\,6.5,  Proposition 6.17]{Wak8}).
Hence,  by  Proposition \ref{P050}, (iii), the deformation functor  $\mr{Def}_{(X, \mcE^\sss)}^{^\mr{Zzz...}}$ of $(X, \mcE^\sss)$ may be pro-representable by   $k[\! [x_1, \cdots, x_{3g-3} ]\! ]$.
This result  is consistent with the  fact  that the moduli stack classifying  such  curves equipped with a dormant indigenous bundle 
 is 
smooth  and of dimension $3g-3$ (cf. ~\cite[Chap.\,II, \S\,2.3, Theorem 2.8]{Mzk2}).
\end{exa}
%-----------------------------------------------------------------

\LSP
%-----------------------------------------------------------------------------------------
\subsection{Ordinariness of dormant indigenous  $(G, H)$-bundles} \label{SSfru002}

In this subsection, we introduce the notion of ordinariness for dormant indigenous $(G, H)$-bundles (of level $1$) and prove a related lifting property  of the underlying space. 
(However, the discussion and results here  will not be used after this subsection.)
Suppose that $G$  is affine, and 
 we shall fix 
  a dormant indigenous $(G, H)$-bundle  $\mcE^\sss := (\mcE_H, \nabla)$  on $X/k$.
  Write $\mcE_G := \mcE_H \times^H G$, and denote by $\mcE_G^\nabla$ the $G$-bundle on $X^{(1)}$ corresponding to the $p$-flat $G$-bundle $(\mcE_G, \nabla)$ via (\ref{EEjj32}).
  In particular, there exists an isomorphism of $G$-bundles $F_{X/k}^{*}(\mcE^\nabla_G) \isom \mcE_G$, which induces an isomorphism of $\mcO_X$-Lie algebras $F^{*}_{X/k}(\mfg_{\mcE_G^\nabla}) \isom \mfg_{\mcE_G}$.
For each $l \geq 0$,     
we obtain the $k$-linear morphism
\begin{align} \label{GJidf}
H^l (\varphi ) : H^l (X^{(1)}, \mfg_{\mcE_G^\nabla}) \migi H^l (X, \mfg_{\mcE_G}/\mfh_{\mcE_H})
\end{align}
induced by the  composite of natural morphisms
$\varphi : F_{X/k}^{-1}(\mfg_{\mcE_G^\nabla}) \left(= \mr{Ker}(\nabla^\mr{ad}_G) \right)\migiincl \mfg_{\mcE_G} \migisurj \mfg_{\mcE_G}/\mfh_{\mcE_H}$.
\SSP
%--------------------------------------------------------------------------------------
\bde \label{ereao283}
We shall say that 
$\mcE^\sss$ (or $(X, \mcE^\sss)$) is {\bf ordinary}
if 
$H^1 (\varphi )$ is an isomorphism and $H^2 (\varphi )$ is injective. 
\ede
%--------------------------------------------------------------------------------------

\SSP

The ordinariness just defined is closed related to the existence of a deformation of $\mcE^\sss$, as described in the following proposition.

\SSP
%--------------------------------------------------------------------------------------
\bpr \label{Etouaojm}
Suppose that 
 $\mcE^\sss$ is ordinary.
Then, 
the restriction
\begin{align}
\delta_{\mcE^\sss}^{^\mr{Zzz...}} : \mr{Def}_{(X, \mcE^\sss)}^{^\mr{Zzz...}} \migi \mr{Def}_X
 \end{align}
of $\delta_{\mcE^\sss}$ (cf. (\ref{L093f5})) to $\mr{Def}_{(X, \mcE^\sss)}^{^\mr{Zzz...}}$
is a semi-isomorphism.
In particular, for a small extension  in $\mfA \mfr \mft_{/k}$  of the form $0 \migi M \migi R \migi k \migi 0$ and a deformation $X_R$ of $X$ over $R$, there exists a unique (up to isomorphism)  deformation  of  $\mcE^\sss$ defined over $X_R$.
If, moreover, $\mr{Def}_X$ is  rigid, then 
this morphism is an isomorphism.
\epr
%--------------------------------------------------------------------------------------
\begin{proof}
Observe that the following diagram is commutative (cf. the commutativity of (\ref{QW70fg0})):
\begin{align} \label{Dig01ddf}
\vcenter{\xymatrix@C=46pt@R=36pt{
& \mcK^\bullet [\widetilde{\nabla}_{H, B}^\mr{ad}] \ar[r]^-{\mr{inclusion}} \ar[d]_-{\wr}^-{\eta^\bullet_{H, B}} & \mcK^\bullet [\widetilde{\nabla}_H^\mr{ad}]\ar[r]^{(d_{\mcE_H}, 0)} \ar[d]_-{\wr}^-{\eta^\bullet_{H}} & \mcT_X[0] \ar[d]_-{\wr}^-{\mr{KS}_{(\mcE_H, \nabla)}}
\\
F_{X/k}^{-1}(\mfg_{\mcE_G^\nabla}) [0]
\ar[r]_-{\mr{inclusion}} & \mcK^\bullet [\nabla_{G, B}^\mr{ad}] \ar[r]_-{\mr{inclusion}}& \mcK^\bullet [\nabla^\mr{ad}_G] \ar[r] & \mfg_{\mcE_G}/\mfh_{\mcE_H}[0],
}}
\end{align}
where the rightmost lower horizontal arrow is given  by the pair of the natural quotient $\mfg_{\mcE_G} \migisurj \mfg_{\mcE_G}/\mfh_{\mcE_H}$ and the zero map $\Omega_{X} \otimes \mfg_{\mcE_G} \migi 0$. 
For each $i \geq 0$, denote by $\delta_{H, B}^i$ the morphism $\widetilde{\mbH}^i_{H, B} \migi H^i$ induced by the composite of the two upper horizontal arrows in this diagram.
The composite of all the lower horizontal arrows coincides with  $\varphi$.
Hence,  since $\mcE^\sss$  is ordinary and (\ref{Ty67}) is an isomorphism,
the commutativity of (\ref{Dig01ddf}) implies that 
$\delta_{H, B}^1$ is an isomorphism and $\delta_{H, B}^2$ is injective.
On the other hand, 
by the definition of $\delta_{\mcE^\sss}^{^\mr{Zzz...}}\!$,
 the triple  
$(\delta_{\mcE^\sss}^{^\mr{Zzz...}}, \delta^1_{H, B}, \delta^2_{H, B})$ forms a morphism of deformation triples $(\mr{Def}_{(X, \mcE^\sss)}^{^\mr{Zzz...}}, \widetilde{\mbH}^1_{H, B}, \widetilde{\mbH}^2_{H, B}) \migi (\mr{Def}_X, H^1, H^2)$ (cf. Proposition \ref{P050}, (i)).
Thus,  the assertion follows from   Propositions \ref{pyp90s4}, (i), and  \ref{E2Rgju2}.
\end{proof}
%--------------------------------------------------------------------------------------
\SSP

%\SSP

%---------------------------
\begin{exa} \label{eso0w}
In the case where $(G, H) = (\mr{PGL}_2, \mr{PGL}_2^\circledcirc)$ and $X$ is a smooth proper curve over $k$,
essentially the same  notion of ordinariness as above  can be found in ~\cite[\S\,3, Definition 3.2]{Wak}.
 To verify the ordinariness of a given dormant  indigenous $(\mr{PGL}_2, \mr{PGL}_2^\circledcirc)$-bundle $\mcE^\sss$, 
we only have to consider whether  $H^1 (\varphi)$ is an isomorphism or not because of  the assumption   $\mr{dim} (X) =1$ (which implies  $H^2 (\varphi) =0$).
Hence, it follows from the discussion in the proof of Proposition \ref{Etouaojm} that  $\mcE^\sss$ is ordinary if and only if the differential $d \delta_{\mcE^\sss}^{^\mr{Zzz...}}$ of $\delta_{\mcE^\sss}^{^\mr{Zzz...}}$ is an isomorphism.

Now, denote by $\mfM_{g}$ the moduli stack classifying smooth proper  curves over $k$ of genus $g$ and by
$\mfM_g^{^\mr{Zzz...}}\!$ the moduli stack classifying pairs $(X, \mcE^\sss)$
consisting of such a curve $X$ 
and a dormant indigenous $ (\mr{PGL}_2, \mr{PGL}_2^\circledcirc)$-bundle $\mcE^\sss$ on $X$.
According to a result of S. Mochizuki  (cf.  ~\cite[Chap.\,II, \S\,2.3, Theorem 2.8]{Mzk2}),
the natural projection $\mfM_g^{^\mr{Zzz...}} \! \migi \mfM_g$ is finite and faithfully flat, and  moreover  the \'{e}tale locus in $\mfM_g^{^\mr{Zzz...}}\!$ relative to this projection   forms a dense open substack.
On the other hand, the above discussion implies that a pair $(X, \mcE^\sss)$ classified by $\mfM_g^{^\mr{Zzz...}}$
is ordinary if and only if the natural projection $\mfM_g^{^\mr{Zzz...}} \! \migi \mfM_g$ is \'{e}tale at the classifying point of  $(X, \mcE^\sss)$.
 Consequencely, every general pair $(X, \mcE^\sss)$ turns out to be  ordinary.
\end{exa}
\SSP

The following assertion gives another example of an  ordinary dormant indigenous  $(G, H)$-bundle.

\SSP

%-----------------------------------------------------------------------------------------
\bpr \label{MMM2fhg}
Let $X$ be 
 an ordinary Abelian variety of dimension $n > 0$.
 Then, any  dormant indigenous $(\mbG_m^{\times n}, \mbG_m^{\times n})$-bundle 
    on $X/k$ is ordinary.
\epr
%-----------------------------------------------------------------------------------------
\begin{proof}
For convenience, we write $G := \mbG_m^{\times n}$ and $H = \{ e \}$.
Let $\mcE^\sss := (\mcE_H, \nabla)$ be a dormant  indigenous $(G, G)$-bundle 
    on $X/k$, and write $\mcE_G := \mcE_H \times^H G$.
    By Proposition \ref{L0045}, (i) (and the bijection $\zeta^{\hh \Rightarrow \sss}$ asserted in Theorem \ref{TheoremA1}),
    $\mcE^\sss$ specifies an element $(\omega_1^+, \cdots, \omega_n^+)$ of $B_{X/k}^{(1), n, 0}$, where $\omega_i^+ := (\omega_i, \mcL_i, \eta_i)$ ($i=1, \cdots, n$).
 If we regard $\nabla_{\mcL_i}^{\mr{can}(0)}$ (for each $i$) as a $k$-connection on $\mcO_X$ via $\eta_i$,
 then the flat $G$-bundle $(\mcE_G, \nabla)$  arises from the flat vector bundle $\bigoplus_{i=1}^n (\mcO_X, \nabla_{\mcL_i}^{\mr{can}(0)})$.
 In particular, we obtain a natural identification $\mfg_{\mcE_G} = \mcO_X^{\oplus n}$.
  Next,   since the sections $\omega_1, \cdots, \omega_n$ are mutually distinct, 
    any locally defined endomorphism $h$  of $\bigoplus_{i=1}^n (\mcO_X, \nabla_{\mcL_i}^{\mr{can}(0)})$
    may be expressed as $\bigoplus_{i=1}^n \mr{mult}_{u_i}$, where each $\mr{mult}_{u_i}$ denotes the locally defined  endomorphism of $\mcO_X$ given by multiplication by a local section $u_i$ of $\mcO_{X^{(1)}}$ ($\subseteq \mcO_X$ via $F_{X/k}$).
    The assignment $h \mapsto (u_1, \cdots, u_n)$ gives 
  an identification  $\left(\mr{Ker}(\nabla_G^\mr{ad}) =  \right)\mfg_{\mcE_G^\nabla} = \mcO_{X^{(1)}}^{\oplus n}$.
  Under the identifications $\mfg_{\mcE_G} = \mcO_X^{\oplus n}$, $\mfg_{\mcE_G^\nabla} = \mcO_{X^{(1)}}^{\oplus n}$ obtained so far, 
$H^l (\varphi)$ (for each $l$) coincides  with the direct sum of $n$ copies of
  the morphism $H^l (X^{(1)}, \mcO_{X^{(1)}}) \migi H^l (X, \mcO_X)$  induced by $F_{X/k}$.
  Hence, since $X$ is assumed to be ordinary,
  this morphism turns out to be an isomorphism,  meaning  that $\mcE^\sss$ is ordinary.
  This completes the proof of the assertion.
    \end{proof}
%----------------------------------------------------------------------------------------

\LSP
%---------------------------[begin subsection]-------------
\subsection{Deformation of dormant indigenous $(G, H)$-bundles of higher level} \label{SS0102}

In the rest of this section, we suppose that  $G$ is affine and fix a positive integer  $\N$. 
Suppose that the base field $k$ is equipped with the trivial $(\N-1)$-PD structure.

Denote by $\mfA \mfr \mft_{/k}^{\N \text{-}\mr{PD}}$ the category  of  
$(\N-1)$-PD rings $(R, \mfa, \mfb, \gamma)$ such that $R$ is a local Artinian $k$-algebras with residue field $k$.
Also, denote by $\mfA \mfr \mft_{/k}^{\N \text{-}\mr{PD}, \circledast}$
the full subcategory of $\mfA \mfr \mft_{/k}^{\N \text{-}\mr{PD}}$
consisting of  $(\N-1)$-PD rings   of the form 
$(R, \mfm_R, \mfm_R, \gamma_R)$ for some PD structure $\gamma_R$ on $\mfm_R$. 
(In particular,  $\mfA \mfr \mft_{/k}^{\N \text{-}\mr{PD}, \circledast}$ does not depend on $\N$.)
If there is no fear of confusion, 
we shall often abuse notation by just writing $R$ for $(R, \mfa, \mfb, \gamma)$ or $(R, \mfm_R, \mfm_R, \gamma_R)$.
%we will  omit ``$\mfa$", ``$\mfb$", and ``$\gamma$"  in the notation.
%Also, we always consider each $R \in \mr{Ob}(\mfA \mfr \mft_{/k}^{\mr{PD}})$  as being equipped with the $(\N-1)$-PD structure $(\mfm_R, \mfm_R, \gamma_R)$.
The discussions conducted  so far will  work even when we replace all ``$\mfA \mfr \mft_{/k}$"'s  appearing there by either  
$\mfA \mfr \mft_{/k}^{\N \text{-}\mr{PD}}$ or $\mfA \mfr \mft_{/k}^{\N \text{-}\mr{PD}, \circledast}$.
To verify this, the following three observations  should be noted:
\begin{itemize}
\item[(i)]
  The $k$-algebra $k_\epsilon$ lies in $\mfA \mfr \mft_{/k}^{\N \text{-}\mr{PD}, \circledast}$ (hence, in $\mfA \mfr \mft_{/k}^{\N \text{-}\mr{PD}}$) because  its maximal  ideal $(\epsilon) \left(=\mfm_{k_\epsilon} \right)$   has evidently an $(\N-1)$-PD structure (and it is uniquely determined);
 \item[(ii)]
 The categories $\mfA \mfr \mft_{/k}^{\N \text{-}\mr{PD}}$, $\mfA \mfr \mft_{/k}^{\N \text{-}\mr{PD}, \circledast}$
admit fiber products.
To be precise, 
  for any morphisms $R' \migi R$, $R'' \migi R$ in  either $\mfA \mfr \mft_{/k}^{\N \text{-}\mr{PD}}$ or
 $\mfA \mfr \mft_{/k}^{\N \text{-}\mr{PD}, \circledast}$,
   % $\mfA \mfr \mft_{/k}^{\N\text{-}\mr{PD}}$, 
 the fiber product $R' \times_R R''$ belongs to the same category because it  has the  $(\N-1)$-PD ideal induced from 
 the $(\N-1)$-PD ideals of 
% the maximal ideal of $R' \times_R R''$ has, in a natural manner,  a PD structure induced from the PD structures 
 $R$, $R'$, and $R''$;
 \item[(iii)]
 If  $(R, \mfa, \mfb, \gamma)$ is an object 
 of  either $\mfA \mfr \mft_{/k}^{\N \text{-}\mr{PD}}$ or
 $\mfA \mfr \mft_{/k}^{\N \text{-}\mr{PD}, \circledast}$,
 then  $\gamma$ extends to any deformation of $X$ over $R$ (cf.  ~\cite[\S\,1.3, Proposition 1.3.8, (iv)]{Ber1}).
\end{itemize}
By  these facts,  it makes sense to speak of  various deformation functors on $\mfA \mfr \mft_{/k}^{\N\text{-}\mr{PD}}$ or $\mfA \mfr \mft_{/k}^{\N \text{-}\mr{PD}, \circledast}$.
In what follows, let $\Box$ denote either the absence or presence of ``$\circledast$". 

Now, let
 $\mcE_\N^\sss := (\mcE_H,  \varepsilon^\ST)$ be  a dormant indigenous $(G, H)$-bundle of level $\N$ on $X/k$.
 Denote by $\nabla$ the $k$-connection on $\mcE_G := \mcE_H \times^H G$ defining $\mcE_\N^\sss |^{\langle 1 \rangle}$ via (\ref{Efjj2}).
 In particular, $\mcE^\sss := (\mcE_H, \nabla)$ specifies a dormant indigenous $(G, H)$-bundle, and various  discussions conducted so far can be applied to this $\mcE^\sss$.
Just as in the case of indigenous $(G, H)$-bundles,
we have  the notion of a (dormant) deformation of $\mcE_\N^\sss$ (resp., $(X, \mcE_\N^\sss)$) over each $R \in \mr{Ob}(\mfA \mfr \mft_{/k}^{\N\text{-}\mr{PD}, \Box})$.
Also, an isomorphism  between such objects  can be defined.
For each $R \in \mr{Ob}(\mfA \mfr \mft^{\N\text{-}\mr{PD}, \Box}_{/k})$, we shall set
\begin{align} \label{gaoeos}
\mr{Def}_{\mcE_\N^\sss}^{^\mr{Zzz...}} (R) := 
\left(\begin{matrix} \text{the set of isomorphism classes of} \\ \text{dormant deformations of $\mcE_\N^\sss$ over $R$} \end{matrix} \right) \hspace{14mm}\\
\left(\text{resp.,} \  \mr{Def}_{(X, \mcE_\N^\sss)}^{^\mr{Zzz...}} (R) := 
\left(\begin{matrix} \text{the set of isomorphism classes of} \\ \text{dormant deformations of $(X, \mcE_\N^\sss)$ over $R$} \end{matrix} \right)\right).\notag
\end{align} 
The assignments $R \mapsto  \mr{Def}_{\mcE_\N^\sss}^{^\mr{Zzz...}} (R)$  (resp., $R \mapsto  \mr{Def}_{(X, \mcE_\N^\sss)}^{^\mr{Zzz...}} (R)$)
defines  a  functor  with g.d.t
\begin{align} \label{Ef009}
\mr{Def}_{\mcE_\N^\sss}^{^\mr{Zzz...}} \ 
\left(\text{resp.,} \ \mr{Def}_{(X, \mcE_\N^\sss)}^{^\mr{Zzz...}} \right) : \mfA \mfr \mft^{\N\text{-}\mr{PD}, \Box}_{/k} \migi \mfS \mfe \mft.
\end{align}

We shall  describe the tangent and obstruction spaces  of $\mr{Def}_{\mcE_\N^\sss}^{^\mr{Zzz...}}$ and $\mr{Def}_{(X, \mcE_\N^\sss)}^{^\mr{Zzz...}}$.
Denote by $\mcE^\nabla_G$ (cf. (\ref{LLKS3})) the $G$-bundle on $X^{(\N)}$
corresponding to $\varepsilon^\ST$ via
(\ref{E3erf002}).
 In particular, there exists  a natural isomorphism of $G$-bundles 
 $F^{(N)*}_{X/k}(\mcE^\nabla_G) \isom \mcE_G$, which induces an isomorphism of $\mcO_X$-Lie algebras 
 $F^{(N)*}_{X/k}(\mfg_{\mcE_G^\nabla}) \isom \mfg_{\mcE_G}$.
 Both $\mfh_{\mcE_H}$ (resp., $\widetilde{\mcT}_{\mcE_H}$) and
 $\mfg_{\mcE_G^\nabla}$ may be regarded as subsheaves of $\mfg_{\mcE_G}$ (resp., $\widetilde{\mcT}_{\mcE_G}$).

\SSP
%-------------------------------------------------------------
\ble \label{Er55k}
The composite  $\widetilde{\nabla}^\mr{ad}_G \circ \nabla : \mcT_X \migi \mfg_{\mcE_G}$ coincides with the zero map.
\ele
%-------------------------------------------------------------
\begin{proof}
Let us take local sections $\partial$, $\partial' \in \mcT_{X}$.
Then, we have
\begin{align}
\langle \partial,  (\widetilde{\nabla}^\mr{ad}_G \circ \nabla) (\partial') \rangle = & \ 
\widetilde{\nabla}^\mr{ad}_{\partial} (\nabla(\partial'))  \\
= & \   [\nabla (\partial), \nabla (\partial')] - \nabla ([\partial, d_{\mcE_G}(\nabla (\partial'))])  \notag \\
= & \   [\nabla (\partial), \nabla (\partial')] - \nabla ([\partial, \partial'])  \notag  \\
= & \  0, \notag
\end{align}
where  $\langle -, - \rangle$ denotes the natural pairing $\mcT_X \times (\Omega_X \otimes \mfg_{\mcE_G}) \migi \mfg_{\mcE_G}$ and  the last equality follows from the flatness of $\nabla$.
This completes the proof of this lemma.
\end{proof}
%-------------------------------------------------------------
\SSP

For each $i \geq 0$, we shall set $\mbH^i_{H, B, \N}$ (resp., $\widetilde{\mbH}^i_{H, B, \N}$) to be the $k$-vector space making  the following square diagram cartesian:
\begin{align}\label{Diag04}
\vcenter{\xymatrix@C=26pt@R=36pt{
\mbH^i_{H, B, \N} \ar[r] \ar[d] &  \mbH^i_{H, B}  \ar[d]
\\
H^i (X^{(\N)}, \mfg_{\mcE_G^\nabla})  \ar[r]& \mbH^i_G
}}
\ \left(\text{resp.,} \ \vcenter{\xymatrix@C=26pt@R=36pt{
\widetilde{\mbH}^i_{H, B, \N} \ar[d] \ar[r]&  \widetilde{\mbH}^i_{H, B}    \ar[d]
\\
 H^i (X^{(\N)}, \mfg_{\mcE_G^\nabla}) \oplus H^i  \ar[r]& \widetilde{\mbH}^i_G
}} \right),
\end{align}
where 
the right-hand vertical arrow arises from the inclusion $\mcK^\bullet [\nabla_{H, B}^\mr{ad}] \migi \mcK^\bullet [\nabla_G^\mr{ad}]$ (resp., $\mcK^\bullet [\widetilde{\nabla}_{H, B}^\mr{ad}] \migi \mcK^\bullet [\widetilde{\nabla}_G^\mr{ad}]$) and 
the lower  horizontal arrow arises from the inclusion $(F^{(\N)}_{X/k})^{-1} (\mfg_{\mcE_G^\nabla}) \migi \mfg_{\mcE_G}$ (resp., both $\nabla$ and the  inclusion $(F^{(\N)}_{X/k})^{-1} (\mfg_{\mcE_G^\nabla}) \migi \widetilde{\mcT}_{\mcE_G}$).
Then, the following assertion holds.

 %------------------------------------------------------------
 \bpr \label{P051}
 Recall the assumption that  $G$ is affine and $\mcE_\N^\sss$ is   dormant.
 \begin{itemize}
 \item[(i)]
 The collection $(\mr{Def}^{^\mr{Zzz...}}_{\mcE_\N^\sss}, \mbH_{H, B, \N}^1, \mbH_{H, B, \N}^2)$ (resp.,  $(\mr{Def}^{^\mr{Zzz...}}_{(X, \mcE_\N^\sss)}, \widetilde{\mbH}_{H, B, \N}^1, \widetilde{\mbH}_{H, B, \N}^2)$)
 form a deformation triple.
  Moreover, the group of infinitesimal automorphism of $\mcE_\N^\sss$ (resp.,  $(X, \mcE_\N^\sss)$)
 is isomorphic to $\mbH_{H, B, \N}^0$ (resp., $\widetilde{\mbH}_{H, B, \N}^0$).
% The tangent space of $\mr{Def}^{^\mr{Zzz...}}_{\mcE_\N^\sss}$
 %is isomorphic to $\mbH_{H, B, \N}^1$, and this functor has an obstruction theory
% by putting $\mbH_{H, B, \N}^2$ as the obstruction space.
% Moreover, the group of infinitesimal automorphism of $\mcE_\N^\sss$
% is isomorphic to $\mbH_{H, B, \N}^0$.
% \item[(ii)]
  %The tangent space of $\mr{Def}^{^\mr{Zzz...}}_{(X, \mcE_\N^\sss)}$
 %is isomorphic to $\widetilde{\mbH}_{H, B, \N}^1$, and this functor has an obstruction theory
 %by putting $\widetilde{\mbH}_{H, B, \N}^2$ as the obstruction space.
 %Moreover, the group of infinitesimal automorphism of $(X, \mcE_\N^\sss)$
 %is isomorphic to $\widetilde{\mbH}_{H, B, \N}^0$.
\item[(ii)]
Suppose that $X$ is proper over $k$.
Then, 
there exist pro-representable hulls for 
$\mr{Def}^{^\mr{Zzz...}}_{\mcE_\N^\sss}$ and $\mr{Def}^{^\mr{Zzz...}}_{(X, \mcE_\N^\sss)}$.
 \item[(iii)]
Regarding the pro-representability of $\mr{Def}^{^\mr{Zzz...}}_{\mcE_\N^\sss}$
 (resp., $\mr{Def}^{^\mr{Zzz...}}_{(X, \mcE_\N^\sss)}$),
the statement of  Corollary \ref{Cf02k2} with $\mbH^i_H$'s (resp., $\widetilde{\mbH}_H^i$'s) replaced by $\mbH_{H, B, \N}^i$ (resp., $\widetilde{\mbH}_{H, B, \N}^i$)
 holds for $\mr{Def}^{^\mr{Zzz...}}_{\mcE_\N^\sss}$
 (resp., $\mr{Def}^{^\mr{Zzz...}}_{(X, \mcE_\N^\sss)}$).
  \end{itemize}
 \epr
 %------------------------------------------------------------
 \begin{proof}
 We only consider the resp'd assertion of (i) because the proofs of  the others follow from this assertion (and its proof).
  For each $R \in \mr{Ob}(\mfA \mfr \mft_{/k}^{\N\text{-}\mr{PD}, \Box})$,  denote by $\mr{Def}_{(X, \mcE_G^\nabla)}(R)$  the set of isomorphism classes of pairs $(X_R, \mcE_{G, R}^\nabla)$ consisting of a deformation $X_R$ of $X$ over $R$ and a deformation $\mcE_{R, G}^\nabla$ of the $G$-bundle $\mcE_G^\nabla$ over $X_R^{(1)}$.
 The  assignment $R \mapsto \mr{Def}_{(X, \mcE_G^\nabla)}(R)$ is verified to define a functor with gdt $\mr{Def}_{(X, \mcE_G^\nabla)} : \mfA \mfr \mft^{\N\text{-}\mr{PD}, \Box}_{/k} \migi \mfS \mfe \mft$.
 If $F_{X/k}^{(\N)\sharp}$ denotes   the morphism of functors   given by $X_R \mapsto X_R^{(\N)}$, then  
 the following  square  diagram of  functors is cartesian:
\begin{align}\label{Diaggh4}
\vcenter{\xymatrix@C=86pt@R=36pt{
\mr{Def}_{(X, \mcE_G^\nabla) } \ar[r]^-{(X_R, \mcE_{R, G}^\nabla) \mapsto X_R} \ar[d]_-{(X_R, \mcE_{R, G}^\nabla) \mapsto (X_R^{(\N)}, \mcE_{R, G}^\nabla)}&  \mr{Def}_{X}  \ar[d]^-{F_{X/k}^{(\N)\sharp}}
\\
\mr{Def}_{(X^{(\N)}, \mcE_G^\nabla)} \ar[r]_-{(X_R^{(\N)}, \mcE_{R, G}^\nabla) \mapsto X_R^{(\N)}}&  \mr{Def}_{X^{(\N)}}
}}
\end{align}
(cf. (\ref{RwwqaW}) for the definition of $\mr{Def}_{(X^{(\N)}, \mcE_G^\nabla)}$).
For each $i =0,1,2$, 
let us write  $\mbH_{G, \nabla}^i$  for  the limit of the diagram
 \begin{align}\label{Diagf04}
\vcenter{\xymatrix@C=46pt@R=36pt{
&  H^i  \ar[d]^-{}
\\
H^i (X^{(\N)}, \widetilde{\mcT}_{\mcE_G^\nabla}) \ar[r]_{H^i (d_{\mcE_G^\nabla})}& H^i (X^{(\N)}, \mcT_{X^{(\N)}}),
}}
\end{align}
where the right-hand vertical arrow  is the morphism arising  from $F_{X/k}^{(\N)\sharp}$.
But, 
  the right-hand vertical arrow
  becomes the zero map  by  the definition  of  Frobenius twist, so
 %  induces the zero map between their tangent spaces  $H^i \migi H^i (X^{(\N)}, \mcT_{X^{(\N)}})$.
this diagram induces  a natural isomorphism
\begin{align} \label{Geiji4}
\mbH_{G, \nabla}^i \isom H^i (X^{(\N)}, \mfg_{\mcE_G^\nabla}) \oplus H^i \left(= \mr{Ker} (H^i (d_{\mcE_G})) \oplus H^i\right).
\end{align}
It follows from the definitions of $\mbH^i_{G, \nabla}$  and (\ref{Diaggh4}) that 
the tangent space of 
$\mr{Def}_{(X, \mcE_G^\nabla)}$ is isomorphic to $\mbH_{G, \nabla}^1$, and this functor has an obstruction theory 
by putting  $\mbH_{G, \nabla}^2$ as the obstruction space.
Moreover, the space of infinitesimal  automorphisms of $(X, \mcE_G^\nabla)$ is isomorphic to $\mbH_{G, \nabla}^0$.
 On the other hand,   the equivalence of categories (\ref{E3erf002}) gives rise to  the following cartesian diagram:
  \begin{align}\label{Diugh4}
\vcenter{\xymatrix@C=136pt@R=36pt{
\mr{Def}^{^\mr{Zzz...}}_{(X, \mcE_\N^\sss)} \ar[r]^-{(X_R, \mcE_{R, \N}^\sss) \mapsto (X_R, \mcE_{R, \N}^\sss |^{\langle 1 \rangle})} \ar[d]_-{(X_R, (\mcE_{R, H}, \varepsilon^\ST_R)) \mapsto (X_R, \mcE^\nabla_{R, G})} &\mr{Def}_{(X, \mcE^\sss)} \ar[d]^-{(X_R, (\mcE_{R, H}, \nabla_R)) \mapsto (X_R, \mcE_{R, G}, \nabla)}
\\
\mr{Def}_{(X, \mcE_G^\nabla)}  \ar[r]_-{(X_R, \mcG) \mapsto  (X_R, F^{(\N)*}_{X_R/R}(\mcG), \nabla_{\mcG}^{\mr{can}(\N)})} &\mr{Def}_{(X, \mcE_G, \nabla)}
 }}
\end{align}
(cf. (\ref{RwwqaW}) for the definition of $\mr{Def}_{(X, \mcE_G, \nabla)}$).
Under the identification $\mbH^i_{G, \nabla} = H^i (X^{(\N)}, \mfg_{\mcE^\nabla_G}) \oplus H^i$ given by  (\ref{Geiji4}),
the differential of the lower horizontal arrow coincides with   the morphism 
$H^i (X^{(\N)}, \mfg_{\mcE^\nabla_G}) \oplus H^i \migi \mbH_G^i$ arising from 
 $(F_{X/k}^{(\N)})^{-1}(\mfg_{\mcE^\nabla_G}) \migiincl  \widetilde{\mcT}_{\mcE_G}$ and
$\nabla : \mcT_{X} \migi \widetilde{\mcT}_{\mcE_G}$ (cf. Lemma \ref{Er55k}).
Thus, we can conclude  the required assertion  by combining  (\ref{Diag04}) and  (\ref{Diugh4}).
 \end{proof}
%------------------------------------------------------------

\SSP

\begin{rema} \label{EREo2}
Similar to  the fact mentioned in Remark \ref{FFmK0},
the tangent and obstruction spaces of 
 the   deformation functor  of  higher level
 % classifying dormant deformations of $(\mcE_G, \varepsilon^\ST)$ 
 can be described  in terms of cohomology groups. 
With 
 the above notation,
we shall denote by
\begin{align} \label{Erqq21}
\mr{Def}_{(\mcE_G, \varepsilon^\ST)}^{^\mr{Zzz...}}
\end{align}
the functor classifying dormant deformations of the  $(\N-1)$-PD stratified  $G$-bundle  $(\mcE_G, \varepsilon^\ST)$.
Then, 
 the proof of the above proposition shows
 that 
 $(\mr{Def}_{(\mcE_G, \varepsilon^\ST)}^{^\mr{Zzz...}}, H^1 (X^{(\N)}, \mfg_{\mcE_G^\nabla}), H^2 (X^{(\N)}, \mfg_{\mcE_G^\nabla}))$ forms a deformation triple, and 
 %the tangent space of 
%$\mr{Def}_{(\mcE_G, \varepsilon^\ST)}^{^\mr{Zzz...}}$ 
%is isomorphic to $H^1 (X^{(\N)}, \mfg_{\mcE_G^\nabla})$, and this functor has an obstruction theory by putting  $H^2 (X^{(\N)}, \mfg_{\mcE_G^\nabla})$
%as the obstruction space.
 the group of infinitesimal automorphisms of $(\mcE_G, \varepsilon^\ST)$ is  isomorphic to $H^0 (X^{(\N)}, \mfg_{\mcE_G^\nabla})$.

If $\N =1$, then $\mr{Def}_{(\mcE_G, \varepsilon^\ST)}^{^\mr{Zzz...}}$ is naturally isomorphic to the  subfunctor $\mr{Def}_{(\mcE_G, \nabla)}^{^\mr{Zzz...}}$ of $\mr{Def}_{(\mcE_G, \nabla)}$  classifying deformations with vanishing $p$-curvature.
Since $F_{X/k}^{-1}(\mfg_{\mcE^\nabla_G}) =  \mr{Ker}(\nabla^\mr{ad}_G)$,
the tangent and obstruction spaces of this functor are isomorphic to 
$H^i (X, \mr{Ker}(\nabla_G^\mr{ad}))$ ($\cong \mbH_{G, B}^i$ by (\ref{Ty67})) for $i=1$ and $2$ respectively.
This fact will be used in the proof of Theorem \ref{Wroauom}.
\end{rema}

 \LSP
%---------------------------[begin subsection]-------------
\subsection{The Ehresmann-Weil-Thurston principle in positive characteristic} \label{SS048}

%In this final  subsection, we formulate and prove  the {\it Ehresmann-Weil-Thurston principle} for $F^\N$-$(G, P)$-structures.

 Let us keep the assumption in the previous subsection, and let
 $\mcS^\hh$  be an $F^\N$-$(G, P)$-structures  on $X/k$.
Denote by 
\begin{align} \label{W23}
\mr{Def}_{(X, \mcS^\hh)} : \mfA \mfr \mft^{\N\text{-}\mr{PD}, \circledast}_{/k} \migi \mfS \mfe \mft
\end{align}
the deformation functor of  the pair $(X, \mcS^\hh)$.
That is to say, for each $R \in \mr{Ob}(\mfA \mfr \mft^{\N\text{-}\mr{PD}, \circledast}_{/k})$,
$\mr{Def}_{(X, \mcS^\hh)} (R)$ is defined as  the set of  pairs $(X_R, \mcS^\hh_R)$ consisting of (an isomorphism class of) a deformation  $X_R$ of $X$ over $R$ and an $F^\N$-$(G, P)$-structure $\mcS^\hh_R$ on $X_R/R$ inducing $\mcS^\hh$ via reduction modulo $\mfm_R$.
The pull-back $\mr{Def}_{(X, \mcS^\hh)} (R') \migi \mr{Def}_{(X, \mcS^\hh)} (R)$ induced by every  morphism $R' \migi R$ can be formulated by applying the discussion in \S\,\ref{SS066}, and thus  we obtain the functor $\mr{Def}_{(X, \mcS^\hh)}$.

On the other hand, 
denote by $\mcE_\N^\sss := (\mcE_H,  \varepsilon^\ST)$
 the dormant indigenous $(G, H)$-bundle of level $\N$ corresponding to $\mcS^\hh$, and moreover, denote by 
$\mcE^\diamondsuit$ the dormant $(\N-1)$-crystal  of $G$-bundles associated with $(\mcE_G, \varepsilon^\ST)$ (cf. Proposition \ref{P019}, Definition \ref{rtyu5}, (ii)).
Also. we shall denote by 
\begin{align} \label{W24}
\mr{Def}_{\mcE^\diamondsuit}^{^\mr{Zzz...}}  : \mfA \mfr \mft_{/k}^{\N\text{-}\mr{PD}, \circledast} \migi \mfS \mfe \mft
\end{align}
  the  functor 
 classifying, for each  $R \in \mr{Ob}(\mfA \mfr \mft_{/k}^{\N\text{-}\mr{PD}, \circledast})$,  dormant $(\N-1)$-crystals on $\mr{Cris}^{(\N-1)}(X/R)$ inducing $\mcE^\dd$ via reduction modulo $\mfm_R$.
Then, the bijection 
$\zeta_\N^{\hh \Rightarrow \sss}$ obtained in Theorem \ref{TheoremA} (resp., the equivalence of categories (\ref{E3GH2}))
 induces an isomorphism  (resp., a morphism) of  functors
\begin{align}
\mu_\N^{\heartsuit \Rightarrow \spadesuit} : \mr{Def}_{(X, \mcS^\hh)} \isom \mr{Def}_{(X, \mcE_\N^\sss)}^{^\mr{Zzz...}} \ \left(\text{resp.,} \
\mu_\N^{\spadesuit \Rightarrow \diamondsuit} : \mr{Def}_{(X, \mcE^\sss_\N)}^{^\mr{Zzz...}} \migi  \mr{Def}^{^\mr{Zzz...}}_{\mcE^\dd}   \right).
\end{align}
In particular,
$\mr{Def}_{(X, \mcS^\hh)}$ is a functor with gdt and has  the same  tangent and obstruction spaces as 
$\mr{Def}_{(X, \mcE_\N^\sss)}^{^\mr{Zzz...}}$ (cf. Proposition \ref{P051}, (i)).
The compose of   $\mu_\N^{\heartsuit \Rightarrow \spadesuit}$ and   $\mu_\N^{\spadesuit \Rightarrow \diamondsuit}$
defines a morphism of functors
\begin{align}
\mu_\N^{\heartsuit \Rightarrow \diamondsuit} : \mr{Def}_{(X, \mcS^\hh)} \migi \mr{Def}_{\mcE^\diamondsuit}^{^\mr{Zzz...}}.
\end{align}
%and  the formation of this morphism commutes with truncation to lower levels.
In other words,  $\mu_\N^{\heartsuit \Rightarrow \diamondsuit} $  is the morphism given  by taking  the monodromy crystals of $F^\N$-$(G, P)$-structures
  (cf. Definition \ref{Efopwd}).
As explained  in the Introduction, the following theorem, i.e., Theorem \ref{TheoremD}, may be regarded as a positive characteristic version of the Ehresmann-Weil-Thurston principle.

\SSP
%---------------------------
\bt  \label{Wroauom}
The morphism $\mu_\N^{\hh \Rightarrow \dd}$ is a semi-isomorphism.
If, moreover, 
$\mr{Def}_{\mcE^\diamondsuit}^{^\mr{Zzz...}}$ is rigid, then
$\mu_\N^{\hh \Rightarrow \dd}$ is an isomorphism.
\et
%---------------------------
\begin{proof}
Since $\mu_\N^{\hh \Rightarrow \sss}$ is an isomorphism,
it suffices to prove the same assertion with $\mu_\N^{\hh \Rightarrow \dd}$ replaced by $\mu_\N^{\sss \Rightarrow \dd}$.
Moreover, we only consider the former assertion because the latter  assertion follows directly from the former assertion and    Proposition \ref{pyp90s4}, (i).

First, we shall prove  
the case of $\N =1$.
Denote by 
$\nabla$ the $k$-connection on $\mcE_G := \mcE_H \times^H G$ corresponding to $\varepsilon^\ST$ (cf. (\ref{Efjj32})),   which determines   a dormant  indigenous $(G, H)$-bundle $\mcE^\sss := (\mcE_H, \nabla)$.
By the equivalence of categories (\ref{Efjj32}), $\mr{Def}_{\mcE^\diamondsuit}^{^\mr{Zzz...}}$ may be identified with
the deformation functor    $\mr{Def}_{(\mcE_G, \nabla)}^{^\mr{Zzz...}}$
introduced in Remark \ref{EREo2}.
 Also, there exists a canonical  identification $\mr{Def}_{(X, \mcE_1^\sss)}^{^\mr{Zzz...}} = \mr{Def}_{(X, \mcE^\sss)}^{^\mr{Zzz...}}$ arising from   (\ref{QQWv4}).
Under these identifications,
 $\mu_1^{\spadesuit \Rightarrow \diamondsuit}$ defines  a morphism of functors  $\mr{Def}_{(X, \mcE^\sss)}^{^\mr{Zzz...}} \migi \mr{Def}_{(\mcE_G, \nabla)}^{^\mr{Zzz...}}$.
By combining  Lemmas \ref{LLODMM2}, (ii), and \ref{GGHD1} described below,
we see that if  
the differential $d \mu_1^{\spadesuit \Rightarrow \diamondsuit}$ 
is regarded as 
 a morphism $\widetilde{\mbH}^1_{H, B} \migi \mbH^1_{G, B}$ according to the results of  Proposition \ref{P050}, (i), and Remark \ref{EREo2},
then it  
coincides with $\mbH^1 (\eta_{H, B}^\bullet)$ (cf. (\ref{Ty67kk})).
Moreover, it follows from  various definitions involved and  the discussion in  the proof of Lemma \ref{LLODMM2}, (i),  that 
 the triple
\begin{align}
(\mu_1^{\sss \Rightarrow \dd}, \mbH^1 (\eta_{H, B}^\bullet), \mbH^2 (\eta_{H, B}^\bullet)) : (\mr{Def}_{(X, \mcE^\sss)}^{^\mr{Zzz...}}, \widetilde{\mbH}^1_{H, B}, \widetilde{\mbH}^2_{H, B}) \migi (\mr{Def}_{(\mcE_G, \nabla)}^{^\mr{Zzz...}}, \mbH^1_{G, B}, \mbH_{G, B}^2)
\end{align}
forms an isomorphism of deformation triples.
Since $\mbH^i (\eta^\bullet_{H, B})$ ($i \geq 0$) are isomorphisms,    
 the morphism $\mu_1^{\sss \Rightarrow \dd}$
 turns out to be a semi-isomorphism (cf. Proposition \ref{E2Rgju2}).
 
Next, we shall consider the case of an arbitrary  positive integer $\N$.
Let us verify the formal smoothness of $\mu_\N^{\sss \Rightarrow \dd}$.
Take a small extension $e := (0 \migi M \migi R \xrightarrow{\sigma} R_0 \migi 0)$ in $\mfA \mfr \mft_{/k}^{\N\text{-}\mr{PD}, \circledast}$.
In the following, we shall prove the surjectivity of  the map
\begin{align} \label{KKo9}
\mr{Def}_{(X, \mcE_\N^\sss)}^{^\mr{Zzz...}} (R) \migi \mr{Def}_{(X, \mcE_\N^\sss)}^{^\mr{Zzz...}} (R_0) \times_{\mr{Def}_{\mcE^\dd}^{^\mr{Zzz...}} (R_0)}\mr{Def}_{\mcE^\dd}^{^\mr{Zzz...}} (R)
\end{align}
arising from   $\mr{Def}_{(X, \mcE_\N^\sss)}^{^\mr{Zzz...}} (\sigma)$ and $\mu_\N^{\sss \Rightarrow \dd}(R_0)$.
Suppose that we are given   
 a dormant deformation  $(X_{R_0}, \mcE^\sss_{R_0, \N})$  (where $\mcE^\sss_{R_0, \N} := (\mcE_{R_0, H}, \varepsilon^\ST_{R_0})$ and we write $\mcE_{R_0, G}:= \mcE_{R_0, H}\times^H G$) of $(X, \mcE^\sss_\N)$ over $R_0$
 and an $(\N-1)$-crystal of $G$-bundles $\mcE_{R}^\dd$ on $X/R$ whose reduction modulo $M$
   is isomorphic to
   the $(\N-1)$-crystal on $X/R_0$ corresponding to 
    $(\mcE_{R_0, G}, \varepsilon^\ST_{R_0})$ (cf. Proposition \ref{P019}).
By the  formal smoothness of $\mu_1^{\sss\Rightarrow\dd}$ proved above,
there exists a deformation $(X_{R}, \mcE_{R, 1}^{\sss})$ (where $\mcE_{R, 1}^{\sss} := (\mcE_{R, H}, \varepsilon_{R, 1}^\ST)$) of $(X, \mcE^\sss_\N |^{\langle 1 \rangle})$ over $R$ which induces  $\mcE_R^\dd |^{(0)}$ via $\mu_1^{\sss \Rightarrow \dd}$ and induces $(X_{R_0}, \mcE_{R_0, \N}^\sss |^{\langle 1 \rangle})$ via reduction modulo $M$.
Then,  $\mcE_R^\dd$ determines   a unique (up to isomorphism) $(\N-1)$-PD stratification $\varepsilon_{R, \N}^\ST$ on $\mcE_{R, G} := \mcE_{R, H} \times^H G$ 
extending $\varepsilon_{R, 1}^\ST$ (cf. Proposition \ref{P019}).
The resulting collection   $(X_R, \mcE_{R, H}, \varepsilon_{R, \N}^\ST)$ forms a deformation of $(X, \mcE^\sss_\N)$ over $R$ inducing $\mcE_R^\dd$ via $\mu_\N^{\sss \Rightarrow \dd}$ and  $(X_{R_0}, \mcE_{R_0, \N}^\sss)$ via reduction modulo $M$.
This implies that (\ref{KKo9}) is surjective.
 
Also, the bijectivity  of $d \mu_\N^{\sss \Rightarrow \dd}$ follows from an entirely  similar argument together with the bijectivity of $d \mu_1^{\sss \Rightarrow \dd}$.
Thus, we conclude that $\mu_\N^{\sss \Rightarrow \dd}$ is a semi-isomorphism, and this completes the proof of the theorem.
 \end{proof}
%---------------------------
\SSP

The following corollary follows  from  the above theorem and the comment in Remark \ref{EREo2}.

\SSP
%---------------------------
\bco \label{FGHFH}
Let us keep the above notation.
Suppose further  that $X$ is proper over $k$ and  $H^0 (X, \mfg_{\mcE_G^\nabla}) =0$.
Then, the functor $\mr{Def}_{(X, \mcS^\hh)}$ is pro-representable.
If, moreover, $H^2 (X, \mfg_{\mcE_G^\nabla}) =0$, then
$\mr{Def}_{(X, \mcS^\hh)}$ may be pro-represented  by $k [\! [ x_1, \cdots, x_d ]\! ]$, where $d := \mr{dim}_k (H^1 (X, \mfg_{\mcE_G^\nabla}))$.
\eco
%---------------------------
\SSP

Now, let us describe  two lemmas, that were 
 used to conclude  Theorem \ref{Wroauom}.
Let $(\mcE_G, \nabla)$ be as in the proof of that theorem, but we write $\mcE := \mcE_G$ for simplicity. 
For each deformation 
$(X_R, \mcE_{R}, \nabla_R)$ 
of $(X, \mcE, \nabla)$ over  
$R \in \mr{Ob}(\mfA \mfr \mft_{/k}^{\N\text{-}\mr{PD}, \circledast})$,
we shall denote by 
\begin{align}
\Xi (X_R, \mcE_{R}, \nabla_R)
\end{align}
the element of $\mr{Def}_{(X, \mcE, \nabla)}(R)$ (cf. (\ref{RwwqaW})) represented by this deformation.
Also, 
for each quasi-nilpotent flat $G$-bundle $(\mcE_R, \nabla_R)$ on a deformation of $X$ over $R \in \mr{Ob}(\mfA \mfr \mft_{/k}^{\N\text{-}\mr{PD}, \circledast})$, we shall denote by  $(\mcE_R, \nabla_R)^\dd$
 (cf. (\ref{EWQQ21}))  the $0$-crystal of $G$-bundles on $X/R$ corresponding to $(\mcE_R, \nabla_R)$ via (\ref{E3GH2}) and (\ref{DFjj2}).

\SSP
%---------------------------------------------------------------------
\ble \label{LLODMM2}
 Let us fix 
a small extension
$e := (0 \migi M \migi R \migi R_0 \migi 0)$ in $\mfA \mfr \mft_{/k}^{\N\text{-}\mr{PD}, \circledast}$ and 
 a flat  deformation $(X_R, \mcE_{R}, \nabla_R)$ of $(X, \mcE, \nabla)$ over $R$. (As mentioned in Remark \ref{koRirt}, the connection $\nabla_R$ is automatically quasi-nilpotent because  of the $p$-flatness of $\nabla$.)
\begin{itemize}
\item[(i)]
Let $v$ be an element of $M \otimes H^1$ and denote by
$X'_R$ the deformation of $X$ over $R$ defined as  the translation by $v$ of the element in  $\mr{Def}_X(R)$ represented by $X_R$.
Also, denote by $(\mcE'_{R}, \nabla'_R)$ a  flat $G$-bundle on $X'_R/R$ with
$(\mcE_{R}, \nabla_R)^\dd \cong (\mcE'_{R}, \nabla'_R)^\dd$ (which always   exists  and  is uniquely determined up to isomorphism according  Proposition \ref{eoqo84939}).
Then, the following equality in $\mr{Def}_{(X, \mcE, \nabla)} (R)$ holds:
\begin{align} \label{EBN23}
\Xi (X'_R, \mcE'_R, \nabla'_R) = 
\Xi (X_R, \mcE_R, \nabla_R)   + \mbH^1 (\nabla) (v),
\end{align}
where $\mbH^1 (\nabla)$ denotes the morphism $H^1 \migi \widetilde{\mbH}^1_G$ induced by $\nabla$ (cf. Lemma \ref{Er55k}).
\item[(ii)]
Let $u$ be an element of $M \otimes \widetilde{\mbH}_{G, Z}^1 \left(\subseteq M \otimes \widetilde{\mbH}_{G}^1 \right)$. 
Denote by   $(\mcE''_{R}, \nabla''_R)$ a unique (up to isomorphism) quasi-nilpotent flat $G$-bundle on $X_R/R$ such that
$(\mcE''_{R}, \nabla''_R)^\dd$
is isomorphic to the $0$-crystal associated with the quasi-nilpotent  flat $G$-bundle representing
the element $\Xi (X_R, \mcE_{R}, \nabla_R) + u$ of $\mr{Def}_{(X, \mcE, \nabla)} (R)$.
Then, the following equality holds:
\begin{align}
\Xi (X_R, \mcE''_{R}, \nabla''_R) = \Xi (X_R, \mcE_R, \nabla_R) + \mbH^1 (\eta^\bullet_{G, Z}) (u)
\end{align} 
(cf. (\ref{Fer4}) for the definition of $\mbH^1 (\eta^\bullet_{G, Z})$).
\end{itemize}
\ele
%---------------------------------------------------------------------
\begin{proof}
First, we shall consider assertion (i).
We can find an open  covering $\mcU := \{ U_\alpha \}_{\alpha \in I}$ of $X$ and $v$ may be represented by a \v{C}ech $1$-cocycle $\{ \partial_{\alpha \beta} \}_{(\alpha, \beta) \in I_2}$ in $\check{C}^1 (\mcU, M \otimes \mcT_X)$ (cf. the proof of Proposition \ref{P022} for the related notations).
Write  $U_{R, \alpha} := X_R |_{U_\alpha}$,  $\mcE_{R, \alpha} := \mcE_R |_{U_\alpha}$ (for each $\alpha \in I$) and
 $U_{R, \alpha \beta} := X_R |_{U_{\alpha \beta}}$,  $\mcE_{R, \alpha \beta} := \mcE_R |_{U_{\alpha \beta}}$ (for each $(\alpha, \beta) \in I_2$).
 Given each $(\alpha, \beta) \in I_2$,
 we shall write $\partial_{\alpha \beta}^\sharp$  (resp., $\widetilde{\partial}^\sharp_{\alpha \beta}$) for the $R$-automorphism
 of $U_{R, \alpha \beta}$ (resp., $\mcE_{R, \alpha \beta}$) determined by the automorphism $\mr{id}_{\mcO_{U_{R, \alpha \beta}}} + \partial_{\alpha \beta}$
 (resp.,  $\mr{id}_{\mcO_{\mcE_{R, \alpha \beta}}} + \nabla (\partial_{\alpha \beta})$)
 of  $\mcO_{U_{R, \alpha \beta}}$  (resp., $\mcO_{\mcE_{R, \alpha \beta}}$).
Let  us regard  $\partial_{\alpha \beta}^\sharp$  as a morphism of $0$-PD thickenings 
in $\mr{Cris}^{(0)}(\overline{X}/S)$.
Then, according to the latter half of the discussion in Remark \ref{EFPPWk} and the discussion at the beginning of \S\,\ref{SS0125},
 the induced isomorphism $\rho_{\partial^\sharp_{\alpha \beta}} : (\partial^\sharp_{\alpha \beta})^*(\mcE_{R, \alpha \beta}) \isom \mcE_{R, \alpha \beta}$ constituting  $(\mcE_R, \nabla_R)^\dd$ is given by $\widetilde{\partial}^\sharp_{\alpha \beta}$.
 It follows that 
 $(\mcE'_R, \nabla'_R)$ may be obtained by gluing together $(\mcE_{R, \alpha}, \nabla_R |_{\mcE_{R, \alpha}})$'s by means of $\widetilde{\partial}^\sharp_{\alpha \beta}$'s.
This implies the required equality (\ref{EBN23}) because the collection   $(\{ \nabla (\partial_{\alpha \beta}) \}_{(\alpha, \beta) \in I_2}, \{ \delta_\alpha \}_{\alpha \in I})$ (cf. (\ref{W105})) with $\delta_\alpha = 0$ for every $\alpha$ represents $\mbH^1 (\nabla) (v) \in \widetilde{\mbH}^1_G$.
We finish the proof of assertion (i).

Next, let us consider assertion (ii).
Denote by $\overline{u}$ the element of $M \otimes H^1$ defined as the image of $u$ via $\delta_G^1 : \left(\widetilde{\mbH}_{G, Z}^1 \subseteq \right) \widetilde{\mbH}_{G}^1 \migi H^1$ (cf. (\ref{QW7011})).
Also, denote by $(X'_R, \mcE'_{R}, \nabla'_R)$ the deformation of $(X, \mcE, \nabla)$ over $R$ represented by  $\Xi (X_R, \mcE_{R}, \nabla_R) + u$.
Then, the assertion can be proved   by applying assertion (i) to the case where $(X_R, \mcE_{R}, \nabla_R, v)$ is replaced by $(X'_R, \mcE'_{R}, \nabla'_R, - \overline{u})$.
This completes the proof of this lemma.
\end{proof}
%---------------------------------------------------------------------

\SSP

%---------------------------------------------------------------------
\ble \label{GGHD1}
Let $v$ be an element of $\mbH^1_G$.
Also, for each $i=1,2$,
let $(X_{R, i}, \mcE_{R, i}, \nabla_{R, i})$
be a deformation of $(X, \mcE, \nabla)$ over $R$, and denote by
$(\mcE^v_{R, i}, \nabla^v_{R, i})$ a unique (up to isomorphism)
flat $G$-bundle on $X_{R, i}/R$
 satisfying the equality
$\Xi (X_{R, i}, \mcE^v_{R, i}, \nabla^v_{R, i}) = \Xi (X_{R, i}, \mcE_{R, i}, \nabla_{R, i}) +v$.
Then, $(\mcE_{R, 1}, \nabla_{R, 1})^\dd \cong (\mcE_{R, 2}, \nabla_{R, 2})^\dd$
implies $(\mcE^v_{R, 1}, \nabla^v_{R, 1})^\dd \cong (\mcE^v_{R, 2}, \nabla^v_{R, 2})^\dd$.
(In the proof of Theorem \ref{Wroauom}, we apply this lemma in the case where one of the deformations $X_{R,1}$, $X_{R, 2}$  is trivial.)
\ele
%---------------------------------------------------------------------
\begin{proof}
The assertion follows from the construction of the assignment
$(\mcE, \nabla) \mapsto (\mcE', \nabla')$ asserted in Proposition \ref{eoqo84939}
(cf.  Remark \ref{GGGh6}).
\end{proof}
%---------------------------------------------------------------------
%\SSP

%%%%%%%%%%%%%%%%%%%%%%%%%%%%%%%---[begin section]---%%%%%%%%%%%%%%
\vspace{10mm}
\section{Appendix: Connections and PD stratifications on a $G$-bundle}\label{S003} \vspace{3mm}

In this Appendix, we recall some definitions concerning connections  on a $G$-bundle and discuss their generalizations to higher-level PD stratifications.
As a consequence of our discussion, we will prove a generalization of  Cartier's theorem (cf. Proposition \ref{Ertwq}, Corollary \ref{EQkji}), asserting  that a principal bundle equipped with a flat  connection (or more generally,  a higher-level PD stratification)  with vanishing $p$-curvature comes, via Frobenius pull-back,  from a principal bundle on the Frobenius twist of the underlying space.

Let $S$ be a scheme  over a ring $R$
and $f: X \migi S$  a smooth $S$-scheme. 
Also, let $G$ be a connected  smooth  algebraic group over $R$, and write $\mfg := \mr{Lie}(G)$.

\LSP
%-----------------------------------------------------------------------------------------
\subsection{Connections on a $G$-bundle} \label{SS002}

First, we briefly  recall  the notion of a connection on a $G$-bundle.
We refer the reader to ~\cite[Chap.\,1 and Chap.\,3]{Wak8}  for detailed  discussions related to  the topic of this subsection. 

Denote by $\mr{GL}(\mfg)$ the algebraic group over $R$ classifying automorphism groups of  $\mfg$ over $R$. 
The {\it adjoint representation}  of $G$ is, by definition,   the morphism of algebraic $R$-groups 
\begin{equation} 
\label{adjrep} \mr{Ad}_G :G \migi \mr{GL}(\mfg) 
\end{equation}
  which, to any
  $h \in G$,  assigns the automorphism $d ( \mr{L}_h \circ \mr{R}_{h^{-1}}) |_e$ of $\mfg$, i.e., the differential of  $\mr{L}_h \circ \mr{R}_{h^{-1}}$ at the identity element $e$.
For an $R$-scheme $T$ and a $T$-rational point $h : T \migi G$, 
we occasionally   identify $\mr{Ad}_G (h)$ with the corresponding $\mcO_T$-linear automorphism of $\mcO_T \otimes \mfg$.

Let $\mcE$ be a (right) $G$-bundle on $X$ and denote by  
  $\pi : \mcE \migi X$ the structure morphism of $\mcE$ over $X$.
Denote by $\mfg_\mcE$  the adjoint vector bundle associated with $\mcE$, i.e., the vector bundle on $X$ corresponding to the $\mr{GL}(\mfg)$-bundle obtained from $\mcE$ via change of structure group by $\mr{Ad}_G$.
More specifically, 
$\mfg_\mcE$  is obtained as the sheaf  of $G$-invariant sections $\pi_*(\mcO_{\mcE} \otimes \mfg)^G$ in $\pi_*(\mcO_{\mcE} \otimes \mfg) \left(= \pi_*(\mcO_\mcE) \otimes  \mfg \right)$ with respect to the right  $G$-action given by $(a \otimes v) \cdot h = \mr{R}^*_h (a) \otimes \mr{Ad}_G (h^{-1})(v)$ for any $a \in \pi_*(\mcO_\mcE)$, $v \in \mfg$, and $h \in G$.
In particular, $\mfg_\mcE$ may be naturally identified with $\pi_*(\mcT_{\mcE/X})^G$ via  the natural  isomorphism $\omega_{\mcE/X}^\triangleright : \mcT_{\mcE/X} \isom \mcO_\mcE \otimes \mfg$ (cf. (\ref{L0116})).

Let
$\widetilde{\mcT}_{\mcE/S}$ denote the subsheaf  $\pi_* (\mcT_{\mcE/S})^G$ of $G$-invariant sections of
$\pi_* (\mcT_{\mcE/S})$.
 The  differential of $\pi$ gives rise to the following short exact sequence:
\begin{align}  \label{Ex0}
 0 \longmigi  \mfg_\mcE \longmigi  \widetilde{\mcT}_{\mcE/S} \xrightarrow{d_\mcE}  \mcT_{X/S} \longmigi 0.
 \end{align}
  Both $\mcT_{X/S}$ and $\widetilde{\mcT}_{\mcE/S}$ have natural 
  Lie bracket operations $[-,-]$ and the surjection $d_\mcE$ preserves these structures.
  The sheaf $\widetilde{\mcT}_{\mcE/S}$ together with  the surjection $d_{\mcE}$ forms a Lie algebroid on $X/S$ (cf. ~\cite[Chap.\,1, \S\, 1.2, Definition 1.11]{Wak8}).

By 
an {\it $S$-connection}
 on  $\mcE$, we mean a split injection of the above sequence, i.e.,   an $\mcO_X$-linear morphism
$\nabla : \mcT_{X/S} \migi \widetilde{\mcT}_{\mcE/S}$ with 
$d_\mcE\circ \nabla = \mr{id}_{\mcT_{X/S}}$.
Recall that, for an $S$-connection $\nabla$ on $\mcE$, 
the {\it curvature} of $\nabla$ is defined as  the global section
 \begin{align} \label{L054}
 \psi_\nabla 
 \in \Gamma (X, (\bigwedge^2 \Omega_{X/S}) \otimes_{\mcO_X} \widetilde{\mcT}_{\mcE/S})
 \end{align}
  determined uniquely by $\psi_\nabla (\partial_1,  \partial_2) = [\nabla (\partial_1), \nabla (\partial_2)] - \nabla  ([\partial_1, \partial_2])$ for any local sections $\partial_1$, $\partial_2 \in \mcT_{X/S}$.
Since the equality $(\mr{id} \otimes d_\mcE) (\psi_\nabla) =0$ holds,
$\psi_\nabla$ may be regarded as an element of $\Gamma  (X, (\bigwedge^2 \Omega_{X/S}) \otimes_{\mcO_X} \mfg_{\mcE})$.
We shall say that  $\nabla$ is   {\it flat} 
 if $\psi^{}_{\nabla} =0$.
 
 By a {\it flat $G$-bundle on $X/S$},
 we shall  means a pair $(\mcF, \nabla)$ consisting of  a $G$-bundle $\mcF$ on $X$ and a flat $S$-connection $\nabla$ on $\mcF$.

 \SSP
 %------------------------------------------------------------------------
 \bde \label{koirt}
Let $\nabla$ be a flat $S$-connection on $\mcE$.
 We shall say that $\nabla$  (or $(\mcE, \nabla)$) is {\bf quasi-nilpotent} if
  for any local section $v \in \widetilde{\mcT}_{\mcE/S}$ and any coordinate system $x_1, \cdots, x_n$  of  $X/S$ (such that $v, x_1, \cdots, x_n$ are defined on a common open subset of $X$),  there exists  a positive integer $s$ with  $\mr{ad}(\nabla (\frac{\partial}{\partial x_i}))^s (v) = 0$ for  every $i= 1, \cdots, n$, where $\mr{ad}$ denotes the adjoint operator on $\widetilde{\mcT}_{\mcE/S}$.
 \ede
 %------------------------------------------------------------------------
 \SSP

Next, let $\mcE$ be a $\mr{GL}_l$-bundle on $X$ (for $l \in \mbZ_{>0}$) and $\mcV$ the rank $l$ vector bundle corresponding to $\mcE$.
Note that  the notion of an $S$-connection recalled above 
is consistent with that of  ~\cite[(1.0)]{Kat}.
That is to say,
a choice of an $S$-connection on $\mcE$ 
can be interpreted as a choice of an $S$-connection on $\mcV$, which is 
 an $f^{-1}(\mcO_S)$-linear morphism
\begin{align}
\nabla : \mcV \migi \Omega_{X/S} \otimes_{\mcO_X} \mcV
\end{align} 
satisfying
 $\nabla (a \cdot v) = d a \otimes v + a \cdot \nabla (v)$ for local sections $a \in \mcO_X$ and $v \in \mcV$.
Throughout the present paper,  we will not distinguish these two notions of $S$-connection.
If $\partial$ is a local section of $\mcT_{X/S}$, then we shall write
$\nabla_\partial : \mcV \migi \mcV$ for the $f^{-1}(\mcO_S)$-linear endomorphism  of $\mcV$ defined as the composite of  $\nabla$ and  $\mr{id}_\mcV \otimes \partial : \mcV \otimes_{\mcO_X} \Omega_{X/S} \migi \left(\mcV \otimes_{\mcO_X} \mcO_X = \right) \mcV$.

For an $S$-connection $\nabla$ on $\mcV$, we shall set
\begin{align} \label{E0124}
\nabla^{(2)}:  \Omega_{X/S} \otimes_{\mcO_X} \mcV \migi (\bigwedge^2 \Omega_{X/S} ) \otimes_{\mcO_X} \mcV
\end{align}
to be the $f^{-1}(\mcO_S)$-linear morphism   determined uniquely by $\nabla^{(2)}(\omega \otimes v) = d \omega \otimes v - \omega \wedge \nabla (v)$ for  any local sections $\omega \in  \Omega_{X/S}$, $v \in \mcV$.
If  $\nabla$ is flat, then we have $\nabla^{(2)}\circ \nabla =0$, namely, 
the sequence 
\begin{align} \label{EOirzH}
\nabla^\bullet : \mcV \xrightarrow{\nabla} \Omega_{X/S} \otimes_{\mcO_X} \mcV \xrightarrow{\nabla^{(2)}}  (\bigwedge^2 \Omega_{X/S} ) \otimes_{\mcO_X} \mcV
\end{align}
 forms a complex of sheaves.

By a {\it flat (vector) bundle} on $X/S$, we mean a pair $(\mcV, \nabla)$ consisting of a vector bundle $\mcV$ on $X$ and a flat $S$-connection  $\nabla$ on $\mcV$.
The notion of an isomorphism between flat bundles can be defined in a natural manner.

\LSP
%\vspace{5mm}
%----------------------------------------------------------------------[begin subsection]-------------
\subsection{Pull-back and base-change of connections} \label{SS023}

Let $\mcE$ be a $G$-bundle on $X$ and $\nabla$ an $S$-connection on $\mcE$.
First,  suppose that we are given another smooth $S$-scheme   $Y$  and
 a (not necessarily \'{e}tale) morphism $y : Y \migi X$ over $S$.
Denote by   $y^*(\mcE) \left(:= Y \times_X \mcE \right)$ the pull-back of $\mcE$ by $y$, which forms a $G$-bundle on $Y$.
Then, we obtain the following morphism of short exact sequences:
 \begin{align} \label{E0179}
\vcenter{\xymatrix@C=46pt@R=36pt{
0 \ar[r] & \mfg_{y^*(\mcE)}\ar[r] \ar[d]_-{\wr}^{d y_{\mcE}} & \widetilde{\mcT}_{y^*(\mcE)/S} \ar[r]^-{d_{y^*(\mcE)}} \ar[d]^-{d \widetilde{y}_{\mcE}} & \mcT_{Y/S}\ar[d]_-{}^{dy} \ar[r] & 0 \\
0 \ar[r] & y^*(\mfg_\mcE) \ar[r] & y^*(\widetilde{\mcT}_{\mcE/S})\ar[r]_-{y^*(d_{\mcE})} & y^*(\mcT_{X/S}) \ar[r] & 0,
 }}
\end{align}
where
the middle and left-hand vertical arrows $d y_{\mcE}$, $d \widetilde{y}_{\mcE}$ are obtained by differentiating the projection $y^*(\mcE) \migi \mcE$ over $S$ and $X$ respectively. 
Note that $d y_{\mcE}$ is an isomorphism.
If $\nabla'$ denotes the split surjection $\widetilde{\mcT}_{\mcE/S} \migisurj \mfg_\mcE$  of  (\ref{Ex0}) corresponding to $\nabla$,
then the composite $d y_{\mcE}^{-1} \circ y^*(\nabla')\circ d \widetilde{y}_{\mcE}$
specifies a split surjection of the upper horizontal arrow in (\ref{E0179}).
Hence,  the corresponding split injection defines 
 an $S$-connection
\begin{align} \label{E0068}
y^*(\nabla) : \mcT_{Y/S} \migi \widetilde{\mcT}_{y^*(\mcE)}
\end{align}
on $y^*(\mcE)$.
Since   the equality  $y^*(\psi_\nabla) = \psi_{y^*(\nabla)}$
 %(resp., $y^*({^p}\psi_\nabla) = {^p}\psi_{y^*(\nabla)}$)
  holds under the identifications given by $d y_\mcE$ and $d \widetilde{y}_\mcE$,  
 the $S$-connection $y^*(\nabla)$ is flat 
 %(resp.,  $p$-flat)
if  $\nabla$ is flat.
% (resp., $p$-flat).
We  shall refer to $y^*(\nabla)$ as the {\it pull-back} of $\nabla$ by $y$.

Next, let $s: S' \migi S$ be a morphism of $k$-schemes and write $X' := S' \times_S X$.
%Denote by $s^*(\mcE)$ the $G$-bundle on $X'$ obtained from $\mcE$ via base-change    by $s$.
We use the notation $s^*(-)$ to denote the result of base-change by  $s$.
Differentiating the projections $X' \migi X$ and $s^*(\mcE) \migi \mcE$ induces the following commutative square:
 \begin{align} \label{E0180}
\vcenter{\xymatrix@C=46pt@R=36pt{
\widetilde{\mcT}_{s^*(\mcE)/S'} \ar[r]^-{d_{s^*(\mcE)}} \ar[d]_-{\wr} & \mcT_{X'/S'}\ar[d]^-{\wr} \\
s^*(\widetilde{\mcT}_{\mcE/S})\ar[r]_-{s^*(d_{\mcE})} & s^*(\mcT_{X/S}).
 }}
\end{align}
The base-change of $\nabla$ by $s$  specifies, via the vertical arrows in this  square, 
an $S'$-connection
\begin{align} \label{E0067}
s^*(\nabla) : \mcT_{X'/S} \migi \widetilde{\mcT}_{s^*(\mcE)/S'}
\end{align}
on $s^*(\mcE)$.
Since  the equality $s^*(\psi_\nabla) = \psi_{s^*(\nabla)}$
 %(resp., $s^*({^p}\psi_\nabla) = {^p}\psi_{s^*(\nabla)}$)
  holds, the $S$-connection $s^*(\nabla)$ is flat 
  %(resp., $p$-flat)
   if  $\nabla$ is flat.
   % (resp., $p$-flat).
We shall refer to $s^*(\nabla)$ as the {\it base-change} of $\nabla$ by $s$.

\LSP
%---------------------------[begin subsection]-------------
\subsection{$p$-curvature and Cartier operators} \label{SS016}

 In this subsection, suppose that $S$ is a scheme over $\mbF_p$.
Let $\mcE$ be as above.
 Recall from ~\cite[Chap.\,3, \S\,3.3, Definition 3.8]{Wak8} that,   given a flat $S$-connection $\nabla$  on $\mcE$, we can obtain the global  section
 \begin{align} \label{L056}
 {^p}\psi_\nabla \in \Gamma (X, F^*_X (\Omega_{X/S}) \otimes_{\mcO_X} \widetilde{\mcT}_{\mcE/S})
 \end{align}
  uniquely determined  by  $\langle F^{-1}_X (\partial), {^p}\psi_\nabla \rangle = \nabla (\partial)^{[p]} - \nabla (\partial^{[p]})$ for each local section $\partial \in \mcT_{X/S}$ (cf. \S\,\ref{SS0131} for the definition of $p$-power operations $(-)^{[p]}$).
  This global section ${^p}\psi_\nabla$ is called the {\it $p$-curvature} of $\nabla$.
Just as in the case of curvature, we see that ${^p}\psi_\nabla$  lies in $\Gamma (X, F^*_X (\Omega_{X/S}) \otimes_{\mcO_X} \mfg_{\mcE})$.
The $S$-connection $\nabla$ is called  {\it $p$-flat} if  ($\psi_\nabla=0$ and) ${^p}\psi_\nabla =0$.
If $\mcE$ is a $\mr{GL}_l$-bundle and $\mcV$ denotes the corresponding vector bundle, then the notion of $p$-curvature is consistent with that of ~\cite[\S\,5, (5.0.3)]{Kat} under the natural identification  $\mcE nd_{\mcO_X}(\mcV) = \mr{Lie}(\mr{GL}_l)_\mcE$.

   \SSP
 %------------------------------------------------------------------------
 \begin{rema} \label{koRirt}
 Suppose that $\nabla$  is $p$-flat.
 Then, one may verify that  it is locally quasi-nilpotent (cf. ~\cite[\S\,4, Remark 4.11]{Og2}).
 \end{rema}
 %------------------------------------------------------------------------
 \SSP

  %------------------------------------------------------------------------
 \begin{rema} \label{koRiGG}
 Let $y$ and $s$ be as in the previous subsection.
 Then, under certain identifications, we have $y^*({^p}\psi_\nabla) = {^p}\psi_{y^*(\nabla)}$ (resp.,  $s^*({^p}\psi_\nabla) = {^p}\psi_{s^*(\nabla)}$).
 This implies that $y^*(\nabla)$ (resp., $s^*(\nabla)$)  is $p$-flat if  $\nabla$ is $p$-flat.
  \end{rema}
 %------------------------------------------------------------------------
 \SSP

  Let $(\mcV, \nabla)$ be a flat  vector bundle on $X/S$.
  Since both $\nabla$ and $\nabla^{(2)}$ are $F^{-1}_{X/S}(\mcO_{X^{(1)}})$-linear,
 the sheaves  $\mr{Ker}(\nabla)$,  $\mr{Ker}(\nabla^{(2)})$, and $\mcH^1 (\nabla^\bullet) \left(:= \mr{Ker}(\nabla^{(2)})/\mr{Im}(\nabla) \right)$ may be regarded as $\mcO_{X^{(1)}}$-modules via the underlying homeomorphism of $F_{X/S}$.
 Recall (cf. ~\cite[Chap.\, 1, Proposition 1.2.4]{Og}) that 
  the {\it Cartier operator} associated with $(\mcV, \nabla)$
 is  a unique $\mcO_{X^{(1)}}$-linear morphism
 \begin{align} \label{Efty}
 C_{(\mcV, \nabla)} : \mr{Ker}(\nabla^{(2)}) \migi \Omega_{X^{(1)}/S} \otimes_{\mcO_{X^{(1)}}} F_{X/S*}(\mcV)
 \end{align}
 satisfying the equality 
 \begin{align} \label{Wer4}
 \langle F_S^* (\partial), C_{(\mcV, \nabla)} (a) \rangle = \langle \partial^{[p]}, a \rangle - \nabla_\partial^{p-1} (\langle \partial, a\rangle)
 \end{align}
  (cf. Example \ref{Eorkcr0} for the definition of $F_S^*(\partial)$) for any local sections $\partial \in \mcT_{X/S}$,
 $a \in \mr{Ker}(\nabla^{(2)})$.
 Let us  consider   the  composite  isomorphism
\begin{align} \label{projformula1}
\Omega_{X^{(1)}/S}\otimes_{\mcO_{X^{(1)}}} F_{X/S*}(\mcV)  \isom 
F_{X/S*}(F_{X/S}^*(\Omega_{X^{(1)}/S}) \otimes_{\mcO_X} \mcV) \isom F_{X/S*}(F_X^*(\Omega_{X/S}) \otimes_{\mcO_X} \mcV), 
\end{align}
where  the first isomorphism follows from the projection formula.
By using this composite, 
we occasionally regard $C_{(\mcV, \nabla)}$ as a morphism $\mr{Ker}(\nabla^{(2)})\migi F_X^*(\Omega_{X/S}) \otimes_{\mcO_{X}} \mcV$.

Suppose further that 
 $(\mcV, \nabla)$ has vanishing $p$-curvature.
 Then,
according to  ~\cite[Theorem 3.1.1]{Og3} (and the discussion following ~\cite[Chap.\,I, Proposition 1.2.4]{Og}),  
there exists an $\mcO_{X^{(1)}}$-linear  isomorphism 
\begin{align} \label{Qw8345}
\overline{C}_{(\mcV, \nabla)} : \mcH^1 (\nabla^\bullet) \isom \Omega_{X^{(1)}/S} \otimes \mr{Ker}(\nabla)
\end{align}
which makes the following square diagram commute:
\begin{align}\label{Cartier3F}
\vcenter{\xymatrix@C=46pt@R=36pt{
\mr{Ker}(\nabla^{(1)})\ar[r]^-{C_{(\mcV, \nabla)}} 
\ar[d]_-{\mr{quotient}}
 & \Omega_{X^{(1)}/S}\otimes_{\mcO_{X^{(1)}}} F_{X/S*}(\mcV) \\
\mcH^1 (\nabla^\bullet) \ar[r]^-{\sim}_-{\overline{C}_{(\mcV, \nabla)}}& \Omega_{X^{(1)}/S} \otimes_{\mcO_{X^{(1)}}} \mr{Ker}(\nabla). \ar[u]_-{\mr{inclusion}}
}}
\end{align}

\LSP
%---------------------------[begin subsection]-------------
\subsection{Differential operators of  level $m$} \label{SS040}
%{\bf .}

In this subsection, we briefly recall 
the rings  of differential operators on $X/S$ of  higher  level.
The case of finite level was originally  defined  and studied  by Berthelot (see ~\cite{Ber1} and ~\cite{Ber2} for more details).

We here  suppose that
  $p$ is  nilpotent in $R$.
  Let $m$ be a nonnegative integer; whenever  we deal with various objects associated with   this  $m$,
  % (i.e., $m \neq \infty$),
we always suppose that the scheme $S$ is equipped with an $m$-PD ideal $(\mfa, \mfb, \gamma)$ such that $\gamma$ extends to $X$.
That is to say,  $\mfa$ is a quasi-coherent ideal of $\mcO_S$ and $(\mfb, \gamma)$ is an $m$-PD structure 
 on $\mfa$
    (cf. ~\cite[\S\,1.3, Definition 1.3.2, (ii)]{Ber1}).
Note that  $\gamma$ always extends to $X$ in the case $\mfb = (p)$.

It follows from the discussion in ~\cite[\S\,1.4]{Ber1} that 
there exists   the  $m$-PD envelope $P_{X/S, (m)}$ of the diagonal embedding  $X \migi X \times_S X$.  
Denote by $\overline{I}$ the defining ideal of the natural morphism $X \migiincl P_{X/S, (m)}$.
Recall that  $\overline{I}$ admits   the $m$-PD-adic filtration $\{ \overline{I}^{\{ l \}} \}_{l \in \mbZ_{\geq 0}}$ constructed in the manner of ~\cite[Definition A. 3]{Ber2}.
Let $\mcP_{X/S, (m)}$ denote the structure sheaf of $P_{X/S, (m)}$.
For each  $l \in \mbZ_{\geq 0}$, we denote by $\mcP_{X/S, (m)}^l$ the quotient sheaf of $\mcP_{X/S, (m)}$   by $\overline{I}^{\{ l+1 \}}$ and by $P^l_{X/S, (m)}$ the closed subscheme of $P_{X/S, (m)}$ defined by the surjection $\mcP_{X/S, (m)} \migisurj \mcP_{X/S, (m)}^l$.

On the other hand, 
for each $l \in \mbZ_{\geq 0}$, denote by $\mcP_{X/S, (\infty)}^l$ 
the quotient of the structure sheaf of $X \times_S X$ by 
 $I^{l+1}$, where $I$ is  the ideal of the diagonal embedding $X \migi X \times_S X$.
Also, denote by $P_{X/S, (\infty)}^l$ the closed subscheme of $X \times_S X$ defined  by  the surjection  $\mcO_{X \times_S X} \migisurj \mcP_{X/S, (\infty)}^l$.

Hence,  in each case of $m \in \mbZ_{\geq 0} \sqcup \{ \infty \}$,  we have the following  sequence of sheaves on $X$:
\begin{align}
\cdots \migi \mcP^l_{X/S, (m)} \migi \mcP^{l-1}_{X/S, (m)} \migi  \cdots \migi \mcP_{X/S, (m)}^1 \migi \mcP^0_{X/S, (m)} \left(= \mcO_X \right).
\end{align}
Let 
$\mr{pr}_1^l$ and $\mr{pr}_2^l$ denote  the morphisms $P_{X/S, (m)}^l \migi X$ induced by 
the first and second projections  $X \times_S X \migi X$
  respectively.
As asserted in ~\cite[\S\,2.1, Proposition 2.1.3, (ii)]{Ber1} and ~\cite[\S\,16, (16.8.9.2)]{EGA4},
we have the canonical  $\mcO_X$-algebra morphism
\begin{align} \label{EOFOSQ98}
\delta_m^{l, l'}: \mcP_{X/S, (m)}^{l+l'} \migi \mcP^l_{X/S, (m)} \otimes_{\mcO_X} \mcP^{l'}_{X/S, (m)}
\end{align}
for each pair of nonnegative integers $(l, l')$ such that  if $m'$ is an element of $\mbZ_{\geq 0} \sqcup \{ \infty \}$ with $m' \leq m$, then 
 the following diagram  is commutative:
\begin{align} \label{COWQ231}
\vcenter{\xymatrix@C=46pt@R=36pt{
\mcP_{X/S, (m)}^{l+l'} \ar[r]^-{\delta_m^{l, l'}} \ar[d]_-{\varsigma^{l+l'}_{m, m'}}& \mcP_{X/S, (m)}^l \otimes_{\mcO_X} \mcP_{X/S, (m)}^{l'} \ar[d]^-{\varsigma^{l}_{m, m'} \otimes \varsigma^{l'}_{m, m'}}
\\
\mcP_{X/S, (m')}^{l+l'} \ar[r]_-{\delta_{m'}^{l, l'}}& \mcP_{X/S, (m')}^l \otimes_{\mcO_X} \mcP_{X/S, (m')}^{l'},
}}
\end{align}
where $\varsigma^{l''}_{m, m'}$ ($l'' \in \mbZ_{\geq 0}$) denotes 
 the natural morphism. 
 The morphism $\delta_m^{l, l'}$ specifies a morphism of $X$-schemes
 \begin{align} \label{MROD3}
 \delta_m^{l, l' \sharp}:  P_{X/S, (m)}^l  \times_X P_{X/S, (m)}^{l'} \migi P_{X/S, (m)}^{l+l'}.
 \end{align}

For each nonnegative integer $l$,
 write $\mcD_{X/S, l}^{(m)}$ for the sheaf 
 \begin{align} \label{MROD2}
\mcD_{X/S, l}^{(m)} := \mcH om_{\mcO_X} (\mr{pr}_1^{l*}(\mcP^l_{X/S, (m)}), \mcO_X).
\end{align}
In particular, we have $\mcD_{X/S, 0}^{(m)} = \mcO_X$ and $\mcD_{X/S, 1}^{(m)}/\mcD_{X/S, 0}^{(m)} = \mcT_{X/S}$.
The  {\it sheaf of differential operators of level $m$} is defined by 
\begin{align} \label{MROD1}
\mcD_{X/S}^{(m)} := \bigcup_{l \in \mbZ_{\geq 0}} \mcD_{X/S, l}^{(m)}.
\end{align}
The various morphisms  $\delta_{m}^{l, l'}$ (for $l, l' \in \mbZ_{\geq 0}$) induces  a structure of  (possibly noncommutative) $\mcO_X$-algebra $\mcD_{X/S}^{(m)} \otimes_{\mcO_X} \mcD_{X/S}^{(m)} \migi \mcD_{X/S}^{(m)}$ on $\mcD_{X/S}^{(m)}$.
The sheaves  $\{ \mcD_{X/S}^{(m)} \}_{m \in \mbZ_{\geq 0}}$ forms  an inductive system and there exists  a canonical identification $\varinjlim_{m \in \mbZ_{\geq 0}}\mcD_{X/S}^{(m)} = \mcD_{X/S}^{(\infty)}$.

By  a $\mcD_{X/S}^{(m)}$-module, we mean a pair $(\mcV, \nabla)$ consisting of an $\mcO_X$-module $\mcV$ and a left $\mcD_{X/S}^{(m)}$-action $\nabla$ on $\mcV$ extending its $\mcO_X$-module structure.

\vspace{5mm}
%---------------------------[begin subsection]-------------
\subsection{$m$-PD stratifications  on a $G$-bundle} \label{SS0239}

In the rest of this Appendix, we suppose that $G$ is affine.
Denote by $R_G$ the $R$-algebra defined as the coordinate ring of $G$.
Also, let us fix $m \in \mbZ_{\geq 0} \sqcup \{ \infty \}$.

%----------------------------------[begin definition]------------------
\SSP
\bde  \label{D057}
\begin{itemize}
\item[(i)]
Let $\mcE$ be a $G$-bundle on $X$.
An {\bf $m$-PD stratification} on $\mcE$ is
a collection
\begin{align}
\varepsilon^\ST :=  \{\varepsilon_l \}_{l \in \mbZ_{\geq 0}},
\end{align}
where each $\varepsilon_l$ denotes an isomorphism $P_{X/S, (m)}^l  \times_X \mcE\left(=\mr{pr}_2^{l *}(\mcE)\right) \isom \mcE \times_X P^l_{X/S, (m)} \left(=\mr{pr}^{l*}_1 (\mcE)\right)$ of $G$-bundles  on $P_{X/S, (m)}^{l}$
satisfying the following conditions:
\begin{itemize}
\item
The equalities $\varepsilon_0 = \mr{id}_\mcE$ and  $\varepsilon_{l_2} |_{P_{X/S, (m)}^{l_1}} = \varepsilon_{l_1}$  hold  for  any pair of integers 
 $(l_1, l_2)$ with $l_1 \leq l_2$.
\item
The cocycle condition holds:
to be precise, for any pair of nonnegative  integers $(l, l')$,  the following diagram is commutative:
\begin{align} \label{E00235}
\vcenter{\xymatrix@C=6pt@R=36pt{
P_{X/S, (m)}^l \times_X P_{X/S, (m)}^{l'} \times_X  \mcE
\ar[rr]_-{\sim}^{\delta_{m}^{l, l' \sharp *}(\varepsilon_{l+l'})} \ar[rd]^-{\sim}_-{q_2^{l, l'*}(\varepsilon_{l+l'})}&&  
\mcE\times_X P_{X/S, (m)}^l \times_X P_{X/S, (m)}^{l'}
\\
& P_{X/S, (m)}^{l} \times_X \mcE \times_X P_{X/S}^{l'} \ar[ru]^-{\sim}_-{q_1^{l, l' *}(\varepsilon_{l+l'})}, &
 }}
\end{align}
where $q_1^{l, l'}$ and $q_2^{l, l'}$  denote  the morphisms defined, respectively,  as
 \begin{align} \label{Ee1}
 q_1^{l, l'} &: P_{X/S, (m)}^l \times_X P_{X/S, (m)}^{l'} \xrightarrow{\mr{pr}_1} P_{X/S, (m)}^l \xrightarrow{\mr{immersion}} P_{X/S, (m)}^{l+l'}, \\
  q_2^{l, l'} &: P_{X/S, (m)}^l \times_X P_{X/S, (m)}^{l'} \xrightarrow{\mr{pr}_2} P_{X/S, (m)}^{l'} \xrightarrow{\mr{immersion}} P_{X/S, (m)}^{l+l'}. \notag
 \end{align}
\end{itemize}

Also, by an {\bf $m$-PD stratified $G$-bundle}
   on $X/S$, we mean    a pair 
\begin{align}
(\mcE, \varepsilon^\ST)
\end{align}
 consisting of a $G$-bundle $\mcE$ on $X$ and  an $m$-PD stratification $\varepsilon^\ST$ on $\mcE$.
\item[(ii)]
Let $(\mcE, \varepsilon^\ST)$ and 
$(\mcE', \varepsilon'^\ST)$ be $m$-PD stratified   $G$-bundles on $X/S$, where $\varepsilon^\ST := \{ \varepsilon_l \}_l$,  $\varepsilon'^\ST :=\{ \varepsilon'_l \}_l$.
An {\bf isomorphism of $m$-PD stratified   $G$-bundles} from 
$(\mcE, \varepsilon^\ST)$ to  $(\mcE', \varepsilon'^\ST)$
is defined as an isomorphism of $G$-bundles $\alpha : \mcE \isom  \mcE'$   such that, for each nonnegative  integer $l$, the following square diagram is commutative:
\begin{align} \label{E00237}
\vcenter{\xymatrix@C=46pt@R=36pt{
  P_{X/S, (m)}^l \times_X \mcE 
 \ar[r]_-{\sim}^-{\varepsilon_l} \ar[d]^-{\wr}_-{\mr{id} \times\alpha} &
 \mcE \times_X P_{X/S, (m)}^l
\ar[d]_-{\wr}^-{\alpha \times \mr{id}}
\\
  P_{X/S, (m)}^l \times_X \mcE'
 \ar[r]^-{\sim}_{\varepsilon'_l}&
 \mcE' \times_X P_{X/S, (m)}^{l}.
 }}
\end{align}

\end{itemize}
 \ede
%-------------------------[end definition]-------------------
\SSP

In what follows, we shall describe an $m$-PD stratification on a $G$-bundle by using  
the classical notion of an $m$-PD stratification (defined on an $\mcO_X$-module).
Let $\mcE$ be a $G$-bundle on $X$.
Since  $\mcE$ is affine over $X$ because of the affineness assumption on $G$,   it defines an $\mcO_X$-algebra;
we shall  write $\mcO_\mcE$ for   this $\mcO_X$-algebra  by abuse of notation (hence, $\mcS pec (\mcO_\mcE) = \mcE$).
If $\mr{R}_\mcE : \mcE \times G \migi G$ denotes
the $G$-action  on $\mcE$,  then it corresponds to 
an $\mcO_X$-algebra morphism $\mr{R}_\mcE^\sharp : \mcO_\mcE \migi \mcO_\mcE \otimes_R  R_G$. (We shall call $\mr{R}_\mcE^\sharp$ the {\it $G$-coaction} on $\mcO_\mcE$.)
Now, let
  $\varepsilon^\ST := \{ \varepsilon_l \}_{l}$ be an $m$-PD stratification on $\mcE$.
For each nonnegative integer $l$, the isomorphism $\varepsilon_l$ defines a $\mcP_{X/S, (m)}^l$-algebra isomorphism $\varepsilon^\natural_l : \mcP_{X/S, (m)}^{l} \otimes_{\mcO_X} \mcO_{\mcE} \isom \mcO_{\mcE} \otimes_{\mcO_X} \mcP_{X/S, (m)}^{l}$.
The  $G$-equivariance condition of  $\varepsilon_l$ is interpreted as the commutativity of the following square diagram:
\begin{align} \label{E00100}
\vcenter{\xymatrix@C=30pt@R=36pt{
\mcP_{X/S, (m)}^{l} \otimes_{\mcO_X} \mcO_{\mcE}  \ar[r]_-{\sim}^{\varepsilon^\natural_l} \ar[d]_-{\mr{id} \otimes\mr{R}_\mcE^\sharp} & \mcO_{\mcE} \otimes_{\mcO_X} \mcP_{X/S, (m)}^{l}  \ar[d]^-{\mr{R}_\mcE^\sharp \otimes \mr{id}}
\\
(\mcP_{X/S, (m)}^{l} \otimes_{\mcO_X} \mcO_{\mcE}) \otimes_R R_G \ar[r]^-{\sim}_-{\varepsilon_l^\natural \otimes \mr{id}_{R_G}}& (\mcO_{\mcE} \otimes_{\mcO_X} \mcP_{X/S, (m)}^{l}) \otimes R_G \left(=  (\mcO_{\mcE} \otimes_R R_G) \otimes_{\mcO_X} \mcP_{X/S, (m)}^{l}\right).
 }} \hspace{-10mm}
\end{align}
Moreover, the commutativity of (\ref{E00235})  reads the commuativity of the following diagram:
\begin{align} \label{E00121}
\vcenter{\xymatrix@C=-32pt@R=36pt{
\mcP_{X/S, (m)}^l \otimes_{\mcO_X} \mcP_{X/S, (m)}^{l'} \otimes_{\mcO_X} \mcO_\mcE
&& 
\mcO_\mcE \otimes_{\mcO_X} \mcP_{X/S, (m)}^l \otimes_{\mcO_X} \mcP_{X/S, (m)}^{l'}
\ar[ll]^-{\sim}_-{\delta_{m}^{l+l'\sharp *} (\varepsilon^\natural_{l+l'})} \ar[ld]_-{\sim}_-{\sim}^-{q_1^{l, l'  *}( \varepsilon^\natural_{l+l'})}
\\
& \mcP^l_{X/S, (m)} \otimes_{\mcO_X} \mcO_\mcE \otimes_{\mcO_X} \mcP^{l'}_{X/S, (m)} \ar[ul]_-{\sim}^-{q_2^{l,  l'  *}( \varepsilon^\natural_{l+l'})}. &
}}
\end{align}
%where $q_1^{l, l' \sharp}$ (resp.,  $q_2^{l, l' \sharp}$) denotes the morphism $\mcP_{X/S, (m)}^{l+l'} \migi \mcP_{X/S, (m)}^l \otimes_{\mcO_X} \mcP_{X/S, (m)}^{l'}$  corresponding to $q_1^{l, l'}$ (resp.,  $q_2^{l, l'}$).
Thus, 
the resulting collection
\begin{align} \label{Ertfq}
\varepsilon^{\ST \natural} := \{ \varepsilon^\natural_l \}_{l \in \mbZ_{\geq 0}}
\end{align}
 forms an $m$-PD stratification on $\mcO_\mcE$ in the usual sense (cf. ~\cite[\S\,2.3, Definition 2.3.1]{Ber1}).
 (By a  {\it $\infty$-PD stratification}, we shall   mean  a stratification in the  usual sense, i.e.,  in the sense of ~\cite[\S\,2, Definition 2.10]{Og2}).

Conversely, suppose that we are given an $m$-PD stratification  $\varepsilon^{\ST \natural} := \{ \varepsilon_l^\natural\}_{l \geq 0}$ on $\mcO_\mcE$ such that  each $\varepsilon^\natural_l$ is   a  $\mcP_{X/S, (m)}^l$-algebra  isomorphism $\mcP_{X/S, (m)}^l \otimes_{\mcO_X} \mcO_\mcE \isom \mcO_\mcE \otimes_{\mcO_X} \mcP_{X/S, (m)}^l$ and
 (\ref{E00100}) 
  for this collection is  commutative.
Then, by applying the functor $\mcS pec (-)$ to each $\varepsilon_l^\natural$,
we obtain a collection
of isomorphisms $P_{X/S, (m)}^l \times_X \mcE \isom \mcE \times_X P_{X/S, (m)}^l $ ($l \geq 0$) forming 
  an  $m$-PD stratification on $\mcE$.
\SSP

%-------------------------------------------------------------------
\begin{rema} \label{Eruy78}
One may verify that giving a collection as in (\ref{Ertfq}) is equivalent to giving 
a {\it compatible} collection
\begin{align} \label{Ler46}
\overline{\varepsilon}^{\ST \natural} := \{ \overline{\varepsilon}_l^{\natural}\}_{l \in \mbZ_{\geq 0}},
\end{align}
where each $\overline{\varepsilon}_l^{\natural}$ denotes an $\mcO_X$-algebra  morphism $\mcO_\mcE \migi \mcO_\mcE \otimes_{\mcO_X} \mcP_{X/S, (m)}^l$, satisfying the following two conditions:
\begin{itemize}
\item
For each $l \geq 0$, the following square diagram is commutative:
\begin{align} \label{E00df1}
\vcenter{\xymatrix@C=46pt@R=36pt{
\mcO_\mcE \ar[r]^-{\overline{\varepsilon}_l^{\natural}} \ar[d]_-{\mr{R}_\mcE^\sharp} &  \mcO_\mcE \otimes_{\mcO_X} \mcP_{X/S, (m)}^l\ar[d]^-{\mr{R}_\mcE^\sharp \otimes \mr{id}}
\\
\mcO_\mcE \otimes_R R_G \ar[r]_-{\overline{\varepsilon}_l^{\natural} \otimes \mr{id}_{R_G}} & (\mcO_\mcE \otimes_{\mcO_X} \mcP_{X/S, (m)}^l) \otimes_R R_G \left(= (\mcO_\mcE \otimes_R R_G) \otimes_{\mcO_X} \mcP_{X/S, (m)}^l \right).
}}
\end{align}
\item
The equality $\overline{\varepsilon}_0^{\natural} = \mr{id}_{\mcO_\mcE}$ holds, and for each pair of nonnegative integers $(l, l')$, the following square diagram is commutative:
\begin{align} \label{E00df1}
\vcenter{\xymatrix@C=46pt@R=36pt{
\mcO_\mcE \otimes_{\mcO_X} \mcP_{X/S, (m)}^{l+l'} \ar[r]^-{\mr{id} \otimes \delta_m^{l, l' }} &  \mcO_\mcE \otimes_{\mcO_X} \mcP_{X/S, (m)}^{l} \otimes_{\mcO_X} \mcP_{X/S, (m)}^{l'}
\\
\mcO_\mcE\ar[r]_-{\overline{\varepsilon}_{l'}^{\natural}} \ar[u]^-{\overline{\varepsilon}_{l+l'}^{\natural}}& \mcO_\mcE \otimes_{\mcO_X} \mcP_{X/S, (m)}^{l'}\ar[u]_-{\overline{\varepsilon}_l^{\natural}\otimes \mr{id}}.
}}
\end{align}
\end{itemize}
In particular, according to ~\cite[\S\,2.3, Proposition 2.3.2]{Ber1}, this collection of data  induces a structure of left $\mcD_{X/S}^{(m)}$-module on $\mcO_\mcE$ extending its $\mcO_X$-module   structure.

Now, let $\varepsilon^{\ST} := \{ \varepsilon_l \}_{l \in \mbZ_{\geq 0} }$ be an $\M$-PD stratification on $\mcE$ and  $\M'$ an element of $\mbZ_{\geq 0}\sqcup \{ \infty \}$ with  $\M' \leq \M$.
Denote by 
 $\overline{\varepsilon}^{\ST \natural} := \{ \overline{\varepsilon}_l^{\natural}\}_{l \geq 0}$
 the collection as in (\ref{Ler46}) corresponding to $\varepsilon^{\ST}$.
The collection $\overline{\varepsilon}^{\ST \natural} |^{(m' )} := \{ (\mr{id}_{\mcO_\mcE} \otimes q_l) \circ \overline{\varepsilon}_l^{\natural} \}_{l}$, where $q_l$ denotes  the natural morphism $\mcP_{X/S, (m)}^{l} \migi \mcP_{X/S, (m')}^l$, 
satisfies the above two conditions for $\M'$.
Hence, it determines  an $\M'$-PD stratification
\begin{align} \label{Rt567}
\varepsilon^{\ST} |^{( m')}
\end{align}
 on $\mcE$, which will be 
called  the {\bf $m'$-th truncation} of $\varepsilon^{\ST}$.
Since $\varinjlim_{m}\mcD_{X/S}^{(m)} = \mcD_{X/S}^{(\infty)}$,  a $\infty$-PD stratification on $\mcE$ may be interpreted as 
a compatible  collection $\{ \varepsilon^\ST_m \}_{m \in \mbZ_{\geq 0}}$ of    $m$-PD stratifications on $\mcE$ for various $m$.
\end{rema}
%-------------------------------------------------------------------

\SSP

%-------------------------------------------------------------------
\bde \label{Erty23}
Let $\varepsilon^\ST := \{ \varepsilon_l \}_l$ be an $m$-PD stratification on $\mcE$ with $m \in \mbZ_{\geq 0}$.
We shall say that $\varepsilon^\ST$ (or $(\mcE, \varepsilon^\ST)$) is {\bf quasi-nilpotent}
if the  $\mcD_{X/S}^{(m)}$-module $\mcO_\mcE$ corresponding to  $\varepsilon^{\ST \natural}$  (cf. Remark \ref{Eruy78} above) is quasi-nilpotent in the sense of ~\cite[\S\,2.3, Definitions]{Ber1}.
\ede
%------------------------------------------------------------------

\SSP

%------------------------------------------------------------------
\begin{rema} \label{Eroz892}
According to  ~\cite[\S\,2.3, Proposition 2.3.7]{Ber1}, giving a quasi-nilpotent $m$-PD stratification on $\mcE$ is equivalent to giving an isomorphism
$P_{X/S, (m)} \times_X \mcE \isom \mcE \times_X P_{X/S, (m)}$
of $G$-bundles on $P_{X/S, (m)}$ satisfying the usual cocycle condition and whose reductions modulo $\overline{I}^{\{ l+1\}}$
 %to  $P_{X/S, (m)}^l$'s 
  for various $l$ form an $m$-PD stratification.
\end{rema}
%------------------------------------------------------------------

\SSP

%------------------------------------------------------------------
 \begin{rema} \label{Eoaseh3}
 We can define the pull-back (resp., the base-change) of an $m$-PD stratification on $\mcE$  to an  \'{e}tale $X$-scheme (resp., an $S$-scheme) by applying  both the discussion in ~\cite[\S\,2.1]{Ber2} and the description in terms of $\mcD_{X/S}^{(m)}$-modules mentioned above.
Also,  the   formation of pull-back (resp., base-change) commutes with truncation to lower levels and preserves the quasi-nilpotency condition.
 We will  omit the details because the construction is straightforward.
\end{rema}
%------------------------------------------------------------------

\LSP
%---------------------------[begin subsection]-------------
\subsection{$m$-crystals of $G$-bundles} \label{SS0140}
%{\bf .}

Let $m$ be an element of $\mbZ_{\geq 0}$ and  $Y$ an $S$-scheme such that the $m$-PD structure  of $S$ extends to $Y$.
Recall (cf. ~\cite[\S\,4.1, Definition 4.1.1]{LQ}) that, for a scheme $U$ over $Y$,  a {\it  $m$-PD thickening} of $U$ over $(S, \mfa, \mfb, \gamma)$ is a collection of data  $(U, T, J, \delta)$  consisting  of a scheme $T$ over $S$, a closed immersion $U \migiincl T$ over $S$ and an $m$-PD structure $(J, \delta)$ on the defining ideal of $U \migiincl T$ compatible with $(\mfb, \gamma)$.
A morphism of $m$-PD thickenings over $(S, \mfa, \mfb, \gamma)$ can be  defined in an obvious way.
The {\it $m$-crystalline site $\mr{Cris}^{(m)}(Y/S)$} is the category of $m$-PD thickenings $(U, T, J, \delta)$ of an open subscheme  $U$ of $Y$ over $(S, \mfa, \mfb, \gamma)$ endowed with the topology induced by the \'{e}tale topology on $T$.

An {\bf $m$-crystal of  $G$-bundles} on $Y/S$ is
a cartesian section $\mcE^\diamondsuit$ on $\mr{Cris}^{(m)}(Y/S)$.
That is to say, $\mcE^\diamondsuit$ is 
 given by the following data:
For every $m$-PD thickening $(U, T, J, \delta)$, a $G$-bundle $\mcE^\dd_T$ on $T$, and for every morphism $u: (U_1, T_1, J_1, \delta_1) \migi (U_2, T_2, J_2, \delta_2)$ in $\mr{Cris}^{(m)}(X/S)$, an isomorphism of $G$-bundles  $\rho_u : u^*(\mcE^\dd_{T_2}) \isom \mcE^\dd_{T_1}$, satisfying the cocycle condition.

Let $\mcE^\dd$ be an $m$-crystal of $G$-bundles on $Y/S$ and $m'$ a nonnegative integer with $m' \leq m$.
Since 
an $m$-PD ideal can be considered as an $m'$-PD ideal,
the $G$-bundles on various thickenings defining $\mcE^\dd$
also form  an $m'$-crystal  of $G$-bundles 
\begin{align} \label{e9w308}
\mcE^\dd |^{( m )}
\end{align}
on $Y/S$; we shall refer to it as 
the {\bf $m'$-th truncation} of $\mcE^\dd$.

Now, let $X$ be as before and suppose that there exists a closed immersion $\iota: Y \migiincl X$ over $S$.
For each $m$-crystal $\mcE^\diamondsuit$ in $G$-bundles on $X/S$, we can  obtain the restriction $\iota^*(\mcE^\diamondsuit)$ of $\mcE^\diamondsuit$.
Then, the following assertion holds.

\SSP
%----------------------------------[begin lemma]------------------
\bpr \label{Ewr456}
Let $(J_0, \gamma_0)$ be a sub-$m$-PD structure of  $(J, \gamma)$ over $I$ (cf. ~\cite[\S\,1.3, Definition 1.3,1]{Ber1}).
Write  $\overline{S}$ for  the closed subscheme of $S$ defined by the ideal $J_0$, and 
write $\overline{X} := X \times_S \overline{S}$.
(Note that the $m$-PD structure on $\overline{S}$ induced from  that on $S$ extends to $\overline{X}$.)
Then, the restriction by the closed immersion $\overline{X} \migiincl X$
%$\mcE^\dd \mapsto \iota^*(\mcE^\dd)$
determines an equivalence of categories   between the groupoids of $m$-crystals of  $G$-bundles on $X/S$ and $\overline{X}/S$.
\epr
%----------------------------------[begin lemma]------------------
\begin{proof}
The assertion follows from various definitions involved  together with the corresponding fact on   $m$-crystals of sheaves  asserted in   ~\cite[\S\,4, Corollary 4.1.6]{LQ}.
\end{proof}

%----------------------------------[begin lemma]------------------
\SSP
\bpr \label{P019}
(Recall that $G$ is assumed to be affine.)
Let us keep the notation in Proposition \ref{Ewr456}.
 Then, there exists a canonical equivalence of categories
\begin{align} \label{E3GH2}
\begin{pmatrix}
\text{the groupoid  of quasi-nilpotent} \\
\text{$m$-PD stratified $G$-bundles on $X/S$} 
\end{pmatrix}
\isom \begin{pmatrix}
\text{the groupoid of $m$-crystals} \\
\text{of  $G$-bundles on $\overline{X}/S$} \\
\end{pmatrix}.
\end{align}
 \epr
%------------------------------[begin proof]-------------------
\begin{proof}
The assertion follows from   a routine argument  of ($m$-)crystals together with ~\cite[\S\,4, Proposition 4.1.7]{LQ}.
\end{proof}
%-------------------------------------------------------------
\SSP

%-------------------------------------------------------------
\begin{rema} \label{EFPPWk}
We shall describe concretely the bijective correspondence  given by   (\ref{E3GH2}).
 First, given 
 an $m$-PD stratified $G$-bundle  $(\mcE, \varepsilon^\ST)$ (where $\varepsilon^\ST := \{ \varepsilon_l \}_l$) on $X/S$, we can construct an $m$-crystal of $G$-bundles $\mcE^\dd$ on $\overline{X}/S$.
 To this end, it suffices to specify $\mcE^\dd_T$  for sufficiently small $(U, T, J, \delta) \in \mr{Ob}(\mr{Cris}^{(m)}(\overline{X}/S))$, e.g., so that there exists an $S$-morphism
 $h : T \migi X$ over $(S, \mfa, \mfb, \gamma)$ extending the open immersion $U \migiincl \overline{X}$.
 For such a morphism $h$, we shall set $\mcE^\dd_T$ to be the pull-back $h^*(\mcE)$.
 If we are given another choice of a morphism $h' : T \migi X$ such as $h$,
 then the pair $(h, h')$ defines a morphism $T \migi P_{X/S, (m)}^{l}$ for a sufficiently large $l$.
 The pull-back of $\varepsilon_l$ via this morphism yields an isomorphism $(h, h')^*(\varepsilon_l) : h'^*(\mcE) \isom h^* (\mcE)$.
 This implies that $\mcE^\dd_T$ is well-defined up to canonical isomorphism.
 Thus, the $\mcE_T^\dd$'s for various $(U, T, J, \delta)$'s form an $m$-crystal of $G$-bundles $\mcE^\dd$.
 One may verify that the resulting assignment $(\mcE, \varepsilon^\ST) \mapsto \mcE^\dd$ coincides with  (\ref{E3GH2}).
 
Conversely,   let $\mcE^\dd$ be an $m$-crystal of $G$-bundles on $\overline{X}/S$.
For each $l \geq 0$ and $i \in \{1, 2 \}$, 
this $m$-crystal  associates to  the projection  $\mr{pr}_i^l : P_{X/S, (m)}^l \migi X$  an isomorphism of $G$-bundles $\varepsilon_{l, i} : \mr{pr}_i^{l*}(\mcE^\dd_{X}) \isom \mcE^\dd_{P_{X/S, (m)}^l}$.
The composite isomorphisms $\varepsilon_{l}:= \varepsilon_{l, 1}^{-1}\circ \varepsilon_{l, 2} : \mr{pr}_2^{l*}(\mcE^\dd_{X}) \isom \mr{pr}_1^{l*}(\mcE^\dd_{X})$ for various $l$  form an $m$-PD stratification on $\mcE_X^\dd$.
The resulting assignment $\mcE^\dd \mapsto (\mcE^\dd_X, \{ \varepsilon_l \}_l)$ realizes the inverse of  (\ref{E3GH2}).

\end{rema}
%-------------------------------------------------------------

\LSP
%---------------------------[begin subsection]-------------
\subsection{Cartier's theorem for $m$-PD stratified $G$-bundles} \label{SS042}
%{\bf .}

Let $\M$ be a nonnegative integer.
In this subsection, suppose that the base  ring $R$ is defined  over $\mbF_p$.
We shall recall  from ~\cite{LQ} the definition of $p$-$\M$-curvature and Cartier's theorem for $\mcD_{X/S}^{(m)}$-modules.
Denote by $\mcK_{X/S}^{(\M)}$ the kernel of the morphism $\mcD_{X/S}^{(\M)} \migi \mcE nd_{f^{-1}(\mcO_S)} (\mcO_X)$ defining the trivial $\mcD_{X/S}^{(\M)}$-module structure  on $\mcO_X$.
If $(\mcV, \nabla)$ is a $\mcD_{X/S}^{(\M)}$-module, then the composite
\begin{align} \label{MSKO3}
{^p \psi}_{\nabla}^{(\M)} : \mcK_{X/S}^{(\M)} \migiincl \mcD_{X/S}^{(\M)} \xrightarrow{\nabla_{}} \mcE nd_{f^{-1}(\mcO_S)} (\mcV)
\end{align}
is called the {\it $p$-$\M$-curvature} of $\nabla$ (cf. ~\cite[\S\,3,  Definition 3.1.1]{LQ}).
We shall write 
\begin{align} \label{Ew90}
\mcV^\nabla
\end{align}
for the subsheaf of $\mcV$ 
on which $\mcD_{X/S}^{(\M)+}$ acts as zero, where $\mcD_{X/S}^{(\M)+}$ denotes the kernel of the natural projection $\mcD_{X/S}^{(\M)} \migi \mcO_X$.
This sheaf  
$\mcV^\nabla$
 is called 
 the {\it sheaf of horizontal sections} of $(\mcV, \nabla)$ and 
occassionaly  regarded  as an $\mcO_{X^{(\M+1)}}$-module via the underlying homeomorphism of $F_{X/S}^{(\M+1)}$.

On the other hand, let $\mcU$ be  an  $\mcO_{X^{(\M+1)}}$-module.
Then,  its pull-back $F_{X/S}^{(\M+1)*}(\mcU)$ by  $F_{X/S}^{(\M+1)}$ admits  a canonical structure  of ${\mcD}_{X/S}^{(m)}$-action
\begin{align} \label{Ew9090}
\nabla^{\mr{can}(\M +1)}_{\mcU} : \mcD_{X/S}^{(\M)}  \migi \mcE nd_{f^{-1}(\mcO_S)} (F_{X/S}^{(\M +1)*}(\mcU)), 
\end{align}
for which 
the pair  $(F_{X/S}^{(\M+1)*}(\mcU), \nabla^{\mr{can}(\M +1)}_{\mcU})$ forms a $\mcD_{X/S}^{(\M)}$-module with vanishing $p$-$\M$-curvature.
We abuse notation by writing   $\nabla^{\mr{can}(\M +1)}_{\mcU}$ for the $S$-connection
on  $F_{X/S}^{(\M +1)*}(\mcU)$  corresponding to the $\mcD_{X/S}^{(0)}$-module structure induced by $\nabla^{\mr{can}(\M +1)}_{\mcU}$.

According to  ~\cite[\S\,3, Corollary 3.2.4]{LQ}, 
 the assignments $\mcU \mapsto (F_{X/S}^{(m+1)*}(\mcU), \nabla^{\mr{can}(\M +1)}_{\mcU})$ and $(\mcV, \nabla) \mapsto \mcV^\nabla$ determine 
an equivalence of categories
\begin{align} \label{E3002}
\begin{pmatrix}
\text{the category of} \\
\text{$\mcO_{X^{(\M+1)}}$-modules} 
\end{pmatrix}
\isom \begin{pmatrix}
\text{the category of $\mcD_{X/S}^{(\M)}$-modules} \\
\text{with vanishing $p$-$\M$-curvature} \\
\end{pmatrix},
\end{align}
which is compatible with the formations of tensor products and extensions.

Next, let   $(\mcE, \varepsilon^\ST)$ be an $m$-PD stratified $G$-bundle on $X/S$.
As discussed after Definition \ref{D057}, $\varepsilon^\ST$ specifies an $m$-PD stratification on the $\mcO_X$-algebra $\mcO_\mcE$ (cf. (\ref{Ertfq})), and hence, specifies 
 a $\mcD_{X/S}^{(m)}$-module structure   $\nabla$ on $\mcO_\mcE$ (cf. Remark \ref{Eruy78}).
In particular, we obtain the $p$-$m$-curvature ${^p}\psi^{(m)}_\nabla$ of $\nabla$.

\SSP
%----------------------------------------------------------
\bde \label{rtyu5}
\begin{itemize}
\item[(i)]
We refer to ${^p}\psi_{\nabla}^{(m)}$ as the {\bf $p$-$m$-curvature} of $\varepsilon^\ST$ (or,  of $(\mcE, \varepsilon^\ST)$), and write ${^p}\psi_{\varepsilon^\ST}^{(m)} := {^p}\psi^{(m)}_{\nabla}$.
\item[(ii)]
Let $\mcE^\dd$ be an $m$-crystal of $G$-bundles on $X/S$.
We shall say that $\mcE^\dd$ is {\bf dormant} if the $m$-PD stratified $G$-bundle    corresponding to $\mcE^\dd$ via (\ref{E3GH2}) has vanishing $p$-$m$-curvature.
\end{itemize}
\ede
%----------------------------------------------------------
\SSP

%----------------------------------------------------------
\begin{rema} \label{E00943zS}
If an $m$-PD stratification $\varepsilon^\ST$ has vanishing $p$-$m$-curvature, then it is quasi-nilpotent in the sense of Definition \ref{Erty23} (cf. ~\cite[\S\,3, Remark 3.1.2]{LQ}).
\end{rema}
%----------------------------------------------------------
\SSP

Suppose  further that $\varepsilon^\ST$ has vanishing $p$-$m$-curvature.
The $\mcO_X$-algebra structure on $\mcO_\mcE$ is  restricted to an $\mcO_{X^{(\M+1)}}$-algebra structure  on the sheaf of horizontal sections $\mcO_\mcE^\nabla$ (with respect to $\nabla$).
We shall write 
\begin{align} \label{LLKS3}
\mcE^\nabla := \mcS pec (\mcO_\mcE^\nabla),
\end{align}
 which is a relative  affine scheme  over $X$.
The commutativity of  (\ref{E00100}) for    $(\mcE, \varepsilon^\ST)$ implies  that 
the $G$-coaction  $\mr{R}_\mcE^\sharp : \mcO_\mcE \migi \mcO_\mcE \otimes R_G$   is restricted  to a morphism $\mcO_\mcE^\nabla \migi \mcO_\mcE^\nabla \otimes R_G$, which  specifies a $G$-action on $\mcE$. 
By (\ref{E3002}),  there exists an isomorphism $(F^{(\M+1)*}_{X/S}(\mcO_\mcE^\nabla), \nabla^{\mr{can}(\M +1)}_{\mcO_\mcE^\nabla}) \isom  (\mcO_\mcE, \nabla)$ 
compatible with the $G$-coactions,  so we have $F_{X/S}^{(m+1)*}(\mcE^\nabla) \cong \mcE$.
Hence,  it follows from the faithful flatness of $F_{X/S}^{(m+1)}$ that  $\mcE^\nabla$ forms a $G$-bundle with respect to the $G$-action defined above.
Moreover, by reversing the steps in the construction of  $\mcE^\nabla$ using $(\mcE, \varepsilon^\ST)$, we obtain the following assertion as a generalization of  Cartier's theorem (cf. ~\cite[\S\,3.2, Corollary 3.2.4]{LQ}).

\SSP
%----------------------------------------------------------
\bpr \label{Ertwq}
(Recall that $G$ is assumed to be affine.) Let $m$ be a nonnegative integer.
\begin{itemize}
\item[(i)]
The assignment $(\mcE, \varepsilon^\ST) \mapsto \mcE^\nabla$ constructed above defines  an equivalence of categories
\begin{align} \label{E3erf002}
\begin{pmatrix}
\text{the groupoid of $m$-PD stratified $G$-bundles } \\
\text{on $X/S$ with vanishing $p$-$m$-curvature} 
\end{pmatrix}
\isom \begin{pmatrix}
\text{the groupoid of } \\
\text{$G$-bundles on  $X^{(m+1)}$} \\
\end{pmatrix}.
\end{align}
In particular, for each $G$-bundle $\mcG$ on $X^{(\M+1)}$,  there exists a canonical structure of $m$-PD stratification on $F^{(\M+1)*}_{X/S}(\mcG)$ with vanishing $p$-$\M$-curvature.
\item[(ii)]
Let $m'$ be another nonnegative  integers with $m \leq m'$.
Then, the following square diagram is $1$-commutative:
\begin{align} \label{E00153}
\begin{CD}
\begin{pmatrix}
\text{ the groupoid of $m'$-PD stratified $G$-bundles} \\
\text{ on $X/S$ with vanishing $p$-$m'$-curvature} 
\end{pmatrix} 
@> (\ref{E3erf002}) \ \text{for} \ m' > \sim >
\begin{pmatrix}
\text{ the groupoid of } \\
\text{$G$-bundles on $X^{(m'+1)}$} 
\end{pmatrix} 
\\
@V (-) |^{( m )}VV @VV F^{(m'-m)*}_{X^{(m+1)}/S}(-) V
\\
\begin{pmatrix}
\text{ the groupoid of $m$-PD stratified $G$-bundles } \\
\text{on $X/S$ with vanishing $p$-$m$-curvature} 
\end{pmatrix} 
@> \sim > (\ref{E3erf002}) \ \text{for} \ m  >
\begin{pmatrix}
\text{the groupoid of } \\
\text{$G$-bundles on $X^{(m+1)}$} 
\end{pmatrix}.
\end{CD}
\end{align}
\end{itemize}
\epr
\SSP

Here, recall the notion of an $F$-divided $G$-bundle (cf. ~\cite[\S\,2.3, Definition 9]{dS}) to construct an equivalence of categories for $\infty$-PD stratified $G$-bundles.

\SSP
%----------------------------------------------------------
\bde \label{EIFO(79}
An {\bf $F$-divided $G$-bundle} on $X/S$ is  a collection  of data
\begin{align} \label{BFIEK}
\mcG_\infty := \{ (\mcG_l, \varpi_l ) \}_{l \in \mbZ_{\geq 0}},
\end{align}
 where each $\mcG_l$ denotes a $G$-bundle on $X^{(l)}$ and each $\varpi_l$ denotes an isomorphism  of $G$-bundles $\mcG_l \isom F^*_{X^{(l)}/S}(\mcG_{l+1})$. 
 The notion of an isomorphism between $F$-divided $G$-bundles can be defined in a natural manner.
\ede
%----------------------------------------------------------

\SSP

%----------------------------------------------------------
\bco \label{C9865}
(Recall that $G$ is assumed to be affine.)
There exists an equivalence of categories
\begin{align} \label{E3erEr}
\begin{pmatrix}
\text{the groupoid of } \\
\text{$\infty$-PD stratified $G$-bundles on $X/S$} 
\end{pmatrix}
\isom \begin{pmatrix}
\text{the groupoid of} \\
\text{$F$-divided $G$-bundles on  $X/S$} \\
\end{pmatrix}.
\end{align}
\eco
%----------------------------------------------------------
\begin{proof}
The assertion follows from Proposition \ref{Ertwq}, (i),  (ii), and the final  comment in Remark \ref{Eruy78}.
\end{proof}

\LSP
%---------------------------[begin subsection]-------------
\subsection{Flat connections  and $0$-PD stratifications} \label{SS0125}

We shall focus  on the case of $m=0$.
 It is well-known that flat connections on an $\mcO_X$-module correspond to   $0$-stratifications on it.
In what follows, we generalize this correspondence to those on $G$-bundles.

Let $\pi : \mcE \migi X$ be a $G$-bundle on $X$ and  $\varepsilon^\ST := \{ \varepsilon_l \}_l$ a $0$-PD stratification on $\mcE$.
According to the discussion following Definition \ref{D057}, 
the data $\varepsilon^\ST$ corresponds to a collection  $\varepsilon^{\ST \natural} := \{ \varepsilon_l^\natural \}_l$ consisting of $\mcP_{X/S, (0)}^l$-algebra isomorphisms $\varepsilon_l^\natural: \mcP^l_{X/S, (0)} \otimes_{\mcO_X} \mcO_\mcE \isom \mcO_\mcE \otimes_{\mcO_X} \mcP_{X/S, (0)}^l$ ($l \geq 0$) as in (\ref{Ertfq}).
This collection
specifies 
 a flat $S$-connection $\nabla : \mcO_\mcE \migi  \Omega_{X/S}  \otimes_{\mcO_X} \mcO_\mcE$ on 
  $\mcO_\mcE$ (cf. ~\cite[\S\,4, Theorem 4.8]{Og2}).
 The isomorphism $\varepsilon_1^\sharp$ preserves the multiplications, so $\nabla$ satisfies the equality $\nabla (a \cdot  b) = \nabla (a) \cdot b + \nabla (b) \cdot a$ for any local sections $a$, $b \in \mcO_\mcE$.
Hence, for any local section $\partial \in \mcT_{X/S}$, the induced $f^{-1}(\mcO_S)$-linear endomorphism  $\nabla_\partial : \mcO_\mcE \migi \mcO_\mcE$ of $\mcO_\mcE$ 
defines  a local section of $\pi_*(\mcT_{\mcE/S})$.
By the commutativity of  (\ref{E00100}) for  $\varepsilon_l^\natural$, we see that  $\nabla_\partial$ lies in $\widetilde{\mcT}_{\mcE/S} \left(=  \pi_*(\mcT_{\mcE/S})^G \right)$.
Since the  restriction of $\nabla$ to $\mcO_X \left(\subseteq \mcO_\mcE \right)$ coincides with the universal derivation $d : \mcO_X \migi \Omega_{X/S}$, the equality  $d_\mcE  (\nabla_\partial) =\partial$ holds.
The morphism 
\begin{align} \label{eowip0}
\nabla_{\varepsilon^\ST} : \mcT_{X/S} \migi \widetilde{\mcT}_{\mcE/S}
\end{align}
 given by assigning $\partial \mapsto \nabla_\partial$ defines an $S$-connection on $\mcE$, which is flat because of 
the flatness of $\nabla$.
If $(\mcE', \varepsilon'^\ST)$ is another $0$-PD stratified $G$-bundle, then each isomorphism  $(\mcE, \varepsilon^\ST) \isom (\mcE', \varepsilon'^\ST)$  is, by construction,  transformed into an isomorphism of  flat $G$-bundles $(\mcE, \nabla_{\varepsilon^\ST}) \isom (\mcE, \nabla_{\varepsilon'^\ST})$.

Moreover, we may reverse the steps in this construction to obtain a $0$-PD stratification  from a flat $S$-connection.
The observation discussed here shows  the following assertion.

\SSP

%----------------------------------[begin definition]------------------
%\vspace{3mm}
\bpr \label{C071}
(Recall that $G$ is assumed to be affine.)
The assignment $(\mcE, \varepsilon^\ST) \mapsto (\mcE, \nabla_{\varepsilon^\ST})$ constructed above
defines 
 an equivalence of categories 
  \begin{align} \label{Efjj2}
\begin{pmatrix}
\text{ the groupoid  of } \\
\text{$0$-PD stratified $G$-bundles on $X/S$} 
\end{pmatrix}
\isom \begin{pmatrix}
\text{ the groupoid of } \\
\text{flat $G$-bundles  on $X/S$} \\
\end{pmatrix},
\end{align}
which is restricted  to an equivalence of categories
  \begin{align} \label{DFjj2}
\begin{pmatrix}
\text{ the groupoid  of quasi-nilpotent} \\
\text{$0$-PD stratified $G$-bundles on $X/S$} 
\end{pmatrix}
\isom \begin{pmatrix}
\text{the groupoid of quasi-nilpotent} \\
\text{flat $G$-bundles  on $X/S$} \\
\end{pmatrix}.
\end{align}
If, moreover,   the base ring $R$ is defined   over  $\mbF_p$, then 
(\ref{DFjj2})
is restricted  to an equivalence of categories
 \begin{align} \label{Efjj32}
\begin{pmatrix}
\text{ the groupoid  of $0$-PD stratified} \\
\text{$G$-bundles on $X/S$ with vanishing $p$-$0$-curvature} 
\end{pmatrix}
\isom \begin{pmatrix}
\text{ the groupoid of } \\
\text{$p$-flat $G$-bundles  on $X/S$} \\
\end{pmatrix}.
\end{align}
 \epr
%-------------------------[end definition]-------------------
\begin{proof}
The remaining  portions follow from the comments in Remark \ref{E00943zS} and ~\cite[\S\,3, Remark 3.1.2]{LQ}.
\end{proof}
%-------------------------[end definition]-------------------
\SSP

%------------------------
\bco \label{EQkji}
(Recall that $G$ is assumed to be  affine.)
There exists a canonical equivalence of categories
\begin{align} \label{EEjj32}
\begin{pmatrix}
\text{the groupoid  of } \\
\text{$G$-bundles on $X^{(1)}$} 
\end{pmatrix}
\isom \begin{pmatrix}
\text{ the groupoid of } \\
\text{$p$-flat $G$-bundles  on $X/S$} \\
\end{pmatrix}
\end{align}
such that, if  $\mcG$ is a $G$-bundle on $X^{(1)}$, then the underlying $G$-bundle of the corresponding $p$-flat $G$-bundle coincides with $F_{X/S}^*(\mcG)$.
\eco
%------------------------
\begin{proof}
The assertion follows from the equivalences of categories (\ref{E3erf002}) and (\ref{Efjj32}).
\end{proof}
\SSP

For each $G$-bundle  $\mcG$ on $X^{(1)}$, 
we shall denote by 
\begin{align} \label{QQwkko}
\nabla_\mcG^\mr{can}
\end{align}
the  $p$-flat connection on the pull-back $F_{X/S}^*(\mcG)$ corresponding to $\mcG$ via (\ref{EEjj32}).

Also, with the notation in Proposition \ref{Ewr456},
let $(\mcE, \nabla)$ be a quasi-nilpotent flat $G$-bundle on $X/S$.
Then, by composing  the equivalences of categories (\ref{E3GH2}) and (\ref{DFjj2}),
we obtain a $0$-crystal of $G$-bundles
\begin{align}\label{EWQQ21}
(\mcE, \nabla)^\dd
\end{align}
on $\overline{X}/S$ corresponding to $(\mcE, \nabla)$.

Finally, we shall conclude the present paper with the following assertion, which says 
that the category of  quasi-nilpotent flat $G$-bundles  on an infinitesimal  thickening  of a smooth scheme  
does not depend on the choice of the thickening.

\SSP
%-------------------------------------------------------------
\bpr  \label{eoqo84939}
Let us keep the assumption in Proposition \ref{Ewr456}.
Also, let $X'$  be another smooth scheme  over $S$ 
such that $\gamma$ extends to $X'$ and 
whose reduction modulo $J_0$ is isomorphic to $\overline{X}$.
Then, for each quasi-nilpotent flat $G$-bundle $(\mcE, \nabla)$ on $X/S$,
there exists a unique (up to isomorphism) quasi-nilpotent flat $G$-bundle
\begin{align}
(\mcE', \nabla')
\end{align}
on $X'/S$
with $(\mcE, \nabla)^\dd \cong (\mcE', \nabla')^\dd$.
Moreover,  the assignment $(\mcE, \nabla) \mapsto (\mcE', \nabla')$ determines an equivalence of categories between the groupoids of quasi-nilpotent flat $G$-bundles on $X/S$ and  $X'/S$.
\epr
%-------------------------------------------------------------
\begin{proof}
The assertion follows from Propositions \ref{Ewr456},  \ref{P019}, and \ref{C071}.
\end{proof}
%-------------------------------------------------------------
\SSP

%-------------------------------------------------------------
\begin{rema} \label{GGGh6}
Let us describe a down-to-earth construction  of the assignment  $(\mcE, \nabla)  \mapsto (\mcE', \nabla')$ asserted in the above proposition.
Since 
$J_0$ is nilpotent and both $X$, $X'$  are smooth over $S$,
there exists a collection of data
\begin{align}
(\{ U_\alpha \}_{\alpha \in I}, \{ \eta_{\alpha \beta} \}_{(\alpha, \beta) \in I_2})
\end{align}
where 
\begin{itemize}
\item
$I$ denotes a finite set and 
$\{ U_\alpha \}_{\alpha \in I}$ is an open covering   of $X$;
\item
$I_2$ denotes the subset of $I \times I$ consisting of pairs $(\alpha, \beta)$ with $U_{\alpha \beta} := U_\alpha \cap U_\beta \neq \emptyset$ and
each $\eta_{\alpha \beta}$ ($(\alpha, \beta) \in I_2$) is an automorphism of $U_{\alpha \beta}$ inducing the identity morphism of $U_{\alpha \beta} \times_S \overline{S}$ via reduction modulo $J_0$,
\end{itemize}
 such that $X'$ may be obtained by gluing together $U_{\alpha \beta}$'s by means of 
 $\eta_{\alpha \beta}$'s.
 Now, let $\mcE^\dd$ be the $0$-crystal of $G$-bundles corresponding  to $(\mcE, \nabla)$ via (\ref{E3GH2}) (for $m=0$) and (\ref{Efjj2}).
 This crystal  associate to  each  $\eta_{\alpha \beta}$  
  a  $G$-equivariant isomorphism  $\eta_{\mcE, \alpha \beta}$ of  $\mcE |_{U_{\alpha \beta}}$ over $\eta_{\alpha \beta}$.
(According to the discussion in Remark \ref{EFPPWk},  $\eta_{\mcE, \alpha \beta}$ also arises from pulling-back the isomorphism $\varepsilon_l$ for a sufficiently large $l$ constituting the $0$-PD stratification  corresponding to $\mcE^\dd$ by the morphism $U_{\alpha \beta} \migi P_{X/S, (0)}^l$ induced  from  $(\iota, \iota \circ \eta_{\alpha \beta}) : U_{\alpha \beta} \migi X \times_S X$, where $\iota$ denotes the natural inclusion $U_{\alpha \beta} \migiincl X$.)
Then, by  the definition of  the assignment $(\mcE, \nabla) \mapsto (\mcE', \nabla')$,
we see that the resulting flat $G$-bundle $(\mcE', \nabla')$ may be obtained by gluing together
$(\mcE |_{U_{\alpha}}, \nabla |_{U_\alpha})$'s by means of $\eta_{\mcE, \alpha \beta}$'s.
\end{rema}
%-------------------------------------------------------------

%%%%%%%%%%%%%%%%%%%%%%%%%%%%%%%%%%%%%%%%%%%%%%%%%%%

\vspace{10mm}
 \begin{thebibliography}{99}

\bibitem{BD1}
A. A. Beilinson, V. Drinfeld,
``Quantization of Hitchin's integrable system and Hecke eigensheaves."
Available at: %http://math.uchicago.edu/~mitya/langlands/hitchin/BD-hitchin.pdf.
http://math.uchicago.edu/~mitya/langlands/hitchin/BD-hitchin.pdf

\bibitem{BG}
N. Bergeron, T. Gelander,
A note on local rigidity.
{\it Geom. Ded.} {\bf  107} (2004), pp. 111-131.





\bibitem{Ber1}
P. Berthelot,
$\mcD$-modules arithm\'{e}tiques I. Op\'{e}rateurs diff\'{e}rentiels de niveau fini.
\textit{Ann. Sci. \'{E}cole Norm. Sup.} {\bf 29} (1996), pp. 185-272.


\bibitem{Ber2}
P. Berthelot,
{\it $\mcD$-modules arithm\'{e}tiques II. Descente par Frobenius}.
\textit{M\'{e}m. Soc. Math. France} {\bf 81} (2000).



\bibitem{Og2}
P. Berthelot, A. Ogus,
\textit{Notes on Crystalline Cohomology}.
Princeton Univ. Press (1978).

\bibitem{BS}
I. Biswas, 
J. P. P. dos Santos, 
 Triviality criteria for bundles over rationally connected varieties. 
 {\it J. Ramanujan Math. Soc.} {\bf 28} (2013), pp. 423-442.



\bibitem{BD2}
I. Biswas, S. Dumitrescu,
Branched holomorphic Cartan geometries and Calabi-Yau manifolds.
\textit{Int. Math. Res. Not. IMRN} {\bf 23} (2019), pp. 7428-7458.



\bibitem{BD}
I. Biswas, S. Dumitrescu,
Generalized holomorphic Cartan geometries.
\textit{Eur. J. Math.} {\bf 6} (2020), pp. 661-680.

\bibitem{BD3}
I. Biswas, S. Dumitrescu, G. Schumacher,
Deformation theory of holomorphic Cartan geometries.
\textit{Indag. Math. (N. S.)} {\bf 31} (2020), pp. 512-524.


\bibitem{CEG}
R. D. Canary, D. B. A. Epstein, P. Green,
Notes on notes of Thurston, in Analytical and geometric aspects of hyperbolic space (Coventry/Durham, 1984),
 {\it  London Math. Soc. Lecture Note Ser.,} {\bf 111},  Cambridge Univ. Press, Cambridge (1987), pp. 3-92.

\bibitem{CG}
D. Cooper, W. Coldman
A $3$-manifold with no real projective structure.
\textit{Annales de la Facult\'{e} des Sciences de Toulouse} {\bf 24} (2015), pp. 1219-1238.

\bibitem{dS}
J. P. P. dos Santos,
Fundamental group schemes for stratified sheaves.
{\it  J. Algebra} {\bf 317} (2007), pp. 691-713.



\bibitem{Ehr}
C. Ehresmann,
Sur les espaces localement homog\'{e}nes.
\textit{L'Enseign. Math.} {\bf 35} (1936), pp. 317-333.


\bibitem{EM}
H. Esnault, V. Mehta,
Simply connected projective manifolds in characteristic $p>0$ have no nontrivial stratified bundles
{\it Invent Math.} {\bf 181} (2010), pp. 449-465. 


 \bibitem{FGA}
B. Fantechi, L. G\"{u}ttsche, L. Illusie, S. Kleiman, N. Nitsure, A. Vistoli, 
{\it Fundamental algebraic geometry. Grothendieck's FGA explained.}.
Mathematical Surveys and Monographs, {\bf 123}  AMS (2005).

 \bibitem{FaMa}
B. Fantechi,  M. Manetti,
 Obstruction calculus for functors of Artin rings, I.
 {\it Journal of Algebra}  {\bf 202} (1998), pp. 541-576.

\bibitem{Fr}
E. Frenkel,
{\it Langlands Correspondence for Loop Groups}.
Cambridge Studies in Advanced Mathematics {\bf 103} Cambridge Univ. Press (2007).


\bibitem{G}
D. Gieseker,
Flat vector bundles and the fundamental group in non-zero characteristics,
{\it Ann. Scuola Norm. Sup. Pisa Cl. Sci.} {\bf 2} (1975), pp. 1-31.



\bibitem{Gol1}
W. M. Goldman, 
Locally homogeneous geometric manifolds, 
{\it Proceedings of the ICM} (2010).


\bibitem{EGA4}
A. Grothendieck,
{\it \'{E}l\'{e}ments de G\'{e}om\'{e}trie Alg\'{e}brique, IV}:
  {\it \'{E}tude locale des sch\'{e}mas et des morphismes de sch\'{e}mas, Quatri\`{e}me partie}.  
Publ. Math. I.H.E.S. {\bf 32} (1967).

\bibitem{G2}
R. C. Gunning,
Special coordinate covering of Riemann surfaces.
\textit{Math. Ann.}  {\bf 170} (1967),  pp. 67-86.

\bibitem{Hej}
D. Hejhal, 
Monodromy groups and linearly polymorphic functions,
{\it  Acta Math.} {\bf 135} (1975), pp. 1-55.


 \bibitem{Hos2}
Y. Hoshi,
Frobenius-projective structures on curves in positive characteristic.
{\it Publ. Res. Inst. Math. Sci. } {\bf 56} (2020), pp. 401-430.


\bibitem{Hos4}
Y. Hoshi,
Frobenius-affine structures and Tango curves.
RIMS Preprint {\bf 1913} (2020).

\bibitem{Hos3}
Y. Hoshi, 
A note on the existence of Tango curves.
{\it Kodai Math. J.} {\bf 44} (2021), pp. 77-80.


\bibitem{Jo14}
K. Joshi,
The degree of the dormant operatic locus.
{\it Internat. Math. Res. Notices} {\bf 9} (2017), pp. 2599-2613.


\bibitem{Kat}
N. M. Katz,
Nilpotent connections and the monodromy theorem: Applications of a result of Turrittin.
\textit{Inst. Hautes Etudes Sci. Publ. Math.} {\bf 39}, (1970), 175-232.



\bibitem{Katz2}
N. M. Katz,
Algebraic solutions of differential equations ($p$-curvature and the Hodge filtration).
\textit{Invent. Math.} {\bf 18} (1972),  pp. 1-118.


\bibitem{KO1}
S. Kobayashi, T. Ochiai,
Holomorphic projective structures on compact complex surfaces I.
\textit{Math. Ann.} {\bf 249} (1980), pp. 75-94.

\bibitem{KO2}
S. Kobayashi, T. Ochiai,
Holomorphic projective structures on compact complex surfaces II.
\textit{Math. Ann.} {\bf 255} (1981), pp. 519-521.


\bibitem{Lan}
A. Langer,
Generic positivity and foliations in positive characteristic.
{\it Adv. Math.} {\bf  277} (2015), pp. 1-23. 


\bibitem{Lau}
B. Laurent,
Almost homogenous curves over an arbitrary field.
{\it Transformation Groups} {\bf 24} (2019), pp. 845-886.

\bibitem{LQ}
B. Le Stum, A. Quir\'{o}s,
Transversal crystals of finite level,
{\it  Ann. Inst. Fourier (Grenoble)} {\bf  47} (1997), pp. 69-100.

\bibitem{Lok}
W. L. Lok,
Deformations of locally homogeneous spaces and Kleinian groups, 
 Doctoral Thesis, Columbia University (1984).

\bibitem{Mil4}
J. S. Milne,
{\it \'{E}tale cohomology}.
Princeton Univ. Press, Princeton, New Jersey (1980)


\bibitem{MIL3}
J. Milne,
 Abelian Varieties.
\textit{Arithmetic Geometry}, Springer, New York, (1986),  pp. 103-150.


\bibitem{Mil2}
J. S. Milne,
{\it Arithmetic Duality Theorem}.
Perspectives in Mathematics, {\bf 1}  
Academic Press Inc., Boston, MA (1986).



\bibitem{Mil}
J. S. Milne,
{\it Algebraic groups. The theory of group schemes of finite type over a field}.
Cambridge Studies in Advanced Mathematics, {\bf 170}  Cambridge University
Press, Cambridge (2017).



\bibitem{Mzk1}
S. Mochizuki,
A theory of ordinary $p$-adic curves.
\textit{Publ. RIMS} {\bf 32} (1996),  pp. 957-1151.

\bibitem{Mzk2}
S. Mochizuki,
{\it Foundations of $p$-adic Teichm\"{u}ller theory}.
American Mathematical Society,  (1999).


\bibitem{Muk}
S. Mukai, 
Counterexamples to Kodaira's vanishing and Yau's inequality in positive
characteristics.
{\it  Kyoto J. Math.} {\bf  53} (2013), pp. 515-532.

\bibitem{Nit}
N. Nitsure,
{\it Notes on Deformation theory.}
Lecture note, (2006).

\bibitem{Og}
A. Ogus,
\textit{$F$-Crystals, Griffiths Transversality, and the Hodge Decomposition}.
Ast\'{e}risque {\bf 221}, Soc. Math. de France, (1994).

\bibitem{Og3}
A. Ogus,
Higgs cohomology, $p$-curvature, and the Cartier isomorphism.
\textit{Compositio. Math.} {\bf 140} (2004),  pp. 145-164.



 \bibitem{Ogu3}
A. Ogus,
\textit{Lectures on Logarithmic Algebraic Geometry}.
Cambridge Univ. Press, (2018).

\bibitem{Ray}
M. Raynaud, 
Contre-example au ``vanishing theorem'' en caracteristique $p>0$. 
{\it C.P. Ramanujan-a tribute, Tata Inst. Fund. Res. Studies in Math.} {\bf 8} Springer, Berlin-New York (1978), pp. 273-278.



\bibitem{Sha}
R. W. Sharpe,
{\it Differential geometry. Cartan's generalization of Klein's Erlangen program}.
Graduate Texts in Mathematics {\bf 166} Springer-Verlag, New York   (1997).


\bibitem{Sch}
M. Schlessinger,
Functors of Artin rings. 
{\it Trans. Amer. Math. Soc.} {\bf 130} (1968), pp. 208-222. 


\bibitem{TV}
M. Talpo, A. Vistoli,
Deformation theory from the point of view of fibered categories,
{\it Handbook of moduli, Vol. III, Adv. Lect. Math.} {\bf 26}, Int. Press, Somerville, MA  (2013), pp. 281-397.


\bibitem{Thu}
W. P. Thurston,
The geometry and topology of three-manifolds,
{\it lecture notes} (1980).




\bibitem{Wak}
Y. Wakabayashi,
An explicit formula for the generic number of dormant indigenous bundles.
\textit{Publ. Res. Inst. Math. Sci.} {\bf 50}  (2014),  pp. 383-409.


\bibitem{Wak7}
Y. Wakabayashi,
Moduli of Tango structures and dormant Miura opers,
\textit{Moscow Math. J.} {\bf 20}  (2020),  pp. 575-636.


\bibitem{Wak6}
Y. Wakabayashi,
Frobenius projective and affine geometry of varieties in positive characteristic.
\textit{arXiv: math. AG/2011.04846} (2020).

\bibitem{Wak8}
Y. Wakabayashi,
A theory of dormant opers on pointed stable curves ---a proof of Joshi's conjecture---,
\textit{arXiv: math. AG/1411.1208v4}, (2021).


\bibitem{Wak9}
Y. Wakabayashi,
Dormant Miura opers, Tango structures, and the Bethe ansatz equations modulo $p$,
\textit{arXiv: math. AG/1905.03364v2}, (2020).
%\textit{arXiv: math. AG/1709.04241}, (2017).

\bibitem{Wei1}
A. Weil, 
On discrete subgroups of Lie groups I,
{\it  Ann. Math.} {\bf  72} (1960), pp. 369-384.

\bibitem{Wei2}
A. Weil,
On discrete subgroups of Lie groups II,
{\it  Ann. Math.}  {\bf  75} (1962), pp. 578-602.

\bibitem{Wei3}
A. Weil,
Remarks on the cohomology of groups, 
{\it Ann. Math.} {\bf 80} (1964), pp. 149-157.




\end{thebibliography}
\end{document}